%% file: 20607.tex
\documentclass[11pt]{article}

\usepackage{amsmath,amssymb}
\usepackage{latexsym}
\usepackage{graphicx}
\usepackage{amscd}
\usepackage{theorem}
\input xy
\xyoption{all}

\newcommand{\R}{{\mathbb{R}}}
\newcommand{\Z}{{\mathbb{Z}}}
\newcommand{\C}{{\mathbb{C}}}
\newcommand{\N}{{\mathbb{N}}}
\newcommand{\T}{{\mathbb{T}}}

\newcommand{\bea}{\begin{eqnarray}}
\newcommand{\eea}{\end{eqnarray}}

\newcommand{\nn}{\nonumber}

\newcommand{\bp}{\begin{pmatrix}}
\newcommand{\ep}{\end{pmatrix}}

\newcommand{\bps}{\begin{smallmatrix}}
\newcommand{\eps}{\end{smallmatrix}}

\newcommand{\ti}{\tilde}
\newcommand{\la}{\langle}
\newcommand{\ra}{\rangle}

\def \cH{{\cal H}}
\def \V{{\cal V}}
\def \T{{\cal T}}

\def \qed{\hfill $\blacksquare$}

\def \lglraw{\longleftrightarrow}
\def \lgraw{\longrightarrow}

\def \hraw{\hookrightarrow}

\def \Id{\mathrm{Id}}
\def \Im{\mathrm{Im}}

\def \deg{\mathrm{deg}}

\def \gf{\mathit{gf}}
\def \rank{\mathrm{rank}}

\def \mes{{\cal D}}

\def \fb{{{\cal F}^*}}

\def \half{\frac{1}{2}}
\def \ov#1{\frac{1}{#1}}

\def \cA{{\cal A}}

\def \fpart#1#2{\frac{\partial #1}{\partial #2}}
\def \flpartial#1{\frac{\overleftarrow{\partial}}{\partial #1}}
\def \frpartial#1{\frac{\overrightarrow{\partial}}{\partial #1}}
\def \flpart#1#2{\frac{#1 \overleftarrow{\partial}}{\partial #2}}
\def \frpart#1#2{\frac{\overrightarrow{\partial} #1}{\partial #2}}
\def \fd#1{\frac{d}{d{#1}}}

\def \g{{\frak g}}
\def \l{{\frak l}}
\def \m{{\frak m}}

\def \cF{{\cal F}}
\def \cG{{\cal G}}

\def \cO{{\cal O}}
\def \cM{{\cal M}}

\def \cMC{{\cal MC}}

\def \cP{{\cal P}}
\def \tcM{{\widetilde {\cal M}}}
\def \tOmega{{\tilde \Omega}}

\def \l({\left(}
\def \r){\right)}

\def \0{{\bf 0}}
\def \1{{\bf 1}}

\def \eb{{\bf e}}

\def \m{{\frak m}}
\def \tri{\triangle}

\def \-{-\hspace*{-0.2cm}-}
\def \ie{{\it i.e.}\ }

\def \cb{{\bar c}}
\def \rh{{\check r}}

\def \mm{\supset\hspace*{-0.2cm}-}
\def \D{-\hspace*{-0.05cm}\framebox{$\delta$}\hspace*{-0.05cm}-}
\def \f{-\hspace*{-0.05cm}\framebox{$\fb$}\hspace*{-0.05cm}-}
\def \Q{-\hspace*{-0.05cm}\framebox{$\m$}\hspace*{-0.05cm}-}
\def \F{-\hspace*{-0.05cm}\framebox{$\cF$}\hspace*{-0.05cm}-}

{%\theorembodyfont{\rmfamily}

 \newtheorem{thm}{Theorem}[section]
 \newtheorem{prop}[thm]{Proposition}%[section]
 \newtheorem{lem}[thm]{Lemma}%[section]
 \newtheorem{cor}[thm]{Corollary}%[section]

\theorembodyfont{\rmfamily}

 \newtheorem{defn}[thm]{Definition}%[section]
 \newtheorem{rem}[thm]{Remark}%[section]
}

\numberwithin{equation}{section}

\begin{document}

\begin{titlepage}
\thispagestyle{empty}
\begin{flushleft}
\hfill v1:\,  January,\, 2003\\
\hfill v2: February, 2007\\
\end{flushleft}

\vskip 3.2 cm

\begin{center}
\noindent{\Large \textbf{Noncommutative homotopy algebras 
}}\\

\vspace*{0.5cm}

\noindent{\Large \textbf{associated with open strings
 }}\\
\renewcommand{\thefootnote}{\fnsymbol{footnote}}

\vskip 2cm
{\large 
Hiroshige Kajiura 
%\footnote{e-mail address: kuzzy@ms.u-tokyo.ac.jp}\\
%\footnote{e-mail address: kajiura@yukawa.kyoto-u.ac.jp}\\

\noindent{ \bigskip }\\

\it
%Graduate School of Mathematical Sciences, University of Tokyo \\
%3-8-1 Komaba, Tokyo 153-8914, Japan\\
Yukawa Institute for Theoretical Physics, Kyoto University \\
Kyoto 606-8502, Japan\\
e-mail: kajiura@yukawa.kyoto-u.ac.jp
\noindent{\smallskip  }\\
}

\bigskip
\bigskip
\end{center}



%abst

\begin{abstract}

We discuss general properties of $A_\infty$-algebras and their applications 
to the theory of open strings. 
The properties of cyclicity for $A_\infty$-algebras are examined in detail. 
We prove the decomposition theorem, which is a stronger version of 
the minimal model theorem, for $A_\infty$-algebras and 
cyclic $A_\infty$-algebras and discuss various consequences of it. 
In particular it is applied to classical open string field theories and 
it is shown that all classical open string field 
theories on a fixed conformal background are cyclic $A_\infty$-isomorphic 
to each other. 
The same results hold for classical closed string field theories, 
whose algebraic structure is governed by cyclic $L_\infty$-algebras. 

\end{abstract}
\vfill

\end{titlepage}
\vfill
\setcounter{footnote}{0}
\renewcommand{\thefootnote}{\arabic{footnote}}
\newpage

\tableofcontents

%section1

\section{Introduction and Summary}

This paper is the extended version of \cite{Ka}. 
We shall discuss general properties of homotopy algebras and 
their application to string theory. 
Homotopy algebras and string theory are related to each other. 
General properties of homotopy algebra govern general properties of 
field theory of string, whereas, string theory or field theory gives 
some insight into the theory of homotopy algebras. 
The direct connection between them is realized in terms of 
formal supermanifolds. 
Homotopy algebras are described (by taking their duals) in terms of 
formal supermanifolds, of which their coordinates are just the fields of 
field theories. 
We concentrate on the theory related to tree-level open strings, 
whose relevant homotopy algebraic structures are $A_\infty$-algebras 
\cite{Sta11}. 
$A_\infty$-algebras appearing in open string theory 
have an additional structure, the cyclicity. 
We call them cyclic $A_\infty$-algebras and examine their 
properties in detail. 
We also give a statement of 
formal noncommutative (odd) symplectic supergeometry 
and examine its local properties. 
It serves as a realization of $A_\infty$-algebras equipped with cyclicity. 
The minimal model theorem \cite{kadei1}  plays a key role 
in studying homotopy algebraic properties of $A_\infty$-algebras. 
We prove a stronger version of the minimal model theorem, 
which we call the decomposition theorem, for $A_\infty$-algebras and 
cyclic $A_\infty$-algebras. 
For $A_\infty$-algebras, a similar result is obtained independently 
in \cite{Le-Ha}. 
Various consequences of the decomposition theorem are then discussed. 
In particular it is applied to the classification of 
classical open string field theories.

In this section, 
we shall provide some background and main ideas of the present work. 
In subsection \ref{ssec:Ainftysp}, 
we shall first recall 
some background history of $A_\infty$-algebras. 
The construction of string field theory and 
the relevance of homotopy algebraic structures to them 
are reviewed in subsection \ref{ssec:osft}. 
In subsection \ref{ssec:Introdual}, we present some of our notations 
related to formal supermanifolds, which play a central role in this paper. 
Subsection \ref{ssec:Introphys} consists of additional comments for 
the noncommutativity of formal supermanifolds and 
their connection to physics of open strings. 
Subsection \ref{ssec:Introsym} is devoted 
to showing the idea of the construction 
of formal noncommutative symplectic supergeometry 
inspired from open strings. 
Subsection \ref{ssec:osft}, subsection \ref{ssec:Introdual} 
and subsection \ref{ssec:Introsym} 
include our basic concept and tools leading 
to some of the main results of this paper. 
The contents and the results of this paper are summarized 
in subsection \ref{ssec:contents}. 
Since in later sections we assume no knowledge presented in this section, 
the readers can skip this section 
and begin with section \ref{sec:Ainfty}.

 \subsection{$A_\infty$-spaces and $A_\infty$-algebras}
\label{ssec:Ainftysp}

\begin{figure}[h]
 \hspace*{0.8cm}\includegraphics{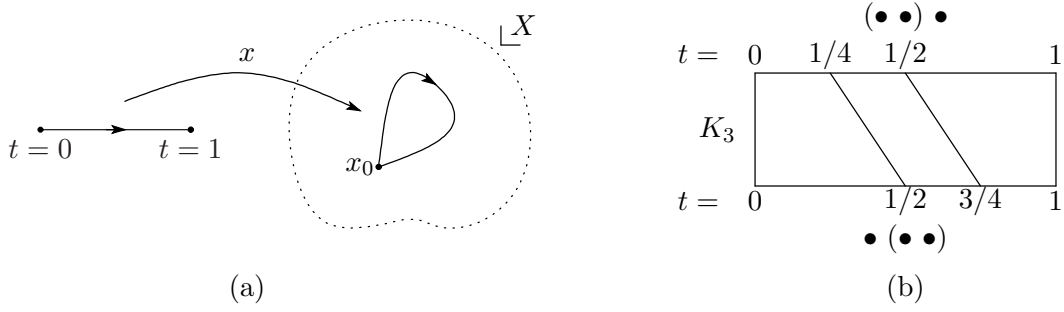}
 \caption{(a). An element in $Y$. (b). A homotopy between 
$m_2(m_2(\bullet,\bullet))$ and $m_2(\bullet,m_2(\bullet,\bullet))$, 
where $\bullet$ symbolizes an element in $Y$. }
 \label{fig:bloop}
\end{figure}

\begin{figure}[h]
 \hspace*{1.3cm}
 \includegraphics{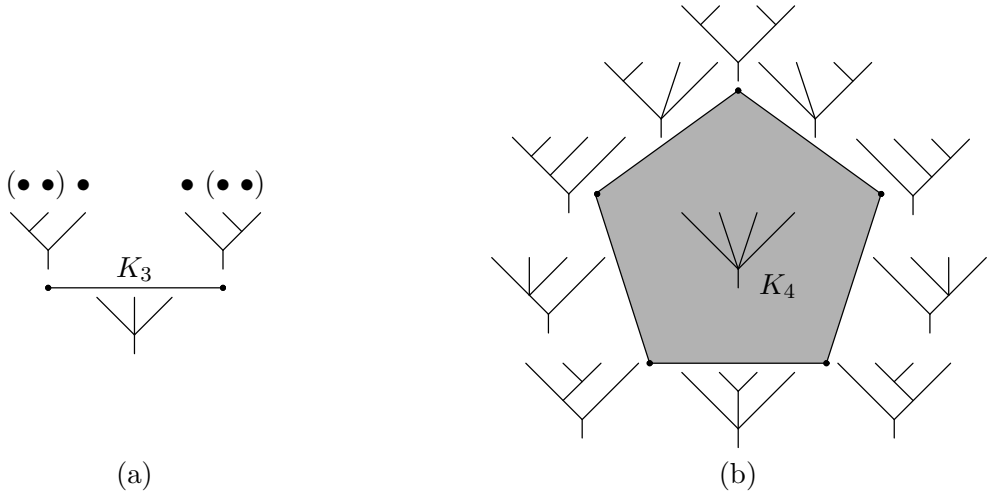}
 \caption{(a). An interval as the associahedron $K_3$. (b). 
A pentagon as the associahedron $K_4$. }
 \label{fig:Ainftysp}
\end{figure}

\begin{figure}[h]
 \hspace*{0.7cm}
 \includegraphics{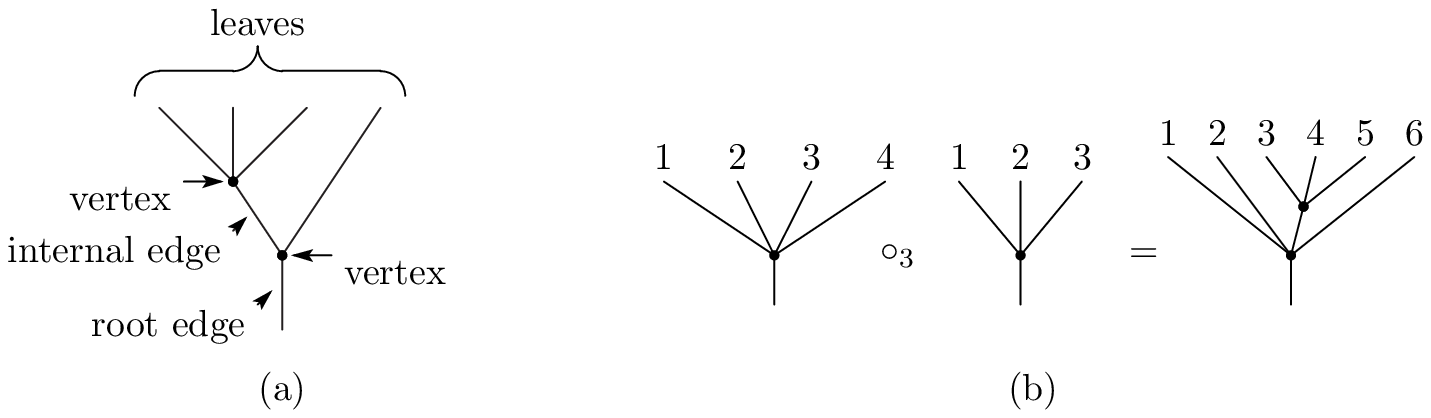}
 \caption{(a). Notation for planar rooted tree. The above one is a 
$4$-tree. (b). An example of grafting, grafting of a $3$-corolla to 
a $4$-corolla along leaf $3$. }
 \label{fig:tree}
\end{figure}

An {\em $A_\infty$-space} was introduced by J.~Stasheff 
as a tool in the study of {\em $H(opf)$-spaces} \cite{Sta1,Stabook}. 
Roughly speaking, $H$-spaces are group-like topological 
spaces. A typical example is a {\em based loop space}. 
Let $Y=\Omega X$ be the space of based loops in $X$. 
For a based point $x_0\in X$, an element of $Y$ is a map 
$x:[0,1]\to X$ where $x(0)=x(1)=x_0$ (Figure \ref{fig:bloop} (a)). 
We have a product as a group-like space 
\begin{equation*}
 m_2 : Y\times Y\to Y\ .
\end{equation*}
It is given by connecting two loops 
as $m_2(x,x')(t)=x(2t)$ for $0\le t\le 1/2$ and 
$m_2(x,x')(t)=x'(2(t-1/2))$ for $1/2\le t\le 1$. 
$m_2$ is not associative but clearly there exists a homotopy 
described by an interval $K_3$ (Figure \ref{fig:bloop} (b))
\begin{equation*}
 m_3 : K_3\times Y\times Y\times Y\lgraw Y\ .
\end{equation*}
When we represent the product by a trivalent planar tree, 
the relation above is characterized pictorially as 
in Figure \ref{fig:Ainftysp}(a). 
Next, when considering possible operations of $(Y)^{\times 4}\to Y$ by 
$m_2$, we have five vertices corresponding to tree graphs 
which consists of trivalent trees. 
Then one gets Figure \ref{fig:Ainftysp}(b) corresponding to the 
`homotopy pentagon relation'. 
One can see that each edge corresponds to $K_3$ and $K_4$ 
bounded by these edges is a pentagon. 
The corresponding homotopy $m_4 :K_4\times (Y)^{\times 4}\to Y$ 
is then defined. 
Repeating this procedure then produces higher homotopies 
\begin{equation*}
 m_n : K_n\times (Y)^{\times n}\lgraw Y\ .
\end{equation*}
For $n\ge 2$, $K_n$ is a polytope of dimension $(n-2)$; 
$K_2$ is a point, $K_3$ is a interval, $K_4$ is a pentagon as above, 
and so on. As indicated in Figure \ref{fig:Ainftysp}(a) or (b), 
$K_n$ is associated with an $n$-corolla, where an
$n$-corolla is an $n$-tree without internal edges and an
$n$-tree is a planar rooted tree with $n$ leaves. 
For a planar rooted $k$-tree, $l$-tree and an integer $1\le i\le k$, 
one can consider the {\it grafting} of $l$-tree to $k$-tree 
along leaf $i$, given by identifying the root edge of the $l$-tree 
with the $i$-th leaf of $k$-tree (see Figure \ref{fig:tree} (b)). 
Associated to the grafting, one can consider the following inclusion 
map 
\begin{equation*}
 K_k\circ_i K_l\hraw K_{k+l-1}\ .
\end{equation*}
Then, by construction, 
the $\{K_n\}_{n\ge 2}$ have the following recursion relation 
\begin{equation}
 \partial K_n=\sum_{\substack{k+l=n+1\\ k,l\ge 2}}
 \sum_{i=1}^k K_k\circ_i K_l
 \label{hedrarec}
\end{equation} 
for the codimension one boundary of $K_n$. 
One can confirm eq.(\ref{hedrarec}) in the case of $n=4$, 
where the summation in the right hand side produces five terms; 
the terms for $k=2,\ i=1,2$ and $k=3,\ i=1,2,3$. 
They corresponds to the five edges of the pentagon in Figure 
\ref{fig:Ainftysp} (b). 
The trees associated to the edges are just the ones 
associated to $K_k\circ_i K_l$. 
There also exist other relations for lower components 
(codimension greater than one boundaries) of $K_n$. 

Generally, a topological space $Y$ equipped with higher homotopies 
$\{m_n\}_{n\ge 2}$ as above is called an 
{$A_\infty$-space} \cite{Sta1}
(for a brief review see \cite{Staoperads}, an origin of this concept 
is M.~Sugawara's work \cite{Sugawara}). 
It is applied to the study of loop
spaces \cite{Adams,BoVo,May}. 
Conversely, it is known that any topological space $Y$ 
that admits the structure of an $A_\infty$-space and 
whose connected components form a group 
is homotopy equivalent to a loop space \cite{Adams}. 
It also appears in the construction of a classical open string 
field theory as will be mentioned in the next subsection. 

The set of associahedra 
$\{K_n\}_{n\ge 2}$ is one of the most typical example of topological 
operads. 
Though we avoid presenting the complicated definition, 
a $(\cdots)$ operad \cite{May} 
is a set of $(\cdots)$-objects that correspond to
corollas and are equipped with natural structures associated with 
trees and their grafting. 
(In the case here, $(\cdots)=$ `topological'. )
The set $\{K_n\}_{n\ge 2}$ is associated to a {\em non-symmetric} operad 
for which the corresponding trees are {\em planar}. 
\footnote{Here non-symmetric corresponds to noncocommutative 
in coalgebra description of $A_\infty$-algebras in section \ref{sec:Ainfty}.
} 
It is known that for any topological operad, the singular chain complex 
forms a differential graded (dg) operad. 
Since the associahedra are presented as cell complexes and the composition 
of trees is cellular, the cellular chains form a dg operad. 
Then the algebra over the dg operad is an
$A_\infty$-algebra \cite{Sta11}, see below. 

\begin{figure}[h]
 \hspace*{0.8cm}
 \includegraphics{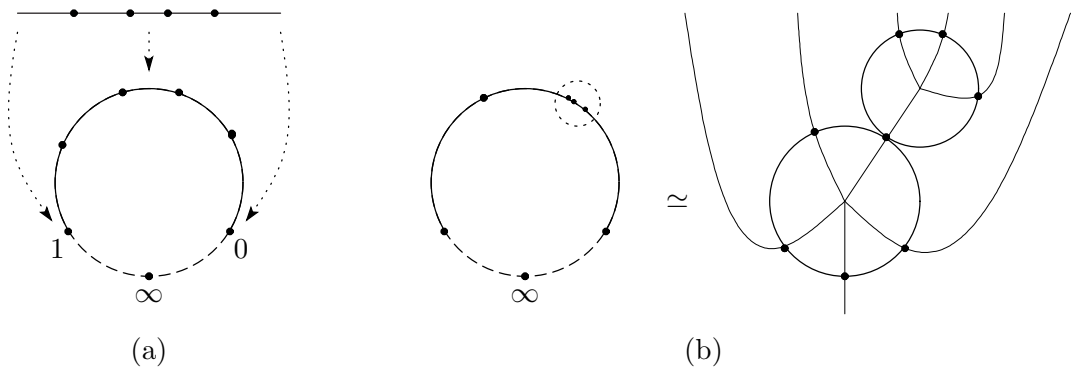}
 \caption{(a). The identification of the interval with $(n-2)$ points on it 
with the boundary of the disk with $(n+1)$ points on the boundary. 
(b). Compactification of moduli space $\cM_7$. The figure above 
represents a boundary component of $\cM_7$. It just corresponds to 
the grafting of trees in Figure \ref{fig:tree} (b).}
 \label{fig:comp1}
\end{figure}

The theory of operads and trees are closely related to 
compactification of configuration spaces. 
It is known that $K_n$ is obtained as the real compactification 
of $(n-2)$ distinct points in an interval 
(cf. the {\em little interval operad}; see \cite{MSS}, p94). 
The configuration space can further be related to the real 
compactification of the moduli space $\cM_{n+1}$ 
of a disk with $(n+1)$ points on the boundary as indicated in 
Figure \ref{fig:comp1}.
The compactified moduli space $\cM_{n+1}$ 
is defined as the configuration space of $(n+1)$-punctures 
on $S^1$ divided by conformal transformations. 
In the case when the Riemann surface is the disk, 
the conformal transformations are elements of 
$SL(2,\R)$. The degree of freedom can be killed by 
fixing three points on the boundary. As usual we denote the three points 
by $0$, $1$ and $\infty$. By choosing $\infty$ as the `root edge', 
the interval is naturally identified with the arc between $0$ and $1$. 
The pattern of the degeneration of points on the boundary is just 
the same as the right hand side of eq.(\ref{hedrarec}), 
which has $2+3+\cdots +(n-1)$ terms 
corresponding to the boundary components. 
Other interesting examples of topological operads and 
their connection to compactifications can be found 
for instance in \cite{GiK,Staconf}. 

Let $\cH$ be a $\Z$-graded vector space and 
$\m:=\{m_n : (\cH)^{\otimes n}\to \cH\}_{n\ge 1}$ a collection of 
multilinear maps. The pair $(\cH,\m)$ is then an $A_\infty$-algebra iff 
$\m$ satisfies the following relations 
(see also Definition \ref{defn:Ainfty})
\begin{equation}
 m_1 m_n
 +\sum_{i=1}^n m_n(\1^{\otimes i-1}\otimes m_1\otimes\1^{\otimes n-i})
=-\sum_{\substack{k+l=n+1\\ k\ge 2, l\ge 2}}\sum_{j=1}^k
 m_k(\1^{\otimes j-1}\otimes\m_l\otimes\1^{\otimes k-j})
 \label{IntroAinfty}
\end{equation}
on $\cH^{\otimes n}$ for each $n\ge 1$. 
The equation for $n=1$ is just $(m_1)^2=0$, which implies that 
$(\cH, m_1)$ forms a complex. 
The equation for $n=2$ is then the Leibniz rule for the action of 
derivation $m_1$ on $m_2$. 
For $n=3$ eq.(\ref{IntroAinfty}) describes 
the associativity of $m_2$ up to homotopy. 
Comparing this with an $A_\infty$-space, one can see that 
a topological space $Y$ corresponds to a graded vector space $\cH$ 
with the 
$m_i$ for $i\ge 2$ on the two sides corresponding to each other, 
where the action of $\partial$ on $K_n$ corresponds to the action of 
$m_1$ on $m_n$ in the left hand side 
of eq.(\ref{IntroAinfty}). 
Namely, the correspondence is in some sense similar to the one between 
singular homology and deRham cohomology. 
This paper deals with this algebra side, 
some `deRham rings up to homotopy'.

Such algebraic treatments of homotopy theory were developed 
in rational homotopy theory by Quillen \cite{Q} and Sullivan \cite{Sul,GM}. 
In particular \cite{Sul} deals with differential forms on a manifold $M$, 
which form a differential graded algebra (dga). 
It is then shown that 
the dga of differential forms on $M$ has the information of 
the rational homotopy type of $M$. 
Note that a dga is an $A_\infty$-algebra $(\cH,\m)$ 
with $m_3=m_4=\cdots=0$. 
In particular, in this situation, the graded vector space $\cH$ is 
the space of differential forms on $M$, 
$m_1$ is the exterior derivative and $m_2$ is the wedge product. 
For $A_\infty$-algebras, there is a notion of homotopy. 
Two homotopy equivalent $A_\infty$-algebras are transformed to each other 
by an $A_\infty$-quasi-isomorphism, 
where quasi-isomorphisms are 
morphisms which preserves the cohomology with respect to $m_1$ 
(Definition \ref{defn:quasiisom}). 
Then it is known that, 
for a given $A_\infty$-algebra $(\cH,\m)$, 
there exists an $A_\infty$-structure on $H(\cH)$, 
the cohomology of the complex $(\cH, m_1)$, 
which is $A_\infty$-quasi-isomorphic to the original $A_\infty$-algebra 
$(\cH,\m)$ \cite{kadei1} (\cite{kadei2} for the case $(\cH,\m)$ is a dga). 
This fact is called the {\it minimal model theorem}. 
The way of constructing minimal models of $A_\infty$-algebras, 
in particular dgas, 
has been developed 
in the framework of {\it homological perturbation theory} 
as an important subject in algebraic topology 
(for example \cite{gugen,GS,gugen-lambe,GLS:chen,GLS,hueb-kadei}, 
and see also \cite{jh-jds} for the dg Lie algebra case). 
The minimal model theorem implies, 
for the case of dga of differential forms, 
that one can recover the rational homotopy type of $M$ by 
considering the $A_\infty$-structure on the deRham cohomology 
instead of the original deRham complex. 
In this case, 
the higher operations $\{m_n\}_{n\ge 3}$ are related to the 
higher Massey-Yoneda products. 
One may also consider the (complex of) modules over $M$ and 
$Ext$ between them. Correspondingly, there are the notion of 
$A_\infty$-modules over $M$ and an $A_\infty$-category on $M$ 
(see \cite{Keller,Le-Ha}). 
It is then known that 
the stories stated above hold in a similar way as for $A_\infty$-algebras. 

Such notions are applied to mathematical physics in many ways. 
One of the application is the {\em homological mirror symmetry 
conjecture} \cite{mirror} 
which states some equivalences between 
an $A_\infty$-category \cite{Fukaya} on Calabi-Yau manifolds 
$M$ ($A$-model side in physical terms) 
and the category of coherent sheaves on the mirror dual manifold ${\hat M}$ 
($B$-model side). 
This conjecture implies that both sides, that is, not only 
the $A$ but also the
$B$-model sides have some $A_\infty$-structures. 
It is known that in some restricted situations both $A$ and $B$ model are 
described by topological Chern-Simons field theories \cite{W2} and 
one can obtain so-called D-brane superpotentials 
from the topological Chern-Simons field theories \cite{Laz,Tom,LazR}. 
This is nothing but the minimal model theorem, where 
a topological Chern-Simons field theory has a structure of dga and a 
D-brane superpotential is regarded as 
the collection of higher Massey-Yoneda products \cite{W2}. 

Furthermore, not only the Chern-Simons field theory above but 
any field theory has a homotopy algebraic structure generally 
only if it satisfies a classical Batalin-Vilkovisky (BV-) master equation
(see subsection \ref{ssec:Introdual}). 
The typical examples are classical string field theories explained 
in the next subsection.

 \subsection{$A_\infty$-structures and classical open string field theory}
\label{ssec:osft}

{\em String field theory} is defined on a fixed conformal background 
of space-time (target space) $M$ to which world sheet of strings 
(Riemann surfaces) are mapped, 
where a conformal background is a background (metric, etc.) of $M$ 
so that the action of a string on $M$ 
has conformal symmetry (see \cite{Sencbg}). 
There exists several classes of string field theories 
corresponding to the classes of Riemann surfaces. 
The most general one is open-closed string field theory \cite{Z2}, 
which associates to the most general class of Riemann surfaces; 
Riemann surfaces with boundaries, genera and punctures. 
It includes various `sub-string field theories'; 
classical open string field theories - associated to 
disks (one boundary and no genus) with punctures only on the boundary, 
classical closed string field theories - associated to 
spheres (no boundary and no genus) with punctures, 
quantum closed string field theories - 
associated to Riemann surfaces with punctures (and genera) and 
without boundary, and so on. 
Genus and multi-boundaries relate to loops of closed strings 
and open strings, respectively. 
We use the term `classical' (resp. quantum) 
for theory without such loop (resp. with such loops). 
There exists an abstract standard way for constructing 
these string field theories \cite{LPP,Z1}. 
We shall review it briefly in the case of 
classical open string field theories below. 
The essence is the same for the other SFTs. 

\begin{figure}[h]
 \hspace*{1.3cm}\includegraphics{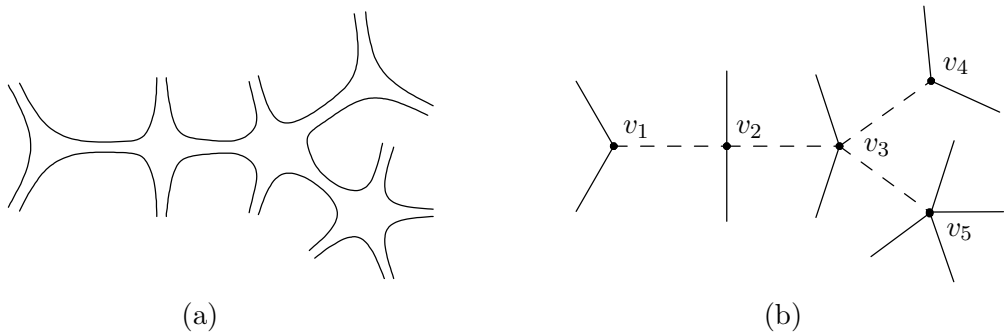}
 \caption{(a). A Riemann surface (disk) 
that describes an open string interaction. 
The corresponding Feynman diagram is the planar graph in Figure (b). 
The dashed lines denote the internal edges that correspond to 
propagators in the string field theory. 
Here the vertices are labeled by $v_1,\dots, v_5$. 
The numbers of legs for the vertices are 
$e_1=3$, $e_2=4$, $e_3=5$, $e_4=3$, $e_5=5$. 
The number of the internal edges equal $I=4$. The graph has 
twelve external edges, and eq.(\ref{dim}) holds because 
$12=3+4+5+3+5-2\cdot 4$.}
 \label{fig:3}
\end{figure}

The open string Hilbert space $\cH$ is a $\Z$-graded vector space. 
The conformal field theory technique gives us a basis system 
$\{\eb_i\}$ of open string states (in terms of the oscillators 
in the mode expansions), 
where the grading of these basis 
is related to the ghost number of string states 
(\cite{Z1,N}). 
For each state $\eb_i$, consider a field $\phi^i$ 
(in the sense of field theory) whose degree is minus the degree of
$\eb_i$ so that the degree of $\Phi:=\eb_i\phi^i$ is set to be zero. 
$\Phi$ is called a {\it string field}. 
\footnote{For the relation to the usual notations in physics 
see also \cite{Ka} subsection 5.2. }
Moreover we have a degree one coboundary operator 
$Q:\cH\to\cH$ and a degree minus one antisymmetric bilinear form 
$\omega(\ ,\ ):\cH\otimes\cH\to\C$ that are also defined canonically 
on the conformal background. 
$Q$ and $\omega$ are called the BRST-operator \cite{KO} and the
BPZ-inner product \cite{BPZ}, respectively. 
They in fact define a degree-zero graded-symmetric bilinear form 
$\V_2:=\omega(\1\otimes Q):\cH\otimes\cH\to\C$, which defines the 
kinetic term (quadratic term with respect to field $\{\phi\}$) 
of the action of a classical open string field theory. 
The action is of the following form, 
\begin{equation*}
 S(\Phi)=\half\omega(\Phi,Q\Phi)+\sum_{k\ge 3}\ov{k}\V_k(\Phi,\dots,\Phi)\ ,
\end{equation*}
where $\V_k:\cH^{\otimes k}\to\C$ is a degree zero cyclic multilinear map. 
We call $\{\V_k\}_{k\ge 3}$ the {\it vertex maps}. 
The term `cyclic' indicates that $\V_k$ satisfies 
$\V_k(\eb_{i_1},\dots,\eb_{i_k})
=(-1)^{\eb_{i_1}(\eb_{i_2}+\cdots +\eb_{i_k})}
\V_k(\eb_{i_2},\dots,\eb_{i_k},\eb_{i_1})$ for any $\eb_i\in\cH$. 
It holds for $k\ge 2$, where the case $k=2$ is equivalent to the fact 
that $\V_2$ is graded-symmetric stated above. 
All the (multi-)linear maps introduced here are 
extended naturally to the polynomials of fields $\phi^i$. 
Thus, the action $S(\Phi)$ is 
a degree zero polynomial function that has the cyclicity. 
To construct a string field theory is then 
to construct vertex maps $\{\V_k\}_{k\ge 3}$
satisfying certain conditions explained below. 
In order to do it, 
some conformal field theory technique provides us with 
the following set-up \cite{V,AGMV,LPP,Z1}. 
Let us consider a disk $D$ with cyclic ordered $n$ 
holomorphic half disks on the boundary $S^1$ for $n\ge 3$. 
Namely, we have $n$ holomorphic maps $f_i$, $i=1,\dots, n$, 
from a half disk $\{z\in\C | \Im(z)\ge 0, |z|\le 1\}$ 
in an upper half plane $H_+=\{z\in\C | \Im(z)\ge 0\}$ 
to the disk $D$ which are injective 
and map the boundary $\Im(z)=0$ of the half
disk to intervals on the boundary of the disk $D$ 
with preserving the orientations. 
Thus, $f_i$ maps the origin $o$ of the half disk to a point 
(puncture) on the boundary of $D$, 
and these $n$ holomorphic maps are defined so that 
$f_1(o),\dots,f_n(o)$ are counterclockwise cyclic ordered and 
the images of the half disks by any two holomorphic maps 
do not overlap with each other. 
In particular, in order to fix the $SL(2,\R)$ automorphisms 
of the disk $D$, we fix three points 
$f_1(o)=:0$, $f_{n-1}(o)=:1$ and $f_n(o)=:\infty$. 
We denote the space of such disks with cyclic ordered $n$ 
holomorphic half disks by $\tcM_n$. 
It forms an infinite dimensional space. 
For a disk $\Sigma_n\in\tcM_n$, the image of the arc 
defined by $|z|=1$, $|\Im(z)|\ge 0$ by each holomorphic map 
$f_i$ is regarded as an open string. 
The disk $\Sigma_n$ thus describes the interaction of such 
$n$ open strings, as in Figure \ref{fig:3} (a), 
with the initial condition of each open string 
being specified by the image of the origin $f_i(o)$. 
An open string state space $\cH$ is associated to each origin $o$ 
of the half disk. 
In particular, since $Q:\cH\to\cH$ is a differential, 
$(\cH,Q)$ forms a complex called the {\em BRST-complex}. 
The kernel (resp. cokernel) of $Q$ is called the {\em on-shell} 
(resp. {\em off-shell}) state space, 
and the cohomology $H(\cH)$ with respect to $Q$, 
the {\em BRST cohomology}, is called the {\em physical state} space.
\footnote{
Here we assume for simplicity 
that the basis $\eb_i$ are taken so that the subbasis of $\{\eb_i\}$ 
can span the on-shell state space or the physical state space. 
} 
For each $\Sigma_n\in\tcM_n$, 
the corresponding {\em correlation function (expectation value)} 
of conformal field theory gives a 
map $\Sigma_n: \cH^{\otimes n}\to\C$ %is given 
(as above we denote the %corresponding 
map also by $\Sigma_n$). 
Moreover, one can consider the tangent space $T\tcM_n$ and 
the space of differential $k$-forms $\tOmega_{diff}^k(\tcM_n)$ on $\tcM_n$ 
for each $k\ge 0$ \cite{V,AGMV, Z1}. 
In particular, associated to the infinitesimal deformations of $\Sigma_n$, 
one can define a map 
$\tOmega_n^k: \cH^{\otimes n}\to\tOmega_{diff}^k(\tcM_n)$ for each $k$. 

Let $\cM_n$, $n\ge 3$, be a suitable compactification of 
the moduli space of disks with $n$ punctures on the boundary. 
The dimension of $\cM_n$ is $(n-3)$. 
There is a projection
$\pi: \tcM_n\to\cM_n$ obtained by forgetting 
the holomorphic maps $f_i$, $i=1,\dots,n$,  
except the image of the origin $f_i(o)$. 
Namely, for $\Sigma_n\in\tcM_n$, 
$\pi(\Sigma_n)$ is the disk with $n$ punctures specified by 
$f_1(o),\dots,f_n(o)$. 

Let us consider a map (section) $\sigma: \cM_n\to\tcM_n$ 
such that $\pi\circ\sigma$ is identity. 
When restricting every $\eb_i\in\cH$ to on-shell, 
the following map 
\begin{equation}
 \ti\V_n(\eb_{i_1},\dots,\eb_{i_n})
 :=\int_{\cM_n}\sigma^*
 \l(\tOmega_n^{n-3}(\eb_{i_1},\dots,\eb_{i_n})\r)\ \in\C 
 \label{corr}
\end{equation}
is in fact independent of the choice of $\sigma$, 
where the degree of the differential form $(n-3)$ 
is the dimension of $\cM_n$. 
These are nothing but 
the tree (on-shell) {\em (scattering) amplitudes} of open strings. 
Since the $n$ insertions (=punctures) 
are on the boundary of the disk, 
$\cA_n$ is a cyclic map. 
Moreover it is known that the on-shell correlation function vanishes 
if one of the external states is $Q$-exact. Thus 
the collection of open string scattering amplitudes 
can be defined on the physical state space $H(\cH)$. 

In this situation, the vertex maps $\{\V_n\}_{n\ge 3}$ 
should be constructed 
so that the perturbation theory reproduces 
the open string %correlation functions 
scattering amplitudes (\ref{corr}). 
In perturbation theory the on-shell scattering amplitudes 
are calculated by Feynman graphs. 
The usual construction of string field theory is then 
to decompose each $\cM_n$ into cells so that 
cells correspond one to one with Feynman graphs. 
The vertex map $\V_n$ is determined by the pair $(\cM^0_n,\sigma)$, 
where $\cM^0_n\in\cM_n$ is a cell of $\cM_n$ 
and $\sigma:\cM^0_n\to\tcM_n$ is a map 
such that $\pi\sigma$ is equal to the identity. 
We give such a pair $(\cM^0_n,\sigma)$ 
so that $\sigma\cM^0_n\subset\tcM_n$ has an $SL(2,\R)$ automorphism 
on the disk $D$ which transforms 
$f_i$ to $f_{i+1}$ for $1\le i\le n-1$ and $f_n$ to $f_1$. 
Then $\V_n$ is given as 
\begin{equation}
 \V_n(\Phi,\dots,\Phi)
 :=
 \int_{\cM_n^0}\sigma^*
 \l(\tOmega_n^{n-3}(\Phi,\dots,\Phi)\r)\ .
 \label{vertex}
\end{equation}
By construction, 
$\V_n(\eb_{i_1},\dots,\eb_{i_n})$ is cyclic.
\footnote{In section 3 in \cite{Ka} the index for the differential 
form $k$ in $\tOmega^k_n$ is omitted. The correspondence of the
notation between \cite{Ka} and this paper is then given by 
$\Omega_n=\sigma^*\tOmega_n$. }

Let us consider an on-shell tree $n$-point %correlation function. 
open string scattering amplitude. 
The corresponding Feynman graphs are tree planar graphs, 
each of which consists of cyclic vertices, internal edges, and
external edges. 
Each internal edge has two distinct vertices. 
Each external edge, called a leaf, 
has a vertex at one end and the other end is free
(see Figure \ref{fig:3} (b)). 
Clearly, by ignoring the distinction between the root edge and the
leaves of a rooted planar tree and regarding the root edge also as a
leaf, one gets a planar tree graph. 
Thus, we have a natural surjection $\rh: G_{n-1}\to G^{cyc}_n$, 
where $G_{n-1}$ is the set of rooted planar $(n-1)$-trees and 
$G^{cyc}_n$ is the set of planar graphs with $n$ leaves. 
Let $G_n^{cyc,I}$ be the set of planar trees with 
$n$ leaves and $I$ internal edges. 
Each element $\Gamma_n^{cyc,I}\in G_n^{cyc,I}$ then has $I+1$
vertices. 
We assign $v_m$, $m=1, 2, \dots, I+1$ to the vertices and let $e_m$
be the number of incident (both internal and external) edges. 
The following identity then holds
\begin{equation}
 n+2I=\sum_{m=1}^{I+1}e_m\label{dim}\ . 
\end{equation} 
For $\Gamma_n^{cyc,I}$, one can consider a number 
$\ti\V_{\Gamma_n^{cyc,I}}(\eb_{i_1},\dots,\eb_{i_n})\in\C$. 
Essentially it is given 
by attaching $\V_{e_m}$ to each vertex $v_m$ and to each internal edge a 
so-called propagator (Definition \ref{defn:BVpropagator}, denoted by
$\V_L^+\in\cH\otimes\cH$) and  $\eb_{i_1},\dots,\eb_{i_n}$
to leaves cyclically (see Definition \ref{defn:cyclicplanar}). 
The tree $n$-point open string scattering amplitude (\ref{corr}) 
is then reproduced by 
\begin{equation}
 \ti\V(\eb_{i_1},\dots,\eb_{i_n})
 :=\sum_{\Gamma^{cyc}_n\in G^{cyc}_n}
 \ov{\sharp\mathrm{Aut}(\Gamma^{cyc}_n)}
\ti\V_{\Gamma^{cyc}_n}(\eb_{i_1},\dots,\eb_{i_n})\ ,
 \label{IntrotiV}
\end{equation}
where each $\eb_i$ is on-shell, 
and $\mathrm{Aut}(\Gamma^{cyc}_n)$ indicates 
the number of the automorphisms acting on $\Gamma^{cyc}_n$. 
The fraction $\ov{\sharp\mathrm{Aut}(\Gamma^{cyc}_n)}$ is called 
the {\em symmetric factor} of the Feynman graph. 
(We shall discuss these Feynman graphs in detail 
in subsection \ref{ssec:PI}. ) 

On the other hand, the propagator $\V_L^+$ is represented by 
an integral over $[0,\infty]$. 
\footnote{This corresponds to the length parameter of 
the open string evolution. } 
Namely, the propagator or the internal edge has modulus 
$\tau\in [0,\infty]$ 
and, in a Riemann surface picture, a strip with fixed width and length
$\tau$ is associated to it. 
Assume that $\{\V_k\}_{k\ge 3}$ are constructed so that 
the associated Riemann surfaces $\{\sigma:\cM^0_k\to\tcM_k\}_{k\ge 3}$ 
can be joined with the strip (propagator) by 
sewing Riemann surfaces. 
Then, each graph $\Gamma^{cyc,I}_n\in G^{cyc,I}_n$ is 
associated with a subspace of $\tcM_n$, 
which we denote by $\tcM_{\Gamma^{cyc,I}_n}\subset\tcM_n$. 
The important point is that 
the compatibility with respect to the sewing of Riemann surface 
is known \cite{LPP} (for classical open string theory, more explicitly
in \cite{SS}), which implies, for instance, 
\begin{equation*}
 \ti\V_{\Gamma_n^{cyc,I}}(\eb_{i_1},\dots,\eb_{i_n})
 =\int_{\tcM_{\Gamma^{cyc,I}_n}}\tOmega_n^{n-3}(\eb_{i_1},\dots,\eb_{i_n})\ .
\end{equation*} 
Here note that the dimension of $\tcM_{\Gamma^{cyc,I}_n}$ is actually 
$n-3$ and independent of $I$. The fact can be confirmed by eq.(\ref{dim}) as 
$(e_1-3)+\cdots+(e_{I+1}-3)+I=(k+2I-3(I+1))+I=k-3$.

Suppose that the vertex maps $\{\V_k\}_{k\ge 3}$ in eq.(\ref{vertex}) 
are constructed consistently up to $k=n-1$ and then 
concentrate on the $n$-point amplitude (\ref{IntrotiV}). 
The Feynman graph without propagator($I=0$) consists only of 
the vertex $\V_n$, which is not determined yet. 
For each $\Gamma^{cyc,I}_n$, $I>0$, we assume the projection 
of $\tcM_{\Gamma^{cyc,I}_n}$ gives an inclusion, 
\begin{equation}
 \pi(\tcM_{\Gamma^{cyc,I}_n})\subset\cM_n\ ,
 \label{pimoduli}
\end{equation}
and for any two distinct elements of 
$G^{cyc}_n$ the images never have a common subspace 
except their boundaries. Let us denote 
\begin{equation}
 \cM_k^I:=\bigcup_{\Gamma^{cyc,I}_n\in
G^{cyc,I}_n}\pi(\tcM_{\Gamma^{cyc,I}_n})\ .
 \label{cMGamma}
\end{equation}
Thus, $\cM^0_n$ is determined as 
\begin{equation}
 \cM_n=\cM_n^0\cup\cM_n^1\cup\cM_n^2\cup\cdots\cup\cM_n^{n-3}\ ,
 \label{decomp}
\end{equation}
where the common subspace $\cM^I\cap\cM^{I'}$, $I\ne I'$, 
has codimension greater than one. 
Furthermore, define $\sigma: \cM^0_n\to\tcM_n$ so that 
$\sigma(\cM^0_n)$ and 
$\bigcup_{I\ge 1}\bigcup_{\Gamma^{cyc,I}_n\in G^{cyc,I}_n}
(\tcM_{\Gamma^{cyc,I}_n})$ 
form a {\it continuous} section of the bundle $\tcM_n\to\cM_n$. 
\begin{figure}[h]
 \hspace*{1.0cm}
 \includegraphics{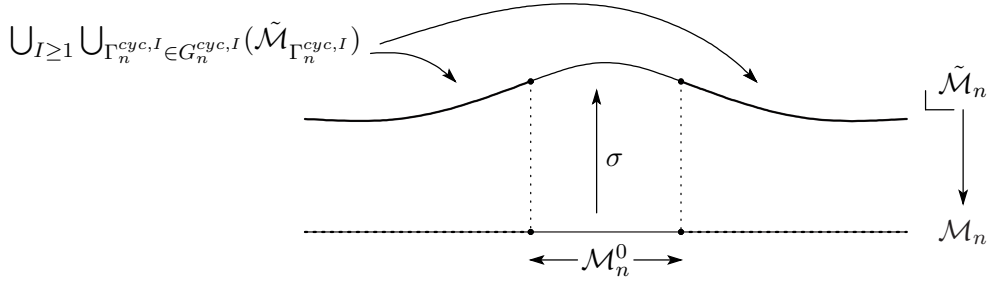}
 \caption{The determination of the pair $(\cM^0_n,\sigma)$. }
 \label{fig:sigma}
\end{figure}
Consequently one obtains $\V_n$ by eq.(\ref{vertex}). 
One can see that the action is obtained by repeating this procedure. 
By construction it is clear that the action reproduces the 
tree open string scattering amplitudes by perturbation theory. 

\def \cg{\hspace*{-0.05cm}{\lglraw\hspace*{-0.56cm}\circ}\hspace*{0.21cm}}

The action constructed as above actually satisfies 
the classical BV-master equation (\ref{Intromeq}). 
First, consider the infinitesimal variation of the decomposition 
of Riemann surfaces, or more precisely, take the boundary $\partial$ of 
eq.(\ref{decomp}). 
Since in eq.(\ref{pimoduli}) we assumed $\pi$ is an inclusion, 
then $\pi$ commute with $\partial$. 
If one takes the boundary of eq.(\ref{cMGamma}), the boundary operator 
$\partial$ acts on each $\tcM_{\Gamma^{cyc,I}_n}$ in the right hand side. 
Here $\tcM_{\Gamma^{cyc,I}_n}$ is a topological space 
$\sigma(\cM^0_{e_1})\times\cdots\times\sigma(\cM^0_{e_{I+1}})
\times [0,\infty]^I$ equipped with the information of the 
planar tree graph $\Gamma^{cyc,I}_n$. 
Then $\partial$ acts by the Leibniz rule on a vertex space $\cM^0_{e_m}$ or 
a propagator $[0,\infty]$. 
Here the boundary of the propagator is $\{0\}-\{\infty\}$. 
It is natural to require that in these construction 
the sum of all the contributions $\{\infty\}$ 
corresponds to the boundary of $\cM_n$. 
Then, acting by $\partial$ on eq.(\ref{decomp}) yields 
\begin{equation}
 0=\partial(\cM^0_n)+\sum_{\substack{k_1+k_2=n+2\\k_1,k_2\ge 3}}
 \half\l(\begin{split}
    &\partial(\cM^0_{k_1})\-(\cM^0_{k_2})\\
   +&(\cM^0_{k_1})\-\partial(\cM^0_{k_2})\\
   +&(\cM^0_{k_1})\cg (\cM^0_{k_2})
 \end{split}\r)
 +\sum_{\substack{k_1+k_2+k_3=n+4\\k_1,k_2,k_3\ge 3}}\l(\cdots\r)+\cdots\ 
 \label{recur1}  
\end{equation}
for $n\ge 3$. We should explain some of the notations used above. 
First, we identify the image of the composition of two maps 
$\sigma: \cM^0_{k_i}\to\tcM_{k_i}$ and 
$\pi:\tcM_{k_i}\subset\tcM_{n}\to\cM_n$ with $\cM^0_{k_i}$ itself and 
wrote $\cM^0_{k_i}$. 
Note that each $\cM^0_{k_i}$ is associated with a vertex. 
We then denoted by $\-$ a topological space $[0,\infty]$ with 
the operation of connecting two vertices with the propagator. 
Alternatively, `$\cg$' indicates the operation of grafting 
two vertices with `$\{0\}$', the contracted propagator.  
%For $n=3, 4$ we get 
%\begin{equation}
% n=3\ :\ 0=\partial\cM^0_3\ ,
% \qquad n=4\ :\  0=\partial\cM^0_4+\half\cM^0_3\partial(\-)\cM^0_3\ .
%\end{equation}
%The first equation ($n=3$) means $\cM^0_3=\cM_3$ has no moduli (a point). 
The equation (\ref{recur1}) is, in fact, equivalent to 
\begin{equation}
 0=\partial(\cM^0_n)+\sum_{\substack{k_1+k_2=n+2\\k_1,k_2\ge 3}}
 \half(\cM^0_{k_1})\cg (\cM^0_{k_2})\ .
 \label{sfeq}
\end{equation}
The right hand side of the identity above 
is the sum of the first term and one of the second term in 
the right hand side of the identity(\ref{recur1}). 
The equivalence holds because 
the other parts of eq.(\ref{recur1}) cancel by induction. 
For example, $\partial(\cM^0_{k_1})\-(\cM^0_{k_2})$ in the second term 
cancels one of the third term $(\cdots)$ of the form 
$\sum_{\substack{k+l=k_1+2\\k,l\ge 3}}\l(\half(\cM^0_k)\cg (\cM^0_l)\r)
\-(\cM^0_{k_2})$. 
The recursion equation (\ref{sfeq}) is called 
the {\em string factorization equation} \cite{SZ}. 

The string factorization equation (\ref{sfeq}) is actually equivalent to 
the BV-master equation. 
In fact, this identity (\ref{sfeq}) 
is an identity between $(n-4)$-dimensional moduli space and
graphically an identity between planar tree graphs with $n$ leaves. 
Thus, let us integrate 
$\sigma^*\l(\tOmega^{n-4}_n(\Phi,\dots,\Phi)\r)$ 
over the identity (\ref{sfeq}). 
The conformal field theory technique leads to the following result 
({\bf cf.} \cite{Z1}): 
\begin{equation}
 0=\delta_1\l(\ov{n}\V_n(\Phi,\dots,\Phi)\r)
 +\half\sum_{\substack{k_1+k_2=n+2\\k_1,k_2\ge 3}}
 \l(\ov{k}\V_k(\Phi,\dots,\Phi), \ov{l}\V_l(\Phi,\dots,\Phi)\r)\ . 
 \label{sfeq2}
\end{equation}
In the equation above, 
$(\ ,\ ):=\flpartial{\phi^i}\omega^{ij}\frpartial{\phi^j}$ and 
$\omega^{ij}$ is the inverse of $\omega_{ij}:=\omega(\eb_i,\eb_j)$. 
This in fact defines an odd Poisson bracket and is called 
{\em the BV-bracket}. Also, $\delta_1$ is defined by 
$\delta_1:= (\ ,\half\omega(\Phi,Q\Phi))$. 
It satisfies $(\delta_1)^2=0$ corresponding to 
the nilpotency $Q^2=0$ (or $\partial^2=0$). 
Summing up eq.(\ref{sfeq2}) for $n\ge 3$ and multiplying 
by two then lead to the {\it classical BV-master equation}
\begin{equation}
 (S(\Phi), S(\Phi))=0\ .\label{Intromeq}
\end{equation}
To summarize, to construct a string field theory is 
to construct $\{\sigma(\cM^0_k)\}_{k\ge 3}$ so that they are 
compatible with the decomposition of the moduli spaces. 
The construction of 
$\{\sigma(\cM^0_k)\}_{k\ge 3}$ is independent of the conformal background 
we choose. Whereas, $\tOmega^\bullet_\bullet$ is determined canonically 
by the conformal background and taking a representation of 
$\{\sigma(\cM^0_k)\}_{k\ge 3}$ by $\{\tOmega^{k-3}_k\}_{k\ge 3}$ 
produces a string field theory action on the conformal background. 
Mathematically, $\{\cM^0_n\}_{n\ge 3}$ forms an operad 
and, by taking its representation, one obtains an algebra $\cH$ 
over the operad, where $\cH$ is called an operad algebra \cite{MSS}.

As seen in the next subsection, 
an action which has cyclic vertices and 
satisfies the classical BV-master equation as above has an 
$A_\infty$-structure. 
The $A_\infty$-algebra in addition possesses an odd symplectic inner 
product and cyclicity. 
Such an algebra is called a cyclic $A_\infty$-algebra 
(see Definition \ref{defn:cycAinfty}). 
The appearance of an $A_\infty$-structure can already be seen 
from eq.(\ref{sfeq2}). This identity is in fact 
a different but an equivalent expression 
of the $A_\infty$-condition (\ref{IntroAinfty}) 
under the situation cyclic symmetry exists. 
The structure of an $A_\infty$-space can also be found 
in an explicit construction of the classical open string field theory 
in \cite{N,Ka}, 
where $\cM^0_{n+1}$ is just the associahedra $K_n$
and string factorization equation (\ref{sfeq}) is just the cyclic 
version of eq.(\ref{hedrarec}) \cite{Ka2}. 
The corresponding operad is called the {\it $A_\infty$-operad} 
\cite{MSS}. 
Similar stories hold for other classical string field theories. 
The underlying operad structure in classical closed string field
theory is the $L_\infty$-operad \cite{KSV,Sta1993}. 
For classical open-closed string case, see \cite{ocha} 
and also a related earlier work \cite{Vo}.

The minimal model theorem appears naturally also in string theory. 
In \cite{KoSo} for any $A_\infty$-algebra an explicit construction of 
the minimal model is given. The construction is just given by 
Feynman graphs. 
For classical open string field theories, these are just the Feynman graphs 
appearing above. 
This implies that the collection of the scattering amplitudes of 
open string theory forms a minimal cyclic $A_\infty$-algebra. 
This statement is essentially already known. 
In \cite{WZ} it is shown that the tree closed string theory 
has the structure 
of the $L_\infty$-algebra (and it is extended to the quantum case in 
\cite{V}). 
Thus, the minimal model theorem implies 
on a fixed conformal background all classical open string field theories 
are $A_\infty$-quasi-isomorphic to each other \cite{Ka}. 
Namely, the difference in the choice of the decomposition of moduli spaces 
leads to homotopy equivalence of $A_\infty$-algebras, 
and the minimal model is obtained 
by homotopic deformation $\cM^0_k\to \cM_k$.
\footnote{This also implies that moduli spaces of open string
correlation functions $\{\cM_k\}$, obtained by a suitable 
compactification, has the structure of an $A_\infty$-space. 
Though in a slightly different context, 
an $A_\infty$-space structure in open string theory is discussed 
in \cite{CoSc}. }
Moreover one can show that these are not only quasi-isomorphic but 
$A_\infty$-isomorphic (Theorem \ref{thm:main}) 
due to Theorem \ref{thm:cycMandC}.

Physically, string field theories have been investigated as a candidate 
for string theory which describes nonperturbative effects. 
This purpose requires the off-shell extension of string theory as above. 
A typical off-shell physics phenomenon is tachyon condensation \cite{tac}. 
Recently, string field theory has been applied in such a direction 
successfully \cite{SZn,Oh} (see also \cite{KoSa}). 
Though we assumed the existence of $\sigma$ which required many 
consistency conditions as above, 
actually there exists many Lorentz-covariant string field theories
(SFTs)
\footnote{String field theories of the type explained here are called 
(Lorentz) covariant string field theories in contrast to 
light cone type string field theories developed earlier. }
; the covariant open or closed SFT with light cone 
type-like vertices(HIKKO's SFT) \cite{HIKKO}, 
a very simple open SFT which consists of only a three-point vertex 
(Witten's open SFT or cubic SFT) \cite{W1}, 
nonpolynomial classical closed SFT constructed by 
`restricted polyhedron' \cite{KKS,KS}, and so on. 
Witten's SFT is treated in the context of 
BV-formalism\cite{Thorn,Bochi}(see \cite{Th}). 
HIKKO's closed SFT is also extended to quantum SFT by 
employing the quantum BV-master equation \cite{H}. 
The quantum master equation is moreover applied to construct quantum closed 
SFT with symmetric vertices by Zwiebach \cite{Z1}. 
Though this theory has infinite sort of 
vertices of higher punctures and higher genus, it has 
a very beautiful algebraic structure. For instance for the classical part, 
the set of the tree vertices has the structure of 
an $L_\infty$-algebra. 
Open-closed SFT is also considered in this direction \cite{Z2}. 
HIKKO-type open-closed SFT is given in \cite{KT,AKT}. 
Recently a one parameter family of classical open string field theories, 
which possess $A_\infty$-structures, has been 
constructed explicitly \cite{N,Ka,Ka2} 
by deforming the Witten's cubic SFT \cite{W1}.

 \subsection{Dual description; 
formal noncommutative supermanifolds}
\label{ssec:Introdual}

For $A_\infty$-algebras, 
we use mainly three descriptions; the coalgebra language, 
its dual language, and superfield description. 
Coalgebras are used to define $A_\infty$-algebras precisely and simply. 
On the other hand, the dual description is geometric and intuitive as 
explained below. 
The superfield description, used in the previous subsection, is 
their mixed version. It is directly equivalent to the dual 
ones but the superfield description uses the notation used in coalgebra. 
This description is convenient to simplify indices.   
The operad structure is implicit in various arguments in this paper, 
but we shall not indicate it explicitly.

As in the previous subsection, given an $A_\infty$-algebra $(\cH,\m)$, 
denote by $\{\eb_i\}$ a basis of $\cH$ and $\{\phi^i\}$ the dual 
coordinates. Reflecting the non(co)commutativity of $\cH$, 
the dual fields are treated as noncommutative as explained in 
section \ref{sec:dual}. 
We call $\Phi=\eb_i\phi^i$ the superfield, which is the string field 
in string field theory. 
Let us describe the $A_\infty$-structure in coordinates as 
\begin{equation*}
 m_k(\eb_{i_1},\dots,\eb_{i_k})=\eb_j c^j_{i_1\cdots i_k}\ .
\end{equation*}
For the collection $c^j_{i_1\cdots i_k}\in\C$ for $k\ge 1$, one can define 
the following degree one vector field, 
called the {\em homological vector field}, 
on a formal noncommutative supermanifold
\begin{equation}
  \delta=\sum_{k=1}^\infty\flpartial{\phi^j}c^j_{i_1\cdots i_k}
 \phi^{i_k}\cdots\phi^{i_1}\ .
 \label{IntroAinftyovf}
\end{equation} 
We often use the Einstein convention 
of summing over repeated indices as above. 
Note that the $A_\infty$-condition is then rewritten as $(\delta)^2=0$. 
We call this $\delta$ an $A_\infty$-odd vector field. 

On the other hand, let us consider a degree zero cyclic function 
of the form
\begin{equation*}
 S=\half\V_{i_1i_2}\phi^{i_2}\phi^{i_1}
 +\sum_{k\ge 3}\ov{k}\V_{i_1\cdots i_k}\phi^{i_k}\cdots\phi^{i_1}
\end{equation*}
where $\V_{i_1\cdots i_k}\in\C$ for $k\ge 2$. 
For a given odd nondegenerate constant Poisson bracket 
$(\ ,\ ):=\flpartial{\phi^i}\omega^{ij}\frpartial{\phi^j}$, 
the Hamiltonian vector field of the Hamiltonian $S$, 
\begin{equation*}
 \delta=(\ ,S)\ , 
\end{equation*}
is nilpotent iff $(S, S)=0$. 
This $\delta$ is nothing 
but an $A_\infty$-odd vector field (\ref{IntroAinftyovf}), 
where the $A_\infty$-structure is written as 
\begin{equation*}
 c^j_{i_1\cdots i_k}=(-1)^{\eb_m}\omega^{jm}\V_{m i_1\cdots i_k}\ .
\end{equation*}
Although one can obtain an $A_\infty$-algebra in such a way, 
it has an additional structure; the $A_\infty$-structure is cyclic 
with respect to the odd Poisson structure. 
Thus, we denote the corresponding algebra by 
$(\cH,\omega,\m)$ or $(\cH,\omega,S)$ and call it a cyclic $A_\infty$-algebra
(see Definition \ref{defn:cycAinfty}). 
Moreover one may notice that the condition $(S, S)=0$ is nothing but 
the classical BV-master equation (\ref{Intromeq}) in the BV-formalism. 
Then one can see that any cyclic field theory equipped with 
a classical BV-structure, including classical open string field theories 
in the previous subsection, has a cyclic $A_\infty$-structure 
(Theorem \ref{thm:ft}).

 \subsection{Noncommutativity, open strings, and D-branes}
\label{ssec:Introphys}

In the explanation above, we set $\{\phi^i\}$ to be 
formally-{\em noncommutative} coordinates. 
Mathematically, it is because, otherwise some informations of 
$A_\infty$-algebras are lost in the dual supermanifold description. 
In the case of field theories equipped with classical BV-structure 
discussed in section \ref{sec:BV}, 
we identify the dual coordinates with the fields of field theory. 
The noncommutativity of fields then implies physically 
the presence of Chan-Paton factor in open string theory. 
Namely, the non(co)commutativity of $\cH$ allows the freedom of the 
choice of the Chan-Paton factor, where 
the fields $\{\phi^i\}$ are described typically 
by $N\times N$ matrices with entries $\C$.
In other words, we have a representation of the theory 
in terms of $N\times N$ matrices. 
\footnote{More precisely when we define the cyclic structure on the 
action we treat the real part and imaginary part of $\C$ separately 
(see subsection \ref{ssec:BV}).}
Note that in the theory of open strings there 
exist D-branes and open strings must end on the D-branes. 
The size of the matrices $N$ then means there exist $N$ (parallel) D-branes. 
The typical gauge structure group is $U(N)$, though the structure group 
depends on the (super)symmetry which the theory has. 
Note that, 
although we can represent classical open string field theory 
with $N\times N$ matrices for any $N$, by definition 
the vertices $\V_{i_1\cdots i_k}\in\C$ are independent of $N$. 
However, just as 
we fix a representation by $N\times N$ matrices, 
the theory reduces to the one equipped with cyclic $L_\infty$-structure. 

For instance, let us represent the noncommutative fields 
by $N\times N$ matrices as
\begin{equation*}
 \phi^i=
\bp \phi^i_{11} & \cdots & \phi^i_{1N}\\
    \vdots      & \ddots & \vdots     \\
    \phi^i_{N1} & \cdots & \phi^i_{NN}
\ep\ .
\end{equation*}
The noncommutative product of $\phi^i$'s are the usual multiplication of 
the matrices. Then the $A_\infty$-odd vector field, 
as in eq.(\ref{IntroAinftyovf}), 
is written in terms of 
the component fields $\phi^i_{pq}, 1\le p,q\le N$ which are 
graded {\em commutative}. 
Correspondingly, the coefficients $c^j_{i_1\cdots i_k}\in\C$ 
of the $A_\infty$-odd vector field 
are graded-symmetrized with respect to the indices $i_1\cdots i_k$ 
and the results turns out to define 
an $L_\infty$-structure (see \cite{LM} for $L_\infty$-algebras 
from symmetrizations of $A_\infty$-algebras without passing through 
the dual supermanifold description). 
Another choice of the structure groups leads to another 
$L_\infty$-algebra as the results of the graded symmetrizations 
of the component fields. 
In particular, if the size of the matrices are one ($N=1$), 
the coefficients $c^j_{i_1\cdots i_k}\in\C$ are completely
symmetrized as the dual supermanifold description of a result 
in \cite{LM}.

Namely, when one fixes a structure group, 
one loses a part of the informations which 
the open string theory has. 
A more familiar example is a gauge theory. 
The action, before being 
treated in the BV-formalism, is of the form
\begin{equation*}
 S(A)=\int F_{\mu\nu}F^{\mu\nu}\ ,\qquad 
F_{\mu\nu}=\partial_\mu A_\nu-\partial_\nu A_\mu+ [A_\mu,A_\nu]\ .
\end{equation*}
If the structure group is $U(N)$, each $A_\mu$ is an antiHermitian matrix. 
However in the case of 
$U(1)$ gauge theory the commutator $[A_\mu,A_\nu]$ vanishes. 
Namely, the structure constants except for the kinetic term 
(quadratic term of the action) are lost. 
For this reason it is reasonable to set the fields fully noncommutative. 
Then we discuss universal structures of open strings independent of 
the choice of Chan-Paton factor. 

The statements above imply 
that many properties which hold for $A_\infty$-algebras 
do also hold for $L_\infty$-algebras. 
That is, at least as far as homotopy algebraic properties are concerned, 
classical closed string theory can 
be understood from that of open string theory.
For this reason we shall discuss only on $A_\infty$-side in this paper.

 \subsection{Formal noncommutative symplectic supergeometry}
\label{ssec:Introsym}

In order to discuss the algebraic properties 
of cyclic $A_\infty$-algebras on a formal noncommutative supermanifold, 
we need some notions of noncommutative symplectic supergeometry, 
where a symplectic structure plays the role of a nondegenerate 
inner product defining the cyclicity of an $A_\infty$-algebra 
(see subsection \ref{ssec:cycAinfty}). 
Such a notion has appeared in \cite{Ko3,Ko4}, where  
a constant symplectic structure is introduced. 
We shall extend it to a nonconstant one in the way inspired 
from the physics of open strings, 
and examine various mathematical properties of them 
such as the Darboux theorem (Theorem \ref{thm:Darboux}) 
in section \ref{sec:sym}. 
Note that 
another nonconstant extension is discussed in \cite{Gi} based on Connes's 
noncommutative differential geometry \cite{connes}. 
Also, a different nonconstant extension of the inner product, called 
the homotopy inner product, is proposed in \cite{tradler}.

When one considers a Poisson algebra 
on a formal noncommutative supermanifold, 
one first needs functions on it. 
We define them so that they can describe linear combinations of 
open string disk correlation functions, 
which are cyclic with respect 
to the open string insertions (punctures) 
on the boundary $S^1$ of the disk. 
Pictorially, such a function is displayed as 
\begin{equation}
 a_{i_1\cdots i_n}\phi^{i_n}\cdots\phi^{i_1}
 =\begin{minipage}[c]{50mm}{\input{A1.tex}}\end{minipage}\ .
% =\begin{minipage}[c]{50mm}{\includegraphics{A1.eps}}\end{minipage}\ .
 \label{S1}
\end{equation}
In order to translate such cyclic objects to purely algebraic terms, one 
needs to cut the boundary of the disk $S^1$ as above. 
The cyclicity is then encoded in the coefficient, that is, 
$a_{i_1\cdots i_n}\in\C$ 
in the left hand side of (\ref{S1}) is graded symmetric 
with respect to the cyclic permutations of $i_1\cdots i_n$. 
When one considers a constant symplectic structure on a formal 
noncommutative supermanifold, 
the corresponding constant odd Poisson bracket is naturally defined so that 
the bracket of two open string disks becomes an open string disk. 
It is then natural to write the odd Poisson bracket as 
the following double lines
\begin{equation}
 (A, B)=\begin{minipage}[c]{50mm}{\input{AB.tex}}\end{minipage}\ .
% (A, B)=\begin{minipage}[c]{50mm}{\includegraphics{AB.eps}}\end{minipage}\ .
 \label{mini1}
\end{equation}
The choice of the place of the cut fixes the ambiguity of the sign
$\pm$ for $(A, B)$. 
The double line notation admits a natural extension to a nonconstant 
odd Poisson structure as follows
\begin{equation}
 (A, B)=\begin{minipage}[c]{50mm}{\input{AB2.tex}}\end{minipage}\ ,
% (A, B)=\begin{minipage}[c]{50mm}{\includegraphics{AB2.eps}}\end{minipage}\ ,
 \label{mini2}
\end{equation}
where $I$ denotes a multiindex and so $\phi^I=\phi^{i_k}\cdots\phi^{i_1}$ if 
$I=[i_k\cdots i_1]$. 
The corresponding equation is 
\begin{equation*}
 (A, B)=\sum_{ij,IJ}\pm\omega^{ij}_{JI}
 \l(\flpart{A}{\phi^i}\phi^I\frpart{B}{\phi^j}\phi^J\r)_c\ , 
\end{equation*}
(see eq.(\ref{poibra}) in Definition \ref{defn:covpoisson}), 
where $_c$ denotes the cyclic symmetrization and 
$\omega^{ij}_{JI}\in\C$ has an appropriate constraint so that 
the bracket satisfies $(B, A)=-(-1)^{AB}(A, B)$ and so on. 
One can define a notion of differential forms on 
formal noncommutative supermanifolds and the class of the odd Poisson 
brackets, which satisfy the Jacobi identity, can naturally be induced 
from closed two-forms (symplectic forms) 
on formal noncommutative supermanifolds (see section \ref{sec:sym}). 

In string theory, 
the nonconstant symplectic structure here is relevant to 
background independent string field theory (recently preferably 
called a boundary string field theory) \cite{W3} (see also \cite{HZ,Ka}). 
Consequently, our definition as above seems to be natural also mathematically.

 \subsection{Plan of this paper}
\label{ssec:contents}

Section \ref{sec:Ainfty} is devoted mostly to fixing our conventions 
for $A_\infty$-algebras. 
The precise definition of cyclic $A_\infty$-algebras is also 
included. 
In subsection \ref{ssec:coalg}, we recall the notion of coalgebras. 
$A_\infty$-algebras are then defined in terms of coalgebras 
(the bar construction) 
in subsection \ref{ssec:Ainfty}. 
The cyclic $A_\infty$-algebras are 
presented in subsection \ref{ssec:cycAinfty}. 
Some basic facts around Maurer-Cartan equations for $A_\infty$-algebras 
are mentioned briefly in subsection \ref{ssec:MCeq}.

In section \ref{sec:dual}, 
$A_\infty$-algebras are realized geometrically in the dual picture. 
In subsection \ref{ssec:dual1}, 
the dual is defined explicitly through an inner product, 
and its graphical realization is also presented. 
The dual picture is used in many papers, but there are 
few where the explicit relation is presented. 
All the tools presented in subsection \ref{ssec:Ainfty} are 
reinterpreted in terms of formal noncommutative supermanifolds 
in subsection \ref{ssec:dual2}. 
We shall then define `superfield' to simplify conventions in the 
dual picture, and mention some mixed description that 
interpolates between the coalgebra side and its dual side 
in subsection \ref{ssec:superfield}. 

In section \ref{sec:sym}, we shall explore local properties 
of symplectic structures on 
the formal noncommutative supermanifolds, which are relevant to the 
dual picture of cyclic $A_\infty$-algebras. 
The notion of formal noncommutative symplectic geometry appears for instance 
in \cite{Ko3,Ko4}. However, nonconstant symplectic structures are not 
explicitly written. 
We first define such covariant symplectic structures inspired by open 
strings. Namely, we consider cyclic formal functions. 
In subsection \ref{ssec:csympoi}, we shall observe some basic properties 
of constant symplectic structures, which serve as 
the starting point of more general cases. 
We then define covariant odd symplectic structures 
in subsection \ref{ssec:sympoi}, 
where we show a key lemma (Lemma \ref{lem:poincare}), the Poincar\'e 
lemma on formal noncommutative supermanifolds. 
Using the lemma, we examine the properties of 
symplectic diffeomorphisms in subsection \ref{ssec:diffeo}, 
and show the Darboux theorem on the formal noncommutative supermanifolds 
(Theorem \ref{thm:Darboux}) in subsection \ref{ssec:Darboux}. 
The study of the formal noncommutative symplectic supergeometry 
is directly related to the notion of cyclic $A_\infty$-algebras. 
We look back over cyclic $A_\infty$-algebras from these dual pictures in 
subsection \ref{ssec:cycAinftyre}.

The purpose of section \ref{sec:MandC} is to understand clearer 
the minimal model theorem 
\cite{kadei1}, one of the key theorem in homotopy algebras. 
For the construction of minimal models of $A_\infty$-structures, 
in particular on the homology of a differential graded algebra (dga), 
various versions of 
homological perturbation theory (HPT) have been developed, for instance,
by \cite{gugen,GS,gugen-lambe,GLS:chen,GLS,hueb-kadei}. 
Alternatively, as mentioned in \cite{Ko1}, 
there exists another stronger version 
of the minimal model theorem. 
It enables us to understand clearly the homotopical structures 
of homotopy algebras. 
We call it the decomposition theorem and prove it explicitly 
(Theorem \ref{thm:MandC}) in subsection \ref{ssec:MandC}. 
The decomposition theorem for cyclic $A_\infty$-algebras is then 
shown in subsection \ref{ssec:cycMandC}. 
The decomposition theorem guarantees the existence of an inverse 
$A_\infty$-quasi-isomorphism of an $A_\infty$-quasi-isomorphism 
(Theorem \ref{thm:inverse}) as stated in \cite{Ko1}. 
We shall explain it in subsection \ref{ssec:inverse}. 
Though the minimal model theorem follows from the decomposition theorem, 
the proof relies on inductive arguments 
and the explicit form of a minimal model is unclear. 
On the other hand, 
it is known that for any $A_\infty$-algebra a minimal model 
can be given explicitly by using some Feynman diagrams \cite{KoSo} 
(see also \cite{GLS,hueb-kadei,mer} and \cite{jh-jds} for $L_\infty$ case). 
We demonstrate in subsection \ref{ssec:minimal} 
that the Feynman diagrams arise naturally from the issue of finding 
the solutions of the Maurer-Cartan equation for an $A_\infty$-algebra 
\cite{Ka}. 
The cyclic $A_\infty$ version of the explicit minimal model is 
discussed in subsection \ref{ssec:cycminimal}, 
which is directly related to section \ref{sec:BV}. 

In section \ref{sec:BV}, these homotopy algebraic structures are 
applied to field theories equipped with classical BV-structures. 
The appearance of cyclic $A_\infty$-structures in 
field theories is explained in subsection \ref{ssec:BV}. 
To consider the perturbative expansion in the BV-formalism, 
we shall review briefly the notion of gauge fixing and propagators 
in our language and examine some properties of propagators 
in subsection \ref{ssec:gf}. 
Subsection \ref{ssec:PI} then shows that 
the tree on-shell correlation functions of a 
classical BV-field theory define just the minimal cyclic $A_\infty$-algebra 
defined in subsection \ref{ssec:cycminimal}
(Corollary \ref{cor:main}, cf.\cite{Ka}). 
Moreover, in subsection \ref{ssec:main}, the arguments in section 
\ref{sec:MandC} are applied to classical open string field theories, 
and it is shown that 
all classical string field theories on a fixed conformal background 
are cyclic $A_\infty$-isomorphic to each other
(Theorem \ref{thm:main}). 
Cyclic $A_\infty$-isomorphic means physically equivalent. 

Finally, in section \ref{sec:homotopy}, 
we shall come back to some basic problems in $A_\infty$-algebras. 
In subsection \ref{ssec:hom}, we shall define homotopy 
between $A_\infty$-morphisms and discuss various 
homotopy invariant algebraic structures of $A_\infty$-algebras. 
In subsection \ref{ssec:homgauge}, the notion of 
gauge equivalence and then the moduli space of $A_\infty$-algebras are 
defined. 
The properties of the moduli spaces are then examined.

% \vspace*{0.1cm}

Throughout this paper, we employ the dual picture, the formal noncommutative 
supermanifolds, in various places. 
To describe the dual of coalgebras by dual coordinates 
has some subtlety when the graded vector space is infinite dimensional. 
For instance, field theory is just such a case. However, 
since field theory is a theory of fields, 
it is well-defined 
as far as assuming that field theory itself is well-defined. 
Moreover, the dual language is used in this paper only for 
intuitive and geometric understanding. Hence 
almost all of the arguments on the dual can be rearranged 
in coalgebra language 
and hold even in the model where 
it is subtle to take a canonical basis system. 
One of the issues we do not discuss is some convergences. 
For instance $A_\infty$-morphisms or the solutions 
of the Maurer-Cartan equations, which are formally preserved 
under the $A_\infty$-morphisms, are defined by polynomials 
of {\em infinite} powers. 
Of course many of the arguments in this paper make sense 
as formal power series. For instance, in the application to field theories, 
each coefficient of the Maurer-Cartan equations defines 
an on-shell S-matrix element.  
However, it is also interesting to examine 
whether the solutions of Maurer-Cartan equations converge. 
This problem of convergence depends on the model 
equipped with an $A_\infty$-structure. 
Thus looking for some `good' models might be a good issue. 
Alternatively, one can also argue these on an appropriate subspace 
due to, for instance, the momentum conservation of the vertices 
in the case of field theory. 
Therefore, some well-defined solutions of the equations of motions 
may be obtained in the subspace.

 \clearpage
 
 \includegraphics{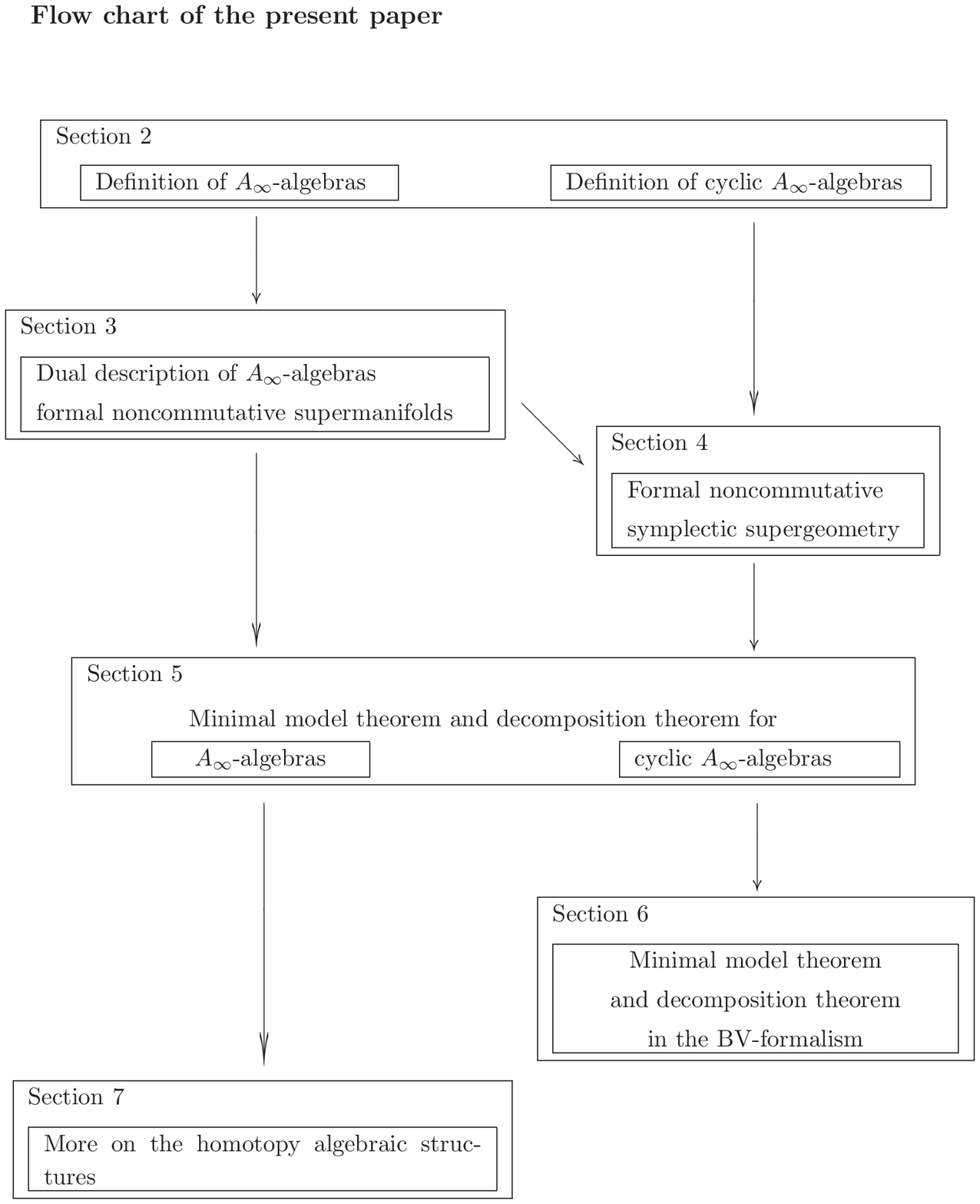}

 \clearpage

%section2

 \section{$A_\infty$-algebras}
\label{sec:Ainfty}

In this section, we shall summarize some basic facts about 
$A_\infty$-algebras ((strong) homotopy associative algebras). 
These facts are applicable in a similar way to 
$L_\infty$-algebras ((strong) homotopy Lie algebras).
We restrict our arguments to $A_\infty$-algebras over a field $k$ of 
characteristic zero. For more simplicity we set $k=\C$. 

$A_\infty$- (and $L_\infty$-) algebras are defined in different ways. 
One way is the operads. 
An $A_\infty$-algebra is obtained by an algebra 
over a non-symmetric dg operad (see \cite{MSS}). 
Another one is the bar construction and then $A_\infty$-algebras are defined 
as coalgebras with some additional structures. 
The bar construction is useful to define $A_\infty$-algebras in a 
simple manner and we take this definition in the present paper. 
For an intuitive or geometric realization of $A_\infty$-algebras, 
the dual picture of coalgebras is suitable. It is the subject of 
the next section. 

First we shall recall the notion of coalgebras in section \ref{ssec:coalg}. 
$A_\infty$-algebras and $A_\infty$-morphisms are then defined 
in terms of coalgebras in subsection \ref{ssec:Ainfty}. 
In subsection \ref{ssec:cycAinfty} we shall give a definition of 
$A_\infty$-algebras with cyclic symmetry. 
For an $A_\infty$-algebra, its Maurer-Cartan equation 
plays some important roles, which 
are explained briefly in subsection \ref{ssec:MCeq}.

 \subsection{Coalgebras, coderivations, and cohomomorphisms}
\label{ssec:coalg}

An element of an $A_\infty$-algebra belong to a 
$\Z$-graded vector space $\cH$. 
In the bar construction, the free tensor coalgebra of $\cH$ is treated 
as a coalgebra. We first provide the notions of coalgebras. 
\begin{defn}[Coalgebra, Coassociativity]
Let $C$ be a (generally infinite dimensional) graded vector space.  
When a {\it coproduct} 
$\triangle: C\lgraw C\otimes C$ is defined on $C$ and 
it is {\it coassociative}, {\it i.e.}
\[
 (\triangle\otimes{\bf 1})\triangle=({\bf 1}\otimes\triangle)\triangle
\]
then $C$ is called a {\it coalgebra}. 
\end{defn}
\begin{defn}[Coderivation]
A linear operator $\m: C\to C$ raising the degree of $C$ by one is 
called {\it coderivation} when   
\begin{equation*}
 \triangle \m= (\m\otimes{\bf 1})\triangle+({\bf 1}\otimes \m)\triangle
\end{equation*}
is satisfied. Here, for $x, y\in C$, 
the sign is defined as 
$({\bf 1}\otimes \m)(x\otimes y)=(-1)^{x}(x\otimes\m(y))$ 
where the $x$ on $(-1)$ denotes the degree of 
$x$. 
\label{defn:coder}
\end{defn}
\begin{defn}[Cohomomorphism]
Given two coalgebras $C$ and $C'$, a {\it cohomomorphism 
(coalgebra homomorphism)} $\cF$ from $C$ to $C'$ is a 
map of degree zero satisfying the condition
\begin{equation}
 \tri\cF=(\cF\otimes\cF)\tri\ .
\end{equation}
\end{defn}
\begin{rem}
Coassociativity of $\tri$, the conditions of coderivations and 
cohomomorphisms imply 
that the following diagrams commute: 
\begin{equation*}
\begin{CD}
 C @>{\triangle}>> C\otimes C\\
 @V{\triangle}VV   @V{\triangle\otimes{\bf 1}}VV\\
 C\otimes C @>{{\bf 1}\otimes\triangle}>> C\otimes C\otimes C
\end{CD}\ ,\quad 
\begin{CD}
 C @>\m>> C\\
 @V{\triangle}VV   @V{\triangle}VV\\
 C\otimes C @>{{\bf 1}\otimes \m+\m\otimes{\bf 1}}>> C\otimes C
\end{CD}\ ,\quad
\begin{CD}
 C @>{\cF}>> C'\\
 @V{\tri}VV   @V{\tri}VV\\
 C\otimes C @>{\cF\otimes\cF}>> C'\otimes C'
\end{CD} \ .
\end{equation*}
If the orientation of these map are reversed and the coproduct is 
replaced by a product, then the coassociativity, the coderivation, and the 
cohomomorphism take place to associativity, a derivation, and 
a homomorphism of the corresponding algebra, respectively. 

Reversing the orientation of the maps corresponds to taking the dual of 
the coalgebra. The precise meaning of the dual in the present 
paper is given in subsection \ref{ssec:dual1}. 
 \label{rem:dual}
\end{rem}
Let $\cH$ be a $\Z$-graded vector space. Namely, 
$\cH=\oplus_{k\in\Z}\cH^k$ where $\cH^k$ is a vector space of degree
$k$. Consider the free tensor coalgebra of $\cH$ 
\begin{equation*}
 C(\cH)=\oplus_{n\geq 0}{\cH^{\otimes n}}
\end{equation*}
as a coalgebra. 
Note that $\cH^{\otimes 0}=\C$, which includes a counit $\1$.
\footnote{
One may or may not include the $\cH^{\otimes 0}$ term for 
defining an $A_\infty$-algebra. 
If includes, one can also define a weak $A_\infty$-algebra uniformly, 
so we use this convention. 
}

Then the coassociative coproduct $\tri:C(\cH)\to C(\cH)\otimes C(\cH)$ 
is uniquely determined. For ${o}_1\cdots {o}_n\in\cH^{\otimes n}$ it is 
given by 
\begin{equation}
 \tri({o}_1\cdots {o}_n)=\sum_{k=0}^{n}
  ( {o}_1\cdots {o}_k)\otimes( {o}_{k+1}\cdots {o}_n)\ ,
\end{equation}
where the term for $k=0$ is $\1\otimes( {o}_1\cdots {o}_n)$ and 
the term for $k=n$ is $( {o}_1\cdots {o}_n)\otimes\1$. 
The form of the coderivation corresponding to this coproduct is 
also given as follows. 
Let $\{m_k: \cH^{\otimes k}\to\cH\}_{k\ge 0}$ 
be a collection of multilinear maps of degree one, that is, 
for any $ {o}_1,\dots, {o}_n\in\cH$ which are homogemeous in degree and 
\begin{equation}
 \begin{array}{cccc}
  m_k\ :&\cH^{\otimes k}&\lgraw& \cH\\
               & {o}_1\otimes\cdots\otimes {o}_k&\mapsto & 
                              m_k( {o}_1,\dots, {o}_k)
 \end{array}\qquad \ ,
\end{equation}
the image $m_k( {o}_1,\dots, {o}_k)\in\cH$ is also homogeneous, 
where its degree is the sum of the degree of $ {o}_i$, $i=1,\dots,k$, 
plus one. Also, 
$m_0 : \C\to\cH$ is defined so that $m_0(\1)$ has degree one. 
The operation on $C(\cH)$ is given as 
 \begin{equation*}
  \m_k( {o}_1\cdots {o}_n)=\sum_{p=1}^{n-k}
 (-1)^{ {o}_1+\cdots+ {o}_{p-1}} {o}_1\cdots {o}_{p-1}
 m_k( {o}_p,\dots, {o}_{p+k-1})
  {o}_{p+k}\cdots {o}_n \ 
%,\quad \eb_i\in\cH\ .
 \end{equation*}
for homogeneous elements $ {o}_1,\dots, {o}_n\in\cH$, 
where $ {o}_1+\cdots+ {o}_{p-1}$ on $(-1)$ denotes the degree of 
$ {o}_1\cdots {o}_{p-1}$. The sign factor appears when $m_k$, which 
has degree one, passes through the $ {o}_1\cdots {o}_{p-1}$. 

Then summing up these $\m_k$ for $k\ge 0$, 
\begin{equation}
 \m=\m_0+\m_1+\m_2+\cdots\ ,
\end{equation}
and this $\m$ is the coderivation. 
The coderivation on the coalgebra $C(\cH)$ is always written in this 
form. 

Moreover, the form of a cohomomorphism $\cF :C(\cH)\to C(\cH')$ 
is determined by a collection of degree zero multilinear maps 
$\{f_k: \cH^{\otimes k}\rightarrow \cH'\}_{k\ge 0}$. 
For homogeneous elements $ {o}_1,\dots, {o}_n\in\cH$, it is given as 
%which are homogeneous of degree zero of the following form
\begin{equation}
 \begin{split}
 \cF( {o}_1\cdots {o}_n)=\sum_{1\leq k_1<k_2\cdots <k_i=n}&
e^{f_0(\1)}f_{k_1}( {o}_1,\dots, {o}_{k_1}) e^{f_0(\1)} 
f_{k_2-k_1}( {o}_{k_1+1},\dots, {o}_{k_2}) e^{f_0(\1)}\\
&\cdots e^{f_0(\1)} f_{n-k_{i-1}}
( {o}_{k_{i-1}+1},\dots, {o}_n)\ ,
 \end{split}
 \label{cohom}
\end{equation}
where each $f(\cdots)$ belongs to $\cH'$ and 
$e^{f_0(\1)}$ is defined by 
\begin{equation*}
 e^{f_0(\1)}=\1+f_0(\1)+f_0(\1)\otimes f_0(\1)+
f_0(\1)\otimes f_0(\1)\otimes f_0(\1)+\cdots\ .
\end{equation*}
If $f_0(\1)=0$, eq.(\ref{cohom}) is simplified since 
$e^{f_0(\1)}=\1$. Note that $\1$ is defined as 
$\cH^{\otimes m}\otimes\1\otimes\cH^{\otimes n}=
\cH^{\otimes (m+n)}$ for $m,n\ge 0$ and $m+n\ge 1$.

 \subsection{$A_\infty$-algebras and $A_\infty$-morphisms}
\label{ssec:Ainfty}

\begin{defn}[$A_\infty$-algebra \cite{Sta1,Sta11}] 
Let $\cH$ be a graded vector space and 
$C(\cH)=\oplus_{k\geq 0}{\cH^{\otimes k}}$ be its tensor coalgebra. 
A {\it weak $A_\infty$-algebra} is a coalgebra $C(\cH)$ with a coderivation 
$\m=\m_0+\m_1+\m_2+\cdots$ satisfying 
\[
 (\m)^2=0\ .
\]
We denote the collection of multilinear maps $\{m_k\}_{k\ge 0}$ also 
by $\m$ and the weak $A_\infty$-algebra by $(\cH,\m)$. 
In particular, $(\cH,\m)$ is called an {\it $A_\infty$-algebra} if $m_0=0$. 
 \label{defn:Ainfty}
\end{defn}
In general, a coderivation $\m:C\to C$ on a coalgebra $C$ satisfying 
$(\m)^2=0$ as above is called 
a {\it codifferential}. Thus, a (weak) $A_\infty$-algebra is a 
{\em differential graded coalgebra} of the tensor coalgebra of a 
graded vector space $\cH$. 

For an $A_\infty$-algebra $(\cH,\m)$, 
if we act $(\m)^2=(\m_1+\m_2+\cdots)^2$ on $ {o}_1\cdots {o}_n\in C(\cH)$ 
for homogeneous elements $ {o}_1,\dots, {o}_n\in\cH$, 
its image belongs to $\cH^{\otimes 1}\oplus\cdots\oplus\cH^{\otimes n}$. 
Then, the condition that the $\cH^{\otimes 1}$ part of the image 
is equal to zero turns out to be 
\begin{equation}
\sum_{\substack{k+l=n+1\\j=0,\dots,k-1}}
{(-1)^{ {o}_1+\cdots+ {o}_j}
 m_k( {o}_1,\dots, {o}_j,m_l( {o}_{j+1},\dots, {o}_{j+l}),
  {o}_{j+l+1},\dots, {o}_n)}=0\ ,
 \label{ainf}
\end{equation}
where $ {o}_i$ on $(-1)$  
denotes the degree of $ {o}_i$. 
The collection of the identities (\ref{ainf}) for $n\geq 1$ is the original 
definition of $A_\infty$-algebras. 
On the other hand, the construction of $A_\infty$-algebras 
using the tensor coalgebra $C(\cH)$ 
is called {\it the bar construction}. 
Actually, eq.(\ref{ainf}) is equivalent (\ie sufficient) to $\m^2=0$ 
due to the anticommutativity of $\m_i$'s. 
This fact is clearer in the dual 
language in the next section, where 
the nilpotent coderivation $\m$ is replaced to a 
nilpotent differential $\delta$ on a formal (noncommutative) supermanifold. 
\footnote{
The dual language fits the field theory. 
In the context of BRST-formalism, physicists usually show 
the nilpotency of BRST-operator $\delta$ on polynomials of fields 
and ghosts (and antifields) by 
confirming the nilpotency on each component field. 
This is just the dual of eq.(\ref{ainf}). 
} 

Let us write down the first three constraints in eq.(\ref{ainf}). 
\bea
\label{a3}
& &m_1^2=0~~\nn\ ,\\
& &m_1(m_2( {o}_1, {o}_2))+m_2(m_1( {o}_1), {o}_2)
+(-1)^{ {o}_1}m_2( {o}_1,m_1( {o}_2))=0~~\ ,\\
& & m_2(m_2( {o}_1,  {o}_2),  {o}_3) 
+ (-1)^{ {o}_1} m_2( {o}_1, m_2( {o}_2,  {o}_3))\nn\\
& &~~~~~~~
+m_1(m_3( {o}_1, {o}_2, {o}_3))+
m_3(m_1( {o}_1),  {o}_2,  {o}_3)+(-1)^{ {o}_1}m_3( {o}_1, m_1( {o}_2),  {o}_3)
\nn\\
& &~~~~~~~+(-1)^{ {o}_1+ {o}_2}m_3( {o}_1,  {o}_2, m_1( {o}_3))=0~~.\nn
\eea
The first equation indicates $m_1$ is nilpotent 
and $(\cH, m_1)$ defines a complex on the $\Z$-graded vector space $\cH$. 
The second equation implies differential $m_1$ satisfies Leibniz rule for 
the product $m_2$. 
The third equation means product $m_2$ is associative up to the 
terms including $m_3$. 
\begin{rem}
In the case $m_n=0$ for $n\geq 3$, an $A_\infty$-algebra 
reduces to a differential graded (associative) algebra (dga). 
The differential $d$ and the product $\bullet$ 
of dga $\g$ correspond to $m_1$ and $m_2$, respectively. However, the product 
$\bullet$ of dga preserves the degree and $m_2$ in $A_\infty$-algebras 
raises the degree by one. For this reason, when a dga $(\g,d,\bullet)$
is considered as 
an $A_\infty$-algebra $(\cH,\m)$ defined in Definition \ref{defn:Ainfty}, 
the degree in the $A_\infty$-algebra is defined 
as the degree of the dga minus one. 
Namely, let $s : \g^k\to(\g[1])^{k-1}=:\cH^{k-1}$ 
be the isomorphism called the {\it suspension}, 
where $\g^k$ is the degree $k$ part of $\g$. 
This $[1]$ `eats' one degree of $\g$, and then the degree 
of the image $(\g[1])^{k-1}$ is $(k-1)$ through the operation. 
Then the following diagram commutes 
\begin{equation*}
 \begin{CD}
  \g^k\otimes\g^l @>{\ \bullet\ }>> \g^{k+l}\\
   @V{s}VV               @V{s}VV\\
  \cH^{k-1}\otimes\cH^{l-1} @>{m_2(\ ,\ )}>> \cH^{(k+l-2)+1}\ .
 \end{CD}
\end{equation*}
There are many literatures where the degree of 
$A_\infty$-algebras are defined with the dga degree. 
Usually Witten's open string field theory \cite{W1}, 
which has the structure of a dga, is also 
defined with the dga degree explained above. 
The degree is certainly natural from the origin of $A_\infty$-algebras 
(see subsection \ref{ssec:Ainftysp}). 
However, when higher products $m_3, m_4,\dots$ are introduced, 
the degree given in Definition \ref{defn:coder} 
is simpler for $A_\infty$-algebras. 
For this reason, we use this convention in the present paper. 
The precise relation between these two conventions can be found in \cite{GJ}.
 \label{rem:degree}
\end{rem}
\begin{defn}[$A_\infty$-morphism]
Given two weak $A_\infty$-algebras 
$(\cH, \m)$ and $(\cH',\m')$, 
{\it a weak $A_\infty$-morphism} $\cF: (\cH, \m)\to (\cH',\m')$ is a 
cohomomorphism from $C(\cH)$ to $C(\cH')$ satisfying 
\begin{equation}
 \cF\m=\m'\cF\ .
 \label{Ainftymorp}
\end{equation}
In particular for two $A_\infty$-algebras $(\cH,\m)$ and $(\cH',\m')$ 
a weak $A_\infty$-morphism $\cF: (\cH, \m)\to (\cH',\m')$ is called 
{\it an $A_\infty$-morphism} iff $f_0=0$. 
\label{defn:amorp}
\end{defn}
For an $A_\infty$-morphism $\cF: (\cH, \m)\to (\cH',\m')$, 
evaluating the condition (\ref{Ainftymorp}) 
with $ {o}_1\cdots {o}_n\in C(\cH)$, $n\ge 1$, 
for homogeneous elements $ {o}_1,\dots, {o}_n\in\cH$, one gets 
a relation in $\oplus_{n'=1}^n{\cH'}^{\otimes n'}$. 
Picking up the ${\cH'}^{\otimes 1}$ part of the equation then yields 
\begin{equation}
 \begin{split}
&\sum_{1\leq k_1<k_2\dots <k_i=n}{m'_i(f_{k_1}( {o}_1,\dots, {o}_{k_1}),
f_{k_2-k_1}( {o}_{k_1+1},\dots, {o}_{k_2}),\dots, 
f_{n-k_{i-1}}( {o}_{k_{i-1}+1},\dots, {o}_n))}\\
&\qquad=\sum_{k+l=n+1}\sum_{j=0}^{k-1}{(-1)^{ {o}_1+\dots + {o}_j}
f_k( {o}_1,\dots, {o}_j,m_l( {o}_{j+1},\dots, {o}_{j+l}), 
 {o}_{j+l+1},\dots, {o}_n)}\ .
\end{split}
\label{amorphism}
\end{equation}
The first two constraints in (\ref{amorphism}) read:
\begin{equation}
 \begin{split}
 m'_1(f_1( {o}_1))&=f_1(m_1( {o}_1))\ ,\\
m'_2(f_1( {o}_1),f_1( {o}_2))&=
f_1(m_2( {o}_1, {o}_2))\\
&+m'_1(f_2( {o}_1, {o}_2))+f_2(m_1( {o}_1), {o}_2)
+(-1)^{ {o}_1}f_2( {o}_1, m_1( {o}_2))\ .\nn
 \end{split}
\end{equation}
In particular, the first equation implies that 
$f_1$ is a chain map 
between the complexes $(\cal{H},m_1)$ and $(\cal{H}',m_1')$. 
In the dual picture explained in the next section, 
$\cF$ is identified with a nonlinear map between two supermanifolds. 
\begin{defn}[$A_\infty$-(quasi)-isomorphism]
An $A_\infty$-morphism $\cF=\{f_1,f_2,\dots\}: (\cH,\m)\to
(\cH',\m')$ is called an {\it $A_\infty$-quasi-isomorphism} if
$f_1$ induces an isomorphism 
between the cohomology spaces $H_{m_1}(\cH)$ and $H_{m_1'}(\cH')$. 
In particular, if $f_1:\cH\to\cH'$ is an isomorphism, $\cF$ is called 
an {\it $A_\infty$-isomorphism}. Moreover, 
if $(\cH',\m')=(\cH,\m)$, 
we call $\cF$ an {\it $A_\infty$-automorphism}. 
 \label{defn:quasiisom}
\end{defn}
These are also defined in weak $A_\infty$ level. 
It is clear that any $A_\infty$-isomorphism $\cF:(\cH,\m)\to
(\cH',\m')$ has its inverse $A_\infty$-isomorphism 
$\cF^{-1}: (\cH',\m')\to (\cH,\m)$. 
Also, if $\cF$ is an $A_\infty$-quasi-isomorphism, there exists an 
inverse quasi-isomorphism (we denote it also by $\cF^{-1}$) 
\cite{K2,Ko1}, which will be discussed in subsection \ref{ssec:inverse}. 
\begin{rem}[Cocommutativity and $L_\infty$-algebras]
$L_\infty$-algebras are obtained by 
imposing cocommutativity upon coalgebra $C(\cH)$. A coalgebra $C$ is
{\it cocommutative} iff there exists an operator 
$\tau: C\otimes C\to C\otimes C$, 
$\tau(x\otimes y)=(-1)^{x y}y\otimes x$ that is compatible with 
the coproduct, 
$$\tau\triangle=\triangle\ . $$
The corresponding tensor coalgebra is $C(\cH)$ divided 
by the ideal generated by 
$ {o}_i\otimes {o}_j-(-1)^{ {o}_i {o}_j} {o}_j\otimes {o}_i$. 
Namely, in this case elements in $\cH$ are set to be graded commutative. 
An $L_\infty$-algebra is then obtained by defining degree one
codifferential (coderivation whose square is zero) 
so that it is compatible with the graded commutativity, 
that is, by graded symmetrizing each multi-linear map $m_k$ \cite{LM} 
(see also \cite{LS,Fukaya2,ocha}, etc.). 
 \label{rem:cocomm}
\end{rem}

 \subsection{Cyclic $A_\infty$-structures}
\label{ssec:cycAinfty}

In this subsection $A_\infty$-structures with cyclic 
symmetry are defined. 
We consider a graded vector space $\cH$ equipped with an odd constant 
symplectic inner product. 
An origin of these definitions is the BV-formalism 
as explained in subsection \ref{ssec:BV}. 
The naturalness of these definitions can be realized 
from the dual picture in section \ref{sec:sym}. 
\begin{defn}[Odd constant symplectic structure]
Let $\cH$ be a graded vector space. 
An {\it odd constant symplectic structure} $\omega :\cH\otimes\cH\to\C$ 
is a nondegenerate skewsymmetric bilinear map of degree minus one. 
%For bases $\eb_i, \eb_j\in\cH$, it is represented as 
%\begin{equation*}
% \omega(\eb_i,\eb_j)=\omega_{ij}\ . 
%\end{equation*}
Namely, 
for any homogeneous elements $ {o},  {o}'\in\cH$, 
$\omega( {o}, {o}')\in\C$ can be nonzero only if 
$\deg( {o})+\deg( {o}')=1$ 
since the degree of $\omega$ is minus one, 
%Moreover $\omega_{ji}=-\omega_{ij}$ since it is skewsymmetric. 
and $\omega( {o}', {o})=-\omega( {o}, {o}')$ 
since it is skewsymmetric. 
Also, for any $ {o}\in\cH$, there exists an element $ {o}'\in\cH$ 
such that $\omega( {o}, {o}')\ne 0$ since $\omega$ is nondegenerate. 
 \label{defn:csym}
\end{defn}
\begin{defn}[Cyclic $A_\infty$-algebra]
Suppose a graded vector space $\cH$ is equipped with 
an odd constant symplectic structure $\omega :\cH\otimes\cH\to\C$ 
in Definition \ref{defn:csym}. 
A triple $(\cH,\omega,\m)$ is called a {\it cyclic $A_\infty$-algebra} 
when $(\cH,\m)$ is an $A_\infty$-algebra and 
$\m$ is {\it cyclic} with respect to $\omega$, that is, 
\begin{equation*}
 \omega( {o}_1,m_k( {o}_2,\dots, {o}_{k+1}))=
 (-1)^{ {o}_2}
\omega( {o}_2,m_k( {o}_3,\dots, {o}_{k+1}, {o}_1)) 
\end{equation*}
holds for any homogeneous elements $ {o}_1,\dots, {o}_{k+1}\in\cH$ 
for each $k\ge 1$. 
 \label{defn:cycAinfty}
\end{defn}
$A_\infty$-algebras with cyclic symmetry as above are considered 
in the context of mathematical physics for instance in \cite{Ko4,Z2,GZ}. 
See \cite{MSS}. Also, a homotopy extension of this cyclicity is
proposed in \cite{tradler}. 
\begin{rem}
Let us define a collection of degree zero multilinear maps 
$S:=\{\V_k :\cH^{\otimes k}\to\C\}_{k\ge 2}$ by 
\begin{equation*}
 \V_{k+1}( {o}_1,\dots, {o}_{k+1})
 :=(-1)^{ {o}_1}\omega( {o}_1,m_k( {o}_2,\dots, {o}_{k+1}))\ .
\end{equation*}
The cyclicity of $\m$ implies 
$\V_{k+1}( {o}_1,\dots, {o}_{k+1})
=(-1)^{ {o}_1}\V_{k+1}( {o}_2,\dots, {o}_{k+1}, {o}_1)$. 
Then a cyclic $A_\infty$-algebra $(\cH,\omega,\m)$ can also be defined 
by triple $(\cH,\omega, S)$. Hereafter we use both notations 
for a cyclic $A_\infty$-algebra. 
Suppose first that $S$ of degree zero is given. 
The degree of the inner product $\omega$ is then determined 
as minus one (odd). 
 \label{rem:Vcyclic}
\end{rem}
\begin{defn}[Cyclic $A_\infty$-morphism]
Let $(\cH,\omega, S)$ and $(\cH',\omega', S')$ be 
two cyclic $A_\infty$-algebras and suppose 
there exists an $A_\infty$-morphism $\cF: (\cH,\m)\to(\cH',\m')$. 
We then call $\cF$ a {\it cyclic $A_\infty$-morphism} when 
\begin{equation}
 \omega'(f_1( {o}),f_1( {o}'))=\omega( {o}, {o}')\ ,
 \label{omegacF1}
\end{equation}
for any $ {o}, {o}'\in\cH$ and for fixed $n\ge 3$, 
\begin{equation}
 \sum_{k,l\ge 1,\ k+l=n}
\omega'(f_k( {o}_1,\dots, {o}_k),
 f_l( {o}_{k+1},\dots, {o}_n))=0\ 
 \label{omegacF2}
\end{equation}
holds for any $ {o}_1,\dots, {o}_n\in\cH$. 
 \label{defn:cycAinftymorp}
\end{defn}
\begin{rem}
We shall explain in subsection \ref{ssec:cycAinftyre} that, 
in the dual language, 
$\cF$ is just the morphism which preserves 
the constant symplectic forms and the actions. 
 \label{rem:cycAinftymorpdual}
\end{rem}
It is clear that cyclic $A_\infty$-algebras and cyclic $A_\infty$-morphisms 
can be defined also at weak level in the same way 
as in the previous subsection.

 \subsection{Maurer-Cartan equations and deformation theory}
\label{ssec:MCeq}

Here we shall define Maurer-Cartan equations for $A_\infty$-algebras 
and their roles in deformation theory 
(see \cite{Fukaya2}). 
Maurer-Cartan equation is an equation for elements of $\cH$. 
Therefore we usually concentrate on its degree-zero part 
so that the equation has an usual meaning. 
Let us consider an element $\Phi\in\cH^0$ where $\cH^0$ is the 
degree-zero subvector space of $\cH$. 
Using the basis of $\cH^0$, $\{\eb_i^0\}$, 
one can express it as $\Phi=\sum_i\eb^0_i\phi^i$ for $\phi^i\in\C$. 
Note that we shall define later a `superfield' $\Phi$ 
as a more extended object than the one here, and consider the Maurer-Cartan 
equation in a similar manner as below. 
Alternatively, we shall discuss the situation here 
in a precise way 
in subsection \ref{ssec:homgauge}. 
In this section we present some abstract definitions without the details. 

Consider formally the following exponential map of $\Phi\in \cH^0$ 
\begin{equation}
 e^\Phi:=\1+\Phi+\Phi\otimes\Phi+\Phi\otimes\Phi\otimes\Phi+\cdots\ .
\end{equation}
$e^\Phi\in C(\cH^0)\subset C(\cH)$ satisfies 
$\tri e^\Phi=e^\Phi\otimes e^\Phi$ and such element 
is called a {\it grouplike element}. 
\begin{defn}[Maurer-Cartan equation]
For an $A_\infty$-algebra $(\cH,\m)$, define 
\begin{equation*}
 \m_*(e^\Phi):=m_1(\Phi)+m_2(\Phi\otimes\Phi)
 +m_3(\Phi\otimes\Phi\otimes\Phi)+\cdots\ .
\end{equation*}
$\m_*(e^\Phi)=0$ is called {\it Maurer-Cartan equation} for $(\cH,\m)$. 
We denote by $\cMC(\cH,\m)$ 
the solution space of the Maurer-Cartan equation for $(\cH,\m)$. 
 \label{defn:mceq}
\end{defn}
Because $\m(e^{\Phi})=e^\Phi\cdot\m_*(e^\Phi)\cdot e^\Phi$, 
$\m_*(e^\Phi)=0$ is equivalent to $\m(e^\Phi)=0$, 
where $\1$ is defined as 
$\cH^{\otimes m}\otimes\1\otimes\cH^{\otimes n}=
\cH^{\otimes (m+n)}$ for $m,n\ge 0$ and $m+n\ge 1$. 
When an $A_\infty$-algebra is a dga, {\it i.e.} $m_3=m_4=\cdots=0$, 
its Maurer-Cartan equation takes the form 
$m_1(\Phi)+m_2(\Phi\otimes\Phi)=0$. It is nothing but the condition of 
a `flat connection'. 
In the case of field theory equipped with a classical BV-structure, 
the theory has a cyclic $A_\infty$-structure, and 
its Maurer-Cartan equation is just the equation of motion 
of the action (see eq.(\ref{eomsft-re})). 
\begin{rem}
The solution space $\cMC(\cH,\m)$ 
of the Maurer-Cartan equation for $(\cH,\m)$ 
parameterizes deformations of original $A_\infty$-algebra $(\cH,\m)$. 

Generally, for a $A_\infty$-algebra $(\cH,\m)$ and a $\Z$-graded
vector space $\ti\cH$, 
suppose a cohomomorphism $\ti\cF:C(\ti\cH)\to C(\cH)$ is given. 
Then there exists the inverse cohomomorphism 
$\ti\cF^{-1}:C(\ti\cH)\to C(\cH)$ such that
$\ti\cF\ti\cF^{-1}=\1$ and $\ti\cF^{-1}\ti\cF=\1$. 
Namely, a weak $A_\infty$-structure
$\ti\m=\ti\cF^{-1}\m\ti\cF:C(\ti\cH)\to C(\ti\cH)$ is induced. 
Actually, it is clear that $\ti\cF\ti\m=\m\ti\cF$ holds. 

In this situation let us consider the cohomomorphism 
$\ti\cF=\{\ti{f}_0, \ti{f}_1,\dots\}$ 
with $\ti{f}_0=\Phi\in\cH^0$, $\ti{f}_1=\Id$ and
$\ti{f}_k=0$ for $k\ge 2$. 
The induced $A_\infty$-structure $\ti\m=\{\ti{m}_0,\ti{m}_1,\dots\}$ 
is given of the form 
\begin{equation*}
 \ti{m}_k( {o}_1,\dots, {o}_k)=
 \sum_{l_0\ge 0,\dots,l_k\ge 0}
 m_{k+l_0+\cdots +l_k}(\Phi^{\otimes l_0}, {o}_1,
 \Phi^{\otimes l_1},\dots, \Phi^{\otimes l_{k-1}}, {o}_k,
 \Phi^{\otimes l_k})\ 
\end{equation*}
for $k\ge 0$. In the equation above, $ {o}_i$ in both sides are 
identified with each other by $\ti{f}_1=\Id :\ti\cH\to\cH$. 
One can easily see, by concentrating on the case $k=0$ in the equation above, 
that the induced weak $A_\infty$-algebra $(\ti\cH,\ti\m)$ is an 
$A_\infty$-algebra if and only if $\Phi\in\cMC(\cH,\m)$.

In the case of cyclic $A_\infty$-algebras, 
this means that, for each solution of the equation of motion, 
another $A_\infty$-algebra is defined. 
Such a property is observed in classical closed string field theory, 
\ie for (cyclic) $L_\infty$-algebras in \cite{Sinfty} 
and is related to the problem of 
background independence of string field theory (see the end of 
subsection \ref{ssec:main}). 

The solution space $\cMC(\cH,\m)$ 
is a subspace of the whole deformation space 
of $A_\infty$-algebra $(\cH,\m)$ that is defined by 
\begin{equation*}
 \{\m_{def}:C(\cH)\to C(\cH);\, \mbox{degree one coderivation}\ | 
 (\m+\m_{def})^2=0\}\ .
\end{equation*}
From string theory point of view, the $A_\infty$-structure $\m$ is 
related to a structure of tree open string interactions 
(for the case of string field theories see subsection \ref{ssec:osft}, 
and for a topological string case see for instance \cite{Fukaya}), 
and the deformation associated to $\cMC(\cH,\m)$ corresponds to 
deformation that comes from condensation of open string fields. 
 \label{rem:MCdeform}
\end{rem}
When we are interested in the space $\cMC(\cH,\m)$, 
it is often convenient to relate $(\cH,\m)$ to another
$A_\infty$-algebra $(\ti\cH,\ti\m)$. 
Suppose that there exists an $A_\infty$-morphism 
$\ti\cF : (\ti\cH,\ti\m)\to (\cH,\m)$. 
The $A_\infty$-morphism $\ti\cF$ then preserves the solutions of 
the Maurer-Cartan equations. 
In the context of field theory, 
this fact means that cyclic $A_\infty$-morphisms 
preserve the equations of motions. 
For $\ti\Phi\in\cMC(\ti\cH,\ti\m)$, 
$\Phi$ is constructed as the pushforward of $\ti\cF$, 
a (nonlinear) coordinate transformation 
between two formal noncommutative supermanifolds (see the next section), 
\begin{equation}
 \Phi=\ti\cF_*(\ti\Phi)=\sum_{n=1}^{\infty}\ti{f}_n(\ti\Phi,\dots,\ti\Phi)\ .
 \label{fredef}
\end{equation}
It by construction satisfies $\ti\cF(e^{\ti\Phi})=e^{\Phi}$. 
The following equality then holds, 
\begin{equation*}
 \m(e^{\Phi})=\m\ti\cF(e^{\ti\Phi})=\ti\cF\ti\m(e^{\ti\Phi})\ ,
\end{equation*}
and one can immediately see that $\Phi\in\cMC(\cH,\m)$ 
(\ie $\m(e^{\Phi})=0$) if $\ti\Phi\in\cMC(\ti\cH,\ti\m)$ 
(\ie $\ti\m(e^{\ti\Phi})=0$). 

More precisely, 
for a given $A_\infty$-algebra $(\cH,\m)$, 
there exists a notion of 
gauge transformation (Definition \ref{defn:gauge}). 
It is, so to speak, an automorphism of the theory generated 
by infinitesimal transformations. 
As seen in subsection \ref{ssec:BV} 
it just corresponds to the gauge transformation in classical BV-field theory. 
By definition it preserves $\cMC(\cH,\m)$, 
\ie , the equation of motion. 
What should be considered is then, instead of $\cMC(\cH,\m)$, 
the moduli space of $A_\infty$-algebra $(\cH,\m)$, 
that is defined by dividing $\cMC(\cH,\m)$ over the gauge transformation 
(see Definition \ref{defn:moduli}). 
An $A_\infty$-morphism, that preserves the solution space 
of Maurer-Cartan equations, in fact induces a well-defined 
morphism between the moduli spaces. 
In particular, it is known that the moduli spaces are isomorphic to 
each other if there exists an $A_\infty$-quasi-isomorphism between 
them (Theorem \ref{thm:quasiisomoduli}). 
Such generality of homotopy algebras is, for instance in $L_\infty$ case, 
applied by M.~Kontsevich \cite{Ko1} to the proof of the existence of 
deformation quantizations \cite{BFFLS} on Poisson manifolds 
and their classification.

%section3

 \section{Dual geometric description of homotopy algebras}
\label{sec:dual}

We shall give the dual description of $A_\infty$-algebras. 
In this picture, $A_\infty$-algebras are understood more geometrically. 
The definition of the dual of a coalgebra in the present paper is 
given in subsection \ref{ssec:dual1}. Its geometric viewpoint 
is explained in subsection \ref{ssec:dual2}, where we deal with a 
formal noncommutative supermanifold. 
In subsection \ref{ssec:superfield} we define 
the notion of superfield which simplify the convention in the dual
picture. It will be used in later discussions. 
These arguments hold similarly for $L_\infty$-algebras.

 \subsection{The dual of coalgebras}
\label{ssec:dual1}

Let $\cH$ be a graded vector space, and 
$C(\cH):=\oplus_{n=0}^\infty\left(\cH^{\otimes n}\right)$ be its tensor 
coalgebra. The basis of $\cH$ is denoted by $\{\eb_i\}$, and 
here we define the dual basis of 
$\eb_{i_1}\cdots\eb_{i_k}\in\cH^{\otimes k}$ 
by introducing a natural pairing as follows. 
At first, denote the dual basis of $\{\eb_i\}$ by $\{\eb^i\}$, and 
define a pairing between $\cH$ and $\cH^*$ as 
\begin{equation}
 \langle \eb^i |\eb_j\rangle =\delta^i_j\ .
 \label{ip1}
\end{equation}
We represent an elements of $C(\cH)$ as 
$g=\sum_{k=1}^\infty g^{i_k\cdots i_1}\eb_{i_1}\cdots\eb_{i_k}$, 
and an element of $C(\cH)^*$, the dual of $C(\cH)$ as  
$a=\sum_{k=1}^\infty a_{i_1\cdots i_k}\eb^{i_k}\cdots\eb^{i_1}$. 
Generalizing the above pairing between $\cH$ and $\cH^*$ (\ref{ip1}), 
here the pairing between $C(\cH)$ and $C(\cH)^*$ is defined as 
\begin{equation}
 \langle \eb^{i_k}\cdots\eb^{i_1} | \eb_{j_1}\cdots\eb_{j_l}\rangle 
 =\epsilon^{i_1\cdots i_k}_{j_1\cdots j_l}\label{ip2}, 
\end{equation}
where $\epsilon^{i_1\cdots i_k}_{j_1\cdots j_l}$ is defined to be zero 
if $k\neq l$ and 
$\epsilon^{i_1\cdots i_k}_{j_1\cdots j_k}
=\delta^{i_1}_{j_1}\cdots\delta^{i_k}_{j_k}$ if $k=l$. 
In addition we define $\la\1|\1\ra=1$. 
Moreover, for $a_1,\dots,a_n\in C(\cH)^*$ and $g_1,\dots,g_n\in C(\cH)$, 
the pairing of $n$-tensor is given by 
\begin{equation*}
 \la a_1\otimes\cdots\otimes a_n| g_1\otimes\cdots\otimes g_n\ra
 =\la a_1| g_1\ra\cdots\la a_n| g_n\ra\ .
 \label{ipn}
\end{equation*}

Since now we have obtained the pairing 
between $C(\cH)$ and its dual $C(\cH)^*$, 
we can translate operations on $C(\cH)$ into those on $C(\cH)^*$. 
For the coproduct $\triangle$ on $C(\cH)$, 
the product $m$ on $C(\cH)^*$ is defined as 
\begin{equation}
 \langle m(a\otimes b) | g\rangle=\langle a\otimes b | \triangle g\rangle\ ,
 \label{defm}
\end{equation}
the derivation $\delta$ corresponding to the coderivation $\m$ is defined as 
\begin{equation}
 \langle \delta(a) | g\rangle=\langle a | \m(g)\rangle\ ,
 \label{defD}
\end{equation}
and homomorphism $\fb$ corresponds 
to the cohomomorphism $\cF$ from $C(\cH)$ 
to another tensor algebra $C(\cH')$ is 
determined as 
\begin{equation}
 \langle \fb(a) | g\rangle=\langle a | \cF(g)\rangle\ .
 \label{deff}
\end{equation}
Because $g\in C(\cH)$ and $a\in C(\cH')^*$, the homomorphism $\fb$ is a 
map from $C(\cH')^*$ to $C(\cH)^*$. 
Therefore $\fb$ is in fact the pullback of $\cF$. 
Here we write the elements of $C(\cH)$ on the left hand side and 
the elements of $C(\cH)^*$ on the right hand side. 
The operations on $C(\cH')$ or $C(\cH')^*$ are distinguished by 
attaching $'$ to them. 
The above definitions of the operations on $C(\cH)^*$ translate 
various conditions for the operations on $C(\cH)$ into those on $C(\cH)^*$ 
as follows. 
The coassociativity of $\tri$ is equivalent to the associativity of $m$ : 
\begin{equation}
 \begin{array}{ccc}
 \langle m(m(a\otimes b)\otimes c) | g\rangle&=&
   \langle a\otimes b\otimes c | (\tri\otimes\1)\tri(g)\rangle\\
                       &&\rotatebox[origin=c]{90}{=}\\
   \langle m(a\otimes m(b\otimes c))
    | g\rangle
 &=&\langle a\otimes b\otimes c|(\1\otimes\tri)\triangle g\rangle
 \end{array}\ .\label{asso}
\end{equation}
The condition that $\m$ is the coderivation is translated 
into the Leibniz rule for $\delta$: 
\begin{equation}
 \begin{array}{ccc}
 \langle \delta\cdot m(a\otimes b) | g\rangle&=&
   \langle a\otimes b | \tri\cdot \m (g)\rangle\\
                       &&\rotatebox[origin=c]{90}{=}\\
   \langle m(\delta\otimes{\bf 1}+{\bf 1}\otimes \delta)(a\otimes b)
    | g\rangle
 &=&\langle a\otimes b|(\m\otimes{\bf 1}+{\bf 1}\otimes \m)\triangle g\rangle
 \end{array}\ .\label{QD}
\end{equation}
The condition that $\cF : C(\cH)\to C(\cH')$ 
is a cohomomorphism is rewritten as the one that 
$\fb : C(\cH')^*\to C(\cH)^*$ is a homomorphism: 
\begin{equation}
 \begin{array}{ccc}
 \langle \fb\cdot m'(a\otimes b) | g\rangle&=&
   \langle a\otimes b | \tri'\cdot \cF(g)\rangle\\
                       &&\rotatebox[origin=c]{90}{=}\\
   \langle m(\fb(a)\otimes \fb(b)) | g\rangle
 &=&\langle a\otimes b|(\cF\otimes\cF)\triangle g\rangle
 \end{array}\ .\label{Ff}
\end{equation}
$(\cH, \m)$ is an $A_\infty$-algebra means 
that $(C(\cH)^*, \delta)$ is a complex on the dual: 
\begin{equation}
 0=\langle \delta\cdot\delta(a) | g\rangle
 =\langle a | \m\cdot\m (g)\rangle=0\ .\label{QQDD}
\end{equation} 
Finally the condition that $\cF$ is an $A_\infty$-morphism is translated 
into the equivariance of $\fb$: 
\begin{equation}
 \begin{array}{ccc}
 \langle \delta\cdot \fb(a) | g\rangle&=&
   \langle a | \cF\cdot \m(g)\rangle\\
                       &&\rotatebox[origin=c]{90}{=}\\
   \langle \fb\cdot \delta'(a) | g\rangle
 &=&\langle a | \m'\cdot\cF(g)\rangle
 \end{array}\ .
 \label{QFDf}
\end{equation}

The above statement will be realized with some graphs.
\footnote{The graphs used below is different from that in the body of this 
paper. In fact, a line denotes the flow of an element of $\cH$ in the 
body of this paper, but the line used below denotes an element of $C(\cH)$. 
} 
In the above explanation, the elements of $C(\cH)$ are written in the left 
hand side of the pairings (ket), and the elements 
of the dual algebra $C(\cH)^*$ are in the right hand side (bra). 
Here, for the algebra on the left hand side, we represent the product $m$, 
the derivation $\delta$, and the homomorphism $\fb$ as  
$m=\mm$, $\delta=\D$, $\fb=\f$. According to the operations of the algebra 
from left, the lines of the graphs are connected to the right direction. 
In other words, the operations on the algebra $C(\cH)^*$ 
in the left hand side from left yields the flow from the left to the right 
on the lines of the graphs. 
Next, for the coalgebra $C(\cH)$ in the right hand side, we represent 
the coproduct $\triangle$, the coderivation $\m$, 
and the cohomomorphism $\cF$ as $\triangle=\mm$, $\m=\Q$, $\cF=\F$, 
and define the orientation of the operation from the right to the left 
on the lines of the graphs. 
Lastly, in order to distinguish the left and right in the pairings, 
we introduce $\langle\ |\ \rangle$ between the algebra $C(\cH)$ and 
the coalgebra $C(\cH)^*$. 

The definition of the algebra $C(\cH)^*$ dual to the coalgebra $C(\cH)$ 
(\ref{defm})(\ref{defD})(\ref{deff}) are written graphically as follows. 
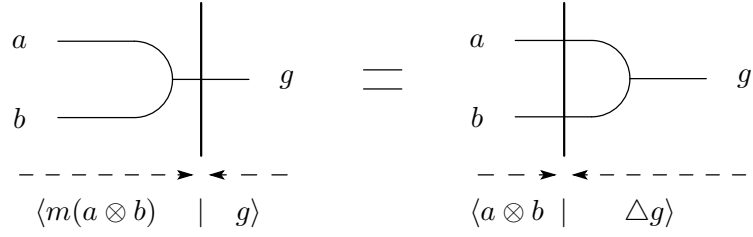
\begin{figure}[h]
 \begin{center}
%WinTpicVersion3.08
\unitlength 0.1in
\begin{picture}( 44.3500, 10.1500)( 11.7000,-18.1500)
% LINE 2 0 3 0
% 2 2005 1000 2405 1000
% 
\special{pn 8}%
\special{pa 2006 1000}%
\special{pa 2406 1000}%
\special{fp}%
% LINE 2 0 3 0
% 2 2005 1400 2405 1400
% 
\special{pn 8}%
\special{pa 2006 1400}%
\special{pa 2406 1400}%
\special{fp}%
% CIRCLE 2 0 3 0
% 4 2405 1200 2405 1400 2405 1400 2405 1000
% 
\special{pn 8}%
\special{ar 2406 1200 200 200  4.7123890 6.2831853}%
\special{ar 2406 1200 200 200  0.0000000 1.5707963}%
% LINE 2 0 3 0
% 2 2605 1200 3005 1200
% 
\special{pn 8}%
\special{pa 2606 1200}%
\special{pa 3006 1200}%
\special{fp}%
% STR 2 0 3 0
% 3 1805 900 1805 1000 5 0
% $a$
\put(18.0500,-10.0000){\makebox(0,0){$a$}}%
% STR 2 0 3 0
% 3 1805 1300 1805 1400 5 0
% $b$
\put(18.0500,-14.0000){\makebox(0,0){$b$}}%
% STR 2 0 3 0
% 3 3205 1100 3205 1200 5 0
% $g$
\put(32.0500,-12.0000){\makebox(0,0){$g$}}%
% LINE 2 0 3 0
% 2 3605 1150 3805 1150
% 
\special{pn 8}%
\special{pa 3606 1150}%
\special{pa 3806 1150}%
\special{fp}%
% LINE 2 0 3 0
% 2 3605 1250 3805 1250
% 
\special{pn 8}%
\special{pa 3606 1250}%
\special{pa 3806 1250}%
\special{fp}%
% LINE 2 0 3 0
% 2 4400 995 4800 995
% 
\special{pn 8}%
\special{pa 4400 996}%
\special{pa 4800 996}%
\special{fp}%
% LINE 2 0 3 0
% 2 4400 1395 4800 1395
% 
\special{pn 8}%
\special{pa 4400 1396}%
\special{pa 4800 1396}%
\special{fp}%
% CIRCLE 2 0 3 0
% 4 4800 1195 4800 1395 4800 1395 4800 995
% 
\special{pn 8}%
\special{ar 4800 1196 200 200  4.7123890 6.2831853}%
\special{ar 4800 1196 200 200  0.0000000 1.5707963}%
% LINE 2 0 3 0
% 2 5000 1195 5400 1195
% 
\special{pn 8}%
\special{pa 5000 1196}%
\special{pa 5400 1196}%
\special{fp}%
% STR 2 0 3 0
% 3 4200 895 4200 995 5 0
% $a$
\put(42.0000,-9.9500){\makebox(0,0){$a$}}%
% STR 2 0 3 0
% 3 4200 1295 4200 1395 5 0
% $b$
\put(42.0000,-13.9500){\makebox(0,0){$b$}}%
% STR 2 0 3 0
% 3 5600 1095 5600 1195 5 0
% $g$
\put(56.0000,-11.9500){\makebox(0,0){$g$}}%
% LINE 1 0 3 0
% 2 2755 800 2755 1600
% 
\special{pn 13}%
\special{pa 2756 800}%
\special{pa 2756 1600}%
\special{fp}%
% LINE 1 0 3 0
% 2 4655 800 4655 1600
% 
\special{pn 13}%
\special{pa 4656 800}%
\special{pa 4656 1600}%
\special{fp}%
% VECTOR 2 1 3 0
% 2 1805 1700 2705 1700
% 
\special{pn 8}%
\special{pa 1806 1700}%
\special{pa 2706 1700}%
\special{da 0.070}%
\special{sh 1}%
\special{pa 2706 1700}%
\special{pa 2638 1680}%
\special{pa 2652 1700}%
\special{pa 2638 1720}%
\special{pa 2706 1700}%
\special{fp}%
% VECTOR 2 1 3 0
% 2 3205 1700 2805 1700
% 
\special{pn 8}%
\special{pa 3206 1700}%
\special{pa 2806 1700}%
\special{da 0.070}%
\special{sh 1}%
\special{pa 2806 1700}%
\special{pa 2872 1720}%
\special{pa 2858 1700}%
\special{pa 2872 1680}%
\special{pa 2806 1700}%
\special{fp}%
% VECTOR 2 1 3 0
% 2 4205 1700 4605 1700
% 
\special{pn 8}%
\special{pa 4206 1700}%
\special{pa 4606 1700}%
\special{da 0.070}%
\special{sh 1}%
\special{pa 4606 1700}%
\special{pa 4538 1680}%
\special{pa 4552 1700}%
\special{pa 4538 1720}%
\special{pa 4606 1700}%
\special{fp}%
% VECTOR 2 1 3 0
% 2 5605 1700 4705 1700
% 
\special{pn 8}%
\special{pa 5606 1700}%
\special{pa 4706 1700}%
\special{da 0.070}%
\special{sh 1}%
\special{pa 4706 1700}%
\special{pa 4772 1720}%
\special{pa 4758 1700}%
\special{pa 4772 1680}%
\special{pa 4706 1700}%
\special{fp}%
% STR 2 0 3 0
% 3 2755 1800 2755 1900 5 0
% $|$
\put(27.5500,-19.0000){\makebox(0,0){$|$}}%
% STR 2 0 3 0
% 3 4655 1800 4655 1900 5 0
% $|$
\put(46.5500,-19.0000){\makebox(0,0){$|$}}%
% STR 2 0 3 0
% 3 2205 1800 2205 1900 5 0
% $\langle m(a\otimes b)$
\put(22.0500,-19.0000){\makebox(0,0){$\langle m(a\otimes b)$}}%
% STR 2 0 3 0
% 3 3005 1800 3005 1900 5 0
% $g\rangle$
\put(30.0500,-19.0000){\makebox(0,0){$g\rangle$}}%
% STR 2 0 3 0
% 3 4355 1800 4355 1900 5 0
% $\langle a\otimes b$
\put(43.5500,-19.0000){\makebox(0,0){$\langle a\otimes b$}}%
% STR 2 0 3 0
% 3 5105 1800 5105 1900 5 0
% $\triangle g\rangle$
\put(51.0500,-19.0000){\makebox(0,0){$\triangle g\rangle$}}%
\end{picture}%
 \end{center}
 \caption[subsection]{$\langle m(a\otimes b) | g\rangle
  =\langle a\otimes b | \triangle g\rangle$\ \ \ eq.(\ref{defm})}
 \label{fig:coalg1}
\end{figure}
\begin{figure}[h]
 \begin{center}
%WinTpicVersion3.08
\unitlength 0.1in
\begin{picture}( 37.9500,  7.0000)( 14.7000,-11.5000)
% LINE 2 0 3 0
% 2 1805 800 2205 800
% 
\special{pn 8}%
\special{pa 1806 800}%
\special{pa 2206 800}%
\special{fp}%
% BOX 2 0 3 0
% 2 2205 650 2505 950
% 
\special{pn 8}%
\special{pa 2206 650}%
\special{pa 2506 650}%
\special{pa 2506 950}%
\special{pa 2206 950}%
\special{pa 2206 650}%
\special{fp}%
% LINE 2 0 3 0
% 2 2505 800 2905 800
% 
\special{pn 8}%
\special{pa 2506 800}%
\special{pa 2906 800}%
\special{fp}%
% STR 2 0 3 0
% 3 1605 700 1605 800 5 0
% $a$
\put(16.0500,-8.0000){\makebox(0,0){$a$}}%
% STR 2 0 3 0
% 3 3105 700 3105 800 5 0
% $g$
\put(31.0500,-8.0000){\makebox(0,0){$g$}}%
% LINE 2 0 3 0
% 2 3405 750 3605 750
% 
\special{pn 8}%
\special{pa 3406 750}%
\special{pa 3606 750}%
\special{fp}%
% LINE 2 0 3 0
% 2 3405 850 3605 850
% 
\special{pn 8}%
\special{pa 3406 850}%
\special{pa 3606 850}%
\special{fp}%
% LINE 2 0 3 0
% 2 4100 800 4500 800
% 
\special{pn 8}%
\special{pa 4100 800}%
\special{pa 4500 800}%
\special{fp}%
% BOX 2 0 3 0
% 2 4500 650 4800 950
% 
\special{pn 8}%
\special{pa 4500 650}%
\special{pa 4800 650}%
\special{pa 4800 950}%
\special{pa 4500 950}%
\special{pa 4500 650}%
\special{fp}%
% LINE 2 0 3 0
% 2 4800 800 5200 800
% 
\special{pn 8}%
\special{pa 4800 800}%
\special{pa 5200 800}%
\special{fp}%
% STR 2 0 3 0
% 3 3900 700 3900 800 5 0
% $a$
\put(39.0000,-8.0000){\makebox(0,0){$a$}}%
% STR 2 0 3 0
% 3 5400 700 5400 800 5 0
% $g$
\put(54.0000,-8.0000){\makebox(0,0){$g$}}%
% LINE 1 0 3 0
% 2 2655 450 2655 1150
% 
\special{pn 13}%
\special{pa 2656 450}%
\special{pa 2656 1150}%
\special{fp}%
% LINE 1 0 3 0
% 2 4355 450 4355 1150
% 
\special{pn 13}%
\special{pa 4356 450}%
\special{pa 4356 1150}%
\special{fp}%
% STR 2 0 3 0
% 3 2360 700 2360 800 5 0
% $\delta$
\put(23.6000,-8.0000){\makebox(0,0){$\delta$}}%
% STR 2 0 3 0
% 3 4650 700 4650 800 5 0
% $\m$
\put(46.5000,-8.0000){\makebox(0,0){$\m$}}%
\end{picture}%
 \end{center}
 \caption[subsection]
{$\langle \delta(a) | g\rangle=\langle a | \m (g)\rangle$
  \ \ \ eq.(\ref{defD})}
 \label{fig:coalg2}
\end{figure}
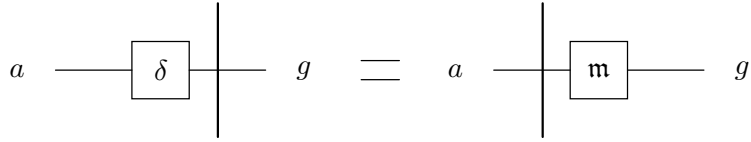
\begin{figure}[h]
 \begin{center}
%WinTpicVersion3.08
\unitlength 0.1in
\begin{picture}( 37.9500,  7.0000)( 14.7000,-11.5000)
% LINE 2 0 3 0
% 2 1805 800 2205 800
% 
\special{pn 8}%
\special{pa 1806 800}%
\special{pa 2206 800}%
\special{fp}%
% BOX 2 0 3 0
% 2 2205 650 2505 950
% 
\special{pn 8}%
\special{pa 2206 650}%
\special{pa 2506 650}%
\special{pa 2506 950}%
\special{pa 2206 950}%
\special{pa 2206 650}%
\special{fp}%
% LINE 2 0 3 0
% 2 2505 800 2905 800
% 
\special{pn 8}%
\special{pa 2506 800}%
\special{pa 2906 800}%
\special{fp}%
% STR 2 0 3 0
% 3 1605 700 1605 800 5 0
% $a$
\put(16.0500,-8.0000){\makebox(0,0){$a$}}%
% STR 2 0 3 0
% 3 3105 700 3105 800 5 0
% $g$
\put(31.0500,-8.0000){\makebox(0,0){$g$}}%
% LINE 2 0 3 0
% 2 3405 750 3605 750
% 
\special{pn 8}%
\special{pa 3406 750}%
\special{pa 3606 750}%
\special{fp}%
% LINE 2 0 3 0
% 2 3405 850 3605 850
% 
\special{pn 8}%
\special{pa 3406 850}%
\special{pa 3606 850}%
\special{fp}%
% LINE 2 0 3 0
% 2 4100 800 4500 800
% 
\special{pn 8}%
\special{pa 4100 800}%
\special{pa 4500 800}%
\special{fp}%
% BOX 2 0 3 0
% 2 4500 650 4800 950
% 
\special{pn 8}%
\special{pa 4500 650}%
\special{pa 4800 650}%
\special{pa 4800 950}%
\special{pa 4500 950}%
\special{pa 4500 650}%
\special{fp}%
% LINE 2 0 3 0
% 2 4800 800 5200 800
% 
\special{pn 8}%
\special{pa 4800 800}%
\special{pa 5200 800}%
\special{fp}%
% STR 2 0 3 0
% 3 3900 700 3900 800 5 0
% $a$
\put(39.0000,-8.0000){\makebox(0,0){$a$}}%
% STR 2 0 3 0
% 3 5400 700 5400 800 5 0
% $g$
\put(54.0000,-8.0000){\makebox(0,0){$g$}}%
% LINE 1 0 3 0
% 2 2655 450 2655 1150
% 
\special{pn 13}%
\special{pa 2656 450}%
\special{pa 2656 1150}%
\special{fp}%
% LINE 1 0 3 0
% 2 4355 450 4355 1150
% 
\special{pn 13}%
\special{pa 4356 450}%
\special{pa 4356 1150}%
\special{fp}%
% STR 2 0 3 0
% 3 2360 700 2360 800 5 0
% $\fb$
\put(23.6000,-8.0000){\makebox(0,0){$\fb$}}%
% STR 2 0 3 0
% 3 4650 700 4650 800 5 0
% $\cF$
\put(46.5000,-8.0000){\makebox(0,0){$\cF$}}%
\end{picture}%
 \end{center}
 \caption[subsection]{$\langle \fb(a) | g\rangle=\langle a | \cF(g)\rangle$
  \ \ \ eq.(\ref{deff})}
 \label{fig:coalg3}
\end{figure}
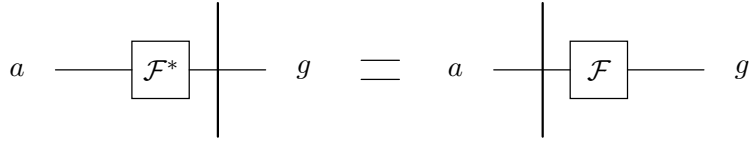
The graphs in both sides of the equations represent the 
$\C$ valued pairings. 
The arrow on the dashed line in Figure \ref{fig:coalg1} 
denotes the orientation of the operations in both sides. 
The $m$ is defined so that the pairing is invariant 
when the $|$ on the right hand side of Figure \ref{fig:coalg1} 
is moved to the left. Then the $\mm$ is $m$ on the left of $|$, 
and it becomes $\tri$ on the right of $|$. 
Similarly, in Figure \ref{fig:coalg2}, 
the $\D$ on the left of $|$ becomes $\Q$ on the right and 
the $\Q$ is $\D$ when is transferred to the left. 
The situation is similar for $\f$ and $\F$ (Figure \ref{fig:coalg3}). 

\begin{figure}[h]
 \begin{center}
%WinTpicVersion3.08
\unitlength 0.1in
\begin{picture}( 54.1500, 22.0000)(  8.6000,-34.0000)
% CIRCLE 2 0 3 0
% 4 2000 3000 2000 3200 2000 3200 2000 2800
% 
\special{pn 8}%
\special{ar 2000 3000 200 200  4.7123890 6.2831853}%
\special{ar 2000 3000 200 200  0.0000000 1.5707963}%
% STR 2 0 3 0
% 3 1000 2700 1000 2800 5 0
% $a$
\put(10.0000,-28.0000){\makebox(0,0){$a$}}%
% STR 2 0 3 0
% 3 1000 3100 1000 3200 5 0
% $b$
\put(10.0000,-32.0000){\makebox(0,0){$b$}}%
% STR 2 0 3 0
% 3 3205 2900 3205 3000 5 0
% $g$
\put(32.0500,-30.0000){\makebox(0,0){$g$}}%
% LINE 2 0 3 0
% 2 3605 2950 3805 2950
% 
\special{pn 8}%
\special{pa 3606 2950}%
\special{pa 3806 2950}%
\special{fp}%
% LINE 2 0 3 0
% 2 3605 3050 3805 3050
% 
\special{pn 8}%
\special{pa 3606 3050}%
\special{pa 3806 3050}%
\special{fp}%
% STR 2 0 3 0
% 3 4200 2695 4200 2795 5 0
% $a$
\put(42.0000,-27.9500){\makebox(0,0){$a$}}%
% STR 2 0 3 0
% 3 4200 3095 4200 3195 5 0
% $b$
\put(42.0000,-31.9500){\makebox(0,0){$b$}}%
% LINE 1 0 3 0
% 2 2760 2600 2760 3400
% 
\special{pn 13}%
\special{pa 2760 2600}%
\special{pa 2760 3400}%
\special{fp}%
% LINE 2 0 3 0
% 2 1895 2800 1995 2800
% 
\special{pn 8}%
\special{pa 1896 2800}%
\special{pa 1996 2800}%
\special{fp}%
% LINE 2 0 3 0
% 2 1895 3200 1995 3200
% 
\special{pn 8}%
\special{pa 1896 3200}%
\special{pa 1996 3200}%
\special{fp}%
% BOX 2 0 3 0
% 2 1595 2650 1895 2950
% 
\special{pn 8}%
\special{pa 1596 2650}%
\special{pa 1896 2650}%
\special{pa 1896 2950}%
\special{pa 1596 2950}%
\special{pa 1596 2650}%
\special{fp}%
% STR 2 0 3 0
% 3 1750 2700 1750 2800 5 0
% $\fb$
\put(17.5000,-28.0000){\makebox(0,0){$\fb$}}%
% BOX 2 0 3 0
% 2 1595 3050 1895 3350
% 
\special{pn 8}%
\special{pa 1596 3050}%
\special{pa 1896 3050}%
\special{pa 1896 3350}%
\special{pa 1596 3350}%
\special{pa 1596 3050}%
\special{fp}%
% STR 2 0 3 0
% 3 1750 3100 1750 3200 5 0
% $\fb$
\put(17.5000,-32.0000){\makebox(0,0){$\fb$}}%
% LINE 2 0 3 0
% 2 2200 3000 3000 3000
% 
\special{pn 8}%
\special{pa 2200 3000}%
\special{pa 3000 3000}%
\special{fp}%
% CIRCLE 2 0 3 0
% 4 5205 3000 5205 3200 5205 3200 5205 2800
% 
\special{pn 8}%
\special{ar 5206 3000 200 200  4.7123890 6.2831853}%
\special{ar 5206 3000 200 200  0.0000000 1.5707963}%
% STR 2 0 3 0
% 3 6410 2900 6410 3000 5 0
% $g$
\put(64.1000,-30.0000){\makebox(0,0){$g$}}%
% LINE 1 0 3 0
% 2 4650 2600 4650 3400
% 
\special{pn 13}%
\special{pa 4650 2600}%
\special{pa 4650 3400}%
\special{fp}%
% LINE 2 0 3 0
% 2 5100 2800 5200 2800
% 
\special{pn 8}%
\special{pa 5100 2800}%
\special{pa 5200 2800}%
\special{fp}%
% LINE 2 0 3 0
% 2 5100 3200 5200 3200
% 
\special{pn 8}%
\special{pa 5100 3200}%
\special{pa 5200 3200}%
\special{fp}%
% BOX 2 0 3 0
% 2 4800 2650 5100 2950
% 
\special{pn 8}%
\special{pa 4800 2650}%
\special{pa 5100 2650}%
\special{pa 5100 2950}%
\special{pa 4800 2950}%
\special{pa 4800 2650}%
\special{fp}%
% STR 2 0 3 0
% 3 4955 2700 4955 2800 5 0
% $\cF$
\put(49.5500,-28.0000){\makebox(0,0){$\cF$}}%
% BOX 2 0 3 0
% 2 4800 3050 5100 3350
% 
\special{pn 8}%
\special{pa 4800 3050}%
\special{pa 5100 3050}%
\special{pa 5100 3350}%
\special{pa 4800 3350}%
\special{pa 4800 3050}%
\special{fp}%
% STR 2 0 3 0
% 3 4955 3100 4955 3200 5 0
% $\cF$
\put(49.5500,-32.0000){\makebox(0,0){$\cF$}}%
% LINE 2 0 3 0
% 2 5405 3000 6205 3000
% 
\special{pn 8}%
\special{pa 5406 3000}%
\special{pa 6206 3000}%
\special{fp}%
% LINE 2 0 3 0
% 2 1200 2800 1600 2800
% 
\special{pn 8}%
\special{pa 1200 2800}%
\special{pa 1600 2800}%
\special{fp}%
% LINE 2 0 3 0
% 2 1200 3200 1600 3200
% 
\special{pn 8}%
\special{pa 1200 3200}%
\special{pa 1600 3200}%
\special{fp}%
% LINE 2 0 3 0
% 2 4400 2800 4800 2800
% 
\special{pn 8}%
\special{pa 4400 2800}%
\special{pa 4800 2800}%
\special{fp}%
% LINE 2 0 3 0
% 2 4400 3200 4800 3200
% 
\special{pn 8}%
\special{pa 4400 3200}%
\special{pa 4800 3200}%
\special{fp}%
% CIRCLE 2 0 3 0
% 4 1995 1600 1995 1800 1995 1800 1995 1400
% 
\special{pn 8}%
\special{ar 1996 1600 200 200  4.7123890 6.2831853}%
\special{ar 1996 1600 200 200  0.0000000 1.5707963}%
% STR 2 0 3 0
% 3 995 1300 995 1400 5 0
% $a$
\put(9.9500,-14.0000){\makebox(0,0){$a$}}%
% STR 2 0 3 0
% 3 995 1700 995 1800 5 0
% $b$
\put(9.9500,-18.0000){\makebox(0,0){$b$}}%
% STR 2 0 3 0
% 3 3200 1500 3200 1600 5 0
% $g$
\put(32.0000,-16.0000){\makebox(0,0){$g$}}%
% LINE 2 0 3 0
% 2 3600 1550 3800 1550
% 
\special{pn 8}%
\special{pa 3600 1550}%
\special{pa 3800 1550}%
\special{fp}%
% LINE 2 0 3 0
% 2 3600 1650 3800 1650
% 
\special{pn 8}%
\special{pa 3600 1650}%
\special{pa 3800 1650}%
\special{fp}%
% STR 2 0 3 0
% 3 4195 1295 4195 1395 5 0
% $a$
\put(41.9500,-13.9500){\makebox(0,0){$a$}}%
% STR 2 0 3 0
% 3 4195 1695 4195 1795 5 0
% $b$
\put(41.9500,-17.9500){\makebox(0,0){$b$}}%
% LINE 1 0 3 0
% 2 2755 1200 2755 2000
% 
\special{pn 13}%
\special{pa 2756 1200}%
\special{pa 2756 2000}%
\special{fp}%
% STR 2 0 3 0
% 3 6405 1500 6405 1600 5 0
% $g$
\put(64.0500,-16.0000){\makebox(0,0){$g$}}%
% LINE 1 0 3 0
% 2 4645 1200 4645 2000
% 
\special{pn 13}%
\special{pa 4646 1200}%
\special{pa 4646 2000}%
\special{fp}%
% LINE 2 0 3 0
% 2 1200 1400 2000 1400
% 
\special{pn 8}%
\special{pa 1200 1400}%
\special{pa 2000 1400}%
\special{fp}%
% LINE 2 0 3 0
% 2 1200 1800 2000 1800
% 
\special{pn 8}%
\special{pa 1200 1800}%
\special{pa 2000 1800}%
\special{fp}%
% LINE 2 0 3 0
% 2 2200 1600 2300 1600
% 
\special{pn 8}%
\special{pa 2200 1600}%
\special{pa 2300 1600}%
\special{fp}%
% BOX 2 0 3 0
% 2 2300 1450 2600 1750
% 
\special{pn 8}%
\special{pa 2300 1450}%
\special{pa 2600 1450}%
\special{pa 2600 1750}%
\special{pa 2300 1750}%
\special{pa 2300 1450}%
\special{fp}%
% STR 2 0 3 0
% 3 2455 1500 2455 1600 5 0
% $\fb$
\put(24.5500,-16.0000){\makebox(0,0){$\fb$}}%
% LINE 2 0 3 0
% 2 2600 1600 3000 1600
% 
\special{pn 8}%
\special{pa 2600 1600}%
\special{pa 3000 1600}%
\special{fp}%
% CIRCLE 2 0 3 0
% 4 5195 1600 5195 1800 5195 1800 5195 1400
% 
\special{pn 8}%
\special{ar 5196 1600 200 200  4.7123890 6.2831853}%
\special{ar 5196 1600 200 200  0.0000000 1.5707963}%
% LINE 2 0 3 0
% 2 4400 1400 5200 1400
% 
\special{pn 8}%
\special{pa 4400 1400}%
\special{pa 5200 1400}%
\special{fp}%
% LINE 2 0 3 0
% 2 4400 1800 5200 1800
% 
\special{pn 8}%
\special{pa 4400 1800}%
\special{pa 5200 1800}%
\special{fp}%
% LINE 2 0 3 0
% 2 5400 1600 5500 1600
% 
\special{pn 8}%
\special{pa 5400 1600}%
\special{pa 5500 1600}%
\special{fp}%
% BOX 2 0 3 0
% 2 5500 1450 5800 1750
% 
\special{pn 8}%
\special{pa 5500 1450}%
\special{pa 5800 1450}%
\special{pa 5800 1750}%
\special{pa 5500 1750}%
\special{pa 5500 1450}%
\special{fp}%
% STR 2 0 3 0
% 3 5655 1500 5655 1600 5 0
% $\cF$
\put(56.5500,-16.0000){\makebox(0,0){$\cF$}}%
% LINE 2 0 3 0
% 2 5800 1600 6200 1600
% 
\special{pn 8}%
\special{pa 5800 1600}%
\special{pa 6200 1600}%
\special{fp}%
% LINE 2 0 3 0
% 2 5250 2200 5250 2400
% 
\special{pn 8}%
\special{pa 5250 2200}%
\special{pa 5250 2400}%
\special{fp}%
% LINE 2 0 3 0
% 2 5350 2200 5350 2400
% 
\special{pn 8}%
\special{pa 5350 2200}%
\special{pa 5350 2400}%
\special{fp}%
% LINE 2 0 3 0
% 2 2050 2200 2050 2400
% 
\special{pn 8}%
\special{pa 2050 2200}%
\special{pa 2050 2400}%
\special{fp}%
% LINE 2 0 3 0
% 2 2150 2200 2150 2400
% 
\special{pn 8}%
\special{pa 2150 2200}%
\special{pa 2150 2400}%
\special{fp}%
\end{picture}%
 \end{center}
 \caption[subsection]{$\langle \fb\cdot m(a\otimes b) | g\rangle
  = \langle m(\fb(a)\otimes \fb(b)) | g\rangle$\ \ \ eq.(\ref{Ff})}
 \label{fig:coalgex}
\end{figure}
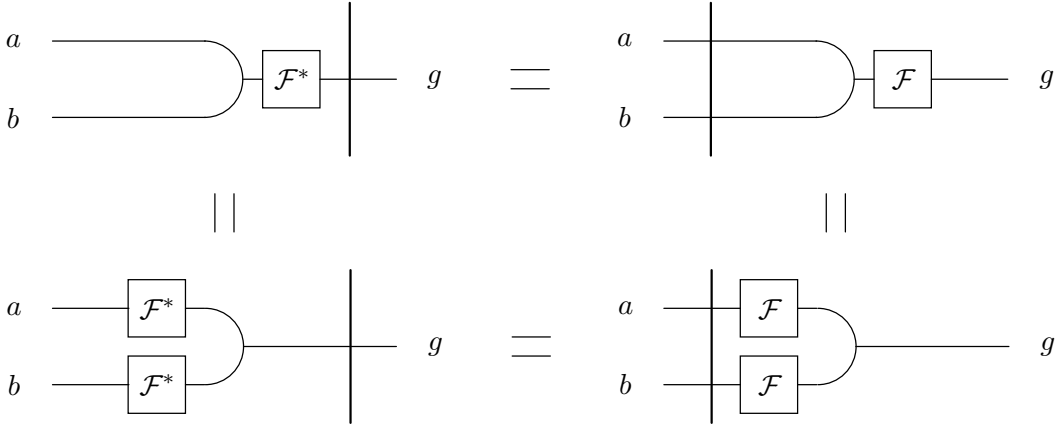

{}From the rewriting above, 
the following dualities can be understood naturally 
by using graphs; 
$\m$ is a coderivation vs. the $\delta$ is a derivation (\ref{QD}), 
the $\cF$ is a cohomomorphism vs. the $\fb$ is a homomorphism (\ref{Ff}), 
the $(\cH, \m)$ is an $A_\infty$-algebra vs. 
the $(C(\cH)^*, \delta)$ is a complex(\ref{QQDD}), 
and the $\cF$ is an $A_\infty$-morphism vs. 
the $\fb$ is $\delta$-equivariant (\ref{QFDf}). 
For instance, eq.(\ref{Ff}) is shown as Figure \ref{fig:coalgex}.  

 \subsection{Formal noncommutative supermanifolds}
\label{ssec:dual2}

In this subsection, 
we represent explicitly $m$, $\delta$ and $\fb$, which correspond to 
$\triangle$, $\m$ and $\cF$, respectively, and realize them geometrically 
on the algebra $C(\cH)^*$ dual to the $C(\cH)$. 
For the coassociative coproduct
\begin{equation*}
 \triangle(\eb_1\cdots \eb_n)=\sum_{k=1}^{n-1}
  (\eb_1\cdots \eb_k)
  \otimes (\eb_{k+1}\cdots\eb_n)\ , 
\end{equation*} 
the corresponding associative product $m$ 
defined in eq.(\ref{defm}) are written as 
\begin{equation}
 m((\eb^{i_k}\cdots\eb^{i_1})\otimes(\eb^{j_l}\cdots\eb^{j_1}) )
 =\eb^{j_l}\cdots\eb^{j_1}\eb^{i_k}\cdots\eb^{i_1}\ .
 \label{Apro}
\end{equation}
For $a=\sum_{k=1}^\infty a_{i_1\cdots i_k}\eb^{i_k}\cdots\eb^{i_1}$ and 
$b=\sum_{l=1}^\infty b_{j_1\cdots j_l}\eb^{j_l}\cdots\eb^{j_1}$, 
$m(a\otimes b)$ becomes 
\begin{equation*} 
 \begin{split}
 &m((\sum_{k=1}^\infty a_{i_1\cdots i_k}\eb^{i_k}\cdots\eb^{i_1})\otimes
  (\sum_{l=1}^\infty b_{j_1\cdots j_l}\eb^{j_l}\cdots\eb^{j_1}) )
 =\sum_{n}(a\cdot b)_{m_1\cdots m_n}\eb^{m_n}\cdots\eb^{m_1}\\
 &(a\cdot b)_{m_1\cdots m_n}=\sum_{p=1}^{n-1}
 \epsilon_{m_1\cdots m_n}^{i_1\cdots i_pj_1\cdots j_{n-p}}
 a_{i_1\cdots i_p}b_{j_1\cdots j_{n-p}}\ .
 \end{split}
\end{equation*}
It is easily seen that by the above definition of $m$, 
$(a\cdot b)_{m_1\cdots m_n}=
\langle m(a\otimes b)| \eb_{m_1}\cdots\eb_{m_n}\rangle
=\langle a\otimes b| \triangle(\eb_{m_1}\cdots\eb_{m_n})\ra$ holds. 

$a, b\in C(\cH)^*$ can be regarded as the polynomial functions 
on the graded vector space $\cH$. 
The dual basis $\{\eb^i\}$ is thought of as the coordinates 
(coordinate functions) of vector space $\cH$, 
though it is graded and noncommutative when the product structure is 
also considered. 
Hereafter we change the notation and denote $\eb^i$ by $\phi^i$. 
$C(\cH)^*$ is also replaced by $C(\phi)$ and its element is
represented as 
\begin{equation*}
  a(\phi)=\sum_{k=1}^\infty a_{i_1\cdots i_k}\phi^{i_k}\cdots\phi^{i_1}\ .
\end{equation*}
The coordinate functions $\{\phi^i\}$ will then be treated as fields 
in the sense of field theory. 
The pair of $\Z$-graded vector space $\cH$ 
and the algebra of formal power series 
of the coordinates $C(\phi)$ on $\cH$ is called 
{\it formal supermanifold} \cite{AKSZ,Ko1}. 
We call so, though this may be an infinitesimal neighborhood or germ of 
a more general global supermanifold. 
Though usually the term `super' indicates $\Z_2$-graded, 
we use it for $\Z$-graded object. 
The term `formal' is that in formal power series or in formal geometry. 
In our situation, 
the coordinates are associative but noncommutative, so 
this is a formal noncommutative supermanifold. 
On the other hand, 
the dual picture of an $L_\infty$-algebra is described by 
a formal (commutative) supermanifold. 
Namely, all the arguments in this paper can be translated into those 
for $L_\infty$-algebras by imposing the graded
commutativity upon $\{\phi^i\}$ 
corresponding to the graded commutativity in the 
$L_\infty$-algebra (Remark \ref{rem:cocomm}).

\vspace*{0.2cm}

 $\bullet$\ \ {\it coderivation}

\vspace*{0.2cm}

Next, for a coderivation $\m=\m_1+\m_2+\cdots$,  
\begin{equation}
 \m_k(\eb_{i_1}\cdots\eb_{i_n})=\sum_{p=1}^{n-k}
 (-1)^{\eb_{i_1}+\cdots+\eb_{i_p}}\eb_{i_1}\cdots\eb_{i_{p-1}}
 m_k(\eb_{i_p},\dots,\eb_{i_{p+k-1}})
 \eb_{i_{p+k}}\cdots\eb_{i_n} \ ,%\quad \eb_i\in\cH\ ,
 \label{Am}
\end{equation}
we construct $\delta$ which is the dual to $\m$. 
By the definition of $\delta$ (\ref{defD}), 
one sees that a derivation corresponding to the coderivation 
may be constructed separately for $k$. 
Let us express $m_k : \cH^{\otimes k}\to \cH$ as 
\begin{equation}
 m_k(\eb_{i_1},\dots,\eb_{i_k})=\eb_jc_{i_1\cdots i_k}^j\ ,\qquad 
 c_{i_1\cdots i_k}^j\in\C
 \label{lc}
\end{equation}
and then $\delta_k : C(\cH)^* \to C(\cH)^*$,   
\[
 \delta_k=\flpartial{\phi^j}c^j_{i_1\cdots i_k}\phi^{i_k}\cdots\phi^{i_1}
\]
is a derivation. Here we identify the span of the 
coordinate $\{\phi^i\}$ with 
$\cH^*$ and replace $\eb^i$ to $\phi^i$. The derivation $\delta$ is 
constructed as 
\begin{equation}
 \delta=\delta_1+\delta_2+\cdots
 =\sum_{k=1}^\infty\flpartial{\phi^j}c^j_{i_1\cdots i_k}
 \phi^{i_k}\cdots\phi^{i_1}\ .
 \label{dualD}
\end{equation}
It is regarded as an 
(odd) {\it formal vector field} on the formal noncommutative pointed 
supermanifold. The formal manifold with such $\delta$ is called 
{\it $Q$-manifold} in \cite{AKSZ}.
\footnote{This $Q$ does not correspond to the BRST operator $Q$ in the 
body of this paper but $\delta$. 
$\delta$ in this paper is written as $Q$ in \cite{AKSZ}. } 
Note that the condition that $\m_k$ is a coderivation is replaced to that 
$\delta_k$ satisfies the Leibniz rule on the polynomials of $\phi^i$'s. 
For instance the operation of $\delta$ on $(\phi^3\phi^2\phi^1)$ is 
defined as 
$$
\delta_k(\phi^3\phi^2\phi^1)=\phi^3\phi^2c^1_{i_1\cdots i_k}
 \phi^{i_k}\cdots\phi^{i_1}+(-1)^{\eb_1}\phi^3c^2_{i_1\cdots i_k}
 \phi^{i_k}\cdots\phi^{i_1}\phi^1+(-1)^{\eb_1+\eb_2}c^3_{i_1\cdots i_k}
 \phi^{i_k}\cdots\phi^{i_1}\phi^2\phi^1\ .
$$ 
The sign arises when the $\delta$ with degree one 
passes through some elements 
which have their degree. 

In field theories equipped with a classical BV-structure, 
this $\delta$ is just the BV-BRST transformation as 
seen in subsection \ref{ssec:BV}. 
Note that in this case $\{\phi^i\}$ consists of both fields and antifields.

\begin{rem}
When $\m$ satisfies $\m\cdot \m=0$, we have relations between $\m_k$. 
Rewriting $m_k$ using eq.(\ref{lc}) yields 
relations between $c^j_{i_1\cdots i_k}$. 
On the other hand, in the dual language, 
the condition $\m\cdot \m=0$ is equivalent to $\delta\cdot\delta=0$. 
Calculating $\delta\cdot\delta$ and concentrating on the terms of 
$n$ powers of $\phi^i$ lead to 
\begin{equation*}
 \begin{split}
 &\sum_{k+l=n+1}\delta_l\cdot\delta_k=
 \l(\flpartial{\phi^i}c^i_{i_1\cdots i_k}\phi^{i_k}\cdots\phi^{i_1}\r)
 \flpartial{\phi^{j}}c^{j}_{j_1\cdots j_l}\phi^{j_l}\cdots\phi^{j_1}\\
 &\qquad=\flpartial{\phi^i}
 \sum_{k+l=n+1}\sum_{m=1}^k(-1)^{\eb_{i_1}+\cdots+\eb_{i_{m-1}}}
 c^i_{i_1\cdots i_k}c^{i_m}_{j_1\cdots j_l}
 \phi^{i_k}\cdots\phi^{i_{m+1}}\l(
 \phi^{j_l}\cdots\phi^{j_1}\r)\phi^{i_{m-1}}\cdots\phi^{i_1}\ .
 \end{split}
\end{equation*}
The coefficient of $\phi^n\cdots\phi^1$ then reads 
\begin{equation}
0= \sum_{\substack{k+l=n+1\\m=0,\dots,k-1}}
 {(-1)^{\eb_1+\cdots+\eb_m}}
 c^i_{1\cdots m,i_m,m+l+1\cdots n}c^{i_m}_{m+1\cdots m+l}\ .
\end{equation}
This is exactly the relation $\m\cdot\m=0$ (or eq.(\ref{ainf})) rewritten 
with $\{c^i_{i_1\cdots i_k}\}$. 
 \label{rem:Qmanifold}
\end{rem} 

\vspace*{0.2cm}

$\bullet$\ \ {\it cohomomorphism}

\vspace*{0.2cm}

In the terminology of the formal supermanifold, a homomorphism corresponding 
to a cohomomorphism $\cF$ are constructed as follows. 
For two graded vector spaces $\cH, \cH'$ and 
a cohomomorphism $\cF: C(\cH)\to C(\cH')$ defined in eq.(\ref{cohom}), 
let us now express $f_n$ as 
\[
 f_n(\eb_{i_1},\dots,\eb_{i_n})=\eb_{j'}f^{j'}_{i_1\cdots i_n}\ ,\qquad 
 f^{j'}_{i_1\cdots i_n}\in\C\ 
\]
for $n\ge 0$. The homomorphism $\fb$ actually gives the pullback 
from $C(\cH')^*$, the formal power series ring on $\cH'$, to $C(\cH)^*$. 
Moreover $\fb : C(\cH')^*\to C(\cH)^*$ is induced from $\cF_*$ below 
\begin{equation}
 \begin{array}{cccc}
  \cF_* :&\cH&\to & \cH' \\
     & \phi&\mapsto& \phi'=\cF_*(\phi)
 \end{array}\ ,\quad 
 \phi^{j'}=\cF_*^{j'}(\phi)
 =f^{j'}+f^{j'}_i\phi^i+f^{j'}_{i_1i_2}\phi^{i_2}\phi^{i_1}
 +\cdots+f^{j'}_{i_1\cdots i_n}\phi^{i_n}\cdots\phi^{i_1}+\cdots\ ,
 \label{coordtransf}
\end{equation}
where $\{\phi^i\}$ and $\{\phi^{i'}\}$ are the coordinates 
on $\cH$ and $\cH'$, respectively. 
Namely, for an element $a(\phi'):=\sum_{k=1}^\infty 
a_{i'_1\cdots i'_k}\phi^{i'_k}\cdots\phi^{i'_1}\in C(\cH')^*$, 
$\fb(a(\phi'))=a(\cF_*(\phi))$ holds. 
One can see that the cohomomorphism $\cF$ is, in the dual geometric 
picture, a nonlinear map $\cF_*$ from a formal supermanifold $\cH$ to 
$\cH'$. If $f^{j'}=0$, the $\cF_*$ preserves the origin.

\vspace*{0.2cm}

$\bullet$\ \ {\it $A_\infty$-morphism}

\vspace*{0.2cm}

The condition that this $\cF$ is an $A_\infty$-morphism 
is equivalent to the statement 
that this map $\cF_*$ between two formal supermanifolds 
is compatible with the actions of $\delta$ and $\delta'$ on both sides, \ie 
$\cF_*$ is a morphism between $Q$-manifolds. 
For any $a(\phi')\in C(\cH')^*$, the condition is 
\begin{equation}
 \fb\delta'(a(\phi'))=\delta\fb a(\phi')\ ,
\end{equation}
and is written explicitly as 
\begin{equation*}
 \fb\l( a(\phi')\flpartial{\phi^{j'}}c^{j'}(\phi') \r)=
 a(\cF_*(\phi))
 \flpartial{\phi^j}c^j(\phi)\ ,
\end{equation*}
where we expressed $\delta=\flpartial{\phi^j}c^j(\phi)$.  
Because $a(\cF_*(\phi))\flpartial{\phi^j}c^j(\phi)=
\fb\l(a(\phi')\flpartial{\phi^{j'}}\r)\flpart{\phi^{j'}}{\phi^j}c^j(\phi)$ 
in the right hand side, we get 
\begin{equation}
 \fb\l(c^{j'}(\phi')\r)=\flpart{\phi^{j'}}{\phi^j}c^j(\phi)\ .
 \label{cc'}
\end{equation}
We can see that when $\delta$ and $\fb$ are given and $\fb$ has its inverse, 
then $\delta'$ is induced as 
$c^{j'}(\phi')=(\cF^{-1})^*\l(\flpart{\phi^{j'}}{\phi^j}c^j(\phi)\r)$. 
This is the dual description of $\m'=\cF\m\cF^{-1}$ 
when $\cF$ is an $A_\infty$-isomorphism. 

\vspace*{0.2cm}

The dual description of cyclicity will be discussed in 
the next section.

 \subsection{Superfield and mixed description}
\label{ssec:superfield}

\begin{defn}[superfield]
For basis $\{\eb_i\}$ of $\cH$ and the dual basis $\{\phi^i\}$, 
we define $\Phi:=\eb_i\phi^i\in\cH\otimes\cH^*$ 
and call it the {\it superfield}. 
 \label{defn:superfield}
\end{defn}
Since the degree of $\phi^i$ is minus the degree of $\eb_i$, 
the superfield $\Phi$ has degree zero. 
The roles are in fact similar 
to those of superfield in supersymmetric field theory, though 
some more extended notions are included. 
Note that the term `superfield' does not mean 
that $\Phi$ is a field (function) on our supermanifold. 
In the case of string field theory it is called the string field. 

It is useful to incorporate the coalgebra description and its dual. 
The multilinear map $m_k:\cH^{\otimes k}\to\cH$ is extended 
to operation on $\Phi^{\otimes k}$ as 
\begin{equation*}
 m_k(\Phi,\dots,\Phi)=\eb_j 
 c^j_{i_1\cdots i_k}\phi^{i_k}\cdots\phi^{i_1}\ .
\end{equation*}
In other words, 
$m_k(\Phi,\dots,\Phi)\in\cH\otimes(\cH^*)^{\otimes k}$ 
is defined by this equation. 
Since $\eb_j$ is identified with $\flpartial{\phi^j}$, 
$m_k(\Phi,\dots,\Phi)$ is identified with $\delta_k$. 
Similarly, for multilinear map $f_k:\cH^{\otimes k}\to\cH'$ 
which defines a cohomomorphism or $A_\infty$-morphism $\cF$, 
its operation on superfields is defined by 
\begin{equation*}
 f_k(\Phi,\dots,\Phi)=e_{j'}f^{j'}_{i_1\cdots i_k}
 \phi^{i_k}\cdots\phi^{i_1}\ .
\end{equation*}
$\cF_*$ in eq.(\ref{coordtransf}) is then represented by 
\begin{equation*}
 \Phi'=\cF_*(\Phi)=f_1(\Phi)+f_2(\Phi,\Phi)+f_3(\Phi,\Phi,\Phi)+\cdots\ .
\end{equation*}
Polynomial functions $C(\phi)$ also have superfield description. 
For $a_{i_1\cdots i_k}\phi^{i_k}\cdots\phi^{i_1}\in C(\phi)$, 
let us define a multilinear map $a_k:\cH^{\otimes k}\to\C$ such that 
$a_k(\eb_{i_1},\dots,\eb_{i_k})=a_{i_1\cdots i_k}$. 
Any element in $C(\phi)$ is then expressed as 
\begin{equation*}
 a_k(\Phi,\dots,\Phi)=a_{i_1\cdots i_k}\phi^{i_k}\cdots\phi^{i_1}\ . 
\end{equation*}
In the next section 
we shall also consider cyclic functions, which is defined by 
imposing cyclic condition upon $a_{i_1\cdots i_k}\in\C$ 
(Definition \ref{defn:func}). 
Namely, they are described by a cyclic multilinear map 
$a_k:\cH^{\otimes k}\to\C$. 

Thus, superfield $\Phi$ can be thought of as an expression of 
general element in $\cH$ 
and also play the role of complete system of $\cH$. 
One can see that in this superfield description 
one can discuss the dual side without the indices. 

Instead of superfield $\Phi\in\cH\otimes\cH^*$, 
one can consider more generally $\cH$ 
over associative graded algebra $C(\ti\phi)$ generated freely by 
$\{\ti\phi^i\}$, 
though we do not use it later in this paper. 
An element of the $C(\ti\phi)$-module $\cH$ is represented as 
$A(\ti\phi)=\eb_i A^i(\ti\phi)$ where $A^i(\ti\phi)\in C(\ti\phi)$. 
The degree of $A(\ti\phi)$ is just the sum
$\deg(\eb_i)+\deg(A^i(\ti\phi))$. 
$C(\cH)$ is then extended to $C(\cH\otimes C(\ti\phi))$. 
The coproduct is defined just in a similar way as that on $C(\cH)$. 
Furthermore, the operations of any coderivation 
$\m: C(\cH)\to C(\cH)$ and any cohomomorphism $\cF: C(\cH)\to C(\cH')$ 
are naturally extended to operations on $C(\ti\phi)$-module. 
It is convenient to rewrite each element in $C(\cH\otimes C(\ti\phi))$ 
into that of the form in $C(\cH)\otimes C(\ti\phi)$ as follows,  
\begin{equation*}
 \eb_{i_1}A_1^{i_1}(\ti\phi)\cdots \eb_{i_n}A_n^{i_n}(\ti\phi)
 =(-1)^{\sum_{k=1}^{n-1} A_k^{i_k}(\sum_{l=k+1}^n (A_l^{i_l}+\eb_{i_l}))}
 \eb_{i_1}\cdots\eb_{i_n}(A_n^{i_n}(\ti\phi)\cdots A_1^{i_1}(\ti\phi))\ , 
\end{equation*}
where $A_k^{i_k}$ in the sign factor in the left hand side 
indicates $\deg(A_k^{i_k}(\ti\phi))$. 
The sign factor is determined as the Kostul sign that 
arises when $A_{i_{k}}^{i_k}$ passes through 
$\eb_{i_{k+1}}\cdots\eb_{i_n}A_{i_{n}}^{i_n}\cdots A_{i_{k+1}}^{i_{k+1}}$ 
for $k=n-1,n-2,\dots, 1$. 
The expression of the right hand side enables us to extend any
operation on $C(\cH)$ to that on $C(\cH\otimes C(\ti\phi))$. 
Though we related elements in $C(\cH\otimes C(\ti\phi))$ 
to that in $C(\cH)\otimes C(\ti\phi)$, 
we should not consider all elements in $C(\cH)\otimes C(\ti\phi)$. 
Only those that comes from $C(\cH\otimes C(\ti\phi))$ have 
a tree structure and can have an $A_\infty$-structure, and so on.

%section 4

 \section{Odd symplectic geometry on formal noncommutative supermanifolds}
\label{sec:sym}

A constant symplectic structure of a cyclic $A_\infty$-algebra 
is a constant symplectic structure on a 
formal noncommutative supermanifold. 
Such noncommutative symplectic geometry is discussed 
in \cite{Ko3}. 
In this section we shall discuss nonconstant symplectic structures. 
They can be defined naturally from the physics of open strings. 
After considering the constant ones 
in subsection \ref{ssec:csympoi}, 
we shall first define general nonconstant ones 
in subsection \ref{ssec:sympoi}. 
Such geometry is discussed by 
A.~Schwarz \cite{Sch} in case of graded commutative supermanifolds, 
where the existence of the Darboux theorem is known. 
In \cite{Sch} the body (even coordinates parts) of the supermanifolds 
are treated global. 
However, in noncommutative case, local properties 
have never been discussed enough. 
Therefore we shall discuss the local properties of nonconstant symplectic 
and Poisson structures. 
We show a key lemma (Lemma \ref{lem:poincare}), the Poincar\'e's 
lemma on the formal noncommutative supermanifolds. 
Due to the lemma, we examine the properties of 
symplectic diffeomorphisms in subsection \ref{ssec:diffeo}, 
and show the Darboux theorem on the formal noncommutative supermanifolds 
(Theorem \ref{thm:Darboux}) in subsection \ref{ssec:Darboux}. 
The study of the formal noncommutative symplectic supergeometry 
is directly related to the notion of cyclic $A_\infty$-algebras. 
We look back over cyclic $A_\infty$-algebras from these dual pictures in 
subsection \ref{ssec:cycAinftyre}.

 \subsection{The constant symplectic structures}
\label{ssec:csympoi}

Consider a $\Z$-graded vector space $\cH$ equipped with a constant 
odd symplectic structure $\omega:\cH\otimes\cH\to\C$ 
%$\omega(\eb_i,\eb_j)=\omega_{ij}$ 
in Definition \ref{defn:csym}. 
For $\{\eb_i\}$ basis of $\cH$, denote 
\begin{equation*}
 \omega_{ij} := \omega(\eb_i,\eb_j)\ . 
\end{equation*}
It defines a symplectic structure on formal noncommutative supermanifold 
of $\cH$. Its inverse is defined by 
\begin{equation*} 
 \omega_{ij}\omega^{jk}=\omega^{kj}\omega_{ji}=\delta_i^k
\end{equation*}
and it defines an odd Poisson structure as seen below. 
The graded antisymmetry implies 
\begin{equation*}
 \omega_{ji}=-\omega_{ij}\ ,\qquad 
\omega^{ji}=-\omega^{ij}=-(-1)^{(i+1)(j+1)}\omega^{ij}\ .
\end{equation*}
As above, in this section we often denote $\deg(\eb_i)$ simply by 
$i$ instead of $\eb_i$. 
The existence of the nondegenerate symplectic structure 
depends on the structure of $\cH$ but 
is natural in the context of the BV-formalism 
in field theory.  

Next, we define an algebra of functions to construct an 
odd Poisson algebra. 
In the previous section we considered an associative noncommutative 
polynomial algebra $C(\cH)^*=C(\phi)$ on a formal noncommutative 
supermanifold. In this section, we need to 
consider additional structure on $C(\phi)$. 
The idea of the following definition is that explained in subsection 
\ref{ssec:Introsym}. 
\begin{defn}[Function on a formal noncommutative supermanifold]
For $\phi^i$ the dual coordinate of $\eb_i$, 
let us consider an element 
$\ov{k}a_{i_1\cdots i_k}\phi^{i_k}\cdots\phi^{i_1}\in C(\phi)$ 
whose coefficient $a_{i_1\cdots i_k}$ is cyclic, 
\begin{equation*}
 a_{i_1\cdots i_k}
 =(-1)^{(i_k+\cdots+i_2)i_1}a_{i_2\cdots i_k i_1}\ .
\end{equation*}
Such elements have the following property 
\begin{equation*}
 \begin{split}
 \ov{k}a_{i_1\cdots i_k}\phi^{i_k}\cdots\phi^{i_1}
 &=\ov{k}a_{i_2\cdots i_k i_1}\phi^{i_1}\phi^{i_k}\cdots\phi^{i_2}\\
 &=(-1)^{(i_k+\cdots+i_2)i_1}
 \ov{k}a_{i_1\cdots i_k}\phi^{i_1}\phi^{i_k}\cdots\phi^{i_2}\ .
 \end{split}
\end{equation*}
That is, the coordinates $\phi^i$ are regarded cyclic. 
We denote by $C(\phi)_c$ the subgroup of $C(\phi)$. 
We also consider the free tensor algebra of $C(\phi)_c$ 
and denote it by $TC(\phi)_c$. 
The free tensor product is defined to be graded commutative and 
denoted by $\bullet$. For instance 
for $A, B\in C(\phi)_c$, 
$B\bullet A=(-1)^{AB}A\bullet B$ holds. 
 \label{defn:func}
\end{defn}
The cyclic symmetry for $C(\phi)_c$ is just the property of open string disk 
as explained in subsection \ref{ssec:Introsym}. 
This is similar to a `trace' in the terminology of noncommutative 
geometry \cite{connes}. 
The restriction $C(\phi)\to C(\phi)_c$ 
(given by cyclic symmetrization) 
is regarded as a trace 
which keeps all other informations of $C(\phi)$. 
Considering $TC(\phi)_c$ are also natural 
from the viewpoints of the BV-formalism. 
As an element in $C(\phi)_c$ is regarded as an $S^1$ as in eq.(\ref{S1}), 
an element in $TC(\phi)_c$ can be characterized by a multiple copy of $S^1$. 
However we do not need $TC(\phi)_c$ essentially in this paper. 
The essential necessity of it is related to `homotopy algebras 
of quantum open strings' if it would be defined. 
\begin{defn}[Odd Poisson structure (Gerstenhaber structure)]
A {\it Gerstenhaber algebra} is an algebra equipped with 
degree zero associative product $\bullet$ and 
degree one bracket $(\ ,\ )$ 
satisfying the following equations: 
\begin{itemize}
 \item[(a)] ${\displaystyle (B, A)=-(-1)^{(A+1)(B+1)}(A, B)}$, 
 \item[(b)] ${\displaystyle (A, B\bullet C)=(A, B)\bullet C
 +(-1)^{(A+1)B}B\bullet (A, C)}$, 
 \item[(c)] ${\displaystyle (-1)^{(A+1)(C+1)}((A, B), C) + cyclic =0}$  
\end{itemize}
where $A, B$ and $C$ are elements of the algebra. 
$A, B$ and $C$ on $(-1)$ denote the degree of $A, B$ and $C$, 
respectively. 
\label{defn:Gersten}
\end{defn}
An odd Poisson algebra can then be constructed 
from the constant symplectic structure in Definition \ref{defn:csym}. 
The odd Poisson bracket can be defined so that it fits the picture in 
eq.(\ref{mini1}). 
\begin{defn}[Constant Poisson structure]
On $TC(\phi)_c$ let us consider a Poisson bracket written as follows, 
\begin{equation}
 (\ ,\ )=\flpartial{\phi^i}\omega^{ij}\frpartial{\phi^j}\ .
 \label{constpoisson}
\end{equation}
For two elements of $C(\phi)_c$, 
the bracket is defined explicitly as 
\begin{equation*}
 \begin{split}
 &\ov{k}\V_{i_1\cdots i_k}\phi^{i_k}\cdots\phi^{i_1}
 \flpartial{\phi^i}\omega^{ij}\frpartial{\phi^j}
  \ov{l}\V_{j_1\cdots j_l}\phi^{j_l}\cdots\phi^{j_1}\\
 &:=\ov{k+l-2}\V_{i_1\cdots i_k}\omega^{i_1j_l}\V_{j_1\cdots j_l}
 \l(\phi^{i_k}\cdots\phi^{i_2}\phi^{j_{l-1}}\cdots\phi^{j_1} +cyclic\r)\\
 &=\ov{k+l-2}\l(\V_{i i_{l}\cdots i_{k+l-2}}\omega^{ij}\V_{i_1\cdots i_{l-1} j}
 +cyclic\r)\phi^{i_{k+l-2}}\cdots\phi^{i_1}\ .
 \end{split}
\end{equation*}
In the second equality, the `$cyclic$' means that 
$\phi^{i_k}\cdots\phi^{i_2}\phi^{j_{l-1}}\cdots\phi^{j_1}$ is moved cyclic 
such as 
\begin{equation*}
 \begin{split}
 &\phi^{i_k}\cdots\phi^{i_2}\phi^{j_{l-1}}\cdots\phi^{j_1}\\
 &\lgraw 
 (-1)^{(i_k+\cdots+i_2+j_{l-1}+\cdots+j_2)j_1}
 \phi^{j_1}\phi^{i_k}\cdots\phi^{i_2}\phi^{j_{l-1}}\cdots\phi^{j_2}\\ 
 &\lgraw \cdots
 \end{split}
\end{equation*}
and these $(k+l-2)$ terms are then summed up in $\l(\cdots\r)$.  
In the third line, `$cyclic$' indicates that 
$i_1\cdots i_{k+l-2}$ is moved 
cyclic with appropriate sign and the resulting $(k+l-2)$ terms are 
then summed up. 
More explicitly, 
for $b_{i_1\cdots i_{k+l-2}}
:=\V_{i i_{l}\cdots i_{k+l-2}}\omega^{ij}\V_{i_1\cdots i_{l-1} j}\in\C$, 
the coefficient of 
$\phi^{i_{k+l-2}}\cdots\phi^{i_1}$ in the third line is written as 
\begin{equation*}
 \ov{k+l-2}\sum_{k'=1}^{k+l-2}
 (-1)^{(i_1+\cdots +i_{k'-1})(i_{k'}+\cdots +i_{k+l-2})}
 b_{i_{k'}\cdots i_{k+l-2}\cdot i_1\cdots i_{k'-1}}\ .
\end{equation*}
By definition, this coefficient is cyclic with respect to 
the indices $i_1\cdots i_{k+l-2}$. 
 \label{defn:cpoi}
\end{defn}
The cyclicity is a natural property of open string, where 
operators $\eb_{i_1}\cdots\eb_{i_{k+l-2}}$ are inserted 
on a boundary of an open string disk $S^1$. 

As in commutative case, this constant Poisson bracket 
actually satisfies the Jacobi identity. 
It follows from Lemma \ref{lem:Jacobi}, but it can also be shown 
by writing pictures as in subsection \ref{ssec:Introsym}. 
Consequently, $( C(\phi)_c, (\ ,\ ) )$ forms a Lie algebra. 
Moreover we can extend the bracket on $TC(\phi)_c$ naturally by using 
Definition \ref{defn:Gersten} (b). 
Thus, $(TC(\phi)_c,(\ ,\ ))$ forms a Gerstenhaber algebra. 
\begin{rem}[Additional structure]
Note that the property of Definition \ref{defn:Gersten}(b) 
holds also for $A\in C(\phi)_c$ and $B, C\in C(\phi)$ 
with the replacement of the product $\bullet$ with 
the usual associative product in $C(\phi)$. 
This implies that $(A,\ )$ and also $(\ ,A)$ act as a derivation on 
$C(\phi)$. 
Pictorially this fact can be seen as follows. For instance 
for $B\in C(\phi)$ and $A\in C(\phi)_c$, the operation of $(\ ,A)$ is 
figured as 
\begin{equation*}
 (\ ,A)=\begin{minipage}[c]{50mm}{\input{Ader.tex}}\end{minipage}
% (\ ,A)=\begin{minipage}[c]{50mm}{\includegraphics{Ader.eps}}\end{minipage}
\end{equation*}
where the double line \ 
%$\begin{minipage}[c]{9mm}{\input{omega.tex}}\end{minipage}$\ 
$\begin{minipage}[c]{9mm}{\includegraphics{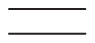}}\end{minipage}$\ 
corresponds to $\omega^{ij}$. 
$(B, A)$ is then expressed as 
\begin{equation*}
 \sum_k\pm\begin{minipage}[c]{85mm}{\input{Ader-l.tex}}\end{minipage}
 =\sum_k\pm\begin{minipage}[c]{50mm}{\input{Ader-r.tex}}\end{minipage}\ .
% \sum_k\pm
% \begin{minipage}[c]{85mm}{\includegraphics{Ader-l.eps}}\end{minipage}
% =\sum_k\pm
% \begin{minipage}[c]{50mm}{\includegraphics{Ader-r.eps}}\end{minipage}\ .
\end{equation*}
Namely, the interval of $B$ becomes also an interval after the insertion 
of $S^1$ of $A$. 
 \label{rem:Gerstender}
\end{rem}
We can also extend $\omega^{ij}$ to be nonconstant 
and define a natural class of nonconstant Poisson structures 
on a formal noncommutative supermanifold. 
They are defined as closed two-forms 
on the formal noncommutative supermanifold 
as seen in the next subsection.

 \subsection{The symplectic and Poisson structures}
\label{ssec:sympoi}

We can extend the constant odd Poisson structure to nonconstant ones. 
\begin{defn}[Covariant odd Poisson bracket]
On $TC(\phi)_c$ we define the following degree one nondegenerate 
bracket
\begin{equation}
 (A, B)=\sum_{ij,IJ}(-1)^{(B-j)J}\omega^{ij}_{JI}
 \l(\flpart{A}{\phi^i}\phi^I\frpart{B}{\phi^j}\phi^J\r)_c\ , 
 \label{poibra}
\end{equation}
where $I,J$ are multi-indices for polynomials of coordinates $\phi^i$ 
and $_c$ denotes the operation of cyclic symmetrization 
as in Definition \ref{defn:cpoi}. 
The nondegeneracy is equivalent to the nondegeneracy of 
$\omega_{ij,\emptyset\emptyset}$. 
The coefficients are required to satisfy 
\begin{equation*}
 \omega^{ji}_{IJ}=-(-1)^{(i+1)(j+1)+IJ}\omega^{ij}_{JI}\ 
\end{equation*}
so that the bracket has the property of Definition \ref{defn:Gersten} (a). 
We then call eq.(\ref{poibra}) an odd Poisson bracket 
if it satisfies the Jacobi identity (Definition \ref{defn:Gersten} (c)). 
 \label{defn:covpoisson}
\end{defn}
The geometric picture of this definition is figured in eq.(\ref{mini2}). 

We would like to know the condition that 
the Jacobi identity holds. 
We shall discuss it by translating these Poisson structure into 
symplectic side. 
For considering symplectic geometry on formal noncommutative 
supermanifolds, we need vector fields, differential forms, 
and so on. 
\begin{defn}[Hamiltonian vector field]
Using the properties of derivation in Definition \ref{defn:Gersten} (b) 
and Remark \ref{rem:Gerstender}, 
we define a Hamiltonian vector field. 
For a Hamiltonian $B\in C(\phi)_c$, we denote 
the corresponding Hamiltonian vector field by 
\begin{equation*}
 \delta_B=(\ ,B)\ .
\end{equation*}
Note that $\delta_B$ has degree $B+1$. 
The operation of $\delta_B: C(\phi)\to C(\phi)$ 
is defined by 
\begin{equation*}
 \begin{split}
 &\delta_B\l( a_{i_1\cdots i_n}\phi^{i_n}\cdots\phi^{i_1}\r) \\
 &\qquad 
 =\sum_{j=1}^n a_{i_1\cdots i_{j-1} j i_{j+1}\cdots i_n} 
 (-1)^{(B+1)(i_1+\cdots +i_{j-1})} 
\phi^{i_n}\cdots\phi^{j+1}
 \l((-1)^{(B-k)J}\omega^{jk}_{JI}
 \phi^I\frpart{B}{\phi^k}\phi^J\r)\phi^{j-1}\cdots\phi^{i_1}
 \end{split}
\end{equation*}
for $a_{i_1\cdots i_n}\phi^{i_n}\cdots\phi^{i_1}\in C(\phi)$. 
One can see that this is a natural extension of 
the derivation in Remark \ref{rem:Gerstender}. 
The operation $\delta_B: C(\phi)_c\to C(\phi)_c$ is then 
defined by cyclic symmetrizing the equation above. 
Namely, one obtains 
\begin{equation*}
 \delta_B A:=(A, B)
\end{equation*}
for $A\in C(\phi)_c$. 
\label{defn:Hamvect}
\end{defn}
\begin{defn}[Exterior derivative and differential form]
We first extend each coordinate $\phi^i$ by exterior derivative $d$ as 
\begin{equation*}
 0 \mapsto \phi^i \mapsto d\phi^i\mapsto 0\ .
\end{equation*}
For $d$ we introduce another degree $\sharp$. 
We assign $\sharp(d)=1$. 
We denote by $\varphi$ the coordinate $\phi$ or its exterior derivative 
$d\phi$. 
A differential form is then defined of the form 
\begin{equation*}
 a(\varphi):=
 \ov{k}a_{\mu_1\cdots\mu_k}\varphi^{\mu_k}\cdots\varphi^{\mu_1}\ .
\end{equation*}
The space of differential forms is then denoted by $C(\varphi)$. 
Denote $\sharp(\varphi_\mu)=:\sharp_\mu\ (=0\ \mbox{or}\ 1)$ and then 
the degree of differential form is simply the sum 
\begin{equation*}
 \sharp(\varphi^{\mu_k}\cdots\varphi^{\mu_1})
 =\sharp_{\mu_k}+\cdots +\sharp_{\mu_1}\ .
\end{equation*}
The space of $l$-th differential forms is denoted by 
$C(\varphi)^{\la l\ra}$. $C(\varphi)$ is then the direct sum of 
$C(\varphi)^{\la l\ra}$ for $l\ge 0$. 
We further define the subspace $C(\varphi)_c\subset C(\varphi)$ 
consisting of {\em cyclic elements} as follows. Recall that, 
for $a(\varphi)\in C(\varphi)$ above, 
$\varphi^{\mu_k}$ is $\phi^{i_k}$ or $d\phi^{i_k}$. 
We then say $a(\varphi)\in C(\varphi)$ is cyclic, 
$a(\varphi)\in C(\varphi)_c$, iff the coefficient satisfies  
\begin{equation*}
 a_{\mu_1\cdots\mu_k}
% =(-1)^{\eb_{i_1}(\eb_{i_2}+\cdots +\eb_{i_k})
 =(-1)^{i_1(i_2+\cdots + i_k)
 +\sharp_{\mu_1}(\sharp_{\mu_2}+\cdots +\sharp_{\mu_k})}
 a_{\mu_2\cdots\mu_k\mu_1}\ . 
\end{equation*}
Namely, we count the degree for $\cH$ and that for differential forms 
independently. 
The space of cyclic $l$-th differential forms is denoted by 
$C(\varphi)_c^{\la l\ra}$. 
The action of exterior derivative is naturally extended 
to that on $C(\varphi)^{\la l\ra}$ and $C(\varphi)_c^{\la l\ra}$. 
Especially, $d: C(\varphi)_c^{\la l\ra}\to C(\varphi)_c^{\la l+1\ra}$ 
is given by 
\begin{equation}
 d a(\varphi)= \ov{k}
 \sum_{i=1}^k (-1)^{\sharp_{\mu_1}+\cdots +\sharp_{\mu_{i-1}}}
 a_{\mu_1\cdots\mu_k}
 \varphi^{\mu_k}\cdots d\varphi^{\mu_i}\cdots\varphi^{\mu_1}\ .
 \label{d}
\end{equation}
Note that $d\varphi^{\mu_i}=0$ iff $\sharp_{\mu_i}=1$. 
In the expression above it is clear that $(d)^2=0$. 
Thus a complex $(C(\varphi)_c^{\la\cdot\ra}, d)$ is given. 
We also consider $(TC(\varphi)_c^{\la\cdot\ra}, d)$ with 
a natural extension. 
 \label{defn:ext}
\end{defn}
One can write eq.(\ref{d}) in a cyclic expression as 
\begin{equation*}
 d a(\varphi)= \l(\ov{k}
 \sum_{i=1}^k (-1)^{\sharp_{\mu_1}+\cdots +\sharp_{\mu_{i-1}}}
 a_{\mu_1\cdots{\check \mu_i}\cdots\mu_k}\r)
 \varphi^{\mu_k}\cdots \varphi^{\mu_1}\ .
\end{equation*}
Here ${\check \mu_i}$ denotes the index corresponding to 
$d^{-1}\varphi_{\mu_i}$. 
That is, 
$a_{\mu_1\cdots{\check \mu_i}\cdots\mu_k}=0$ if $\sharp_{\mu_i}=0$. 
Thus $d$ can be regarded as operation on the coefficient such as 
$(da)_{\mu_1\cdots\mu_k}
=\sum_{i=1}^k (-1)^{\sharp_{\mu_1}+\cdots +\sharp_{\mu_{i-1}}}
a_{\mu_1\cdots{\check \mu_i}\cdots\mu_k}$. 
In this expression the complex is realized as an analogue of 
Cech cohomology. 
Note that in noncommutative case 
$C(\varphi)^{\la l\ra}$ and $C(\varphi)_c^{\la l\ra}$ 
do not vanish trivially for any $l\ge 0$ 
even if the supermanifold is finite dimensional. 
\begin{lem}[Poincar\'e's lemma]
The cohomology with respect to $d$ is trivial. 
\label{lem:poincare}
\end{lem}
\begin{pf}
One can construct a homotopy operator 
$d^{-1}:C(\varphi)_c^{\la l+1\ra}\to C(\varphi)_c^{\la l\ra}$ 
which satisfies 
\begin{equation*}
 d d^{-1}+d^{-1}d =\Id\ .
\end{equation*}
It is constructed explicitly as 
\begin{equation*}
 d^{-1}a(\varphi)= \ov{k}
 \sum_{i=1}^k (-1)^{\sharp_{\mu_1}+\cdots +\sharp_{\mu_{i-1}}}
 a_{\mu_1\cdots\mu_k}
 \varphi^{\mu_k}\cdots (d^{-1}\varphi^{\mu_i})\cdots\varphi^{\mu_1}\ .
\end{equation*}
This fact completes the proof. 
\qed
\end{pf}

This local triviality of the deRham complex implies that this 
noncommutative geometry provides some natural extension of usual commutative 
geometry. 
We shall use this fact for later subsections. 
We shall see that 
the symplectic diffeomorphism discussed in the next subsection 
is related to the first deRham cohomology, whereas, 
the Darboux theorem in subsection \ref{ssec:Darboux} is equivalent to 
the triviality of the second deRham cohomology. 
\begin{rem}[Compatibility with transformations]
Let $(\cH,TC(\phi)_c)$ and $(\cH',TC(\phi')_c)$ be two formal 
noncommutative supermanifolds and 
$\cF^* :TC(\phi')_c\to TC(\phi)_c$ a pullback induced by a 
map $\phi^{i'}=f(\phi)
=f^{i'}_j+f^{i'}_j\phi^j+f^{i'}_{j_1j_2}\phi^{j_2}\phi^{j_1}+\cdots$. 
For $a(\phi')\in C(\phi')_c$, $\cF^*(a(\phi))$ is written explicitly as 
\begin{equation*}
 \cF^*(a(\phi'))=\l(\ a(f(\phi))\ \r)_c
\end{equation*}
where $_c$ is the operation of cyclic symmetrization 
as in Definition \ref{defn:cpoi}. 
An exterior derivative $d$ on $TC(\varphi')_c$ is related to 
that on $TC(\varphi)_c$ by 
\begin{equation*}
 \cF^*(d\phi^{i'})=d\phi^i\frpart{f(\phi)}{\phi^i}\quad 
(= \flpart{f(\phi)}{\phi^i}d\phi^i)\ .
\end{equation*}
By using this relation, 
$\cF^*:TC(\phi')_c\to TC(\phi)_c$ can be extended naturally to 
$\cF^*:TC(\varphi')_c\to TC(\varphi)$. 
$\cF^* : (TC(\varphi')_c, d)\to (TC(\varphi)_c, d)$ is then 
compatible with the exterior derivatives on both sides. That is, 
$$
 \cF^*d=d \cF^*\ .
$$
 \label{rem:dtransf}
\end{rem}

\begin{defn}[Odd symplectic structure]
The {\it odd symplectic structure} 
on a formal noncommutative supermanifold 
is defined by degree minus one closed $2$-form in $C(\varphi)_c$.
\footnote{It is not defined as an closed element of $TC(\varphi)_c$. 
The definition above is natural from the viewpoints of open string physics.
} 
We represent it as 
\begin{equation*}
 \omega=\sum_{ij,IJ}\omega_{ji,JI}\l(\phi^I d\phi^i\phi^J d\phi^j\r)_c\ ,
\end{equation*} 
where $I$ and $J$ are multi-indices. 
Note that, by the definition of $C(\varphi)_c$, 
the polynomial that consists of $\phi$'s and $d\phi$'s are cyclic. 
We then denote by $_c$ 
that $\phi^I d\phi^i\phi^J d\phi^j$ is defined to be cyclic symmetrized. 
{}From this, the coefficient satisfies 
\begin{equation*}
 \omega_{ij,IJ}=-\omega_{ji,JI}\ .
\end{equation*} 
\label{defn:covsym}
\end{defn}

In order to relate 
the odd Poisson structure and the odd symplectic structure, 
we need to define the contraction of forms with vector fields. 
\begin{defn}[Contraction]
For an element in $C(\varphi)^{\la k\ra}$, 
the contraction with $k$ vector fields are defined by 
\begin{equation*}
 \begin{split}
& \phi^{I_1}d\phi^{i_1}\phi^{I_2}d\phi^{i_2}\cdots d\phi^{i_k}\phi^{I_{k+1}}
 (\delta_{A_1},\dots,\delta_{A_k})\\
&\ \ =\phi^{I_1}(-1)^{(A_1+1)(I_2+i_2+\cdots +i_k+I_{k+1})}
(d\phi^{i_1},\delta_{A_1})\phi^{I_2}
(-1)^{(A_2+1)(I_3\cdots +i_k+I_{k+1})}
(d\phi^{i_2},\delta_{A_2})\\
&\hspace*{3.0cm}\cdots\cdots
(-1)^{(A_k+1) I_{k+1}}
(d\phi^{i_k},\delta_{A_k})\phi^{I_{k+1}}\ .
 \end{split}
\end{equation*}
When we write $\delta_A=\flpartial{\phi^i}A^i(\phi)$, 
$(d\phi^i,\delta_{A})=A^i(\phi)$. 

The contraction of an element in $C(\varphi)_c^{\la k\ra}$ with 
$k$ vector fields follows from the formula above. 
Any element in $C(\varphi)_c^{\la k\ra}$ is written in the form 
$a(\varphi):=
a_{I_k,i_k,\dots,I_1,i_1}\phi^{I_1}d\phi^{i_1}\phi^{I_2}d\phi^{i_2}\cdots
d\phi^{i_k}$ where $a_{I_k,i_k,\dots,I_1,i_1}\in\C$ has a condition
of cyclicity. 
The contraction with $k$ vector fields is then given by 
\begin{equation*}
 a(\varphi)(\delta_{A_1},\dots,\delta_{A_k})
 =a_{I_k,i_k,\dots,I_1,i_1} 
 \l(\phi^{I_1}d\phi^{i_1}\phi^{I_2}d\phi^{i_2}\cdots d\phi^{i_k} 
 (\delta_{A_1},\dots,\delta_{A_k}) \r)_c\ .
\end{equation*}
One can see that the definition of contraction is well-defined, 
that is, 
the cyclicity of $C(\varphi)_c$ is compatible with the cyclicity
of $C(\phi)_c$. 
By the cyclic permutation of $k$ vector fields, one gets  
\begin{equation*}
 a(\varphi)(\delta_{A_2},\dots,\delta_{A_k},\delta_{A_1})
 =(-1)^{k-1}(-1)^{(A_1+1)(\sum_{l=2}^k(A_l+1))}
 a(\varphi)(\delta_{A_1},\dots,\delta_{A_k})\ .
\end{equation*}
 \label{defn:contraction}
\end{defn}
The Poisson bracket and the symplectic structure is then related to
each other by 
\begin{equation}
 (A, B)=\omega(\delta_A, \delta_B)\ .
 \label{poisymcorresp}
\end{equation}
The odd Poisson structure and the odd symplectic structure 
are inverse to each other. 
Explicitly in component language, these are related by 
\begin{equation*}
 \sum_{i,I+J=K,I'+J'=K'}\omega_{ji,JI}\omega^{ik}_{J'I'}(-1)^{I'+J'+kJ'}
 =\sum_{i,I+J=K,I'+J'=K'}\omega^{ki}_{I'J'}\omega_{ij,IJ}(-1)^{iI'}
 =\delta_{K,\emptyset}\delta_{K',\emptyset}\delta_j^k\ .
\end{equation*}
Note that $\omega_{ji,JI}$ is uniquely determined 
when $\omega^{ik}_{J'I'}$ is given and vice versa. 
This fact can easily be checked by induction with respect to the powers 
of $\phi$.  
\begin{lem}
An odd bracket $(\ ,\ )$ of the form in eq.(\ref{poibra}) satisfies 
the Jacobi identity iff the corresponding two-form $\omega$ is closed. 
Namely, $(\ ,\ )$ is an odd Poisson bracket iff $\omega$ is an 
odd symplectic form. 
\label{lem:Jacobi}
\end{lem}
\begin{pf}
It follows from 
calculating $d\omega(\delta_A,\delta_B,\delta_C)$ for $A,B,C\in
C(\phi)_c$ directly and 
using the correspondence (\ref{poisymcorresp}). 
\qed
\end{pf}
\begin{rem}
$\delta_{C(\phi)_c}$ is an algebraic homomorphism between the 
Gerstenhaber algebra and the Lie super algebra defined as follows. 
\begin{equation*}
 \delta_{(B, A)}=[\delta_A, \delta_B]\ ,
 \qquad [\delta_A, \delta_B]
 :=\delta_A\delta_B-(-1)^{(A+1)(B+1)}\delta_B\delta_A\ .
\end{equation*}
 \label{rem:alghom}
\end{rem}

 \subsection{The symplectic diffeomorphisms and the Hamiltonian flows}
\label{ssec:diffeo}

It is known in the usual commutative situation 
that the infinitesimal symplectic diffeomorphisms 
are generated by Hamiltonian vector fields. 
Namely, the following holds 
\begin{equation*}
 (A, B)+\delta_\epsilon(A, B)=(A+\delta_\epsilon A, B+\delta_\epsilon B)
\end{equation*}
up to $\epsilon^2$. 
This fact holds also in the noncommutative situation. 
\begin{prop}
Any infinitesimal symplectic diffeomorphisms on a formal noncommutative 
symplectic supermanifold can be expressed as a Hamiltonian vector fields.
 \label{prop:hvf}
\end{prop}
\begin{pf} 
Let us define inner product of 
$a(\varphi)=a_{I_k,i_k,\dots,I_1,i_1}
\phi^{I_1}d\phi^{i_1}\phi^{I_2}d\phi^{i_2}\cdots
d\phi^{i_k}\in C(\varphi)_c^{\la k\ra}$ 
with vector field $\delta_A=\flpartial{\phi^j}A^j(\phi)$ by 
\begin{equation*}
 \iota_{(\delta_A)}a(\varphi)=
a_{I_k,i_k,\dots,I_1,i_1}\phi^{I_1}d\phi^{i_1}\phi^{I_2}d\phi^{i_2}\cdots
A^k(\phi)\ .
\end{equation*} 
It is shown directly that 
the infinitesimal transformation of $a(\varphi)
\in C(\varphi)_c^{\la k\ra}$ by $\delta_A$ is then Lie derivative, 
\begin{equation*}
 (d\iota_{\delta_A}+\iota_{\delta_A}d)a(\varphi)\ .
\end{equation*}

Here assume that an infinitesimal transformation 
$\delta_\epsilon:=\flpartial{\phi^i}\epsilon^i(\phi)$ preserves 
the symplectic form $\omega$, 
\begin{equation*}
 0=2(d\iota_{\delta_\epsilon}+\iota_{\delta_\epsilon}d)\omega\ .
\end{equation*}
Because $d\omega=0$, the second term is dropped out. 
Thus, the condition of the vector field $\delta_\epsilon$ 
to preserve the symplectic form is 
\begin{equation*}
 0=d\iota_{\delta_\epsilon}\omega
 =d\l(\sum_{ij,IJ}\omega_{ji,JI}\phi^I d\phi^i\phi^J\epsilon^i(\phi)\r)_c\ .
\end{equation*}
By Lemma \ref{lem:poincare} $d$-closed forms are $d$-exact, 
so $\sum_{ij,IJ}\omega_{ji,JI}\l(\phi^I d\phi^i\phi^J\epsilon^i(\phi)\r)_c$ 
is written as $d\epsilon$ 
for some degree one element $\epsilon\in C(\phi)_c$. Namely, 
for any $A\in C(\phi)_c$ 
\begin{equation*}
 \omega(\delta_A,\delta_\epsilon)=(d\epsilon)(\delta_A)
\end{equation*}
holds. Since the left hand side is $\delta_\epsilon A$ and 
the right hand side is $(\epsilon, A)=(A,\epsilon)$, one gets 
\begin{equation*}
 \delta_\epsilon=(\ ,\epsilon)\ .
\end{equation*}
\qed
\end{pf}

Note that the integral of this infinitesimal symplectic diffeomorphism 
is written of the form 
\begin{equation}
 e^{(\ ,\epsilon)}=\1+\sum_{k\ge 1}\ov{k!}(\ ,\epsilon)^k
 \label{fsymdiff}
\end{equation}
where $(\ ,\epsilon)^k$ acts on $A\in TC(\phi)_c$ as 
$(\cdots((A,\epsilon),\epsilon),\dots,\epsilon)$. 
In fact, from the properties of the Poisson bracket, 
one can check the transformation satisfies 
\begin{itemize}
 \item[(a)] ${\displaystyle 
e^{(\ , \epsilon)}AB= e^{(\ , \epsilon)}A\cdot e^{(\ , \epsilon)}B\ ,
\qquad A, B\in C(\phi),\ AB\in C(\phi)}$, 
 \item[(a')] ${\displaystyle 
e^{(\ , \epsilon)}A\bullet B= e^{(\ , \epsilon)}A\bullet e^{(\ , \epsilon)}B\ ,
\qquad A, B\in TC(\phi)_c}$, 
 \item[(b)] ${\displaystyle 
e^{(\ , \epsilon)}(A, B)= (e^{(\ , \epsilon)}A, e^{(\ , \epsilon)}B)\
,
\qquad A, B\in TC(\phi)_c}$\ . 
\end{itemize}
Condition $(a)$ implies that the finite transformation is a 
homomorphism induced by a coordinate transformation. 
Condition $(a')$ and $(b)$ just mean that the transformation preserves 
the product $\bullet$ and symplectic form, respectively.

 \subsection{The Darboux theorem for noncommutative odd 
symplectic structures}
\label{ssec:Darboux}

\begin{thm}
Any symplectic form on a formal noncommutative supermanifold can be 
transformed to be constant by a coordinate transformation. 
 \label{thm:Darboux}
\end{thm}
\begin{pf}
Let $\omega$ be any symplectic form 
on a formal noncommutative supermanifold. 
We shall consider to transform it to be constant from the lower power 
of the coordinates $\phi$. 
Suppose now that $\omega$ is transformed to be constant up to 
$k$ powers of $\phi$. 
We then consider the transformation of the form 
\begin{equation*}
 \phi^i\lgraw \phi^i+f^i(\phi)\ ,\qquad 
f^i(\phi):=f^i_{i_1\cdots i_{k+1}}\phi^{i_{k+1}}\cdots\phi^{i_1}\ .
\end{equation*}
By this transformation, $\omega$ is transformed to 
\begin{equation*}
 \begin{split}
 &\omega_{ji,\emptyset\emptyset}d\phi^i d\phi^j
 +\sum_{|I+J|=k}\omega_{ji,JI}\l(\phi^I d\phi^i\phi^J d\phi^j\r)_c+\cdots\\
 &\qquad \lgraw \omega_{ji,\emptyset\emptyset}d\phi^i d\phi^j
 +\sum_{|I+J|=k}\omega_{ji,JI}\l(\phi^I d\phi^i\phi^J d\phi^j\r)_c
 +2 d\l(\omega_{ji,\emptyset\emptyset}f^i(\phi)d\phi^j)\r)_c+\cdots
 \end{split}
\end{equation*}
Here note that since $\omega$ is closed separately with respect to the 
powers of $\phi$, $\omega_{ji,JI}\phi^I d\phi^i\phi^J d\phi^j$ is 
closed and furthermore exact due to Lemma \ref{lem:poincare}. 
Therefore it can be canceled by 
$2 d\l(\omega_{ji,\emptyset\emptyset}f^i(\phi)d\phi^j)\r)$ 
with appropriate $f^i$ since the constant part 
$\omega_{ji,\emptyset\emptyset}$ is nondegenerate. 
Thus, $\omega$ is transformed to be constant up to $(k+1)$ powers of $\phi$. 
Repeating this process completes the proof. 
\qed\end{pf}

Since the Poincar\'e's lemma (Lemma \ref{lem:poincare}) holds, 
one can also extend the above results to a noncommutative version of 
Darboux-Weinstein's theorem \cite{We}. 
In contrast, one can also prove Proposition \ref{prop:hvf} in the previous
subsection by power expansion as in Theorem \ref{thm:Darboux}.

 \subsection{Cyclic $A_\infty$-algebras from the dual pictures}
\label{ssec:cycAinftyre}

Here we shall reconsider the meaning of cyclic $A_\infty$-algebras 
from the viewpoints of odd symplectic geometry on 
formal noncommutative supermanifolds.

Let us consider any degree zero function $S\in C(\phi)_c$ which has 
critical point at the origin $\phi=0$ 
on a formal noncommutative supermanifolds. 
Such a function can be written as 
\begin{equation}
 S=\half\V_{i_1i_2}\phi^{i_2}\phi^{i_1}
 +\sum_{k\ge 3}\ov{k}\V_{i_1\cdots i_k}\phi^{i_k}\cdots\phi^{i_1}
 \label{action-pre}
\end{equation}
where $\V_{i_1\cdots i_k}\in\C$ and $S\in C(\phi)_c$. 
As seen later, this is just the action of field theory. 

For a given covariant odd Poisson structure 
in Definition \ref{defn:covpoisson}, 
let us consider the Hamiltonian vector field of $S$ 
\begin{equation}
 \delta =(\ ,S)\ . 
 \label{HvfS}
\end{equation}
By definition it has degree one. 
It is furthermore nilpotent 
iff $S$ satisfies 
\begin{equation}
 (S, S)=0\ , 
 \label{SS}
\end{equation}
where $(\ ,\ )$ is a covariant odd Poisson structure 
in Definition \ref{defn:covpoisson}. 
This fact follows from the Jacobi identity of $(\ ,\ )$. 
Namely, the Hamiltonian vector field $\delta$ of $S$ 
satisfying $(S, S)=0$ 
defines an $A_\infty$-structure. 
We remark that eq.(\ref{SS}) corresponds to 
the classical BV-master equation explained later. 
Let us call such a triple $(C(\phi)_c,\omega, S)$ or 
equivalently $(C(\phi)_c,\omega, \delta)$ 
a {\it pre cyclic $A_\infty$-supermanifold}, where 
$\omega$ is the (covariant) odd symplectic structure 
(Definition \ref{defn:covsym}) corresponding to the given 
covariant odd Poisson structure (Definition \ref{defn:covpoisson}). 

On the other hand, by Theorem \ref{thm:Darboux} 
one can transform any odd symplectic formal noncommutative supermanifolds 
to the one where the symplectic structure $\omega$ is constant, \ie , 
a skew symmetric bilinear form on $\cH$. 
In this case, the $(k+1)$-power terms of $S$ just correspond to the 
$A_\infty$-structure $m_k$ such as 
\begin{equation}
 c^j_{i_1\cdots i_k}=(-1)^{\eb_m}\omega^{jm}\V_{mi_1\cdots i_k}\ ,
\qquad k\ge 2\ ,
 \label{defc}
\end{equation}
where $\omega^{jm}$ is constant. Note that the equation above is 
equivalent to the one in Remark \ref{rem:Vcyclic}. 
The $A_\infty$-odd vector field is then written as 
\begin{equation*}
 \delta=(\ ,S)=\sum_{k=1}^\infty\flpartial{\phi^j}c^j_{i_1\cdots i_k}
 \phi^{i_k}\cdots\phi^{i_1}\ .
\end{equation*}
Thus, any pre cyclic $A_\infty$-supermanifold $(C(\phi)_c,\omega, S)$ 
is isomorphic to the one with 
a constant symplectic structure where 
the $(k+1)$-power terms of $S$ just correspond to the 
$A_\infty$-structure $m_k$. 
We call such a pre cyclic $A_\infty$-supermanifold $(C(\phi)_c,\omega, S)$ 
a {\it cyclic $A_\infty$-supermanifold}. 
This is exactly the dual description of the cyclic $A_\infty$-algebra 
$(\cH,\omega, S)$ in Definition \ref{defn:cycAinfty}. 

Next, let us consider the relation between 
two pre cyclic $A_\infty$-supermanifolds 
$(C(\phi)_c,\omega, S)$ and $(C(\phi')_c,\omega', S')$. 
Suppose that there exists a morphism $\cF_*$ (eq.(\ref{coordtransf}) )
preserving the origins ($f^{j'}=0$) and the symplectic forms, 
\begin{equation*}
 \cF^*\omega'=\omega\ .
\end{equation*} 
\begin{prop}
The following two statements are equivalent; 
\begin{itemize}
 \item $\cF$ preserves the value of the action $S=\cF^*S'$. 

 \item $\cF:(\cH,\m)\to (\cH',\m')$ is an $A_\infty$-morphism. 
\end{itemize}
 \label{prop:Ap}
\end{prop}
Note that here $f_1$ may not be an isomorphism. 
This equivalence follows from the fact that the symplectic structures 
on both sides are non-degenerate. 
We call $\cF: (C(\phi)_c,\omega,S)\to (C(\phi')_c,\omega',S')$ 
preserving the origins, the symplectic forms $\cF^*\omega'=\omega$ 
and the actions $S=\cF^*S'$ a {\it morphism between 
the (pre) cyclic $A_\infty$-supermanifolds}. 

In this situation, both pre cyclic $A_\infty$-supermanifolds 
$(C(\phi)_c,\omega,S)$ and $(C(\phi')_c,\omega',S')$ are 
isomorphic to some cyclic $A_\infty$-supermanifolds 
$(C(\phi)_c,\ti\omega,\ti{S})$ and $(C(\phi')_c,\ti{\omega'},\ti{S'})$ 
where $\ti\omega$ and $\ti{\omega'}$ are constant. 
Let $\ti\cF :(C(\phi)_c,\ti\omega,\ti{S})\to (C(\phi)_c,\omega,S)$ and 
$\ti{\cF'} :(C(\phi')_c,\ti{\omega'},\ti{S'})\to (C(\phi')_c,\omega',S')$ 
be two isomorphisms. 
The composition map $(\ti{\cF'})^{-1}\cF\ti\cF : 
(C(\phi)_c,\ti\omega,\ti{S})\to (C(\phi')_c,\ti{\omega'},\ti{S'})$ is 
then an $A_\infty$-morphism preserving the constant symplectic structures. 
Thus, there exists a functor from 
the category of pre cyclic $A_\infty$-supermanifolds 
to the category of cyclic $A_\infty$-supermanifolds. 

A morphism between two cyclic $A_\infty$-supermanifolds is 
in fact the dual description of a cyclic $A_\infty$-morphism in 
Definition \ref{defn:cycAinftymorp}. 
For two cyclic $A_\infty$-supermanifolds 
$(C(\phi)_c,\omega, S)$, $(C(\phi')_c,\omega', S')$ and a morphism 
$\cF_*: (C(\phi)_c,\omega, S)\to (C(\phi'),\omega', S')$, 
the condition $\cF^*\omega'=\omega$ is written as 
\begin{equation*}
 \left( \sum_{k',l'}\omega'_{k'l'}
d(\cF_*(\phi^{l'}))d(\cF_*(\phi^{k'}))\right)_c
 =\sum_{i,j}\omega'_{ij}d\phi^j d\phi^i\ .
\end{equation*}
Reading the coefficient of each term of polynomial fields separately 
and dualizing back, one gets 
$\omega'(f_1(\eb_i),f_1(\eb_j))=\omega(\eb_i,\eb_j)$ (eq.\ref{omegacF1}) 
and 
\begin{equation}
 \sum_{k,l\ge 1,\ k+l=n}
 \sum_{cyc(i,j)}
\omega'(f_k(\dots,\eb_{j_q},\eb_i,\eb_{i_1},\dots),
f_l(\dots,\eb_{i_p},\eb_j,\eb_{j_1},\dots))=0\ 
 \label{omegacF2'}
\end{equation}
for fixed $n\ge 3$, $p\ge 0$ and $q=n-p-2\ge 0$, 
where $cyc(i,j)$ indicates all possible 
cyclic permutations for 
$\{ \eb_i, \eb_{i_1}, \dots, \eb_{i_p}, 
\eb_j, \eb_{j_1}, \dots, \eb_{j_q} \}$ but 
keeping $\eb_i$ in $f_k(\cdots)$ and $\eb_j$ in $f_l(\cdots)$. 
However, eq.(\ref{omegacF2'}) includes overlapping equivalent
identities. For instance, when we consider the case of 
order $\{ \eb_i, \eb_j, \eb_{j_1}, \dots, \eb_{j_{n-2}} \}$, 
the summation $cyc(i,j)$ drops out and eq.(\ref{omegacF2'}) simply
gives 
\begin{equation*}
  \sum_{k,l\ge 1,\ k+l=n}
\omega'(f_k(\eb_{j_l},\dots,\eb_{j_{n-2}},\eb_i),
f_l(\eb_j,\eb_{j_1},\dots,\eb_{j_{l-1}}))=0\ . 
\end{equation*}
Namely, the condition (\ref{omegacF2'}) just reduces to 
eq.(\ref{omegacF2}) and the condition of morphisms between 
cyclic $A_\infty$-supermanifolds is actually equivalent to the 
condition of cyclic $A_\infty$-morphisms 
in Definition \ref{defn:cycAinftymorp}.

%section5

 \section{The minimal model theorem}
\label{sec:MandC}

The following theorem is 
one of the key theorems in homotopy algebra. 
\begin{thm}[Minimal model theorem \cite{kadei1}]
Given any $A_\infty$-algebra $(\cH,\m)$, let $\cH^p$ be 
the cohomology of $\cH$ with respect to $m_1$. 
Then there necessarily exists an 
$A_\infty$-algebra $(\cH^p,\ti\m^p)$ and an $A_\infty$-quasi-isomorphism 
from $(\cH^p,\ti\m^p)$ to $(\cH,\m)$. 
 \label{thm:minimal}
\end{thm}
The purpose of this section is to clarify and to develop 
basic properties of $A_\infty$-algebras around this theorem. 
For the construction of minimal models of $A_\infty$-structures, 
in particular of dgas, 
various versions of homological perturbation theory (HPT) 
have been developed, for instance by 
\cite{gugen,GS,gugen-lambe,GLS:chen,GLS,hueb-kadei}. 
On the other hand, 
it was mentioned in \cite{Ko1} that there exists 
another stronger version of the minimal model theorem, 
which we call the decomposition theorem. 
In subsection \ref{ssec:MandC} we shall show explicitly 
the decomposition theorem (Theorem \ref{thm:MandC}). 
A similar result is obtained independently in \cite{Le-Ha} 
in the framework of a closed model category \cite{Q}. 
The decomposition theorem includes the minimal model theorem and 
implies various basic properties of homotopy algebras. 
In subsection \ref{ssec:cycMandC} we show the decomposition theorem 
for cyclic $A_\infty$-algebras. 
The decomposition theorem guarantees the existence of an inverse 
$A_\infty$-quasi-isomorphism of an $A_\infty$-quasi-isomorphism 
(Theorem \ref{thm:inverse}) as stated in \cite{Ko1}. 
We shall explain it in subsection \ref{ssec:inverse}. 
Though the minimal model theorem follows from the decomposition theorem, 
the proof relies on inductive arguments and the form of the 
minimal model is not explicit. On the other hand, 
for any $A_\infty$-algebra, its minimal model 
can be given explicitly and more recently in \cite{KoSo} 
in terms of some Feynman diagrams. 
The Feynman diagram expression provides us with intuitive 
understanding of the minimal model, though it is equivalent to 
minimal models given by formula as in the traditional 
homological perturbation theory 
(see \cite{GLS,hueb-kadei,mer}). 
We demonstrate in subsection \ref{ssec:minimal} 
that it arises naturally from the issue of finding 
the solutions of the Maurer-Cartan equation 
for an $A_\infty$-algebra \cite{Ka}. 
Subsection \ref{ssec:cycminimal} presents the cyclic $A_\infty$ version, 
which is directly related to section \ref{sec:BV}, where 
we shall show that the minimal cyclic algebra is derived by 
semiclassical perturbation theory of field theory 
in the BV-formalism.

 \subsection{The decomposition theorem for $A_\infty$-algebras}
\label{ssec:MandC}

\begin{defn}[Minimal $A_\infty$-algebra]
An $A_\infty$-algebra $(\cH,\m)$ is called {\it minimal} if 
$m_1=0$ on $\cH$. 
 \label{defn:minimalalg}
\end{defn}
\begin{defn}
An $A_\infty$-algebra $(\cH,\m)$ is called {\it linear contractible} 
if $m_k=0$ for $k\ge 2$ and $Q=m_1$ has trivial cohomology. 
 \label{defn:cont}
\end{defn}
The following theorem holds. 
\begin{thm}[Decomposition theorem for $A_\infty$-algebras]
Any $A_\infty$-algebra is 
$A_\infty$-isomorphic to the direct sum of a minimal $A_\infty$-algebra 
and a linear contractible $A_\infty$-algebra. 
 \label{thm:MandC}
\end{thm}
This subsection is devoted to proving this theorem 
based on the strategy presented in \cite{Ko1}. 
We shall prove it 
in terms of formal noncommutative supermanifolds. 
However the result holds 
even in the case that the dual side is not well-defined in a strict sense, 
since the proof can be translated into 
that on the corresponding coalgebra side (see Remark \ref{rem:valid}).

For $Q$ the coboundary operator of the complex $(\cH, Q:=m_1)$, 
we first give the analogue of the Hodge-Kodaira decomposition; 
consider a degree minus one linear map $Q^+:\cH\to\cH$ satisfying 
\begin{equation}
 QQ^++Q^+ Q+P=\1\ 
 \label{HKdecomp}
\end{equation}
on $\cH$ such that $P^2=P$ and $QP=0$ hold. 
Namely, $P$ is the projection onto the cohomology of $(\cH,Q)$. 
By definition $PQ=0$ and $QQ^+Q=Q$ hold, which lead to 
$(QQ^+)^2=QQ^+$, $(Q^+Q)^2=Q^+Q$ and 
\begin{equation*}
 (QQ^+)(Q^+Q)=(Q^+Q)(QQ^+)=0\ , \quad 
 P(Q^+Q)=(Q^+Q)P=0\ ,\quad P(QQ^+)=(QQ^+)P=0 \ .
\end{equation*}
In physical terms, $QQ^+$, $Q^+Q$ and $P$ are the projections onto 
$Q$-trivial states, 
unphysical states and physical states, respectively. 
We denote $QQ^+=P^t$ and $Q^+Q=P^u$ and express the decomposition 
of $\cH$ as 
\begin{equation}
 \cH=\cH^t\oplus\cH^u\oplus\cH^p\ ,\qquad 
 \cH^p:=P\cH\ ,\quad \cH^t:=P^t\cH\ ,\quad \cH^u:=P^u\cH\ .
 \label{cHdecomp}
\end{equation}

In other words, we give a splitting of the complex $(\cH,Q)$: 
\begin{equation*}
 \xymatrix{
  (\cH,Q)\ \ \ar@<0.5ex>[r]^\pi &\ \ (\cH^p, 0) \ar@<0.5ex>[l]^\iota
}
\end{equation*} 
with a contracting homotopy $Q^+:\cH\to\cH$ such that 
$\1-P=QQ^++Q^+Q$ for $P:=\iota\circ\pi$ 
(we shall sometimes omit $\circ$). 
Here $0$ is the zero differential on $\cH^p$. 
As above, we shall often denote the image of $\cH^p$ by $\iota$ also by
$\cH^p$. 
By definition $P^2=P$ holds. Also, since $\iota:(\cH^p,0)\to (\cH,Q)$ 
and $\pi:(\cH,Q)\to (\cH^p,0)$ are chain maps, 
$QP=PQ=0$ hold. Thus, these data give the decomposition
(\ref{cHdecomp}). 
Note that this set-up is a special case (in particular the
differential on $\cH^p$ is zero) of the 
{\em strong deformation retract} or SDR which is the starting point 
of various versions of homological perturbation theory 
(see \cite{gugen,GS,gugen-lambe,GLS:chen,GLS,hueb-kadei}). 

We decompose the dual coordinates $\{\phi^i\}$ of $\cH$ into 
$\{x^i\}$, $\{y^j\}$ and $\{p^k\}$, 
the dual coordinates corresponding to the basis of 
$\cH^t$, $\cH^u$ and $\cH^p$, respectively. 

Our goal for the proof of Theorem \ref{thm:MandC} is to 
construct a local diffeomorphism around the origin of the formal 
noncommutative supermanifold such that 
$x^i$ and $y^j$ span the contractible directions, whereas 
$p^k$ is a coordinate of the minimal part. 

For the first step, we should examine 
the properties of the cohomology with respect to $Q$ on formal 
noncommutative supermanifolds. 
\begin{defn}[$\delta_1$-complex]
Let $C(\phi)^k$ be the space of associative noncommutative 
polynomial functions on a formal noncommutative supermanifold of 
total degree $k$. 
$\delta_1 : C(\phi)^k\to C(\phi)^{k+1}$, which is the dual of $Q$ on $\cH$, 
then defines a complex. 
We denote it by $(C(\phi)^\star,\delta_1)$ and call 
the $\delta_1$-complex. 
 \label{defn:Qcpx}
\end{defn}
\begin{lem}[$\delta_1$-cohomology]
%Let us decompose the coordinates $\phi^i$ into 
%$x^i, y^j, p^k$ where $x^i$, $y^j$ and $p^k$ are 
%the dual coordinates corresponding to the basis of 
%$\cH^t$, $\cH^u$ and $\cH^p$, respectively. 
For any element $a_{i_1\cdots i_k}\phi^{i_k}\cdots\phi^{i_1}
\in C(\phi)^\star$, 
its cohomology part is represented as 
\begin{equation}
 a_{i_1\cdots i_k}p^{i_k}\cdots p^{i_1}
 =\l(a_{i_1\cdots i_k}\phi^{i_k}\cdots\phi^{i_1}\r)\Big|_{x=y=0}\ .
 \label{cohproj}
\end{equation}
We denote by $C(p)^\star$ the space of such elements. 
 \label{lem:Qcoh}
\end{lem}
\begin{pf}
This follows from the fact that one can construct a degree minus one 
homotopy operator $H^*$ on $C(\phi)^\star$ such that 
\begin{equation*}
 \Id-P=\delta_1 H^*+H^*\delta_1\ ,
\end{equation*}
where $P$ is the projection corresponding to eq.(\ref{cohproj}). 
This $H^*$ will be constructed 
explicitly later in Lemma \ref{lem:homIdP}. 
\footnote{More precisely, we will construct an operator $H$ in 
Lemma \ref{lem:homIdP} and $H^*$ is just the dual of $H$.}
\qed
\end{pf}

\noindent
One may notice that this corresponds to 
the dual version (supermanifold description) of so-called 
the {\em tensor trick} 
(see \cite{gugen-lambe,GLS:chen,GLS,hueb-kadei}). 

\vspace*{0.2cm}

\noindent
{\it proof of Theorem \ref{thm:MandC}.}\quad 
The $A_\infty$-odd vector field is then written as 
\begin{equation*}
 \begin{split}
 &\delta=\delta_1+\delta_2+\cdots\ ,\\
 &\delta_1=\flpartial{x^i}c^i_j y^j\ ,\ \cdots,\ \  
 \delta_n=\flpartial{\phi^i}c^i_{k_1\cdots k_n}\phi^{k_n}\cdots\phi^{k_1}\ 
 \end{split}
\end{equation*}
for $n\geq 2$. 
Note that $\delta_1$ acts nontrivially on the contractible part only. 
Now we would like to bring $\delta_n$, $n\ge 2$, to the form 
$\flpartial{p^i}c^i_{k_1\cdots k_n}p^{k_n}\cdots p^{k_1}$ 
by a coordinate transformation. 
Suppose now that $\delta_n$ is such a form for $n\le l$. 
By the $(l+1)$ powers-part of the $A_\infty$-condition, 
we have 
\begin{align}
 &\flpartial{x^i}
 \l({c_x}^i_{j_1\cdots j_{l+1}}\phi^{j_{l+1}}\cdots\phi^{j_1}\r)
 \flpartial{x^k}c^k_j y^j
 +\flpartial{x^i}c^i_j 
 {c_y}^j_{j_1\cdots j_{l+1}}\phi^{j_{l+1}}\cdots\phi^{j_1} =0\ ,
 \label{cx}\\
 &\flpartial{y^i}
 \l({c_y}^i_{j_1\cdots j_{l+1}}\phi^{j_{l+1}}\cdots\phi^{j_1}\r)
 \flpartial{x^k}c^k_j y^j=0\ ,\label{cy}\\
 &\flpartial{p^i}
 \l({c_p}^i_{j_1\cdots j_{l+1}}\phi^{j_{l+1}}\cdots\phi^{j_1}\r)
 \flpartial{x^k}c^k_j y^j
 +\sum_{n\ge 2}\delta_n\delta_{l+2-n}=0\ ,\label{cp}
\end{align}
where 
${c_x}^i_{j_1\cdots j_{l+1}}\phi^{j_{l+1}}\cdots\phi^{j_1}
\in C(\phi)$ 
is the coefficient of $\delta_{l+1}$ with respect to 
$\flpartial{x^i}$, and so on. 
We sometimes ignore the upper indices $i$ and 
indicate by $c_x, c_y$ and $c_p$ these coefficients in $C(\phi)$. 
One can see that $c_y$ is $\delta_1$-closed from eq.(\ref{cy}). 
Moreover, eq.(\ref{cx}) implies $c_y$ does not have $\delta_1$-cohomology 
since the first term includes $y^j$ and the second term does not. 
Therefore $c_y$ is $\delta_1$-exact due to Lemma \ref{lem:Qcoh}. 
Alternatively, 
the first and second terms of eq.(\ref{cp}) are independent of each other 
since the first term includes $y^j$ and second term does not. 
Thus, $c_p$ is $\delta_1$-closed. 

We would like to remove $c_x$, $c_y$ and the $\delta_1$-exact part of $c_p$. 
Note that in fact we know a transformation which removes $c_x$ and $c_y$. 
It is given by the $A_\infty$-quasi-isomorphism 
in subsection \ref{ssec:minimal}. 
However it is instructive to construct the transformation inductively 
without assuming this knowledge. 

Consider the following coordinate transformation
\begin{equation}
 \phi^i={\phi'}^i+f^i(\phi)\ ,\qquad f^i(\phi)=
 f^i_{k_1\cdots k_{l+1}}{\phi'}^{k_{l+1}}\cdots{\phi'}^{k_1}\ .
 \label{MandCtransf}
\end{equation}
Its inverse transformation is ${\phi'}^i=\phi^i-f^i(\phi)+\cdots$. 
We write the transformation (\ref{MandCtransf}) separately with respect to 
$x$, $y$ and $p$ such as 
$$
 x^i={x'}^i+f_x^i(\phi)\ ,\qquad y^j={y'}^j+f_y^j(\phi)\ ,\qquad 
 p^k={p'}^k+f_p^k(\phi)\ .
$$
$\delta$ is then transformed as (see eq.(\ref{cc'}))
\begin{equation*}
 \delta\ \lgraw\ \delta+\flpartial{x^i}c^i_j f_y^j(\phi)
 -\flpartial{\phi^k}\l(f^k(\phi)\flpartial{x^i}c^i_jy^j\r)+\cdots\ .
\end{equation*}
The conditions for the correction terms $f^i(\phi)$ 
to cancel $\delta_{l+1}$ except its minimal part are
\begin{align}
 & c^i_j f_y^j(\phi)-f_x^i(\phi)\flpartial{x^k}c^k_jy^j
 +{c_x}^i_{j_1\cdots j_{l+1}}\phi^{j_{l+1}}\cdots\phi^{j_1}
 =0\ ,\label{cx2}\\
 &-f_y^i(\phi)\flpartial{x^i}c^i_jy^j
 +{c_y}^i_{j_1\cdots j_{l+1}}\phi^{j_{l+1}}\cdots\phi^{j_1}
 =0\ ,\label{cy2}\\
 &-f_p^i(\phi)\flpartial{x^i}c^i_jy^j
 +{c_p}^i_{j_1\cdots j_{l+1}}\phi^{j_{l+1}}\cdots\phi^{j_1}
 =(\mbox{polynomial of $p$'s})\ .\label{cp2}
\end{align}
We shall perform the transformation by dividing it into two steps. 
First, find a solution of eq.(\ref{cx2}) and eq.(\ref{cy2}) 
for $f_x$ and $f_y$ which removes $c_x$ and $c_y$. 
There exist ambiguities for the solution; one 
solution is given by $f_x=0$ and 
$f_y^j=-{\bar c}^j_k c^k_{x,j_1\cdots j_{l+1}}
\phi^{j_{l+1}}\cdots\phi^{j_1}$. 
It is clear that it satisfies eq.(\ref{cx2}). 
In order to confirm eq.(\ref{cy2}) one may employ eq.(\ref{cx}), 
an $A_\infty$-condition for $\delta$ under the induction hypothesis. 
Next, $c_p$ is $\delta_1$-closed 
as stated above and eq.(\ref{cp2}) implies the $\delta_1$-exact part is 
removed by $f_p$. The remaining $\delta_1$-cohomology part of $c_p$ 
then forms the minimal $A_\infty$-structure. 
This completes the induction. 
\qed
\begin{rem}
In the proof above, one can see that 
the coordinate transformation $f^i(\phi)$ was constructed only from 
$Q^+$ and $\m$. Therefore, the proof can be rewritten purely on the coalgebra
side and is independent of 
whether the dual supermanifold 
description is strictly well-defined or not 
(see \cite{KaTe}). 
 \label{rem:valid}
\end{rem}
Given an $A_\infty$-algebra $(\cH,\m)$, let us denote by 
$(\cH_{dc},\m_{dc})$ 
a direct sum of a minimal $A_\infty$-algebra and a 
linear contractible $A_\infty$-algebra obtained by 
Theorem \ref{thm:MandC} and $(\cH_{dc}^p,\m_{dc}^p)$ its minimal part. 
One can obtain an $A_\infty$-isomorphism 
$\cF_{dc}: (\cH_{dc},\m_{dc})\to (\cH,\m)$ by composing inductively 
$A_\infty$-isomorphisms (\ref{MandCtransf}). 
As stated above, the $A_\infty$-isomorphism is not unique so 
the decomposed model $(\cH_{dc},\m_{dc})$ is not unique. 
Moreover, even if $(\cH_{dc},\m_{dc})$ is fixed, the 
$A_\infty$-isomorphism $\cF_{dc}: (\cH_{dc},\m_{dc})\to (\cH,\m)$ 
is not unique since there exists gauge transformations (see 
Definition \ref{defn:gauge}). 

In the situation above, 
both the inclusion $\iota : \cH_{dc}^p\to \cH_{dc}$ and 
the projection $\pi :\cH_{dc}\to \cH_{dc}^p$ 
such that $\iota\circ\pi=P$ 
naturally extend to  
$A_\infty$-quasi-isomorphisms 
by setting the leading linear maps 
of the $A_\infty$-quasi-isomorphisms to be $\iota$ and $\pi$, and 
their higher multilinear maps zero. 
We write them also as 
$\iota :(\cH_{dc}^p,\m_{dc}^p)\to (\cH_{dc},\m_{dc})$ and 
$\pi :(\cH_{dc},\m_{dc})\to (\cH_{dc}^p,\m_{dc}^p)$. 
Then one can see that $\cF_{dc}\circ\iota$ is 
an $A_\infty$-quasi-isomorphism from $(\cH_{dc}^p,\m_{dc}^p)$ to 
$(\cH,\m)$ and an inverse $A_\infty$-quasi-isomorphism is given by 
$\pi\circ (\cF_{dc})^{-1}: (\cH,\m)\to (\cH_{dc}^p,\m_{dc}^p)$. 

Although the existence of a minimal model for any $A_\infty$-algebra 
is guaranteed from Theorem \ref{thm:MandC}, 
the minimal model is not unique as stated above. 
A trivial ambiguity is the one given by an $A_\infty$-isomorphism 
transforming a minimal $A_\infty$-algebra 
to another one. 
On this point, from Theorem \ref{thm:MandC}, one can state the following. 
\begin{cor}
For an $A_\infty$-algebra $(\cH,\m)$, suppose that it has
a minimal model $(\cH^p,\ti\m^p)$ 
and an $A_\infty$-quasi-isomorphism 
$\ti\cF^p :(\cH^p,\ti\m^p)\to (\cH,\m)$ are given. 
Then, there exist a decomposed $A_\infty$-algebra $(\cH_{dc},\m_{dc})$ 
whose minimal part is the minimal $A_\infty$-algebra
$(\cH^p,\ti\m^p)$ and an $A_\infty$-isomorphism 
$\cF_{dc}:(\cH_{dc},\m_{dc})\to (\cH,\m)$. 
Equivalently, the $A_\infty$-quasi-isomorphism 
$\ti\cF^p :(\cH^p,\ti\m^p)\to (\cH,\m)$ 
can be described by the composition 
$\cF_{dc}\circ\iota$, where 
$\iota: (\cH^p,\ti\m^p)\to (\cH_{dc},\m_{dc})$ is the inclusion. 
 \label{cor:comp}
\begin{cor}
Given a minimal model $(\cH^p,\ti\m^p)$ of an $A_\infty$-algebra
$(\cH,\m)$ 
and an $A_\infty$-quasi-isomorphism 
$\ti\cF^p :(\cH^p,\ti\m^p)\to (\cH,\m)$, 
there exists an 
inverse $A_\infty$-quasi-isomorphism 
$(\ti\cF^p)^{-1} : (\cH,\m)\to (\cH^p,\ti\m^p)$. 
 \label{cor:inverse}
\end{cor}
\end{cor}
\begin{cor}[Uniqueness of minimal $A_\infty$-algebras]
For an $A_\infty$-algebra $(\cH,\m)$, 
its minimal $A_\infty$-algebra is unique up to an isomorphism 
on $\cH^p$. 
 \label{cor:uniqueness}
\end{cor}
\begin{pf}
These three corollaries are shown at the same time as follows. 
For an arbitrary decomposed $A_\infty$-algebra 
$(\cH'_{dc},\m'_{dc})$ of an $A_\infty$-algebra $(\cH,\m)$ and 
an $A_\infty$-isomorphism 
$\cF'_{dc}:(\cH'_{dc},\m'_{dc})\to (\cH,\m)$, 
let us consider the following diagram
\begin{equation*}
\xymatrix{
 (\cH,\m) \ar@<0.5ex>[r]^{(\cF'_{dc})^{-1}}\ \  & 
 \ \ (\cH'_{dc},\m'_{dc}) \ar@<0.5ex>[l]^{\cF'_{dc}} \ \ \, 
 \ar@<0.5ex>[d]^{\pi'} \\
 (\cH^p,\ti\m^p)\ \  \ar[u]^{\ti\cF^p} \ar[r]^{\cF^p} 
  &\ \ ({\cH'}_{dc}^p,{\m'}_{dc}^p)\ \ . \ar@<0.5ex>[u]^{\iota'}
}
\end{equation*}
Note that $\cF'_{dc}$ has its inverse $A_\infty$-isomorphism 
$(\cH'_{dc})^{-1}$ and $\iota'$ has its 
inverse $A_\infty$-quasi-isomorphism 
$\pi'$. The composition $\pi'\circ(\cF'_{dc})^{-1}\circ\ti\cF^p$ then gives 
an $A_\infty$-quasi-isomorphism 
from $(\cH^p,\ti\m^p)$ to $({\cH'}_{dc}^p,{\m'}_{dc}^p)$. 
Denote this composition map by 
$\cF^p: (\cH^p,\ti\m^p)\to({\cH'}_{dc}^p,{\m'}_{dc}^p)$. 
Since $\cH^p={\cH'}_{dc}^p$, it is not only 
a quasi-isomorphism but also an isomorphism. 
On the other hand, one can consider 
an $A_\infty$-algebra $(\cH_{dc},\m_{dc})$ 
which is the direct sum of the minimal $A_\infty$-algebra 
$(\cH^p,\ti\m^p)$ 
and the contractible part of $(\cH'_{dc},\m'_{dc})$. 
Clearly, there exists an $A_\infty$-isomorphism 
$\cF: (\cH_{dc},\m_{dc})\to(\cH'_{dc},\m'_{dc})$ which is the natural 
extension of $\cF^p$. 
Thus, the composition of $\iota: (\cH^p,\ti\m^p)\to (\cH_{dc},\m_{dc})$
with $\cF$ yields an $A_\infty$-quasi-isomorphism 
$\cF\circ\iota: (\cH^p,\ti\m^p)\to(\cH'_{dc},\m'_{dc})$. 
This leads to Corollary \ref{cor:comp}; 
$\ti\cF^p=\cF_{dc}\circ\iota$ where $\cF_{dc}:=\cF'_{dc}\circ\cF$. 

It is clear that the existence of inverse quasi-isomorphisms 
(Corollary \ref{cor:inverse}) follows from 
Corollary \ref{cor:comp}. Namely, an inverse quasi-isomorphism of 
$\cF_{dc}\circ\iota$ is given by 
$\pi\circ\cF_{dc}$. 

The uniqueness of the minimal model (Corollary \ref{cor:uniqueness}) then 
follows from Corollary \ref{cor:inverse}. 
Since Corollary \ref{cor:inverse} implies the existence of 
an $A_\infty$-quasi-isomorphism between any two minimal models 
and moreover the quasi-isomorphism is an isomorphism, 
one can see that there exists an $A_\infty$-isomorphism between any two 
minimal models of an $A_\infty$-algebra $(\cH,\m)$. 
 \qed
\end{pf}

We shall construct a minimal model explicitly by using Feynman graphs 
in Definition \ref{defn:minimal} in subsection 
\ref{ssec:minimal}. 
Corollary \ref{cor:uniqueness} then implies that 
the explicit minimal model is unique 
up to $A_\infty$-isomorphisms on $\cH^p$. 

\begin{rem}[Geometric realization]
One can realize the results around the decomposition theorem above 
geometrically as follows. 
\begin{figure}[h] \label{fig:minimal}
 \hspace*{1.3cm}
 \includegraphics{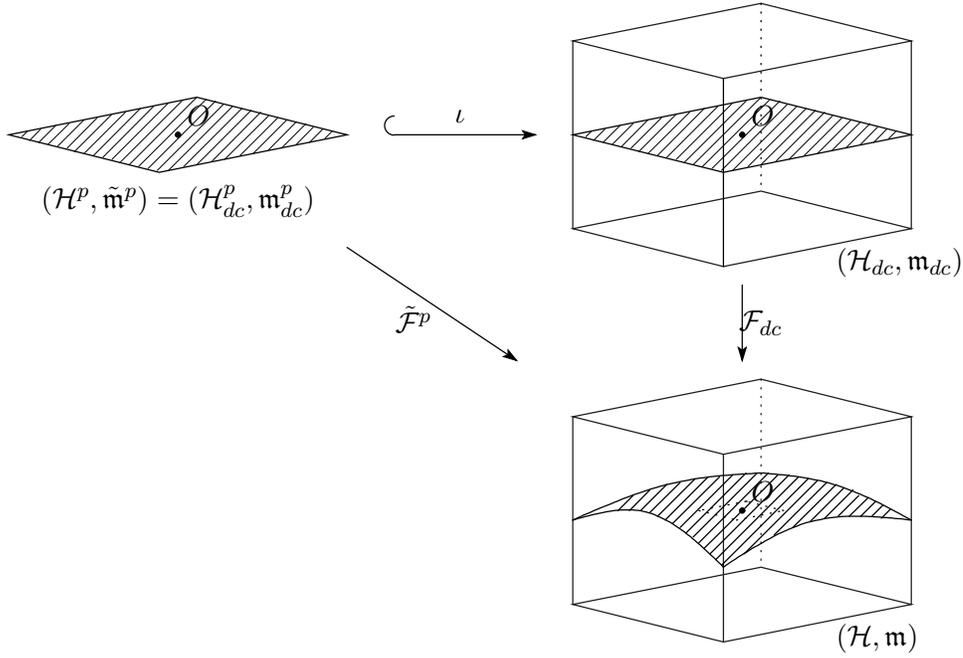}
 \caption{The embedding of the minimal model. }
\end{figure}
An $A_\infty$-algebra $(\cH,\m)$ is equivalent to a formal noncommutative 
supermanifold with an $A_\infty$-odd vector field $\delta$. 
Denote $\delta=\delta_1+\delta_\bullet$ 
where $\delta_1$ is the dual of $m_1$ 
and $\delta_\bullet$ is the rest. 
First we fix the basis of the formal noncommutative supermanifold 
from $\delta$ and its homotopy operator ($Q^+$) and 
decompose $\cH=\cH^p\oplus\cH^t\oplus\cH^u$ as a graded vector space. 
Theorem \ref{thm:MandC} implies that 
there exists a nonlinear coordinate transformation so that 
the vector field $\delta_\bullet$ flows along the $\cH^p$ direction and 
does not depend on $\cH^t\oplus\cH^u$ direction. 
This fact guarantees the existence of the minimal model. 
Moreover for a minimal model $(\cH^p,\ti\m^p)$, 
an $A_\infty$-quasi-isomorphism $\ti\cF^p :(\cH^p,\ti\m^p)\to (\cH,\m)$ 
can be regarded as an embedding of a hypersurface $\cH^p$ into 
$\cH$ as in Figure \ref{fig:minimal}, 
though it is a formal noncommutative graded hypersurface. 
It follows from the condition on the $A_\infty$-morphisms 
$\ti\cF^p\ti\m^p=\m\ti\cF^p$ that on the hypersurface in $\cH$ 
the $A_\infty$-odd vector field $\delta=\delta_\bullet$ is tangent to it.
\footnote{Note that $\delta_1$ vanishes on the hypersurface. } 
Corollary \ref{cor:comp} then implies that such an embedding can necessarily 
be given by the inclusion $\iota$ 
and a nonlinear coordinate transformation $\cF_{dc}$ 
preserving the tangent space at the origin $O$. 
Thus, one can see that many ordinary geometric intuitions are valid 
in formal noncommutative graded situations. 
 \label{rem:geom}
\end{rem}

 \subsection{The decomposition theorem for cyclic $A_\infty$-algebras}
\label{ssec:cycMandC}

In this subsection we prove the decomposition theorem for cyclic 
$A_\infty$-algebras. 
For a given cyclic $A_\infty$-algebra $(\cH,\omega,\m)$, 
the proof requires 
a homotopy operator which gives an orthogonal decomposition of $\cH$ 
with respect to the inner product $\omega$. 
Let us begin with arbitrary homotopy operator $Q^+$ which 
defines a Hodge-Kodaira decomposition of $\cH$
$$ 
QQ^++Q^+Q+P=\1\ .
$$ 
There are ambiguities of the choice of $Q^+$; 
$\cH^t=QQ^+\cH$ is unique, but $\cH^p=P\cH$ is unique 
modulo $\cH^t$ and $\cH^u=Q^+Q\cH$ is unique 
modulo $\cH^t\oplus\cH^p$. 
For the odd symplectic inner product $\omega$ 
of the cyclic $A_\infty$-algebra $(\cH,\omega,\m)$ 
with homogeneous basis $\{\eb_i\}$ of $\cH$, denote 
\begin{equation*}
 \omega_{ij}:=\omega(\eb_i,\eb_j)\ .
\end{equation*}
{}From the cyclicity $\omega(\cH,Q\cH)=-\omega(Q\cH,\cH)$ holds, 
which implies the following properties independent of the ambiguities; 
if $\eb_j\in\cH^t$ then $\omega_{ij}=0$ for $\eb_i\in\cH^t\oplus\cH^p$, 
and if $\eb_j\in\cH^t\oplus\cH^p$ then $\omega_{ij}=0$ for 
$\eb_i\in\cH^t$. If we denote the block element of matrix $\{\omega_{ij}\}$ 
where $\eb_i\in\cH^u$ and $\eb_j\in\cH^p$ as $\omega_{up}$ and similar for the
other eight block elements, the matrix $\{\omega_{ij}\}$ is represented as 
the left hand side of eq.(\ref{omegadecomp}). 
This implies that by basis transformation corresponding 
to the ambiguity in $Q^+$ the inner product $\omega$ is decomposed 
as the right hand side 
\begin{equation}
 \{\omega_{ij}\}=\bp \omega_{uu}&\omega_{up}&\omega_{ut}\\
                     \omega_{pu}&\omega_{pp}&\omega_{pt}\\
                     \omega_{tu}&\omega_{tp}&\omega_{tt}\ep
                =\bp \omega_{uu}&\omega_{up}&\omega_{ut}\\
                     \omega_{pu}&\omega_{pp}& 0         \\
                     \omega_{tu}& 0         & 0         \ep\ 
 \lgraw          \bp  0         & 0         &\omega_{ut}\\
                      0         &\omega_{pp}& 0         \\
                     \omega_{tu}& 0         & 0         \ep\ . 
 \label{omegadecomp}
\end{equation}
Thus, there exists a homotopy operator $Q^+$, 
which defines the state space $\cH^u$, 
so that $\omega$ is decomposed orthogonally. 
\begin{defn}[Homotopy operator of $Q$ compatible with $\omega$]
Given a cyclic $A_\infty$-algebra $(\cH,\omega,\m)$, 
let $QQ^++Q^+Q+P=\1$ be a Hodge-Kodaira decomposition of $\cH$. 
We call $Q^+$ a homotopy operator of $Q$ compatible with $\omega$ 
if $\omega(\cH^p,\cH^u)=\omega(\cH^u,\cH^u)=0$. 
(When ($\cH,\omega,\m$) is a cyclic $A_\infty$-algebra, 
$\omega(\cH^t,\cH^p)=\omega(\cH^t,\cH^t)=0$ is automatically satisfied. )
 \label{defn:homopomega}
\end{defn}
Note that a homotopy operator $Q^+$ of $Q$ compatible with $\omega$ 
satisfies 
\begin{equation}
\omega(\1\otimes Q^+)=\omega(Q^+\otimes\1)\ ,\qquad 
\omega(\1\otimes P)=\omega(P\otimes\1)\ .
 \label{homotopy-compatible}
\end{equation}
In particular, the first equation follows from 
$\omega(\1\otimes Q^+)=\omega(QQ^+\otimes Q^+)=\omega(Q^+\otimes
QQ^+)=\omega(Q^+\otimes\1)$. 
As stated above, any cyclic $A_\infty$-algebra $(\cH,\omega,\m)$ 
has a homotopy operator of $Q$ compatible with $\omega$. 
It satisfies additional conditions compared with 
a homotopy operator 
defining the usual Hodge-Kodaira decomposition; however, 
it is still not unique. The remaining ambiguity is just related 
to the choice of the gauge fixing (Definition \ref{defn:gf})
when the propagator (Definition \ref{defn:BVpropagator}) 
is constructed in the BV-formalism. 

The properties of the $\delta_1$-complex are similar to those 
in the previous non-cyclic case. 
\begin{defn}[Cyclic $\delta_1$-complex]
Let $C(\phi)_c^k$ be the space of associative noncommutative 
polynomial cyclic functions on a formal noncommutative supermanifold of 
total degree $k$. 
The action $\delta_1$ on it satisfies the Leibniz rule and so 
preserves the cyclicity. 
$\delta_1 : C(\phi)_c^k\to C(\phi)_c^{k+1}$ then defines a complex. 
We denote it by $(C(\phi)_c^\star,\delta_1)$ and call it
the cyclic $\delta_1$-complex. 
 \label{defn:cycQcpx}
\end{defn}
\begin{lem}[Cyclic $\delta_1$-cohomology]
For any element $a_{i_1\cdots i_k}\phi^{i_k}\cdots\phi^{i_1}
\in C(\phi)_c^\star$, 
its cohomology part is represented as 
\begin{equation*}
 a_{i_1\cdots i_k}p^{i_k}\cdots p^{i_1}
 =\l(a_{i_1\cdots i_k}\phi^{i_k}\cdots\phi^{i_1}\r)\Big|_{x=y=0}\ .
\end{equation*}
We denote by $C(p)^\star_c$ the space of such elements. 
 \label{lem:cycQcoh}
\end{lem}
\begin{pf}
The corresponding homotopy operator on $C(\phi)^\star_c$ can be obtained 
by the cyclic symmetrization of 
the homotopy operator in Lemma \ref{lem:Qcoh}. 
\qed
\end{pf}
\begin{thm}[Decomposition theorem for cyclic $A_\infty$-algebra]
Any cyclic $A_\infty$-algebra is 
cyclic $A_\infty$-isomorphic to the direct sum of 
a minimal cyclic $A_\infty$-algebra 
and a linear contractible cyclic $A_\infty$-algebra. 
 \label{thm:cycMandC}
\end{thm}
Note that this implies 
the symplectic form is also decomposed into the direct sum. 

\begin{pf}
Let us represent the $A_\infty$-odd vector field 
of a cyclic $A_\infty$ algebra $(\cH,\omega,S)$ as 
\begin{equation*}
 \delta=(\ ,S)\ ,\qquad S=S_2+S_3+\cdots\ ,\qquad 
S_k=\ov{k}\V_{i_1\cdots i_k}\phi^{i_k}\cdots\phi^{i_1}\ . 
\end{equation*}
As in the proof of Theorem \ref{thm:MandC}, 
the strategy of the proof is to construct a cyclic $A_\infty$-isomorphism 
so that the induced cyclic $A_\infty$-structure 
$(\cH,\omega,S')$ is of the form 
$S'=S_2+S'_3+S'_4+\cdots$ with $S'_k\in C(p)^\star_c$ for $k\ge 3$. 
We choose 
the coordinates of linear contractible direction as $x$ and $y$, and 
minimal direction as $p$. 
Suppose that $S_k$ belongs to $C(p)^k_c\subset C(\phi)^k_c$ up to $l$, 
that is, 
$S_k=\ov{k}\V_{i_1\cdots i_k}p^{i_k}\cdots p^{i_1}$ for $3\le k\le l$. 
The $(l+1)$ powers-part of the $A_\infty$-condition $\ov{2}(S, S)=0$ implies 
\begin{equation*}
 0=S_{l+1}\flpartial{x^i}c^i_j y^j+
 {\ov 2} \sum_{k\ge 3}(S_k, S_{l+3-k})\ .
\end{equation*}
The first term and second term are independent because 
the first one includes $y^j$ and the second one does not. 
Thus, $S_{l+1}$ is $\delta_1$-closed.  
Recall that by Proposition \ref{prop:hvf} 
any coordinate transformation preserving a constant symplectic structure 
can be written in the form in eq.(\ref{fsymdiff}). 
The transformation 
\begin{equation*}
 S\ \lgraw\ e^{(\ ,\epsilon(\phi))}S(\phi)
 =S(\phi)-\epsilon(\phi)\flpartial{x^i}c^i_j y^j +\cdots\ ,\qquad 
 \epsilon(\phi)=\epsilon_{i_1\cdots i_{l+1}}\phi^{i_{l+1}}\cdots\phi^{i_1}
\end{equation*}
can cancel the exact part of $S_{l+1}$ by an appropriate $\epsilon(\phi)$ 
and we can transform $S_{l+1}$ to be an element of $C(p)^{l+1}_c$. 
Thus the statement of this theorem can be proved by induction. 
\qed
\end{pf}

This procedure is similar to homological perturbation theory in 
\cite{HT}. 

Note that this result is stronger than Theorem \ref{thm:MandC}. 
Namely, it claims that any cyclic $A_\infty$-algebra can be
transformed to the direct sum of a
minimal and a contractible one with the cyclicity (or equivalently 
the constant symplectic structure) being preserved. 
\begin{cor}
For a cyclic $A_\infty$-quasi-isomorphism $\cF^p$ from 
a minimal cyclic $A_\infty$-algebra 
to another cyclic $A_\infty$-algebra, 
there exists an inverse cyclic $A_\infty$-quasi-isomorphism. 
 \label{cor:cyclicinverse}
\end{cor}
\begin{pf}
The proof is just the same as for Corollary \ref{cor:inverse}. 
\end{pf}

 \subsection{Existence of the inverse of $A_\infty$-quasi-isomorphisms}
\label{ssec:inverse}

The decomposition theorem in the previous subsections 
gives clear understanding for the properties of homotopy algebras. 
One usage is to prove the following. 
\begin{thm}[Existence of the inverse quasi-isomorphism]
Let $(\cH,\m)$ and $(\cH',\m')$ be two $A_\infty$-algebras and assume that 
an $A_\infty$-quasi-isomorphism $\cF$ from $(\cH,\m)$ to $(\cH',\m')$ 
is given. Then there exists an inverse $A_\infty$-quasi-isomorphism 
$(\cF)^{-1} : (\cH',\m')\to (\cH,\m)$. 
 \label{thm:inverse}
\end{thm}
The outline of Theorem \ref{thm:inverse} is presented by M.~Kontsevich 
in \cite{Ko1} for $L_\infty$-case. 
The author was noticed by M.~Akaho \cite{Akaho} 
the necessity of the decomposition theorem, 
instead of minimal model theorem, 
for the proof of Theorem \ref{thm:inverse}. \\
\begin{pf}
First, we transform both $A_\infty$-algebras $(\cH,\m)$ and $(\cH',\m')$ 
to their minimal models $(\cH^p,\ti\m^p)$ and $({\cH'}^p,\ti{\m'}^p)$ by 
$A_\infty$-quasi-isomorphisms $\ti\cF^p$ and $\ti{\cF'}^p$ 
in Theorem \ref{thm:MandC}. 
$\ti\cF^p$ and $\ti{\cF'}^p$ have their inverse quasi-isomorphisms, and 
the $A_\infty$-quasi-isomorphism $\cF^p$ from $(\cH^p,\ti\m^p)$ to 
$({\cH'}^p,\ti{\m'}^p)$ 
is then given by the composition $(\ti{\cF'}^p)^{-1}\circ\cF\circ\ti\cF^p$ 
\begin{equation*}
\xymatrix{
  (\cH,\m)\ \ \ar[r]^\cF \ar@<0.5ex>[d]^{(\ti\cF^p)^{-1}}
  & \ \ (\cH',\m') \ar@<0.5ex>[d]^{(\ti{\cF'}^p)^{-1}}\\
 (\cH^p,\ti\m^p)\ \ \ar@<0.5ex>[u]^{\ti\cF^p} \ar@<0.5ex>[r]^{\cF^p} &  
 \ \ ({\cH'}^p,\ti{\m'}^p) \ar@<0.5ex>[u]^{\ti{\cF'}^p} 
 \ar@<0.5ex>[l]^{(\cF^p)^{-1}}
}
\end{equation*}
so that the diagram commutes. 
Because the quasi-isomorphism $\cF^p$ is isomorphism, 
it has its inverse, and one can obtain 
an $A_\infty$-quasi-isomorphism as 
$\ti\cF^p\circ(\cF^p)^{-1}\circ(\ti{\cF'}^p)^{-1}$. 
 \qed
\end{pf}

In the situation $(\cH,\m)=(\cH^p,\ti\m^p)$, $(\cH',\m')=(\cH,\m)$ 
and $\cF'=\ti\cF^p$ the statement of this theorem reduces to 
Corollary \ref{cor:inverse}. 

It is clear that this theorem holds also for cyclic $A_\infty$-algebras.

 \subsection{Maurer-Cartan equations, Feynman graphs 
and the minimal model theorem}
\label{ssec:minimal}

For a given $A_\infty$-algebra $(\cH,\m)$, 
the existence of its minimal $A_\infty$-algebra was shown 
in subsection \ref{ssec:MandC}. 
The existence was proved inductively and an explicit form of 
the minimal $A_\infty$-algebra is unclear. 
In \cite{KoSo} (see also \cite{GLS,hueb-kadei,mer} 
and \cite{jh-jds} for $L_\infty$ case) 
the explicit form is presented by using Feynman graphs. 
In this subsection we shall discuss the `meaning' of 
the minimal $A_\infty$-algebra given explicitly in \cite{KoSo}. 

Here we construct the $A_\infty$-morphism $\{\ti{f}^p_k\}$ and 
$A_\infty$-structure $\{\ti{m}^p_k\}$ with $k\ge 2$ naturally 
as the problem of finding the solutions for 
the Maurer-Cartan equation. 
This explanation of the minimal model theorem 
is inspired by a lecture by K.~Fukaya \cite{Fukayalec} 
(see \cite{Fukaya2}). 
We shall then 
prove that the $(\cH^p,\ti\m^p)$ and $\ti\cF$ are indeed an 
$A_\infty$-algebra and an $A_\infty$-quasi-isomorphism, respectively, 
for the sake of completeness. 
The procedure of finding the solution is quite natural and standard, 
and so similar procedures can be found in various problems. 
In the context of cyclic $A_\infty$/$L_\infty$-algebras 
for open/closed string field theory, 
the Maurer-Cartan equations are 
the equations of motions of the actions, 
and the procedure relates to the way 
of finding some classical solutions ($L_\infty$ case \cite{MS,KZ}) 
or constructing the tachyon potential 
($A_\infty$ case \cite{MT}, see also \cite{Ka}).

Consider solving the Maurer-Cartan equation (MC-eq.) for
an $A_\infty$-algebra:
\begin{equation}
 \sum_{k\ge 1}m_k(\Phi)=0 \label{eomsft}\ .
\end{equation}
Hereafter we often use a shorthand notation $m_k(\Phi)$ for 
$m_k(\Phi,\dots,\Phi)$ as above. 
This is an extended 
MC-eq. of that in Definition \ref{defn:mceq} in the sense
that here we take $\Phi=e_i\phi^i\in\cH\otimes\cH^*$ the superfield in
Definition \ref{defn:superfield}, that is, 
we include $\eb_i$ of any degree and then $\phi^i$ is the formal
noncommutative coordinate. 
As in the previous subsections, consider the Hodge-Kodaira decomposition 
\begin{equation*}
 QQ^++Q^+ Q+P=\1\ \quad (Q:=m_1)
\end{equation*}
and decompose $\cH$ as $\cH=\cH^t\oplus\cH^u\oplus\cH^p$ where 
$\cH^p:=P\cH$, $\cH^t:=QQ^+\cH$ and $\cH^u:=Q^+Q\cH$. 
As seen in subsection \ref{ssec:BV}, in the case of field theory, 
$Q^+$ can be regarded as the propagator and also 
plays the role of the gauge fixing. $P$ is then the projection 
onto the physical states. 
\footnote{
The arguments in this subsection can be generalized 
to some situation where $QQ^++Q^+ Q+P=\1$ still holds but 
the condition that $QQ^+$, $Q^+Q$ and $P$ are
projections is not satisfied. 
As an application of such a generalization to SFT, 
see section 6 of \cite{Ka}, where 
a one parameter family of quasi-isomorphic $A_\infty$-algebras 
is constructed from a given $A_\infty$-algebra 
in terms of tree graphs as in this subsection. 
}

Here we assume that $\Phi$ is sufficiently `small'; 
smaller than the radius of convergence, or equivalently 
$\Phi$ is assumed to be multiplied by the corresponding 
small parameter $\hbar<<1$, 
or instead the $\hbar$ is treated as a formal parameter. 
Then the solution is almost the solution of $Q(\Phi)=0$, 
in the sense that $Q(\Phi)\sim\cO(\hbar^2)$. 
Since the solutions for eq.(\ref{eomsft}) are preserved 
under the gauge transformation 
$\delta_\alpha\Phi=Q(\alpha)+m_2(\alpha,\Phi)+m_2(\Phi,\alpha)+\cdots
\sim Q(\alpha)$, 
we will find the gauge fixed solutions $Q^+\Phi=0$. 
We express the gauge fixed $\Phi$ as $\Phi|_{gf}=\Phi^p+\Phi^u$ 
where $\Phi^p\in \cH^p$ and $\Phi^u\in \cH^u$. 
As explained below, $\Phi^u$ can be solved recursively 
for the power of $\Phi^p$. 
Because here we regard that $\Phi^p$ is `small', 
one can define a degree by the power of $\Phi^p$. 
By substituting $\Phi|_{gf}=\Phi^p+\Phi^u$ 
into $\Phi$, the MC-eq. (\ref{eomsft}) turns out to be 
\begin{equation}
 Q(\Phi^u)+\sum_{k\ge 2}m_k(\Phi^p+\Phi^u)=0\ ,\label{eom2}
\end{equation}
and acting by $Q^+$ on both sides of this equation yields 
\begin{equation}
  \Phi^u=-\sum_{k\ge 2}Q^+ m_k(\Phi^p+\Phi^u)\ . \label{eom3}
\end{equation}
As will be presented explicitly later, 
$\Phi^u$ is expressed in terms of the powers of $\Phi^p$ 
by substituting the right hand side of eq.(\ref{eom3}) 
into $\Phi^u$ in the right hand side of eq.(\ref{eom3}) itself recursively. 
However not all $\Phi|_{gf}=\Phi^p+\Phi^u$ 
expressed in terms of $\Phi^p$ give 
the solution of eq.(\ref{eomsft}) because eq.(\ref{eom3}) is derived from 
`$Q^+$ acting via eq.(\ref{eomsft})'. 
In order to find $\Phi^u$ which is the solution of eq.(\ref{eomsft}), 
we substitute eq.(\ref{eom3}) in the MC-eq.(\ref{eomsft}) 
once again, 
\begin{equation}
 \begin{split}
 0&=Q(\Phi^p+\Phi^u)+ \sum_{k\ge 2}m_k(\Phi)\\
  &=(Q^+ Q+P-1)\sum_{k\ge 2}m_k(\Phi)+\sum_{k\ge 2}m_k(\Phi)\\
  &=Q^+ Q\sum_{k\ge 2}m_k(\Phi) +\sum_{k\ge 2}P m_k(\Phi)
 \end{split}
\label{eom4}
\end{equation}
and we can get a condition (obstruction) for $\Phi^p$. 
The first term in the third line of eq.(\ref{eom4}) vanishes due to 
MC-eq.(\ref{eomsft}) 
because $Q\sum_{k\ge 2}m_k(\Phi)=-QQ(\Phi)=0$, 
and a condition for $\Phi^p$ is derived as 
\begin{equation}
 \sum_{k\ge 2}P m_k(\Phi^p+\Phi^u)=0\ .\label{mceq1}
\end{equation}

The minimal $A_\infty$-algebra $(\cH^p,\ti\m^p)$ and the 
$A_\infty$-quasi-isomorphism $\ti\cF^p: (\cH^p,\ti\m^p)\to (\cH,\m)$ 
are given by eq.(\ref{mceq1}) and eq.(\ref{eom3}), respectively. 
As mentioned above, 
a nonlinear map from $\cH^p$ to $\cH$ is obtained 
by substituting the right hand side of eq.(\ref{eom3}) into the right hand 
side of the equation $\Phi|_{gf}=\Phi^p+\Phi^u$ recursively. 
Here we want to distinguish the element of $\cH^p$ from that of $\cH$, 
so we rewrite $\Phi^p\in\cH^p$ as $\ti\Phi^p\in\cH^p$. 
Let us represent the nonlinear map by a collection of multilinear maps 
$\ti{f}^p_l:(\cH^p)^{\otimes l}\to\cH$ as 
\begin{equation}
\Phi|_{gf}=\ti{f}^p_1(\ti\Phi^p)+\ti{f}^p_2(\ti\Phi^p,\ti\Phi^p)
+\ti{f}^p_3(\ti\Phi^p,\ti\Phi^p,\ti\Phi^p)+\cdots
 \label{eom33}
\end{equation}
where $\ti{f}^p_1:\cH^p\to\cH$ the inclusion map. 
Alternatively, eq.(\ref{mceq1}) is also 
expressed as an equation for $\ti\Phi^p$ 
by substituting eq.(\ref{eom33}) into eq.(\ref{mceq1}). 
Let us write it also by a collection of multilinear maps 
$\ti{m}^p_l:(\cH^p)^{\otimes l}\to\cH^p$ as 
\begin{equation}
 \iota\ \sum_{k\ge 2}\ti{m}^p_k(\ti\Phi^p)=0\ ,\label{mceq2}
\end{equation}
where $\iota:\cH^p\to\cH$ is the inclusion. 
Once $\ti{m}^p_k(\ti\Phi^p)$ and $\ti{f}^p_k(\ti\Phi^p)$ 
are obtained, 
$\ti{m}^p_k(\eb^p_{i_1},\dots,\eb^p_{i_k})$ and 
$\ti{f}^p_k(\eb^p_{i_1},\dots,\eb^p_{i_k})$ can be obtained immediately 
as the coefficient of $\phi^{i_k}\cdots\phi^{i_1}$, 
where $\eb_i^p$ are bases of $\cH^p$ and 
$\ti\Phi^p=\eb^p_i\phi^i$.
As shown in the end of this subsection, 
$(\cH^p,\ti\m^p:=\{\ti{m}^p_k\}_{k\ge 2})$ forms a minimal 
$A_\infty$-algebra. 
The equation (\ref{mceq2}) is just the Maurer-Cartan equation 
for $(\cH^p,\ti\m^p)$. 
$\ti\cF^p:=\{\ti{f}^p_l\}_{l\ge 1}$ is then an $A_\infty$-quasi-isomorphism 
from $(\cH^p,\ti\m^p)$ to $(\cH,\m)$.

{}From the field theory point of view the above result means that 
if the expectation value of physical states satisfying 
the Maurer-Cartan equation (\ref{mceq2}) is given, 
the solution of the equations of motions for field theory (\ref{eomsft}) 
is obtained by the $A_\infty$-quasi-isomorphism (\ref{eom33}). 

Here we summarize the arguments above and define 
the $A_\infty$-structure and $A_\infty$-quasi-isomorphism precisely. 
\begin{defn}[An explicit minimal model]
Given an $A_\infty$-algebra $(\cH,\m)$, 
assume that we have a splitting of the complex $(\cH,Q)$ such that 
$\iota:\cH^p\to\cH$ is the inclusion, $\pi:\cH\to\cH^p$ is the 
projection and $Q^+:\cH\to\cH$ is the contracting homotopy 
$\1-P=QQ^++Q^+Q$ for $P=\iota\circ\pi$. 
Then, 
an $A_\infty$-structure $\ti\m^p$ on the cohomology $\cH^p$ and an 
$A_\infty$-quasi-isomorphism $\ti\cF^p :(\cH^p,\ti\m^p)\to(\cH,\m)$ 
are constructed as follows. 
$\ti\cF^p=\{\ti{f}^p_l:(\cH^p)^{\otimes l}\to\cH\}_{l\ge 1}$ 
is defined recursively 
with respect to $k$ as 
 \begin{equation*}
 \ti{f}^p_k =-Q^+\sum_{i\ge 2}\sum_{1\le k_1<k_2\cdots<k_i=k}
 m_i(\ti{f}^p_{k_1}\otimes\ti{f}^p_{k_2-k_1}\otimes
 \cdots\otimes\ti{f}^p_{k-k_{i-1}})
 \end{equation*}
with $\ti{f}^p_1:=\iota:\cH^p\to\cH$. 
 $\ti\m^p=\{\ti{m}^p_k:(\cH^p)^{\otimes k}\to\cH\}_{k\ge 2}$ 
is then given by 
 \begin{equation*}
  \ti{m}^p_k =\pi\sum_{i\ge 2}\sum_{1\le k_1<k_2\cdots <k_i=k}
 m_i(\ti{f}^p_{k_1}\otimes\ti{f}^p_{k_2-k_1}\otimes
 \cdots\otimes\ti{f}^p_{k-k_{i-1}})\ .
 \end{equation*}
 \label{defn:minimal}
\end{defn}
It is also convenient 
to present an alternative description of 
Definition \ref{defn:minimal} in terms of rooted planar tree graphs. 
It is related to Feynman graphs in field theory 
in section \ref{sec:BV}. 
\begin{defn}[Rooted planar tree graph]
A rooted planar tree graph is a simply connected 
rooted planar tree without loops. 
It consists of vertices, internal edges and external edges. 
Both ends of an internal edge are on two vertices. 
An external edge has one end on a vertex and another end is free. 
The number of incident edges at a vertex is greater than three. 
The term `planar' means the cyclic order of edges at each vertex is 
distinguished. 
A rooted planar tree has a {\it root} that is a free end of an 
external edge. The external edge is called the {\it root edge}. 
The vertex on which the root edge ends is the {\it root vertex}. 
The free ends of the remaining external edges are called the {\it leaves}. 
\begin{figure}[h]
 \hspace*{0.7cm}
 \includegraphics{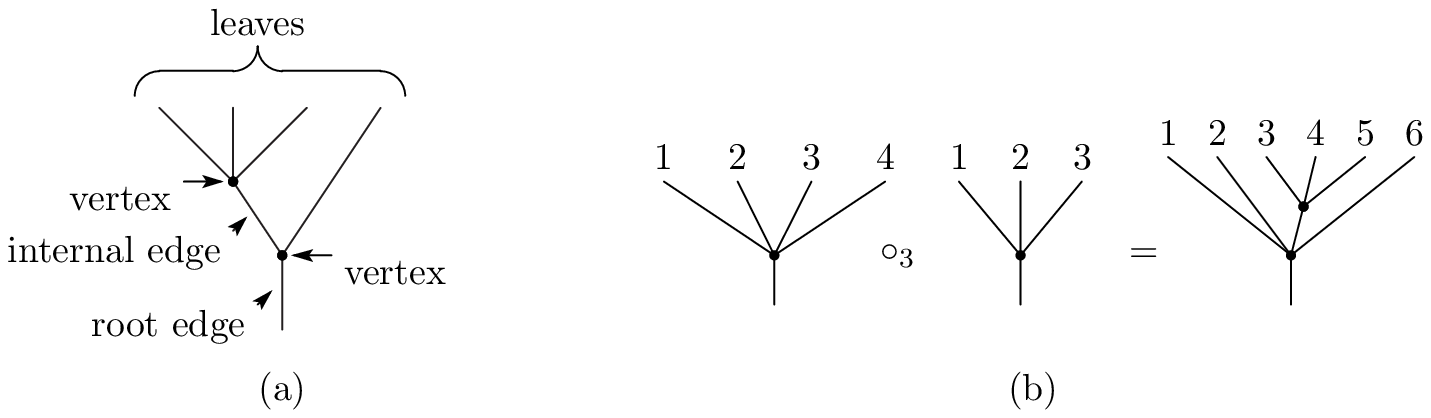}
 \caption{(a). Notation for planar rooted tree. The above one is a 
$4$-tree. (b). An example of grafting, grafting of $3$-corolla to 
$4$-corolla along leaf $3$.}
 \label{fig:tree2}
\end{figure}
We call a rooted planar tree that has $k$ leaves a $k$-tree and 
denote it by $\Gamma_k$. 
We denote by $G_k$ the set of $k$-trees. 

For $k$-tree $\Gamma_k$, $l$-tree $\Gamma_l$ and an integer $1\le i\le k$, 
the {\it grafting} of $l$-tree to $k$-tree 
along leaf $i$ is given by identifying the root edge of the $\Gamma_l$ 
with the $i$-th leaf of $\Gamma_k$ (see Figure \ref{fig:tree2} (b)). 
The resulting $(k+l-1)$-tree is denoted by $\Gamma_k\circ_i\Gamma_l$. 

$G_1$ has only one element $|$ that has no vertex. 
A $k$-tree that has only one vertex is called a {\it $k$-corolla}. 
Any other tree, that has more than one vertices, 
is obtained by grafting corollas. 
 \label{defn:planar}
\end{defn}
\begin{defn}[Minimal model; an alternative description]
Let us define a map $\ti{f}: G_k\to (\cH^{\otimes k}\to\cH)$ 
as follows. 
To the $1$-tree (that has no vertex) we associate the identity operator 
$\Id:\cH\to\cH$. 
To a $k$-corolla, we associate $-Q^+m_k:\cH^{\otimes k}\to\cH$. 
For any $k$-tree $\Gamma_k$ and $l$-tree $\Gamma_l$, 
denote the associated endomorphisms by 
$\ti{f}_{\Gamma_k}\in (\cH^{\otimes k}\to\cH)$ and 
$\ti{f}_{\Gamma_l}\in (\cH^{\otimes l}\to\cH)$. 
$\ti{f}$ is then defined so that it is compatible with the grafting of 
the trees. 
Namely, to $\Gamma_k\circ_i\Gamma_l$ we associate 
\begin{equation*}
 \ti{f}_{\Gamma_k\circ_i\Gamma_l}=
 \ti{f}_{\Gamma_k}\circ
\l(\1^{\otimes (i-1)}\otimes\ti{f}_{\Gamma_l}\otimes\1^{\otimes (k-i)}\r)
: \cH^{\otimes (k+l-1)}\to\cH\ .
\end{equation*}
Thus, for any $k$-tree, $\ti{f}_{\Gamma_k}$ is defined. 
The $A_\infty$-quasi-isomorphism $\ti\cF^p=\{\ti{f}^p_k\}_{k\ge 1}$ 
is then defined by 
\begin{equation*}
 \ti{f}^p_k=\sum_{\Gamma_k\in G_k}\ti{f}_{\Gamma_k}\circ 
(\iota)^{\otimes k}\ ,
\end{equation*}
where $\iota :\cH^p\to\cH$ is the inclusion. 

The minimal $A_\infty$-structure $\ti\m^p$ is defined by 
using another map 
$\ti{m}: G_k\to (\cH^{\otimes k}\to\cH)$. 

To the $1$-tree we associate the differential $Q:\cH\to\cH$. 
To a $k$-corolla, we associate 
$m_k:\cH^{\otimes k}\to\cH$. 
To any $k$-tree $\Gamma_k$ denote the associated endomorphisms by 
$\ti{m}_{\Gamma_k}\in (\cH^{\otimes k}\to\cH)$. 
For a grafting $\Gamma_k\circ_i\Gamma_l$ we associate 
\begin{equation*}
 \ti{m}_{\Gamma_k\circ_i\Gamma_l}=
 \ti{m}_{\Gamma_k}\circ
\l(\1^{\otimes (i-1)}\otimes\ti{f}_{\Gamma_l}\otimes\1^{\otimes (k-i)}\r)
 : \cH^{\otimes (k+l-1)}\to\cH\ .
\end{equation*}
Thus, for any tree graph, 
$\ti{m}$ is defined so that it is compatible with this grafting. 
$\ti\m^p$ is then given by 
\begin{equation*}
 \ti{m}^p_k=\pi\circ\sum_{\Gamma_k\in G_k}\ti{m}_{\Gamma_k}\circ 
(\iota)^{\otimes k}\ ,
\end{equation*}
where $\pi:\cH\to\cH^p$ is the projection. 
Note that $\ti{m}^p_1$ automatically vanishes. 

As a result, for a given $l$-tree $\Gamma_l$, $\ti{m}^p_{\Gamma_l}$ 
is given by attaching $m_k$ to each vertex that has $(k+1)$-incident edges, 
$-Q^+$ to each internal edge, $\iota$ to each leave and 
$\pi$ to the root edge. 
$\ti{f}^p_{\Gamma_l}$ is also given in the same way but replacing 
$\pi$ on the root edge to $-Q^+$. 
 \label{defn:minimal2}
\end{defn}

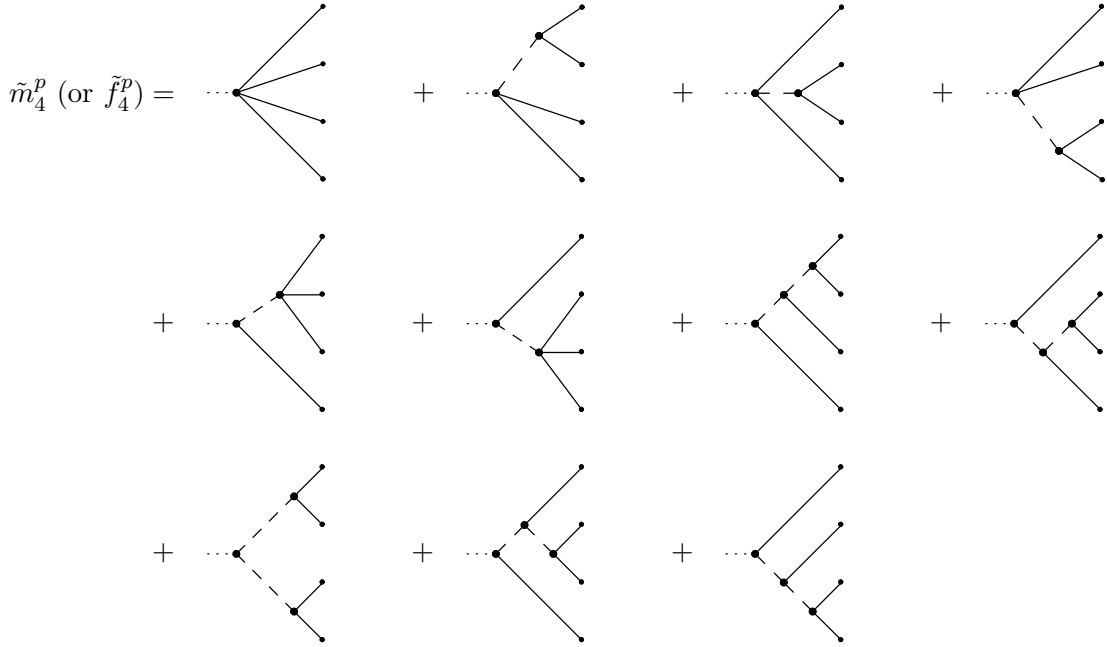
\begin{figure}[h]
 \hspace*{-2.5cm}
%WinTpicVersion3.08
\unitlength 0.1in
\begin{picture}( 69.6850, 33.2087)(-13.6220,-36.8996)
% STR 2 0 3 0
% 3 326 762 326 838 5 0
% $\ti{m}^p_4\ (\mbox{or}\ \ti{f}^p_4)=$
\put(3.2087,-8.2480){\makebox(0,0){$\ti{m}^p_4\ (\mbox{or}\ \ti{f}^p_4)=$}}%
% LINE 2 0 3 0
% 2 4309 2523 3850 2064
% 
\special{pn 8}%
\special{pa 4242 2484}%
\special{pa 3790 2032}%
\special{fp}%
% LINE 2 0 3 0
% 2 4309 2830 3850 3289
% 
\special{pn 8}%
\special{pa 4242 2786}%
\special{pa 3790 3238}%
\special{fp}%
% DOT 2 0 3 0
% 5 1552 375 1552 681 1552 988 1552 1294 1552 1294
% 
\special{pn 8}%
\special{sh 1}%
\special{ar 1528 370 10 10 0  6.28318530717959E+0000}%
\special{sh 1}%
\special{ar 1528 671 10 10 0  6.28318530717959E+0000}%
\special{sh 1}%
\special{ar 1528 973 10 10 0  6.28318530717959E+0000}%
\special{sh 1}%
\special{ar 1528 1274 10 10 0  6.28318530717959E+0000}%
\special{sh 1}%
\special{ar 1528 1274 10 10 0  6.28318530717959E+0000}%
% LINE 2 0 3 0
% 2 1552 375 1092 835
% 
\special{pn 8}%
\special{pa 1528 370}%
\special{pa 1075 822}%
\special{fp}%
% LINE 2 0 3 0
% 2 1552 1294 1092 835
% 
\special{pn 8}%
\special{pa 1528 1274}%
\special{pa 1075 822}%
\special{fp}%
% LINE 2 0 3 0
% 2 4309 1911 4156 1757
% 
\special{pn 8}%
\special{pa 4242 1881}%
\special{pa 4091 1730}%
\special{fp}%
% LINE 2 0 3 0
% 2 1552 3136 1399 2983
% 
\special{pn 8}%
\special{pa 1528 3087}%
\special{pa 1377 2937}%
\special{fp}%
% LINE 2 0 3 0
% 2 4309 2217 4003 1911
% 
\special{pn 8}%
\special{pa 4242 2183}%
\special{pa 3940 1881}%
\special{fp}%
% LINE 2 0 3 0
% 2 1552 3442 1399 3596
% 
\special{pn 8}%
\special{pa 1528 3388}%
\special{pa 1377 3540}%
\special{fp}%
% LINE 2 0 3 0
% 2 4309 3136 4003 3442
% 
\special{pn 8}%
\special{pa 4242 3087}%
\special{pa 3940 3388}%
\special{fp}%
% LINE 2 0 3 0
% 2 4309 3442 4156 3596
% 
\special{pn 8}%
\special{pa 4242 3388}%
\special{pa 4091 3540}%
\special{fp}%
% LINE 2 0 3 0
% 2 1552 681 1092 835
% 
\special{pn 8}%
\special{pa 1528 671}%
\special{pa 1075 822}%
\special{fp}%
% LINE 2 0 3 0
% 2 1552 988 1092 835
% 
\special{pn 8}%
\special{pa 1528 973}%
\special{pa 1075 822}%
\special{fp}%
% LINE 2 0 3 0
% 2 1552 2217 1322 1911
% 
\special{pn 8}%
\special{pa 1528 2183}%
\special{pa 1302 1881}%
\special{fp}%
% LINE 2 0 3 0
% 2 1552 1604 1322 1911
% 
\special{pn 8}%
\special{pa 1528 1579}%
\special{pa 1302 1881}%
\special{fp}%
% LINE 2 0 3 0
% 2 1552 1911 1322 1911
% 
\special{pn 8}%
\special{pa 1528 1881}%
\special{pa 1302 1881}%
\special{fp}%
% LINE 2 0 3 0
% 2 1092 2064 1552 2523
% 
\special{pn 8}%
\special{pa 1075 2032}%
\special{pa 1528 2484}%
\special{fp}%
% LINE 2 1 3 0
% 2 1092 2064 1322 1911
% 
\special{pn 8}%
\special{pa 1075 2032}%
\special{pa 1302 1881}%
\special{da 0.070}%
% LINE 2 0 3 0
% 2 2930 1911 2701 2217
% 
\special{pn 8}%
\special{pa 2884 1881}%
\special{pa 2659 2183}%
\special{fp}%
% LINE 2 0 3 0
% 2 2930 2523 2701 2217
% 
\special{pn 8}%
\special{pa 2884 2484}%
\special{pa 2659 2183}%
\special{fp}%
% LINE 2 0 3 0
% 2 2930 2217 2701 2217
% 
\special{pn 8}%
\special{pa 2884 2183}%
\special{pa 2659 2183}%
\special{fp}%
% LINE 2 0 3 0
% 2 2471 2064 2930 1604
% 
\special{pn 8}%
\special{pa 2433 2032}%
\special{pa 2884 1579}%
\special{fp}%
% LINE 2 1 3 0
% 2 2471 2064 2701 2217
% 
\special{pn 8}%
\special{pa 2433 2032}%
\special{pa 2659 2183}%
\special{da 0.070}%
% LINE 2 2 3 0
% 2 2471 2064 2318 2064
% 
\special{pn 8}%
\special{pa 2433 2032}%
\special{pa 2282 2032}%
\special{dt 0.045}%
% LINE 2 2 3 0
% 2 1092 2064 939 2064
% 
\special{pn 8}%
\special{pa 1075 2032}%
\special{pa 925 2032}%
\special{dt 0.045}%
% LINE 2 2 3 0
% 2 1092 3289 939 3289
% 
\special{pn 8}%
\special{pa 1075 3238}%
\special{pa 925 3238}%
\special{dt 0.045}%
% LINE 2 2 3 0
% 2 1092 835 939 835
% 
\special{pn 8}%
\special{pa 1075 822}%
\special{pa 925 822}%
\special{dt 0.045}%
% LINE 2 2 3 0
% 2 3850 3289 3696 3289
% 
\special{pn 8}%
\special{pa 3790 3238}%
\special{pa 3638 3238}%
\special{dt 0.045}%
% LINE 2 2 3 0
% 2 3850 2064 3696 2064
% 
\special{pn 8}%
\special{pa 3790 2032}%
\special{pa 3638 2032}%
\special{dt 0.045}%
% STR 2 0 3 0
% 3 2088 762 2088 838 5 0
% $+$
\put(20.5512,-8.2480){\makebox(0,0){$+$}}%
% STR 2 0 3 0
% 3 3467 762 3467 838 5 0
% $+$
\put(34.1240,-8.2480){\makebox(0,0){$+$}}%
% LINE 2 0 3 0
% 2 4309 1604 4156 1757
% 
\special{pn 8}%
\special{pa 4242 1579}%
\special{pa 4091 1730}%
\special{fp}%
% LINE 2 1 3 0
% 2 4156 1757 4003 1911
% 
\special{pn 8}%
\special{pa 4091 1730}%
\special{pa 3940 1881}%
\special{da 0.070}%
% LINE 2 1 3 0
% 2 4003 1911 3850 2064
% 
\special{pn 8}%
\special{pa 3940 1881}%
\special{pa 3790 2032}%
\special{da 0.070}%
% LINE 2 0 3 0
% 2 1552 2830 1399 2983
% 
\special{pn 8}%
\special{pa 1528 2786}%
\special{pa 1377 2937}%
\special{fp}%
% LINE 2 0 3 0
% 2 1552 3749 1399 3596
% 
\special{pn 8}%
\special{pa 1528 3690}%
\special{pa 1377 3540}%
\special{fp}%
% LINE 2 0 3 0
% 2 4309 3749 4156 3596
% 
\special{pn 8}%
\special{pa 4242 3690}%
\special{pa 4091 3540}%
\special{fp}%
% LINE 2 1 3 0
% 2 1399 3596 1092 3289
% 
\special{pn 8}%
\special{pa 1377 3540}%
\special{pa 1075 3238}%
\special{da 0.070}%
% LINE 2 1 3 0
% 2 1399 2983 1092 3289
% 
\special{pn 8}%
\special{pa 1377 2937}%
\special{pa 1075 3238}%
\special{da 0.070}%
% LINE 2 1 3 0
% 2 3850 3289 4003 3442
% 
\special{pn 8}%
\special{pa 3790 3238}%
\special{pa 3940 3388}%
\special{da 0.070}%
% LINE 2 1 3 0
% 2 4003 3442 4156 3596
% 
\special{pn 8}%
\special{pa 3940 3388}%
\special{pa 4091 3540}%
\special{da 0.070}%
% DOT 0 1 3 0
% 2 1092 835 1092 835
% 
\special{pn 20}%
\special{sh 1}%
\special{ar 1075 822 10 10 0  6.28318530717959E+0000}%
\special{sh 1}%
\special{ar 1075 822 10 10 0  6.28318530717959E+0000}%
% DOT 0 0 3 0
% 2 1322 1911 1322 1911
% 
\special{pn 20}%
\special{sh 1}%
\special{ar 1302 1881 10 10 0  6.28318530717959E+0000}%
\special{sh 1}%
\special{ar 1302 1881 10 10 0  6.28318530717959E+0000}%
% DOT 0 0 3 0
% 2 1092 2064 1092 2064
% 
\special{pn 20}%
\special{sh 1}%
\special{ar 1075 2032 10 10 0  6.28318530717959E+0000}%
\special{sh 1}%
\special{ar 1075 2032 10 10 0  6.28318530717959E+0000}%
% DOT 0 0 3 0
% 2 2471 2064 2471 2064
% 
\special{pn 20}%
\special{sh 1}%
\special{ar 2433 2032 10 10 0  6.28318530717959E+0000}%
\special{sh 1}%
\special{ar 2433 2032 10 10 0  6.28318530717959E+0000}%
% DOT 0 0 3 0
% 2 2701 2217 2701 2217
% 
\special{pn 20}%
\special{sh 1}%
\special{ar 2659 2183 10 10 0  6.28318530717959E+0000}%
\special{sh 1}%
\special{ar 2659 2183 10 10 0  6.28318530717959E+0000}%
% DOT 0 0 3 0
% 2 3850 3289 3850 3289
% 
\special{pn 20}%
\special{sh 1}%
\special{ar 3790 3238 10 10 0  6.28318530717959E+0000}%
\special{sh 1}%
\special{ar 3790 3238 10 10 0  6.28318530717959E+0000}%
% DOT 0 0 3 0
% 2 1092 3289 1092 3289
% 
\special{pn 20}%
\special{sh 1}%
\special{ar 1075 3238 10 10 0  6.28318530717959E+0000}%
\special{sh 1}%
\special{ar 1075 3238 10 10 0  6.28318530717959E+0000}%
% DOT 0 0 3 0
% 2 3850 2064 3850 2064
% 
\special{pn 20}%
\special{sh 1}%
\special{ar 3790 2032 10 10 0  6.28318530717959E+0000}%
\special{sh 1}%
\special{ar 3790 2032 10 10 0  6.28318530717959E+0000}%
% DOT 0 0 3 0
% 2 4156 1757 4156 1757
% 
\special{pn 20}%
\special{sh 1}%
\special{ar 4091 1730 10 10 0  6.28318530717959E+0000}%
\special{sh 1}%
\special{ar 4091 1730 10 10 0  6.28318530717959E+0000}%
% DOT 0 0 3 0
% 2 4003 1911 4003 1911
% 
\special{pn 20}%
\special{sh 1}%
\special{ar 3940 1881 10 10 0  6.28318530717959E+0000}%
\special{sh 1}%
\special{ar 3940 1881 10 10 0  6.28318530717959E+0000}%
% DOT 0 0 3 0
% 2 1399 2983 1399 2983
% 
\special{pn 20}%
\special{sh 1}%
\special{ar 1377 2937 10 10 0  6.28318530717959E+0000}%
\special{sh 1}%
\special{ar 1377 2937 10 10 0  6.28318530717959E+0000}%
% DOT 0 0 3 0
% 2 1399 3596 1399 3596
% 
\special{pn 20}%
\special{sh 1}%
\special{ar 1377 3540 10 10 0  6.28318530717959E+0000}%
\special{sh 1}%
\special{ar 1377 3540 10 10 0  6.28318530717959E+0000}%
% DOT 0 0 3 0
% 2 4003 3442 4003 3442
% 
\special{pn 20}%
\special{sh 1}%
\special{ar 3940 3388 10 10 0  6.28318530717959E+0000}%
\special{sh 1}%
\special{ar 3940 3388 10 10 0  6.28318530717959E+0000}%
% DOT 0 0 3 0
% 2 4156 3596 4156 3596
% 
\special{pn 20}%
\special{sh 1}%
\special{ar 4091 3540 10 10 0  6.28318530717959E+0000}%
\special{sh 1}%
\special{ar 4091 3540 10 10 0  6.28318530717959E+0000}%
% LINE 2 2 3 0
% 2 2318 838 2471 838
% 
\special{pn 8}%
\special{pa 2282 825}%
\special{pa 2433 825}%
\special{dt 0.045}%
% LINE 2 2 3 0
% 2 3696 838 3850 838
% 
\special{pn 8}%
\special{pa 3638 825}%
\special{pa 3790 825}%
\special{dt 0.045}%
% LINE 2 0 3 0
% 2 2930 379 2701 532
% 
\special{pn 8}%
\special{pa 2884 374}%
\special{pa 2659 524}%
\special{fp}%
% LINE 2 0 3 0
% 2 2930 685 2701 532
% 
\special{pn 8}%
\special{pa 2884 675}%
\special{pa 2659 524}%
\special{fp}%
% LINE 2 0 3 0
% 2 2930 992 2471 838
% 
\special{pn 8}%
\special{pa 2884 977}%
\special{pa 2433 825}%
\special{fp}%
% LINE 2 0 3 0
% 2 2930 1298 2471 838
% 
\special{pn 8}%
\special{pa 2884 1278}%
\special{pa 2433 825}%
\special{fp}%
% LINE 2 1 3 0
% 2 2701 532 2471 838
% 
\special{pn 8}%
\special{pa 2659 524}%
\special{pa 2433 825}%
\special{da 0.070}%
% DOT 0 0 3 0
% 2 2471 838 2471 838
% 
\special{pn 20}%
\special{sh 1}%
\special{ar 2433 825 10 10 0  6.28318530717959E+0000}%
\special{sh 1}%
\special{ar 2433 825 10 10 0  6.28318530717959E+0000}%
% DOT 0 0 3 0
% 2 2701 532 2701 532
% 
\special{pn 20}%
\special{sh 1}%
\special{ar 2659 524 10 10 0  6.28318530717959E+0000}%
\special{sh 1}%
\special{ar 2659 524 10 10 0  6.28318530717959E+0000}%
% DOT 2 0 3 0
% 2 2930 379 2930 379
% 
\special{pn 8}%
\special{sh 1}%
\special{ar 2884 374 10 10 0  6.28318530717959E+0000}%
\special{sh 1}%
\special{ar 2884 374 10 10 0  6.28318530717959E+0000}%
% DOT 2 0 3 0
% 2 2930 685 2930 685
% 
\special{pn 8}%
\special{sh 1}%
\special{ar 2884 675 10 10 0  6.28318530717959E+0000}%
\special{sh 1}%
\special{ar 2884 675 10 10 0  6.28318530717959E+0000}%
% DOT 2 0 3 0
% 2 2930 992 2930 992
% 
\special{pn 8}%
\special{sh 1}%
\special{ar 2884 977 10 10 0  6.28318530717959E+0000}%
\special{sh 1}%
\special{ar 2884 977 10 10 0  6.28318530717959E+0000}%
% DOT 2 0 3 0
% 2 2930 1298 2930 1298
% 
\special{pn 8}%
\special{sh 1}%
\special{ar 2884 1278 10 10 0  6.28318530717959E+0000}%
\special{sh 1}%
\special{ar 2884 1278 10 10 0  6.28318530717959E+0000}%
% LINE 2 0 3 0
% 2 3850 838 4309 1298
% 
\special{pn 8}%
\special{pa 3790 825}%
\special{pa 4242 1278}%
\special{fp}%
% LINE 2 0 3 0
% 2 3850 838 4309 379
% 
\special{pn 8}%
\special{pa 3790 825}%
\special{pa 4242 374}%
\special{fp}%
% LINE 2 1 3 0
% 2 3850 838 4079 838
% 
\special{pn 8}%
\special{pa 3790 825}%
\special{pa 4015 825}%
\special{da 0.070}%
% LINE 2 0 3 0
% 2 4079 838 4309 685
% 
\special{pn 8}%
\special{pa 4015 825}%
\special{pa 4242 675}%
\special{fp}%
% LINE 2 0 3 0
% 2 4079 838 4309 992
% 
\special{pn 8}%
\special{pa 4015 825}%
\special{pa 4242 977}%
\special{fp}%
% DOT 2 0 3 0
% 2 4309 379 4309 379
% 
\special{pn 8}%
\special{sh 1}%
\special{ar 4242 374 10 10 0  6.28318530717959E+0000}%
\special{sh 1}%
\special{ar 4242 374 10 10 0  6.28318530717959E+0000}%
% DOT 2 0 3 0
% 2 4309 685 4309 685
% 
\special{pn 8}%
\special{sh 1}%
\special{ar 4242 675 10 10 0  6.28318530717959E+0000}%
\special{sh 1}%
\special{ar 4242 675 10 10 0  6.28318530717959E+0000}%
% DOT 2 0 3 0
% 2 4309 992 4309 992
% 
\special{pn 8}%
\special{sh 1}%
\special{ar 4242 977 10 10 0  6.28318530717959E+0000}%
\special{sh 1}%
\special{ar 4242 977 10 10 0  6.28318530717959E+0000}%
% DOT 2 0 3 0
% 2 4309 1298 4309 1298
% 
\special{pn 8}%
\special{sh 1}%
\special{ar 4242 1278 10 10 0  6.28318530717959E+0000}%
\special{sh 1}%
\special{ar 4242 1278 10 10 0  6.28318530717959E+0000}%
% DOT 0 0 3 0
% 2 3850 838 3850 838
% 
\special{pn 20}%
\special{sh 1}%
\special{ar 3790 825 10 10 0  6.28318530717959E+0000}%
\special{sh 1}%
\special{ar 3790 825 10 10 0  6.28318530717959E+0000}%
% DOT 0 0 3 0
% 2 4079 838 4079 838
% 
\special{pn 20}%
\special{sh 1}%
\special{ar 4015 825 10 10 0  6.28318530717959E+0000}%
\special{sh 1}%
\special{ar 4015 825 10 10 0  6.28318530717959E+0000}%
% STR 2 0 3 0
% 3 4849 758 4849 835 5 0
% $+$
\put(47.7264,-8.2185){\makebox(0,0){$+$}}%
% DOT 2 0 3 0
% 2 5692 375 5692 375
% 
\special{pn 8}%
\special{sh 1}%
\special{ar 5603 370 10 10 0  6.28318530717959E+0000}%
\special{sh 1}%
\special{ar 5603 370 10 10 0  6.28318530717959E+0000}%
% DOT 2 0 3 0
% 2 5692 681 5692 681
% 
\special{pn 8}%
\special{sh 1}%
\special{ar 5603 671 10 10 0  6.28318530717959E+0000}%
\special{sh 1}%
\special{ar 5603 671 10 10 0  6.28318530717959E+0000}%
% DOT 2 0 3 0
% 2 5692 988 5692 988
% 
\special{pn 8}%
\special{sh 1}%
\special{ar 5603 973 10 10 0  6.28318530717959E+0000}%
\special{sh 1}%
\special{ar 5603 973 10 10 0  6.28318530717959E+0000}%
% DOT 2 0 3 0
% 2 5696 1298 5696 1298
% 
\special{pn 8}%
\special{sh 1}%
\special{ar 5607 1278 10 10 0  6.28318530717959E+0000}%
\special{sh 1}%
\special{ar 5607 1278 10 10 0  6.28318530717959E+0000}%
% STR 2 0 3 0
% 3 705 1983 705 2060 5 0
% $+$
\put(6.9390,-20.2756){\makebox(0,0){$+$}}%
% DOT 2 0 3 0
% 2 1548 1600 1548 1600
% 
\special{pn 8}%
\special{sh 1}%
\special{ar 1524 1575 10 10 0  6.28318530717959E+0000}%
\special{sh 1}%
\special{ar 1524 1575 10 10 0  6.28318530717959E+0000}%
% DOT 2 0 3 0
% 2 1548 1907 1548 1907
% 
\special{pn 8}%
\special{sh 1}%
\special{ar 1524 1877 10 10 0  6.28318530717959E+0000}%
\special{sh 1}%
\special{ar 1524 1877 10 10 0  6.28318530717959E+0000}%
% DOT 2 0 3 0
% 2 1548 2213 1548 2213
% 
\special{pn 8}%
\special{sh 1}%
\special{ar 1524 2179 10 10 0  6.28318530717959E+0000}%
\special{sh 1}%
\special{ar 1524 2179 10 10 0  6.28318530717959E+0000}%
% DOT 2 0 3 0
% 2 1548 2520 1548 2520
% 
\special{pn 8}%
\special{sh 1}%
\special{ar 1524 2481 10 10 0  6.28318530717959E+0000}%
\special{sh 1}%
\special{ar 1524 2481 10 10 0  6.28318530717959E+0000}%
% STR 2 0 3 0
% 3 2084 1983 2084 2060 5 0
% $+$
\put(20.5118,-20.2756){\makebox(0,0){$+$}}%
% DOT 2 0 3 0
% 2 2927 1600 2927 1600
% 
\special{pn 8}%
\special{sh 1}%
\special{ar 2881 1575 10 10 0  6.28318530717959E+0000}%
\special{sh 1}%
\special{ar 2881 1575 10 10 0  6.28318530717959E+0000}%
% DOT 2 0 3 0
% 2 2927 1907 2927 1907
% 
\special{pn 8}%
\special{sh 1}%
\special{ar 2881 1877 10 10 0  6.28318530717959E+0000}%
\special{sh 1}%
\special{ar 2881 1877 10 10 0  6.28318530717959E+0000}%
% DOT 2 0 3 0
% 2 2927 2213 2927 2213
% 
\special{pn 8}%
\special{sh 1}%
\special{ar 2881 2179 10 10 0  6.28318530717959E+0000}%
\special{sh 1}%
\special{ar 2881 2179 10 10 0  6.28318530717959E+0000}%
% DOT 2 0 3 0
% 2 2927 2520 2927 2520
% 
\special{pn 8}%
\special{sh 1}%
\special{ar 2881 2481 10 10 0  6.28318530717959E+0000}%
\special{sh 1}%
\special{ar 2881 2481 10 10 0  6.28318530717959E+0000}%
% LINE 2 2 3 0
% 2 5083 838 5236 838
% 
\special{pn 8}%
\special{pa 5003 825}%
\special{pa 5154 825}%
\special{dt 0.045}%
% LINE 2 0 3 0
% 2 5696 1298 5466 1145
% 
\special{pn 8}%
\special{pa 5607 1278}%
\special{pa 5380 1127}%
\special{fp}%
% LINE 2 0 3 0
% 2 5696 992 5466 1145
% 
\special{pn 8}%
\special{pa 5607 977}%
\special{pa 5380 1127}%
\special{fp}%
% LINE 2 0 3 0
% 2 5696 685 5236 838
% 
\special{pn 8}%
\special{pa 5607 675}%
\special{pa 5154 825}%
\special{fp}%
% LINE 2 0 3 0
% 2 5696 379 5236 838
% 
\special{pn 8}%
\special{pa 5607 374}%
\special{pa 5154 825}%
\special{fp}%
% LINE 2 1 3 0
% 2 5466 1145 5236 838
% 
\special{pn 8}%
\special{pa 5380 1127}%
\special{pa 5154 825}%
\special{da 0.070}%
% DOT 0 0 3 0
% 2 5236 838 5236 838
% 
\special{pn 20}%
\special{sh 1}%
\special{ar 5154 825 10 10 0  6.28318530717959E+0000}%
\special{sh 1}%
\special{ar 5154 825 10 10 0  6.28318530717959E+0000}%
% DOT 0 0 3 0
% 2 5466 1145 5466 1145
% 
\special{pn 20}%
\special{sh 1}%
\special{ar 5380 1127 10 10 0  6.28318530717959E+0000}%
\special{sh 1}%
\special{ar 5380 1127 10 10 0  6.28318530717959E+0000}%
% STR 2 0 3 0
% 3 3463 1983 3463 2060 5 0
% $+$
\put(34.0846,-20.2756){\makebox(0,0){$+$}}%
% DOT 2 0 3 0
% 2 4305 1600 4305 1600
% 
\special{pn 8}%
\special{sh 1}%
\special{ar 4238 1575 10 10 0  6.28318530717959E+0000}%
\special{sh 1}%
\special{ar 4238 1575 10 10 0  6.28318530717959E+0000}%
% DOT 2 0 3 0
% 2 4305 1907 4305 1907
% 
\special{pn 8}%
\special{sh 1}%
\special{ar 4238 1877 10 10 0  6.28318530717959E+0000}%
\special{sh 1}%
\special{ar 4238 1877 10 10 0  6.28318530717959E+0000}%
% DOT 2 0 3 0
% 2 4305 2213 4305 2213
% 
\special{pn 8}%
\special{sh 1}%
\special{ar 4238 2179 10 10 0  6.28318530717959E+0000}%
\special{sh 1}%
\special{ar 4238 2179 10 10 0  6.28318530717959E+0000}%
% DOT 2 0 3 0
% 2 4305 2520 4305 2520
% 
\special{pn 8}%
\special{sh 1}%
\special{ar 4238 2481 10 10 0  6.28318530717959E+0000}%
\special{sh 1}%
\special{ar 4238 2481 10 10 0  6.28318530717959E+0000}%
% STR 2 0 3 0
% 3 705 3209 705 3285 5 0
% $+$
\put(6.9390,-32.3327){\makebox(0,0){$+$}}%
% DOT 2 0 3 0
% 2 1548 2826 1548 2826
% 
\special{pn 8}%
\special{sh 1}%
\special{ar 1524 2782 10 10 0  6.28318530717959E+0000}%
\special{sh 1}%
\special{ar 1524 2782 10 10 0  6.28318530717959E+0000}%
% DOT 2 0 3 0
% 2 1548 3132 1548 3132
% 
\special{pn 8}%
\special{sh 1}%
\special{ar 1524 3083 10 10 0  6.28318530717959E+0000}%
\special{sh 1}%
\special{ar 1524 3083 10 10 0  6.28318530717959E+0000}%
% DOT 2 0 3 0
% 2 1548 3439 1548 3439
% 
\special{pn 8}%
\special{sh 1}%
\special{ar 1524 3385 10 10 0  6.28318530717959E+0000}%
\special{sh 1}%
\special{ar 1524 3385 10 10 0  6.28318530717959E+0000}%
% DOT 2 0 3 0
% 2 1548 3745 1548 3745
% 
\special{pn 8}%
\special{sh 1}%
\special{ar 1524 3687 10 10 0  6.28318530717959E+0000}%
\special{sh 1}%
\special{ar 1524 3687 10 10 0  6.28318530717959E+0000}%
% STR 2 0 3 0
% 3 3463 3209 3463 3285 5 0
% $+$
\put(34.0846,-32.3327){\makebox(0,0){$+$}}%
% DOT 2 0 3 0
% 2 4305 2826 4305 2826
% 
\special{pn 8}%
\special{sh 1}%
\special{ar 4238 2782 10 10 0  6.28318530717959E+0000}%
\special{sh 1}%
\special{ar 4238 2782 10 10 0  6.28318530717959E+0000}%
% DOT 2 0 3 0
% 2 4305 3132 4305 3132
% 
\special{pn 8}%
\special{sh 1}%
\special{ar 4238 3083 10 10 0  6.28318530717959E+0000}%
\special{sh 1}%
\special{ar 4238 3083 10 10 0  6.28318530717959E+0000}%
% DOT 2 0 3 0
% 2 4305 3439 4305 3439
% 
\special{pn 8}%
\special{sh 1}%
\special{ar 4238 3385 10 10 0  6.28318530717959E+0000}%
\special{sh 1}%
\special{ar 4238 3385 10 10 0  6.28318530717959E+0000}%
% DOT 2 0 3 0
% 2 4305 3745 4305 3745
% 
\special{pn 8}%
\special{sh 1}%
\special{ar 4238 3687 10 10 0  6.28318530717959E+0000}%
\special{sh 1}%
\special{ar 4238 3687 10 10 0  6.28318530717959E+0000}%
% LINE 2 2 3 0
% 2 5228 2064 5075 2064
% 
\special{pn 8}%
\special{pa 5146 2032}%
\special{pa 4996 2032}%
\special{dt 0.045}%
% DOT 0 0 3 0
% 2 5228 2064 5228 2064
% 
\special{pn 20}%
\special{sh 1}%
\special{ar 5146 2032 10 10 0  6.28318530717959E+0000}%
\special{sh 1}%
\special{ar 5146 2032 10 10 0  6.28318530717959E+0000}%
% STR 2 0 3 0
% 3 4841 1983 4841 2060 5 0
% $+$
\put(47.6476,-20.2756){\makebox(0,0){$+$}}%
% DOT 2 0 3 0
% 2 5684 1600 5684 1600
% 
\special{pn 8}%
\special{sh 1}%
\special{ar 5595 1575 10 10 0  6.28318530717959E+0000}%
\special{sh 1}%
\special{ar 5595 1575 10 10 0  6.28318530717959E+0000}%
% DOT 2 0 3 0
% 2 5684 1907 5684 1907
% 
\special{pn 8}%
\special{sh 1}%
\special{ar 5595 1877 10 10 0  6.28318530717959E+0000}%
\special{sh 1}%
\special{ar 5595 1877 10 10 0  6.28318530717959E+0000}%
% DOT 2 0 3 0
% 2 5684 2213 5684 2213
% 
\special{pn 8}%
\special{sh 1}%
\special{ar 5595 2179 10 10 0  6.28318530717959E+0000}%
\special{sh 1}%
\special{ar 5595 2179 10 10 0  6.28318530717959E+0000}%
% DOT 2 0 3 0
% 2 5684 2520 5684 2520
% 
\special{pn 8}%
\special{sh 1}%
\special{ar 5595 2481 10 10 0  6.28318530717959E+0000}%
\special{sh 1}%
\special{ar 5595 2481 10 10 0  6.28318530717959E+0000}%
% LINE 2 0 3 0
% 2 5688 1604 5228 2064
% 
\special{pn 8}%
\special{pa 5599 1579}%
\special{pa 5146 2032}%
\special{fp}%
% LINE 2 0 3 0
% 2 5381 2217 5688 2523
% 
\special{pn 8}%
\special{pa 5297 2183}%
\special{pa 5599 2484}%
\special{fp}%
% LINE 2 0 3 0
% 2 5535 2064 5688 1911
% 
\special{pn 8}%
\special{pa 5448 2032}%
\special{pa 5599 1881}%
\special{fp}%
% LINE 2 0 3 0
% 2 5535 2064 5688 2217
% 
\special{pn 8}%
\special{pa 5448 2032}%
\special{pa 5599 2183}%
\special{fp}%
% LINE 2 1 3 0
% 2 5228 2064 5381 2217
% 
\special{pn 8}%
\special{pa 5146 2032}%
\special{pa 5297 2183}%
\special{da 0.070}%
% LINE 2 1 3 0
% 2 5535 2064 5381 2217
% 
\special{pn 8}%
\special{pa 5448 2032}%
\special{pa 5297 2183}%
\special{da 0.070}%
% DOT 0 0 3 0
% 2 5381 2217 5381 2217
% 
\special{pn 20}%
\special{sh 1}%
\special{ar 5297 2183 10 10 0  6.28318530717959E+0000}%
\special{sh 1}%
\special{ar 5297 2183 10 10 0  6.28318530717959E+0000}%
% DOT 0 0 3 0
% 2 5535 2064 5535 2064
% 
\special{pn 20}%
\special{sh 1}%
\special{ar 5448 2032 10 10 0  6.28318530717959E+0000}%
\special{sh 1}%
\special{ar 5448 2032 10 10 0  6.28318530717959E+0000}%
% LINE 2 2 3 0
% 2 2471 3289 2318 3289
% 
\special{pn 8}%
\special{pa 2433 3238}%
\special{pa 2282 3238}%
\special{dt 0.045}%
% DOT 0 0 3 0
% 2 2471 3289 2471 3289
% 
\special{pn 20}%
\special{sh 1}%
\special{ar 2433 3238 10 10 0  6.28318530717959E+0000}%
\special{sh 1}%
\special{ar 2433 3238 10 10 0  6.28318530717959E+0000}%
% STR 2 0 3 0
% 3 2084 3209 2084 3285 5 0
% $+$
\put(20.5118,-32.3327){\makebox(0,0){$+$}}%
% DOT 2 0 3 0
% 2 2927 2826 2927 2826
% 
\special{pn 8}%
\special{sh 1}%
\special{ar 2881 2782 10 10 0  6.28318530717959E+0000}%
\special{sh 1}%
\special{ar 2881 2782 10 10 0  6.28318530717959E+0000}%
% DOT 2 0 3 0
% 2 2927 3132 2927 3132
% 
\special{pn 8}%
\special{sh 1}%
\special{ar 2881 3083 10 10 0  6.28318530717959E+0000}%
\special{sh 1}%
\special{ar 2881 3083 10 10 0  6.28318530717959E+0000}%
% DOT 2 0 3 0
% 2 2927 3439 2927 3439
% 
\special{pn 8}%
\special{sh 1}%
\special{ar 2881 3385 10 10 0  6.28318530717959E+0000}%
\special{sh 1}%
\special{ar 2881 3385 10 10 0  6.28318530717959E+0000}%
% DOT 2 0 3 0
% 2 2927 3745 2927 3745
% 
\special{pn 8}%
\special{sh 1}%
\special{ar 2881 3687 10 10 0  6.28318530717959E+0000}%
\special{sh 1}%
\special{ar 2881 3687 10 10 0  6.28318530717959E+0000}%
% LINE 2 0 3 0
% 2 2930 3749 2471 3289
% 
\special{pn 8}%
\special{pa 2884 3690}%
\special{pa 2433 3238}%
\special{fp}%
% LINE 2 0 3 0
% 2 2624 3136 2930 2830
% 
\special{pn 8}%
\special{pa 2583 3087}%
\special{pa 2884 2786}%
\special{fp}%
% LINE 2 0 3 0
% 2 2777 3289 2930 3442
% 
\special{pn 8}%
\special{pa 2734 3238}%
\special{pa 2884 3388}%
\special{fp}%
% LINE 2 0 3 0
% 2 2777 3289 2930 3136
% 
\special{pn 8}%
\special{pa 2734 3238}%
\special{pa 2884 3087}%
\special{fp}%
% LINE 2 1 3 0
% 2 2471 3289 2624 3136
% 
\special{pn 8}%
\special{pa 2433 3238}%
\special{pa 2583 3087}%
\special{da 0.070}%
% LINE 2 1 3 0
% 2 2777 3289 2624 3136
% 
\special{pn 8}%
\special{pa 2734 3238}%
\special{pa 2583 3087}%
\special{da 0.070}%
% DOT 0 0 3 0
% 2 2624 3136 2624 3136
% 
\special{pn 20}%
\special{sh 1}%
\special{ar 2583 3087 10 10 0  6.28318530717959E+0000}%
\special{sh 1}%
\special{ar 2583 3087 10 10 0  6.28318530717959E+0000}%
% DOT 0 0 3 0
% 2 2777 3289 2777 3289
% 
\special{pn 20}%
\special{sh 1}%
\special{ar 2734 3238 10 10 0  6.28318530717959E+0000}%
\special{sh 1}%
\special{ar 2734 3238 10 10 0  6.28318530717959E+0000}%
\end{picture}%
 \caption{For example $\ti{m}^p_4$ and $\ti{f}^p_4$ are given. 
The large dots represent the vertices $\{m_k\}$. 
The dashed lines denote the internal edges and we attach 
$-Q^+$ on them. The dotted line on each graph is the root edge, 
on which we attach 
$\pi$ for $\ti{m}^p_k$ and $-Q^+$ for $\ti{f}^p_k$. 
For $\ti{m}^p_4$ and $\ti{f}^p_4$, all such $4$-trees 
are summed up with weight $+1$. }
 \label{fig:m4f4}
\end{figure}

An explicit example is given 
in Figure \ref{fig:m4f4}. 
In the order of the graphs in Figure \ref{fig:m4f4}, we have 
\begin{equation*}
 \begin{split}
\ti{m}_4^p( {o}^p_1, {o}^p_2, {o}^p_3, {o}^p_4)&
=\pi\circ m_4( {o}^p_1, {o}^p_2, {o}^p_3, {o}^p_4)+
\pi\circ m_3(-Q^+ m_2( {o}^p_1, {o}^p_2), {o}^p_3, {o}^p_4)\\&+
\pi\circ m_3( {o}^p_1,-Q^+ m_2( {o}^p_2, {o}^p_3), {o}^p_4)+
\pi\circ m_3( {o}^p_1, {o}^p_2,-Q^+ m_2( {o}^p_3, {o}^p_4))\\&+
\pi\circ m_2(-Q^+ m_3( {o}^p_1, {o}^p_2, {o}^p_3), {o}^p_4)+
\pi\circ m_2( {o}^p_1,-Q^+ m_3( {o}^p_2, {o}^p_3, {o}^p_4))\\&+
\pi\circ m_2(-Q^+ m_2(-Q^+ m_2( {o}^p_1, {o}^p_2), {o}^p_3), {o}^p_4)\\&+
\pi\circ m_2( {o}^p_1,-Q^+ m_2(-Q^+ m_2( {o}^p_2, {o}^p_3), {o}^p_4))\\&+
\pi\circ m_2(-Q^+ m_2( {o}^p_1, {o}^p_2),-Q^+ m_2( {o}^p_3, {o}^p_4))\\&+
\pi\circ m_2(-Q^+ m_2( {o}^p_1,-Q^+ m_2( {o}^p_2, {o}^p_3)), {o}^p_4)\\&+
\pi\circ m_2( {o}^p_1,-Q^+ m_2( {o}^p_2,-Q^+ m_2( {o}^p_3, {o}^p_4)))\ ,
 \end{split}
\end{equation*} 
for $ {o}^p_1, {o}^p_2, {o}^p_3, {o}^p_4\in\cH^p$ 
and  $\ti{f}^p_4$ is obtained similarly but replaced each 
$\pi$ on the outgoing line to $-Q^+$.

We emphasize that the definition above is derived essentially by 
using eq.(\ref{eom3}) recursively. 
\begin{equation*}
 \Phi|_{gf}=\ti\Phi^p-Q^+\sum_{k\ge 2}m_k(\Phi|_{gf})\ .
\end{equation*} 
Let us extend the equation above to 
\begin{equation}
 \Phi=\ti\Phi-Q^+\sum_{k\ge 2}m_k(\Phi)\ , 
 \label{rec-extend}
\end{equation}
where $\ti\Phi\in\cH$. 
One can construct an $A_\infty$-isomorphism $\ti\Phi\mapsto\Phi$ 
recursively in a similar way as in Definition \ref{defn:minimal}. 
Then one can get another $A_\infty$-algebra 
as the pullback of $(\cH,\m)$ by this isomorphism. 
Let us denote this $A_\infty$-algebra by $(\cH,\ti\m)$ and the 
$A_\infty$-isomorphism by $\ti\cF: (\cH,\ti\m)\to (\cH,\m)$. 
The explicit forms of $\ti\m:=\{\ti{m}_k\}_{k\ge 1}$ and 
$\ti\cF:=\{\ti{f}_k\}_{k\ge 1}$ are given 
in the terminology in Definition \ref{defn:minimal2} as 
\begin{equation*}
 \ti{f}_k=\sum_{\Gamma_k\in G_k}\ti{f}_{\Gamma_k}\ \ (k\ge 1)\ ,\qquad 
 \ti{m}_1=m_1=Q\ ,\quad 
 \ti{m}_k=P\circ\sum_{\Gamma_k\in G_k}\ti{m}_{\Gamma_k}\ \ (k\ge 2)\ .
\end{equation*}
In particular, one has $\ti{f}_1=\Id$. 
Note that these are the operations not on $\cH^p$ but on $\cH$. 
That is, $\ti{m}_k: \cH^{\otimes k}\to\cH^p\subset\cH$ and 
$\ti{f}_k:\cH^{\otimes k}\to\cH$ do not vanish generally 
even if one of the $\cH$ in $\cH^{\otimes k}$ includes an element 
in $\cH^t\oplus\cH^u$. 
For $\iota :\cH^p\to\cH$ the inclusion map, 
we have the relations 
$\ti{m}^p_k=\pi\circ\ti{m}_k\circ(\iota)^{\otimes k}$ 
and $\ti{f}^p_k=\ti{f}_k\circ(\iota)^{\otimes k}$ for $k\ge 2$. 
\begin{rem}
Since $\ti\cF:(\cH,\ti\m)\to (\cH,\m)$ is not only 
an $A_\infty$-quasi-isomorphism but also an $A_\infty$-isomorphism, 
it has its inverse $A_\infty$-isomorphism 
$(\ti\cF)^{-1} :(\cH,\m)\to (\cH,\ti\m)$, which is given by 
\begin{equation*}
 \begin{array}{ccccl}
 (\ti\cF)^{-1}_* &:& \cH &\lgraw &\cH \\
              & & \Phi&\mapsto&\ti\Phi:=\Phi+Q^+\sum_{k\ge 2} m_k(\Phi)\ .
 \end{array}
\end{equation*}
 \label{rem:inverse}
\end{rem} 
In the rest of this subsection we shall give a proof of the statement below. 
\begin{lem}
$(\cH^p,\ti\m^p)$ and $(\cH,\ti\m)$ are in fact $A_\infty$-algebras and 
$\ti\cF^{(p)}$ 
$($by $\ti\cF^{(p)}$ we denote $\ti\cF$ or $\ti\cF^p$ $)$ 
is an $A_\infty$-morphism. 
 \label{lem:minimal}
\end{lem}
\begin{pf}
The fact that $(\cH^p,\ti\m^p)$ is an $A_\infty$-algebra and $\ti\cF^p$ is 
an $A_\infty$-quasi-isomorphism between $(\cH^p,\ti\m^p)$ to $(\cH,\m)$ 
immediately follows from the fact that 
$(\cH,\ti\m)$ is an $A_\infty$-algebra and $\ti\cF$ is 
an $A_\infty$-quasi-isomorphism between $(\cH,\ti\m)$ to $(\cH,\m)$. 
The latter fact is explicitly described by the equations on $\cH$ as 
\begin{equation*}
 \m\ti\cF=\ti\cF\ti\m\ ,\qquad (\ti\m)^2=0\ .
\end{equation*}
The former is obtained by restricting 
these equations to $\cH^p$, that is, 
$\m\, \ti{\cF}\, \iota=\ti\cF\, \ti{\m}\, \iota$ 
and $(\ti\m)^2\, \iota=0$ on $\cH^p$. 
Therefore we will prove the latter fact. 
In order to see this, it is enough to confirm the following 
two facts : $\m\ti\cF=\ti\cF\ti\m$ and 
$(\ti\m)^2=0$ on $\cH$. 
These are shown at the same time by 
checking
\begin{equation}
\ti\m=\ti\cF^{-1}\m\ti\cF
 \label{tmfmf}
\end{equation} 
since $\ti\cF^{-1}\ti\cF=\1$ and $\ti\cF\ti\cF^{-1}=\1$. 
We shall show it below in the dual picture. 
Let $\delta$ be $A_\infty$-odd vector fields corresponding 
to $\m$. 
Recall that 
$(\ti\cF)^{-1}_*(\ti\Phi)=\Phi+Q^+\sum_{k\ge 2}m_k(\Phi,\dots,\Phi)$, 
and we write it as 
\begin{equation}
 \ti\phi^i=\phi^i+{\bar c}^i_j\sum_{k\ge 2}
 c^j_{i_1\cdots i_k}\phi^{i_k}\cdots\phi^{i_1}\ ,
 \label{phitiphi}
\end{equation}
where $Q^+\eb_j=\eb_i{\bar c}^i_j$. 
The coordinate transformation corresponding to $\ti\cF_*:\ti\Phi\to\Phi$ 
is obtained by using 
\begin{equation}
 \phi^i=\ti\phi^i-{\bar c}^i_j\sum_{k\ge 2}
 c^j_{i_1\cdots i_k}\phi^{i_k}\cdots\phi^{i_1}
 \label{tiphiphi}
\end{equation} 
recursively. 

Let $\ti\delta$ be the 
$A_\infty$-odd vector field dual to $\ti\cF^{-1}\m\ti\cF$ in the 
right hand side of eq.(\ref{tmfmf}). 
{}From now we rewrite $\delta$ as this $\ti\delta$ 
using eq.(\ref{phitiphi}) and eq.(\ref{tiphiphi}), 
and show that the $\ti\delta$ is in fact 
the $A_\infty$-odd vector field dual to $\ti\m$ defined by 
eq.(\ref{rec-extend}). 
For 
$$
 \delta=\flpartial{\phi^i}\l(c^i_j\phi^j+\sum_{k\ge 2}
 c^i_{i_1\cdots i_k}\phi^{i_k}\cdots\phi^{i_1}\r)\ ,
$$
$\ti\delta$ is obtained by 
\begin{equation}
 \ti\delta=\flpartial{\ti\phi^i}\flpart{\ti\phi^i}{\phi^l}
 \l( c^l_j\phi^j+\sum_{k\ge 2}
 c^l_{i_1\cdots i_k}\phi^{i_k}\cdots\phi^{i_1}\r)\ .
 \label{tidelta1}
\end{equation}
Here $\flpart{\ti\phi^i}{\phi^l}=\delta^i_k+{\bar c}^i_j\sum_{k\ge 2}
c^j_{i_1\cdots i_k}\l(\phi^{i_k}\cdots\phi^{i_1}\r)\flpartial{\phi^l}$ 
by using eq.(\ref{phitiphi}) 
and all $\phi$'s are supposed to be substituted into by eq.(\ref{tiphiphi}). 
Thus eq.(\ref{tidelta1}) is further rewritten as 
\begin{equation}
 \begin{split}
 \ti\delta=&\flpartial{\ti\phi^i}
 c^i_j\l(\ti\phi^j-{\bar c}^j_l\sum_{n\ge 2}
 c^l_{j_1\cdots j_n}\phi^{j_n}\cdots\phi^{j_1}\r)\\
  &+\flpartial{\ti\phi^i}\sum_{k\ge 2}
 c^i_{i_1\cdots i_k}\phi^{i_k}\cdots\phi^{i_1}\\
  &+\flpartial{\ti\phi^i}{\bar c}^i_j\sum_{k\ge 2}
c^j_{i_1\cdots i_k}\l(\phi^{i_k}\cdots\phi^{i_1}\r)\flpartial{\phi^l}
\sum_{n\ge 1}c^l_{j_1\cdots j_n}\phi^{j_n}\cdots\phi^{j_1} \ .
 \end{split}
 \label{tidelta2}
\end{equation}
Note that, by the $A_\infty$-condition for $\delta$, the equation of the 
third line is replaced by 
$$
\flpartial{\ti\phi^i}{\bar c}^i_j\l(-c^j_l\r)
\sum_{n\ge 1}c^l_{j_1\cdots j_n}\phi^{j_n}\cdots\phi^{j_1}\ .
$$ 
Thus one can get 
\begin{equation*}
  \ti\delta=\flpartial{\ti\phi^i}
 c^i_j\ti\phi^j
 +\flpartial{\ti\phi^i}\sum_{k\ge 2}
 (\delta^i_l-c^i_j{\bar c}^j_l-{\bar c}^i_jc^j_l)
 \sum_{k\ge 2}
 c^l_{i_1\cdots i_k}\phi^{i_k}\cdots\phi^{i_1}\ .
\end{equation*}
The second term of the first line, the term of the second line 
and the one of the third line in eq.(\ref{tidelta2}) are gathered as 
the second term above. 
$(\delta^i_l-c^i_j{\bar c}^j_l-{\bar c}^i_jc^j_l)$ restricts the indices 
$l$ to those in $\cH^p$, \ie, the dual description of $P$. 
Therefore by substituting eq.(\ref{tiphiphi}) recursively 
in the equation above 
one can see that $\ti\delta$ above is just the dual expression of $\ti\m$. 
Thus the proof is completed. \qed
\end{pf} 

Another proof in terms of the `superfield' $\Phi$ is presented in \cite{Ka}.

 \subsection{Minimal cyclic $A_\infty$-algebras and Feynman graphs}
\label{ssec:cycminimal}

There exists an explicit construction of the minimal model also for 
cyclic $A_\infty$-algebras. 
In this subsection we shall see 
that the arguments in the previous subsection 
are applicable directly to cyclic version. 

Our starting point is the (extended) Maurer-Cartan equation (\ref{eomsft}) 
in the previous subsection. 
One can see that the equation is identified with $\delta=0$. 
As stated in subsection \ref{ssec:cycAinftyre}, 
a cyclic $A_\infty$-algebra $(\cH,\omega,\m)$ has 
a degree zero cyclic function, 
the action $S\in C(\phi)_c$ (eq.(\ref{action-pre})). 
The $A_\infty$-odd vector field is then given by $\delta=(\ ,S)$. 
Since the odd constant Poisson bracket $(\ ,\ )$ is nondegenerate, 
it is just equivalent to the equation of motion of the action, 
\begin{equation}
 \frpartial{\phi^j}S=0\ .\label{eomsft-re}
\end{equation}
The action is expressed in terms of the superfield $\Phi$ as 
\begin{equation}
 S(\Phi)=\sum_{k\ge 2}\ov{k}\V_k(\Phi,\dots,\Phi)
  =\half\omega(\Phi,Q\Phi)
 +\sum_{k\ge 3}\ov{k}\omega(\Phi,m_{k-1}(\Phi,\dots,\Phi))\ , 
 \label{action-superfield}
\end{equation}
where $\V_k(\eb_{i_1},\dots,\eb_{i_k})=\V_{i_1\cdots i_k}$. 
\begin{lem}[A minimal cyclic $A_\infty$-algebra]
For a given cyclic $A_\infty$-algebra $(\cH,\omega,\m)$, 
suppose we have a Hodge-Kodaira decomposition of $\cH$ with 
a homotopy operator $Q^+$ of $Q$ compatible with $\omega$. 
Then, construct the explicit minimal $A_\infty$-algebra $(\cH^p,\ti\m^p)$ 
of the $A_\infty$-algebra $(\cH,\m)$ 
and the $A_\infty$-quasi-isomorphism 
$\ti\cF^p : (\cH^p,\ti\m^p)\to (\cH,\m)$ in Definition \ref{defn:minimal}. 
Next, define a symplectic structure $\ti\omega^p$ on $\cH^p$ by 
restricting the symplectic structure $\omega$ on $\cH$ to $\cH^p$, 
that is, $\ti\omega^p:=\omega(\iota\otimes\iota)$ for $\iota:\cH^p\to\cH$. 
\begin{itemize}
 \item[(a)]
The $A_\infty$-structure $\ti\m^p$ is cyclic with respect to
$\ti\omega^p$, that is, $(\cH^p,\ti\omega^p,\ti\m^p)$ defines 
a minimal cyclic $A_\infty$-algebra. 

 \item[(b)]
$\ti\cF^p : (\cH^p,\ti\omega^p,\ti\m^p)\to (\cH,\omega,\m)$ 
is a cyclic $A_\infty$-quasi-isomorphism. 
\end{itemize}
 \label{lem:cycminimal}
\end{lem}
\begin{pf}
Recall that, for each $k\ge 2$, $\ti{m}^p_k$ is defined so that 
it is associated to the summation over all rooted planar $k$-trees. 
The sum of all $k$-trees is cyclic, that is, 
invariant with respect to the cyclic permutations of the root and 
the leaves (see also subsection \ref{ssec:PI}). 
Statement (a) follows from this fact together with the conditions 
(\ref{homotopy-compatible}). 
For statement (b), by the definition 
of cyclic $A_\infty$-morphisms (Definition \ref{defn:cycAinftymorp}),
it is sufficient to show that 
$\ti\cF^p$ preserves the symplectic structures $\ti\omega^p$ and $\omega$. 
$(\ti\cF^p)^*\omega$ is written as 
\footnote{The equation below is actually equivalent to 
the one for odd Poisson structure. 
The equivalence follows from $\omega_{ij}\omega^{jk}=\delta_i^k$ and 
$\ti\omega^p_{ij}\ti\omega^{p,jk}=\delta_i^k$. }
\begin{equation*}
 \l((\ti\cF^p)^*\omega\r)_{ij}
 =\frpart{\phi^k}{\ti{p}^i}\omega_{kl}\flpart{\phi^l}{\ti{p}^j}
 =(-1)^{\eb^p_i}
 \omega\l(\frpartial{\ti{p}^i}\Phi,\Phi\flpartial{\ti{p}^j}\r)\ ,
\end{equation*}
where $\ti{p}^i$ is the dual coordinate of the basis vector 
$\eb^p_i\in\cH^p$. 
Since $\Phi=\ti\Phi^p-Q^+\sum_{k\ge 2}m_k(\Phi)$ and 
the image of $Q^+$ vanishes in the symplectic inner product in the right 
hand side of the above equation, 
the right hand side becomes 
$(-1)^{\eb^p_i}\omega(\frpartial{\ti{p}^i}\ti\Phi^p,
\ti\Phi^p\flpartial{\ti{p}^j})=\ti\omega^p_{ij}$. 
Thus it has been shown that 
the map $\ti\cF^p$ preserves the symplectic structures. 
 \qed
\end{pf}

These facts together with Proposition \ref{prop:Ap} further imply that 
the action of the corresponding minimal cyclic $A_\infty$-algebra 
\begin{equation}
 \ti{S}(\ti\Phi^p)=\sum_{k\ge 2}\ov{k+1}
 \omega^p(\ti\Phi^p,\ti{m}^p_k(\ti\Phi^p))\ 
 \label{Sp}
\end{equation}
can also be obtained by the pullback of $S(\Phi)$ by $\ti\cF^p$
\begin{equation}
 \ti{S}(\ti\Phi^p)=(\ti\cF^p)^*S(\Phi)
 =\l(\ S((\ti\cF^p)_*(\ti\Phi^p))\ \r)_c\ .
 \label{fppreserve}
\end{equation}
It is also interesting to show the equality 
$\ti{S}(\ti\Phi^p)=(\ti\cF^p)^*S(\Phi)$ directly; 
$(\ti\cF^p)^*S(\Phi)$ coincides with $\ti{S}(\ti\Phi^p)$ 
due to some nontrivial combinatorial cancellations 
(see subsection 5.3 of \cite{Ka}).

As stated previously, the homotopy types of cyclic $A_\infty$-algebras 
are classified by their minimal cyclic $A_\infty$-algebras. 
For a given cyclic $A_\infty$-algebra, 
its minimal one is unique up to cyclic $A_\infty$-isomorphisms. 
Namely, we have at least a decomposed cyclic $A_\infty$-algebra 
whose minimal part is the one given explicitly in this subsection. 
Note however that a cyclic $A_\infty$-isomorphism from the decomposed
one to the original one is not given explicitly. 
It might be interesting to explore the relation 
between the explicit construction of a minimal model 
and the way of the proof of the decomposition theorem in 
subsection \ref{ssec:cycMandC}.

%section6

 \section{The minimal model theorem in the BV-formalism}
\label{sec:BV}

In this section we shall apply the homotopy algebraic 
structures discussed in the previous section 
to field theory equipped with classical BV-structures. 
In subsection \ref{ssec:BV}, 
it is shown that 
any cyclic field theory equipped with a classical BV-structure has a 
cyclic $A_\infty$-structure (Theorem \ref{thm:ft}). 
Subsection \ref{ssec:gf} is devoted to a brief review of 
perturbative expansion in the BV-formalism and to 
translating it into our language. 
The perturbative expansion is necessary 
for computing correlation functions in the subsequent subsections. 
Also, we shall show that the propagator in the BV-formalism, 
which is defined in subsection \ref{ssec:gf}, 
is a homotopy operator $Q^+$ in the previous 
section (Proposition \ref{prop:BVp}). 
In subsection \ref{ssec:PI} it is then shown that 
the tree on-shell correlation functions of a 
cyclic field theory equipped with a classical BV-structure 
define just the minimal cyclic $A_\infty$-algebra 
defined in subsection \ref{ssec:cycminimal} 
(Corollary \ref{cor:main} (cf.\cite{Ka})). 
Moreover in subsection \ref{ssec:main} the arguments in section 
\ref{sec:MandC} are applied to the classification of 
classical open string field theories, 
and it is shown that 
all classical string field theories on a fixed conformal 
background are cyclic $A_\infty$-isomorphic to each other
(Theorem \ref{thm:main}). 
The theorem means all classical string field theories on a fixed conformal 
background are related by field redefinitions and so 
physically equivalent to each other.

 \subsection{Cyclic $A_\infty$-structures in the BV-formalism}
\label{ssec:BV}

The BV-formalism is formulated on formal supermanifolds equipped with 
odd symplectic forms, where the coordinates of supermanifolds are just the
fields. 
Since we discuss field theories related to open string theory, 
we let the fields noncommutative as explained in subsection 
\ref{ssec:Introphys}. 
We shall first explain that our noncommutative symplectic
supergeometry just fits the BV-formalism. 

For any odd symplectic form one can take a Darboux coordinate 
due to Theorem \ref{thm:Darboux}. 
We denote the odd symplectic form on $\Z$-graded vector space $\cH$ 
in a Darboux coordinate by 
\begin{equation}
  \{\omega_{ij}\}=\bp \0 & -\1\\ \1 & \0\ep\ .
 \label{symD}
\end{equation}
In this coordinate, let us decompose the basis $\{\eb_i\}$ into 
$\{\eb_a\}$ and $\{\eb_a^*\}$ such that 
$-\omega(\eb_a,\eb_b^*)=\omega(\eb_b^*, \eb_a)=\delta_{ab}$ 
and $\omega(\eb_a,\eb_b)=\omega(\eb_a^*,\eb_b^*)=0$. 
Then we have $\cH=\cH_+\oplus\cH_-$, where $\cH_+$ and $\cH_-$ are 
the $\Z$-graded vector spaces spanned by $\{\eb_a\}$ and $\{\eb_a^*\}$, 
respectively. As above, 
we use indices $a, b, \dots$ for the basis of $\cH_+$ or $\cH_-$. 
For bases $\eb_a$ and $\eb_a^*$, we denote the associated dual fields 
by $\phi^a$ and $\phi^{a,*}$, respectively. Also, it is more convenient to 
prepare the notation $\phi^*_a:=\omega_{ab}\phi^{b,*}$ 
where $\omega_{ab}$ is the one in eq.(\ref{symD}). 
The odd Poisson bracket associated to eq.(\ref{symD}) is then written as 
\begin{equation}
 (\ ,\ )=\flpartial{\phi^i}\omega^{ij}\frpartial{\phi^j}
        =\flpartial{\phi^a}\frpartial{\phi^*_a}
        -\flpartial{\phi^*_a}\frpartial{\phi^a}\ .
 \label{BVbracket}
\end{equation}
This is just the situation in the BV-formalism \cite{BV1,BV2,HT,GPS}. 
In the context of the BV-formalism, 
$\{\phi^a\}$ consists of usual fields of degree zero, 
ghost fields of degree one, and ghosts of ghosts with degree two, ...
and also so-called antighosts whose degree is defined negative 
(so that the gauge fixing (Definition \ref{defn:gf}) 
can be performed in the BV-formalism). 
For each $\phi^a$, its {\it antifield} 
$\phi^*_a$ is then introduced, where 
its degree is defined as 
\begin{equation*}
 \deg(\phi^*_a)=-1-\deg(\phi^a)\ .
\end{equation*}
This is equivalent to $\deg(\eb_a)+\deg(\eb_a^*)=1$. 
For a constant symplectic form $\omega$, 
this fact determines the degree of $\omega$ to be minus one. 
Thus, our definition for the degree of $\omega$ 
(Definition \ref{defn:csym}) is natural also 
from the viewpoint of the BV-formalism.

The BV-formalism is applied mainly to two cases. 
Originally \cite{BV1,BV2}, it is a general method to 
quantize gauge-invariant actions consistently. 
In such a usage one begins with the gauge invariant action 
which does not include antifields, one 
adds the terms including antifields to the original action 
so that the action satisfies the master equation and 
is proper (Definition \ref{defn:proper}). 
The BV-quantization of Poisson-$\sigma$ model in \cite{CF1} 
is a good example. 
On the other hand, it is used to determine higher terms of actions. 
Namely, starting from an action which consists only a kinetic term, 
when one includes higher interaction terms as deformation of the
action preserving its symmetry, 
the BV-master equation becomes constraints for the determination of
the higher interaction terms. 
String field theory as reviewed in subsection \ref{ssec:osft} 
is just the latter case. 
A similar application to topological field theories is 
given by \cite{BH} and developed 
for example in \cite{AKSZ,I,Park,CF2,Ba}. 

In any case the action in the BV-formalism 
is, by power series of fields, written as 
\begin{equation}
 S=\half\V_{i_1i_2}\phi^{i_2}\phi^{i_1}
 +\sum_{k\ge 3}\ov{k}\V_{i_1\cdots i_k}\phi^{i_k}\cdots\phi^{i_1}\ ,
\qquad \V_{i_1\cdots i_k}\in\C\ ,
 \label{action}
\end{equation}
where $\{\phi^i\}$ consists of both fields and antifields. 
We consider the case when the action is cyclic $S\in C(\phi)_c$. 
We call such a field theory {\it a cyclic field theory}. 
Moreover suppose that $S$ satisfies the classical BV-master equation 
\begin{equation}
 (S, S)=0\ ,
 \label{meq}
\end{equation}
where $(\ ,\ )$ is the odd Poisson bracket in eq.(\ref{BVbracket}). 
In this situation we say that the action is equipped with a 
{\it classical BV-structure}. 

The BV-BRST transformation is then defined 
by the Hamiltonian vector field of $S$
\begin{equation}
 \delta =(\ ,S)\ .
 \label{BVt}
\end{equation}
With this $\delta$, the classical master equation (\ref{meq}) is written 
as $\delta S=0$. Moreover, by using the Jacobi identity of the BV-bracket 
and the classical BV-master equation (\ref{meq}), 
$(\delta)^2=0$ holds. 
Namely, the following 
three statements are equivalent; the action $S$ satisfies the BV-master 
equation (\ref{meq}), the action $S$ is invariant under the BV-BRST 
transformation (\ref{BVt}), and the BV-BRST transformation $\delta$ is 
nilpotent. 

As stated in subsection \ref{ssec:cycAinftyre}, 
this $\delta$ is nothing but the $A_\infty$-odd vector field. 
It is written in the form 
\begin{equation}
 \delta=(\ ,S)=\sum_{k=1}^\infty\flpartial{\phi^j}c^j_{i_1\cdots i_k}
 \phi^{i_k}\cdots\phi^{i_1}\ , 
\end{equation}
where 
$c^j_{i_1\cdots i_k}=(-1)^{\eb_l}\omega^{jl}\V_{li_1\cdots i_k}$ 
(eq.(\ref{defc})). 
Note that the $\{c^j_{i_1\cdots i_k}\}$ not only defines an 
$A_\infty$-structure but also has a cyclic structure 
since the set $\{\V_{ii_1\cdots i_k}\}$ defines the cyclic field theory. 
The algebraic structures of these field theories are 
the cyclic $A_\infty$-structure in Definition \ref{defn:cycAinfty}. 

Generally the following fact holds. 
\begin{thm}
Any field theory has a cyclic $A_\infty$-structure 
if the action is cyclic, $S\in C(\phi)_c$, and 
satisfies a classical BV-master equation. 
 \label{thm:ft}
\end{thm}
In Theorem \ref{thm:ft} we assumed the fields $\{\phi^i\}$ are associative. 
An example of field theory with a cyclic action is nonabelian 
single trace gauge theory. 
Here nonabelian means each field is $N\times N$ matrix for 
some $N\in\N$ (see subsection \ref{ssec:Introphys}). 
Single trace then corresponds to $S\in C(\phi)_c$, not $TC(\phi)_c$. 
If we assume the fields (graded) commutative, the $A_\infty$-structure 
reduces to an $L_\infty$-structure. 
Namely, Theorem \ref{thm:ft} implies that any field theory which consists of 
usual graded commutative fields has a (cyclic) $L_\infty$-structure 
if it is equipped with a classical BV-structure 
(such as in \cite{BH,AKSZ,I,Park,CF2,Ba}). 
Of course this fits also the case 
when the matrix fields are decomposed into 
graded commutative component fields as mentioned 
in subsection \ref{ssec:Introphys}.

Usually field theory deals with fields $\{\phi^i\}$, 
the dual side, where $\Z$-graded vector space $\cH$ is implicit. 
However we can now express the $A_\infty$-structure explicitly 
by employing the arguments in subsection \ref{ssec:superfield}. 
For each field $\phi^i$ one can define its dual base $\eb_i\in\cH$ 
whose degree is minus the degree of $\phi^i$. 
Then one can associate to $\V_{i_1\cdots i_k}$ for $k\ge 2$ 
a cyclic multilinear map $\V_k :\cH\otimes\cdots\otimes\cH\to\C$ as 
\begin{equation*}
 \V_k(\eb_{i_1},\dots,\eb_{i_k})=\V_{i_1\cdots i_k}\ .
\end{equation*}
It can also be written as 
$\V(\eb_{i_1},\dots,\eb_{i_k})
=(-1)^{\eb_{i_1}}
\omega(\eb_{i_1},m_{k-1}(\eb_{i_2},\dots,\eb_{i_k}))$ 
(Remark \ref{rem:Vcyclic})
since $\omega$ is nondegenerate. 
Also, for the operation $\omega(\ ,\ ):\cH\otimes\cH\to\C$ 
we count the degree minus one by not $`,'$ but $\omega$, and 
extend naturally the operation to the one over $C(\phi)$. 
Each term of the action (\ref{action}) is then represented as 
\begin{equation*}
 \begin{split}
 \V_{i_1\cdots i_{k+1}}\phi^{i_{k+1}}\cdots\phi^{i_1} 
 &=(-1)^{\eb_{i_1}}\omega_{i_1 j}c^j_{i_2\cdots i_{k+1}}
 \phi^{i_{k+1}}\cdots\phi^{i_2}\cdot\phi^{i_1}\\
 &=\phi^{i_1}\omega(\eb_{i_1},\eb_{j}c^j_{i_2\cdots i_{k+1}})
 \phi^{i_{k+1}}\cdots\phi^{i_2}\\
 &=\omega(\eb_{i_1}\phi^{i_1},m_k(\eb_{i_2}\phi^{i_2},
 \dots,\eb_{i_{k+1}}\phi^{i_{k+1}}))\\
 &=\omega(\Phi,m_k(\Phi,\dots,\Phi))\ ,
 \end{split}
\end{equation*}
and the action is given by eq.(\ref{action-superfield}), 
\begin{equation}
 S=\half\omega(\Phi,Q\Phi)+\sum_{k\ge 2}\ov{k+1}\omega(\Phi,m_k(\Phi))\ .
 \label{actionf}
\end{equation}
{}From the form of the kinetic term one can see that 
the degrees of $Q$, $\omega$ and
$\Phi$ are assigned consistently. 
Of course classical open string field theory 
constructed as in subsection \ref{ssec:osft} also takes this form. 
In this case, 
$\Phi=\eb_i\phi^i$ is just the string field 
where $\{\eb_i\}$ is the basis of the string Hilbert space $\cH$ 
and $\phi^i$ are its coordinate whose degree is minus 
the degree of $\eb_i$. 
The inner product $\omega$ is known as BPZ-inner product.

 \subsection{Gauge fixing in the BV-formalism}
\label{ssec:gf}

To proceed the path-integral in the BV-formalism, 
it is necessary to fix a gauge. 
Roughly speaking, the gauge fixing corresponds to killing the 
degree of freedom of gauge transformations. 
When choose a gauge fixing, one can construct a propagator canonically. 
The aim of this subsection is to explain these facts and 
to show that the propagator is just a homotopy operator $Q^+$ 
as in the previous section (Proposition \ref{prop:BVp}). 
The statement of the path-integral (perturbative expansion) given 
in this subsection is formal. We shall write down the definition 
explicitly in the next subsection.

Formally, given an action $S\in C(\phi)_c$ as in eq.(\ref{actionf}), 
the starting point of the path-integral 
is the partition function of the field theory, 
\begin{equation}
 Z=\int\mes\Phi\ e^{-S}\ .
 \label{Z}
\end{equation}
Let us separate the action into the quadratic term and others. 
We denote the quadratic term by $S_2=\half\omega(\Phi,Q\Phi)$ and 
the rest terms by $S_{int}:=s_3+s_4+\cdots$. 
Then we have $e^{-S}=e^{-S_{int}}e^{-S_2}$ in eq.(\ref{Z}). 
In perturbation theory, the partition function (\ref{Z}) 
is computed by perturbative expansion, which is essentially 
the Gaussian integral of $e^{-S_2}$ where 
$e^{-S_{int}}$ is Taylor expanded. 
Here we should define the integration $\int\mes\Phi$ in some sense. 
In fact, the Hessian of the kinetic term $S_2$ 
is degenerate and one cannot integrate 
over the whole space $\Phi$ by perturbative expansion. 

For this purpose, 
the properties of the kinetic term should be examined. 
The degeneracy is directly related to the gauge transformation. 
For the BV-BRST transformation of the string field, 
$\delta\Phi=\sum_{k\ge 1}m_k(\Phi)=\m_*(e^\Phi)$, 
the gauge transformation is defined correspondingly as 
a degree zero transformation $\delta_\alpha$ given by 
\begin{equation}
 \begin{split}
 \delta_\alpha\Phi&=\m_*(e^\Phi\alpha e^\Phi)\\
 &:=Q\alpha+m_2(\alpha,\Phi)+m_2(\Phi,\alpha)
 +m_3(\alpha,\Phi,\Phi)+m_3(\Phi,\alpha,\Phi)+m_3(\Phi,\Phi,\alpha)+\cdots\ ,
 \end{split}
 \label{gauge}
\end{equation}
where $\alpha=\eb_i\alpha^i$ is a gauge parameter of degree minus
one. More precisely, $\alpha^i$ is treated as a graded parameter that 
belongs to $\cH^*[1]$ and can be identified with
$\alpha^i=\phi^i[1]$. 
The degree of $\alpha^i$ is then minus the degree of $\eb_i$ minus one. 
One may notice that this is the gauge transformation for 
$A_\infty$-algebras discussed in subsection \ref{ssec:homgauge}.
We shall discuss more on the properties of gauge transformations there. 
The gauge transformation is written as 
$\delta_{\alpha}\Phi=m_*(e^\Phi)\flpartial{\phi^i}\alpha^i$, and in the 
language of the component fields, it is 
\begin{equation}
 \delta_\alpha=\sum_{k=1}^\infty\flpartial{\phi^j}c^j_{i_1\cdots i_k}
 \l(\phi^{i_k}\cdots\phi^{i_1}\flpartial{\phi^i}\alpha^i\r)\ .
 \label{gaugedual}
\end{equation}
By standard arguments in the BV-formalism \cite{BV2,HT}, one sees 
the following facts. 
First, the action is invariant under this $\delta_\alpha$ because 
$0=(S,S)\flpartial{\phi^i}\alpha^i$ implies $\delta_\alpha S=0$. 
Moreover, 
\begin{equation}
 0=\left.\omega_{kj}\frpartial{\phi^j}(S,S)\flpartial{\phi^i}\alpha^i
 \right|_{\fpart{S}{\phi}=0}
  =2\left.\l(\omega_{kj}\frpartial{\phi^j}S\flpartial{\phi^l}\r)
 \l(\omega^{lm}\frpartial{\phi^m}S\flpartial{\phi^i}\r)\alpha^i
 \right|_{\fpart{S}{\phi}=0}
 \label{hessian}
\end{equation}
indicates that the generator of the gauge transformation is degenerate 
and the rank of the Hessian for the quadratic part 
of the action $S_2$ is less than half of the number of the basis 
$\{\eb_i\}$ on the space $\{\phi |\fpart{S}{\phi}=0\}$, 
though the number of the basis is infinity. 
The origin $\phi=0$ is also the solution for $\{\phi |\fpart{S}{\phi}=0\}$, 
and eq.(\ref{hessian}) at the origin 
is nothing but the condition $(Q)^2=0$. 
Hereafter for simplicity we choose the origin for the solution 
for $\{\phi |\fpart{S}{\phi}=0\}$ (without loss of generality). 
\begin{defn}[Proper]
At the origin $\phi=0$ if the ratio of the rank of the Hessian 
over the number of the basis is just half, 
the action is called {\it proper} (see \cite{BV1,BV2,HT,GPS}). 
(This is a traditional definition, but 
when we say an action is proper, we assume 
an additional condition stated later in eq.(\ref{Ttransf}).) 
 \label{defn:proper}
\end{defn}
The Hessian at the origin is $\V_{i_1 i_2}$ 
in eq.(\ref{action}), which is determined by $Q$. 
Let us consider the decomposition (\ref{cHdecomp}) 
$\cH=\cH^t\oplus\cH^u\oplus\cH^p$. 
The rank of the Hessian is equal to 
the rank of unphysical states $\cH^u$ which generate the gauge 
transformations. $\rank (\cH^u)$ is equal to $\rank (\cH^t)$, where 
$\cH^t$ is $Q$-trivial states. 
The condition of the proper is then equivalent to 
the condition that $\rank (\cH^u)/ \rank (\cH)=\half$. 
Note that it does not imply that $\rank (\cH^p)=0$. 
When an action is proper, $\cH^p$ corresponds to so-called 
the Green kernel and $\rank (\cH^p)/\rank (\cH^u)=0$ holds. 
String field theory is also proper 
where $Q$ is the BRST-operator \cite{KO} of conformal 
field theory on a fixed background. 
(The reducibility of 
the gauge group of string field theory action 
then comes from the Virasoro symmetry of $Q$.)

Given a proper action, the gauge fixing and the path-integral measure 
$\mes\Phi$ are defined as  follows. 
\begin{defn}[Gauge fixing]
Let us consider a degree minus one element $\Psi$ in $C(\phi)_c$, 
which is called {\it the gauge fixing fermion}. 
The power of fields for $\Psi$ 
is assumed to be greater than or equal to two. 
A {\it BV-gauge fixing} is defined as restriction of antifields to 
the lagrangian submanifold $\phi^*_a=\fpart{\Psi}{\phi^a}$. 
The path-integral measure in eq.(\ref{Z}) is then the 
integration over the space of fields $\{\phi^a\}$, 
the dual of the graded vector space $\cH_+$. 
We denote it by $\mes\Phi_{gf}$. 
 \label{defn:gf}
\end{defn}
Thus, choosing $\Psi$ determines the gauge fixing. 
In the original context of the BV-formalism, 
restricting the antifields to zero ($\Psi=0$) recovers 
the original action that consists of only degree zero fields, 
where the rank of the Hessian is possibly less 
than the rank of the fields. 
We call it the trivial gauge. 
The gauge fixing is then performed by shifting the trivial gauge
$\phi^*_a=0$ to $\phi^*_a=\fpart{\Psi}{\phi^a}$ 
so that the rank of the Hessian is 
equal to the rank of the fields, \ie half of the rank of the total
space $\cH$. 
In case of string field theories, 
however, the antifields are originally included 
in the quadratic term $S_2$, 
and BV-master equation is used in order to determine the form of higher 
vertices. 
Therefore the trivial gauge fixing 
can also be a candidate for consistent gauge fixing. 
This trivial gauge is called {\it Siegel gauge} in string field theory.

Let us now examine some properties of the kinetic term of a proper action. 
We write $Q\eb_j=\eb_kc^k_j$ and in matrix expression 
\begin{equation*}
 c:=\{c^k_j\}=\bp c_1 & c_2\\ c_3 & c_4\ep\ .
\end{equation*}
The kinetic term $\V_{ij}$ is then written as 
\begin{equation}
 \omega c=\bp -c_3 & -c_4 \\ c_1 & c_2\ep\ .
 \label{omegac}
\end{equation}
$\V_{ij}$ is graded symmetric since vertices are defined to be cyclic. 
This implies that $c$ satisfies 
\begin{equation*}
 (\omega c)_{ji}=(-1)^{\eb_i}(\omega c)_{ij}
\end{equation*}
where $\omega$ is the one in eq.(\ref{symD}). 
When we write $\dag$ for the transpose with the sign factor, 
the above equation becomes 
\begin{equation*}
\omega c=(\omega c)^\dag \ .
\end{equation*}
One then obtains that 
$c_2=c_2^\dag$, $c_3=c_3^\dag$ and $c_4=-c_1^\dag$.

Alternatively, when fixing a gauge, 
one can bring the gauge fixing condition $\phi^*_a=\fpart{\Psi}{\phi^a}$ 
to the form $\phi^*_a=0$ by a coordinate transformation 
of the form 
\begin{equation*}
 \Phi'=\Phi+\Phi\flpartial{\phi^*_a}\frpart{\Psi}{\phi^a}\ .
\end{equation*}
The second term in the right hand side 
is written as $-(\Phi,\Psi)$. 
This transformation preserves the symplectic form, since 
it is a special case of the transformation eq.(\ref{fsymdiff}), 
that is, 
\begin{equation}
 \Phi'=e^{(\ ,\Psi)}\Phi\ . \label{gauge-fsymdiff}
\end{equation}
For the kinetic term $\omega c$, only the linear part of the 
coordinate transformation is relevant. 
When we represent the gauge fixing fermion as 
$\Psi_a=\phi^a\psi_{ab}\phi^b+\cdots$, the linear part of the 
transformation is as follows
\begin{equation*}
 \bp \1 & \0 \\ \psi & \1 \ep\ .
\end{equation*}
In this coordinate the kinetic term is just $-c_3$ 
since the gauge fixing is $\phi^*_a=0$. 
The rank of $c_3$ in eq.(\ref{omegac}) is then half of the 
rank of total space. 

Though $c_3$ is degenerate because of the Green kernel generally, 
let us first consider the case $c_3$ is nondegenerate, where 
the cohomology with respect to $Q=m_1$ is trivial. 
In this situation, the condition $c^2=0$ implies that $c$ can be written as 
\begin{equation*}
 \omega c=T^\dag
   \bp -c_3 & \0 \\ \0 & \0\ep
   T\ ,\qquad 
 T:=\bp \1 & (c_3)^{-1}c_1 \\ \0 & \1\ep\ .
\end{equation*}
Note that $T$ preserves the BV-symplectic form in the Darboux
coordinate
\begin{equation}
 T^\dag \bp \0 & -\1 \\ \1 & \0\ep  T =\bp \0 & -\1 \\ \1 & \0\ep\ .
 \label{Tpres}
\end{equation}

By the definition of proper (Definition \ref{defn:proper}) 
it is natural that the properties of the kinetic term above holds 
also in general situations. 
Namely, when we say an action is proper, we assume that, 
in the coordinate where the gauge fixing is $\phi^*_a=0$, 
there exists a linear transformation 
\begin{equation*}
 T=\bp \1 & t \\ \0 & \1\ep\ ,\qquad t-t^\dag=0
\end{equation*}
which preserves the symplectic form (\ref{Tpres}) and 
transforms the kinetic term $\{c^k_j\}$ of the form 
\begin{equation}
 \bp c_1 & c_2\\ c_3 & c_4\ep 
 =T^{-1}\bp 0 & 0 \\ c_3 & 0 \ep T\ .
 \label{Ttransf}
\end{equation}
One can see that $\cH^p$ is just linearly isomorphic to two copies 
of the kernel of $c_3$. 
$t=(c_3)^{-1}c_1$ is a solution of eq.(\ref{Ttransf}). 
On the other hand, $c^2=0$ leads to $P_L c_4=c_4$ and $c_1P_L=c_1$, and 
in addition $P_L c_1=c_1$ and $c_4P_L=c_4$ hold 
if eq.(\ref{Ttransf}) is satisfied.

Given a gauge fixing $\Psi$, 
a propagator is constructed canonically as follows. 
\begin{defn}[Propagator in BV-formalism]
Let $\Psi\in C(\phi)_c$ be a gauge fixing fermion. 
By this gauge fixing, the quadratic term of the action is 
written in terms of fields only (not antifields) as 
$\half\phi^a(\V_{gf})_{ab}\phi^b$ 
for some graded symmetric (\ie cyclic) matrix $\V_{gf}$. 
It is also regarded as a bilinear map 
$\V_{gf}:\cH_+\otimes\cH_+\to\C$ such that 
$\V_{gf}(\eb_a,\eb_b)=(\V_{gf})_{ab}$. 
The gauge fixing fermion is taken so that the rank of $\V_{gf}$ is 
maximal. 
Note that $\V_{gf}$ is degenerate only for on-shell states. 
The degeneracy corresponding to gauge orbits 
(orbits of the gauge transformations) is killed 
by the gauge fixing. 
Let $P_{gf}:\cH_+\to\cH_+$ be the projection onto the kernel of $\V_{gf}$. 
The {\it BV-propagator} $\V_{gf}^+$ is then given 
by the inverse of $\V_{gf}$ such that 
\begin{equation}
 \V_{gf}^+\V_{gf}=\V_{gf}\V_{gf}^+=\1-P_{gf} 
  \label{propagator}
\end{equation}
on $\cH_+$ in the matrix expression. 
 \label{defn:BVpropagator}
\end{defn}
\begin{prop}
A propagator in the BV-formalism is a homotopy operator of $(\cH, Q)$.  
 \label{prop:BVp}
\end{prop}
\begin{pf}
As stated previously, 
in the coordinate where the gauge fixing is $\phi^*=0$, 
the gauge fixed kinetic term is $\V_{gf}=-c_3$. 
The propagator $\V_{gf}^+\in\cH_+\otimes\cH_+$ given in 
Definition \ref{defn:BVpropagator} is naturally extended to 
the one in $\cH\otimes\cH$. We denote it by $\V_L^+$. 
Let us define a degree one operator $Q^+:\cH\to\cH$, 
$Q^+(\eb_i)=\eb_j\cb^j_i$ in matrix expression by 
\begin{equation}
 \cb:=\V_L^+\omega=\bp \V_{gf}^+ & \0 \\ \0 & \0 \ep\omega
 =\bp \0 & \V_{gf}^+\\ \0 & \0 \ep\ .
 \label{propagatormatrix}
\end{equation}
$Q^+$ then satisfies the Hodge-Kodaira decomposition (\ref{HKdecomp})
\begin{equation}
 QQ^+ +Q^+Q+P=\1\ ,
 \label{HKdecomp2}
\end{equation}
where $P:\cH\to\cH$ is the projection onto the on-shell state. That is, 
$Q^+$ is a homotopy operator. 
This is because the action is proper. 
More precisely, 
one can make a coordinate transformation $T$ 
in eq.(\ref{Ttransf}) where the kinetic term $\omega c$ is transformed to be 
of the form 
\begin{equation*}
 \bp -c_3 & \0 \\ \0 & \0\ep\ .
\end{equation*}
Note that $T$ preserves the matrix (\ref{propagatormatrix}) corresponding 
to the homotopy operator $Q^+$. 
In this coordinate it is clear that eq.(\ref{HKdecomp2}) holds. 
\qed
\end{pf}

 \subsection{Path integral, Feynman diagram and the minimal model theorem}
\label{ssec:PI}

In subsection \ref{ssec:minimal}, 
we obtained an explicit form of the minimal model 
$(\cH^p,\ti\omega^p,\ti\m^p)$ 
of cyclic $A_\infty$-algebra $(\cH,\omega,\m)$. 
There, a cyclic $A_\infty$-morphism $\ti\cF^p$ from 
minimal cyclic $A_\infty$-algebra 
$(\cH^p,\ti\omega^p,\ti\m^p)$ to $(\cH,\omega,\m)$ was also constructed. 
Here let $(\cH,\omega,\m)$ be the cyclic $A_\infty$-algebra 
of a field theory whose action is proper (Definition \ref{defn:proper}). 
Then 
we can see that the $n$-point vertex defined by $A_\infty$-structure 
$\ti\m^p$ is nothing but the tree level $n$-point correlation function 
of the field theory(Lemma \ref{lem:main}). 

This subsection is devoted to show 
that the scattering amplitudes 
of the field theory computed by the Feynman rule coincides with 
$\pm\ov{n}\omega(\eb^p_{i_1},
\ti{m}^p_{n-1}(\eb^p_{i_2},\dots,\eb^p_{i_n}))$ where 
$\eb^p_{i_1},\dots,\eb^p_{i_n}\in\cH^p$ 
are the external states of the amplitude.

Let us first write down the explicit definition of the 
perturbative expansion. 
It is obtained by fixing the gauge, constructing the propagator and 
eq.(\ref{Z}). 
Let $\cO(\Phi)\in TC(\phi)_c$ be any operators. 
Then its path-integral is formally given by 
\begin{equation}
 \la\cO(\Phi)\ra\sim\int\mes\Phi\ \cO(\Phi)e^{-S}\Big|_{\gf}
 =\int\mes\Phi_{gf}\ \cO(\Phi_{\gf})
  e^{-\sum_{k\ge 3}\ov{k}\V_k(\Phi_{\gf},\dots,\Phi_{\gf})}
  e^{-S_2|_{\gf}}\ ,
 \label{ev}
\end{equation}
where $|_{\gf}$ denotes a gauge fixing (Definition \ref{defn:gf}) 
in the previous subsection, 
$\Phi_{\gf}$ denotes the gauge fixed $\Phi$, and 
$S_2|_{\gf}=\half\omega(\Phi_{\gf},Q\Phi_{\gf})$. 
The path integral (\ref{ev}) is not well-defined precisely 
since we do not define the integral $\mes\Phi_{gf}$ on a formal 
noncommutative supermanifold. 
Instead, we shall define the path integral precisely at the level of
the perturbative expansion below. 
Since $S_2|_{\gf}$ is quadratic with respect to $\Phi_{\gf}$, 
$e^{-S_2|_{\gf}}$ is a gaussian. Then 
the perturbative expansion is essentially the gaussian integral 
where $e^{-\sum_{k\ge 3}\ov{k}\V_k(\Phi_{\gf},\dots,\Phi_{\gf})}$ 
is Taylor expanded. 
In field theory it is calculated by so-called Wick contraction. 
One can rewrite it in a purely algebraic manner as follows. 
We take it as the starting point of our definition. 
\begin{defn}[Perturbative expansion]
For $\cO(\Phi)\in TC(\phi)_c$, 
The perturbative expansion of $\cO(\Phi)$ is defined by 
\begin{equation}
 \la\cO(\Phi)\ra=
 \left.
 \l(\cO(\Phi)\cdot
 e^{-\sum_{k\ge 3}\ov{k}\V_k(\Phi,\dots,\Phi)}\r)
 e^{\l(\half\V_L^{+,ij}\flpartial{\phi^i}\flpartial{\phi^j}\r)}
\right|_{\Phi=0}\ ,
 \label{ev2}
\end{equation}
where $\V_L^+$ is the propagator constructed in the gauge $|_{gf}$ as in 
Definition \ref{defn:BVpropagator} and the below. 
Here the derivation $\flpartial{\phi^i}$ acts 
on $\l(\cO(\Phi)\cdot
 e^{-\sum_{k\ge 3}\ov{k}\V_k(\Phi,\dots,\Phi)}\r)$ from right 
with an appropriate Kostul sign. 
 \label{defn:feynman-pre}
\end{defn}
The above expression is the one derived directly from eq.(\ref{ev}) 
only if the gauge fixing fermion is quadratic
$\Psi=\phi^a\psi_{ab}\phi^b$ 
(since otherwise $\V_k(\Phi_{gf},\dots,\Phi_{gf})$ in eq.(\ref{ev}) 
includes terms higher than $k$ powers of fields $\{\phi^a\}$). 
However, as stated in the previous subsection, 
a gauge fixing is equivalent to a particular field transformation 
preserving the symplectic form; 
perturbative expansion with gauge fixing $\Psi$ is 
performed by trivial gauge $\Psi=0$ after 
field redefinition (\ref{gauge-fsymdiff}). 
Thus, we indicate by $\V_k(\Phi,\dots,\Phi)$ in eq.(\ref{ev2}) 
the term of $k$ powers after the such a field redefinition associated
to $\Psi$. 
The dependence or independence of the choice of the gauge fixing 
will be stated in the end of this subsection. 

The perturbative expansion (\ref{ev2}) gives 
a well-defined linear map $\la\cdots\ra: TC(\phi)_c\to\C$. 
Especially the path integral above reduces to the ordinary one 
if the fields are (graded) commutative. 
\footnote{For the Feynman rule 
in operator language, see \cite{Ka}. 
Comparing it from the argument here in component field theory picture, 
one can see the equivalence between them. 
The reference \cite{Th} also provides us with useful informations. } 
The value (\ref{ev2}) is calculated by {\it Feynman rules} as follows. 
Let us consider eq.(\ref{ev2}) in the case when 
$\cO(\Phi)=a_1(\phi)\bullet\cdots\bullet a_n(\phi)$ 
where $a_r(\phi)
:=\ov{k_r}a_{i^r_1\cdots i^r_{k_r}}\phi^{i^r_{k_r}}\cdots\phi^{i^r_1}
\in C(\phi)_c$ for each $1\le r\le n$. 
We call each $a_r(\phi)$ an {\it observable vertex} and assign 
$o_r$ to it. Moreover we associate $o^r_a$, $1\le a\le k_r$ 
to each field in $a_r(\Phi)$. 
On the other hand,  
$e^{-\sum_{k\ge 3}\ov{k}\V_k(\Phi,\dots,\Phi)}$ is Taylor expanded. 
Each term is characterized by a multiindex 
$\Lambda:=\{\lambda_3\in\Z_{\ge 0},\lambda_4\in\Z_{\ge 0},\dots\}$ 
and given by 
\begin{equation*}
 \Pi_{l\ge 3}\ov{\lambda_l!}\l(-\ov{l}\V_l(\Phi)\r)^{\lambda_l}
 =\ov{\lambda_3!}\l(-\ov{3}\V_3(\Phi)\r)^{\lambda_3}\bullet
\ov{\lambda_4!}\l(-\ov{4}\V_4(\Phi)\r)^{\lambda_4}\bullet\cdots\ .
\end{equation*}
It is a product of $m:=(\sum_{l\ge 3}\lambda_l)$ numbers of vertices. 
Namely, the equation above is 
\begin{equation*}
 \V_{e_1}(\Phi)\bullet\cdots\bullet\V_{e_m}(\Phi)
\end{equation*}
up to an appropriate coefficient, 
where $3\le e_1\le\cdots\le e_m$. 
To each $\V_{e_q}(\Phi)$, $1\le q\le m$, we assign a {\it vertex} $v_q$. 
Moreover we associate $v^q_a$, $1\le a\le e_q$ to each field 
in $\V_{e_q}(\Phi)$. 
We indicate both observable vertices and vertices by 
(observable) vertices. 
The set of fields of (observable) vertices is denoted by 
\begin{equation*}
 Vert(\cO,\Lambda)
 :=\{
 o^1_1,..,o^1_{k_1},\dots .\, , o^n_1,..,o^n_{k_n}
 ,v^1_1,..,v^1_{e_1},\dots .\, ,v^m_1,..,v^m_{e_m}
 \}\ .
\end{equation*}

In this situation, consider the value 
\begin{equation*}
 \left.
\l(\cO(\Phi)\bullet \Pi_{l\ge 3}\ov{\lambda_l!}
 \l(-\ov{l}\V_l(\Phi)\r)^{\lambda_l}\r)
 e^{\l(\half\V_L^{+,ij}\flpartial{\phi^i}\flpartial{\phi^j}\r)}
\right|_{\Phi=0}\ .
\end{equation*}
The value vanishes if $(k_1+\cdots+k_n)+(e_1+\cdots +e_m)$ is odd. 
When it is even, the equation above is equal to 
\begin{equation}
 \l(\cO(\Phi)\bullet \Pi_{l\ge 3}\ov{\lambda_l!}
 \l(-\ov{l}\V_l(\Phi)\r)^{\lambda_l}\r)
 \ov{2^{I'} I'!}\l(\V_L^{+,ij}\flpartial{\phi^i}\flpartial{\phi^j}\r)^{I'}
 \label{ev4}
\end{equation}
where $2I'=(k_1+\cdots+k_n)+(e_1+\cdots +e_m)$. 
We have $2I'$ differentials that act on (observable) vertices. 
We assign $ed^p_a$, $1\le p\le I'$, $a=1, 2$ to 
each differential and denote the set by 
\begin{equation*}
 Edge(I'):=\{ed^1_1,ed^1_2, ed^2_1,ed^2_2,\dots,
 ed^{I'}_1,ed^{I'}_2\}\ . 
\end{equation*} 
Consider isomorphisms from the set $Edge(I')$ to 
the set $Vert(\cO,\Lambda)$. 
We have $(2I')!$ such isomorphisms. 
Denote the set of the isomorphisms by $\ti{F}(\cO,\Lambda)$.

Then eq.(\ref{ev2}) is represented in the following form 
\begin{equation}
 \Pi_{l\ge 3}\ov{\lambda_l!}\ov{l^{\lambda_l}}\ov{2^{I'} I'!}\cdot
 \sum_{\ti{\Upsilon}(\cO,\Lambda)\in\ti{F}(\cO,\Lambda)} 
 N_{\ti{\Upsilon}(\cO,\Lambda)}\ ,
 \label{ev5}
\end{equation}
where $N_{\ti{\Upsilon}(\cO,\Lambda)}\in\C$ is given by the product 
of $a_{i^r_1\cdots i^r_{k_r}}$, $\V_{i^q_1\cdots i^q_{e_q}}$, 
$\V_L^{+,ij}$ and the Kostul sign factor. 
One can consider two actions on the set $Edge(I')$; 
the exchange between $ed^p_1$ and $ed^p_2$ for each $p$, 
and the exchange between the pair $(ed^p_1,ed^p_2)$ and 
$(ed^{p'}_1,ed^{p'}_2)$ for any $p\ne p'$. 
One can also consider an action on $Vert(\cO,\Lambda)$, 
the cyclic permutations in $o^r_1,\dots,o^r_{k_r}$ for each $r$ or 
in $v^q_1,\dots,v^q_{e_q}$ for each $q$. 
These actions induce automorphisms on $\ti{F}(\cO,\Lambda)$ 
and $N_{\ti{\Upsilon}(\cO,\Lambda)}$ is independent of the
automorphisms. 
Let us denote by $F(\cO,\Lambda)$ the set of 
isomorphisms $\ti{F}(\cO,\Lambda)$ over these automorphisms. 
Moreover we introduce the direct sum, 
${\displaystyle F(\cO):=\bigoplus_{\Lambda} F(\cO,\Lambda)}$. 

The perturbative expansion eq.(\ref{ev2}) is the sum of 
eq.(\ref{ev5}) with respect to $\Lambda$. It is rewritten as 
\begin{equation}
 \la\cO(\Phi)\ra=
 \sum_{\Upsilon(\cO)\in F(\cO)}\ov{\Lambda !} 
  N_{\Upsilon(\cO)}\ , 
 \label{ev6}
\end{equation}
where $\Lambda !:=\Pi_{l\ge 3}\ov{\lambda_l!}$. 
Note that ${\displaystyle \ov{l^{\lambda_l}}\ov{2^{I'} I'!}}$ in
eq.(\ref{ev5}) is canceled by the automorphisms. 
Associated to each $\Upsilon(\cO)\in F(\cO,\Lambda)$, we define 
Feynman graphs for a cyclic field theory $(\cH,\omega,S)$ as follows. 
\begin{defn}[Feynman graph]
In the situation above, 
let us arrange the observable vertices $o_r$ and 
vertices $v_q$ from left to right on a plane 
such as 
\begin{equation*}
 \bullet_{o_1}\ \ \ \bullet_{o_2}\ \ \ \cdots \ \ \ \bullet_{o_n}\ \ \ 
 \bullet_{v_1}\ \ \ \cdots\ \ \ \bullet_{v_{m-1}}\ \ \ \bullet_{v_m}\ \ \ .
\end{equation*}
Moreover, connect any two (observable) vertices to each other 
by edges so that 
the number of incident edges is $k_r$ at observable vertex $o_r$ 
and $e_q$ at vertex $v_q$. 
We call such graphs the {\it Feynman graphs} for 
a cyclic field theory $(\cH,\omega, S)$. 
Hereafter we identify an element $\Upsilon(\cO)\in F(\cO)$ 
with a Feynman graph. 
The cyclic order of the edges around each (observable) vertex 
is distinguished. 
Two edges that intersect on the plane can pass through each other 
and two graphs before and after this process are not distinguished. 
 \label{defn:feynman}
\end{defn}
Connecting two (observable) vertices by an edge is called 
the {\it Wick contraction}. The contraction indicates that 
two differentials in $\half\V_L^{+,ij}\flpartial{\phi^i}\flpartial{\phi^j}$ 
act on the fields in the two (observable) vertices. 

Now we are interested in tree on-shell amputated amplitudes. 
\begin{defn}[Feynman graphs for $n$-point amplitudes]
An $n$-point amplitude is calculated by eq.(\ref{ev6}) 
with $a_r(\Phi)=\phi^{j_r}$, that is, 
$\cO(\Phi)=\phi^{j_1}\bullet\cdots\bullet\phi^{j_n}$. 
Each $\phi^{j_r}$ is called an {\it external field}. 
We call the corresponding vertex a {\it external vertex}. 
A Feynman graph $\Upsilon(\cO)\in F(\cO)$ is called 
{\it connected} if any two (observable) vertices of $\Upsilon(\cO)$ 
is connected to each other by the edges. 
When $\Upsilon(\cO)$ includes a circle $S^1$ 
consisting of the edges, the $S^1$ is called a {\it loop}. 
The set of connected Feynman graphs is 
denoted by $F^{conn}(\cO)\subset F(\cO)$. 
Moreover we denote the set of every connected tree 
Feynman graphs by $\T(\cO)\subset F^{conn}(\cO)$, 
where a {\it tree} Feynman graph means a Feynman graph without loops. 
The value 
\begin{equation*}
 \left.\la\cO(\Phi)\ra\right|_{conn}:=\sum_{\Upsilon(\cO)\in F^{conn}(\cO)}
 \ov{\Lambda !} N_{\Upsilon(\cO)}
\end{equation*}
is then the {\it $n$-point amplitude}. Moreover, 
restricting the Feynman graphs in $F^{conn}(\cO)$ to those 
in $\T(\cO)$ one gets the {\it $n$-point tree amplitude}, 
\begin{equation*}
 \left.\la\cO(\Phi)\ra\right|_{\substack{conn \\ tree}}
 :=\sum_{\Upsilon(\cO)\in\T(\cO)}\ov{\Lambda !} N_{\Upsilon(\cO)}\ .
\end{equation*}
 \label{defn:npt}
\end{defn}
For a connected tree Feynman graph, 
the Wick contraction by the edges can be divided 
into two processes; the contractions between $n$ external fields 
$\phi^{j_1}\bullet\cdots\bullet\phi^{j_n}$ and 
the vertices, and the contractions between the vertices. 
As explained below, 
the latter process produces some function of $n$ powers of 
$\Phi$ that associates to planar tree graphs 
(Definition \ref{defn:planar} below). 
We shall define cyclic functions associated to the planar tree
graphs in Definition \ref{defn:cyclicplanar} 
as we defined the multilinear map $\ti{m}_{\Gamma_k}$ associated to 
rooted planar $k$-tree. 
After the latter process, the former one, 
the contractions of the cyclic function with $n$ external fields, 
finishes the calculation for the value associated to the 
Feynman graph. 
\begin{defn}[Cyclic function associated to planar graphs]
Let $G^{cyc}_n$ be the set of planar tree graphs with $n$-leaves. 
An element $\Gamma^{cyc}_n\in G^{cyc}_n$ is a rooted planar $(n-1)$-tree 
without the distinction of the root edge and the leaves, that is, 
the root edge is regarded as a leaf.
\footnote{Here a notation is changed compared to that in \cite{Ka}. 
$\Gamma^{cyc}_n$ here is denoted by $\Gamma^{cyc}_{n-1}$ in 
\cite{Ka}. } 
Denote the natural surjection by $\rh: G_{n-1}\to G^{cyc}_n$. 
Now we regard $G_{n-1}$ and $G^{cyc}_n$ as vector spaces and 
denote them also by themselves. 
Namely, for the set $G_{n-1}$ (resp. $G^{cyc}_n$), 
their elements are regarded as the bases 
of the vector space $G_{n-1}$ (resp. $G^{cyc}_n$). 
$\rh: G_{n-1}\to G^{cyc}_n$ is then regarded as a vector bundle. 
For an element $\Gamma^{cyc}_n\in G^{cyc}_n$, 
there exist $n$ choices to pick up one of the leaves as the root edges. 
Summing over these $n$ numbers of $(n-1)$-trees and dividing by $n$ 
defines a section $s: G^{cyc}_n\to G_{n-1}$. 

First we define a 
$(\cH)^{\otimes k}\to\cH$ valued linear function $\ti{m}^{cyc}_k$ on 
vector space $G_k$ in a similar way 
as in Definition \ref{defn:minimal2}. 
To an elementary $k$-tree we associate $m_k:\cH^{\otimes k}\to\cH$. 
For any $k$-tree $\Gamma_k$ denote the associated endomorphisms by 
$\ti{m}^{cyc}_{\Gamma_k}: (\cH)^{\otimes k}\to\cH$. 
For a grafting $\Gamma_k\circ_i\Gamma_l$ we associate 
\begin{equation*}
 \ti{m}^{cyc}_{\Gamma_k\circ_i\Gamma_l}=
 \ti{m}^{cyc}_{\Gamma_k}\circ
\l(\1^{\otimes (i-1)}\otimes\ti{f}_{\Gamma_l}\otimes\1^{\otimes (k-i)}\r)
: \cH^{\otimes (k+l-1)}\to\cH\ ,
\end{equation*}
where $\ti{f}_{\Gamma_l}:(\cH)^{\otimes l}\to\cH$ are those in 
Definition \ref{defn:minimal2}. 
$\ti{m}^{cyc}_{\Gamma_k}$ are then defined 
so that they are compatible with this grafting. 
Note that 
$\ti{m}_{\Gamma_k}=\pi\circ\ti{m}^{cyc}_{\Gamma_k}$, that is, 
removing $\pi$ 
on the root edge in Definition \ref{defn:minimal2} leads 
to $\ti{m}^{cyc}$. 

$(\cH)^{\otimes n}\to\C$ valued cyclic function on vector space
$G^{cyc}_n$ is then defined by 
\begin{equation*}
 \ov{n}\ti\V_{\Gamma^{cyc}_n}=(s^*\ti{m}^{cyc})(\Gamma^{cyc}_n)
 :=\omega\circ\l(\1\otimes\l(\ti{m}^{cyc}(s(\Gamma^{cyc}_n))\r)\r)\ .
\end{equation*}
 \label{defn:cyclicplanar}
\end{defn}
\begin{defn}[(Amputated) tree on-shell scattering amplitude]
For each element in $\T(\cO)$,
$\cO(\Phi)=\phi^{j_1}\bullet\cdots\bullet\phi^{j_n}$, 
remove the external vertices 
together with the edges whose one end is the external vertices. 
One gets a tree graph with $n$ free ends. 
Denote by $\T_n$ the set of the graphs obtained in such a way. 
We denote the surjection by $Ampu: \T(\cO)\to \T_n$ and call it 
the {\it amputation map}. 
It is an $n!$-to-one map. 
We regard $\T(\cO)$ and $\T_n$ also as vector spaces. 
Any element in $\T_n$ is isomorphic to an element in $G^{cyc}_n$ 
as a planar tree graph. 
Thus we have a surjection $t : \T_n\to G^{cyc}_n$. 
\begin{equation*}
 \xymatrix{
   &  & G_{n-1}\ \ar@<0.7ex>[d]^{\rh} \\
  F(\cO)\supset\T(\cO)\ \ar[r]^{\ \, \qquad Ampu} & \T_n \ar[r]^{t} & 
  G^{cyc}_n \ar@<0.7ex>[u]^{s}
}
\end{equation*}
The {\it (amputated) tree correlation functions} 
for a cyclic field theory $(\cH,\omega, S)$ is 
the collection of the following cyclic functions $\{\ti\V_v\}_{n\ge 3}$, 
\begin{equation}
 \ti\V_n:=
 \sum_{\Upsilon_n\in \T_n}\ov{\Lambda !} \ti\V_{t(\Upsilon_n)}\ .
 \label{Smatrix}
\end{equation}
The {\it tree on-shell correlation functions} is given by 
\begin{equation*}
 \ti\V^p_n:=\ti\V_n\circ (\iota)^{\otimes n}\ 
 : (\cH^p)^{\otimes n}\to\C\ .
\end{equation*}
It is called also the {\it tree on-shell scattering amplitudes} or 
the {\it tree on-shell S(cattering)-matrices}. 
 \label{defn:Smatrix}
\end{defn}
\begin{rem}
The scattering amplitudes are usually defined on the gauge fixed
subspace of $\cH$. In this sense the scattering amplitude in the 
Definition above is an extended one that is defined on the whole 
graded vector space $\cH$. 
 \label{rem:Smatrixgf}
\end{rem}
\begin{lem}
For $\cO(\Phi)=\phi^{j_1}\bullet\cdots\bullet\phi^{j_n}$ and 
for a fixed $\Upsilon_n\in \T_n$, we have 
\begin{equation}
\V_L^{+,j_1 j'_1}\frpartial{\phi^{j'_1}}\cdots
\V_L^{+,j_n j'_n}\frpartial{\phi^{j'_n}}\cdot
\l(-\ov{n}\ti\V_{t(\Upsilon_n)}(\Phi,\dots,\Phi)\r)
=\sum_{\Upsilon(\cO)\in Ampu^{-1}(\Upsilon_n)}N_{\Upsilon(\cO)}\ .
 \label{tiVampu}
\end{equation}
Consequently, 
\begin{equation*}
 \sum_{\Upsilon(\cO)\in \T(\cO)} 
 \ov{\Lambda !} N_{\Upsilon(\cO)}
=\V_L^{+,j_1 j'_1}\frpartial{\phi^{j'_1}}\cdots
 \V_L^{+,j_n j'_n}\frpartial{\phi^{j'_n}}\cdot 
 \l(-\ov{n}\ti\V_n(\Phi,\dots,\Phi)\r)
\end{equation*}
holds \cite{Ka}.
 \label{lem:main1}
\qed\end{lem} 
This implies that the Definition \ref{defn:Smatrix} 
actually gives the amputated $n$ point tree amplitudes 
in a usual sense in field theory. 
\begin{lem}[\cite{Ka}]
$\{\ti\V_n\}_{n\ge 3}$ is given by 
\begin{equation}
 \ti\V_n
 =\sum_{\Gamma_{n-1}\in G_{n-1}}
\omega\circ\l(\1\otimes\ti{m}^{cyc}_{\Gamma_{n-1}}\r)
=\sum_{\Gamma^{cyc}_n\in G^{cyc}_n}
 \ov{\sharp\l(\rh^{-1}(\Gamma^{cyc}_n)\r)}
 \ti\V_{\Gamma^{cyc}_n}\ . 
 \label{tiV}
\end{equation}
Here $\sharp\l(\rh^{-1}(\Gamma^{cyc}_n)\r)$ indicates the 
number of the elements of the set $\rh^{-1}(\Gamma^{cyc}_n)$. 
 \label{lem:main}
\qed\end{lem}
In the terminology of Feynman graphs, 
$1/\sharp\l(\rh^{-1}(\Gamma^{cyc}_n)\r)$ is the {\it symmetric factor} 
of graph $\Gamma^{cyc}_n$.
\footnote{The number $\sharp\l(\rh^{-1}(\Gamma^{cyc}_n)\r)$ coincides 
with the number of the automorphisms acting on $\Gamma^{cyc}_n$ 
which we denoted by $\sharp\mathrm{Aut}(\Gamma^{cyc}_n)$ 
in eq.(\ref{IntrotiV}) in the Introduction. }
\begin{cor}[\cite{Ka}]
For a given cyclic field theory $(\cH,\omega, S)$, 
the tree on-shell correlation functions
are given by 
$$
 \ti\V^p_n=\omega\circ\l(\1\otimes\ti{m}^p_{n-1}\r)
$$ 
and therefore they form 
the minimal cyclic $A_\infty$-algebra $(\cH^p,\ti\omega^p,\ti\m^p)$ 
in eq.(\ref{Sp}). 
 \label{cor:main}
\end{cor}
For the proof of Lemma \ref{lem:main1}, 
one may calculate each $\Upsilon(\cO)\in \T(\cO)$ 
using the correspondence 
\begin{equation*}
 \flpartial{\phi^i}\V_L^{+,ij}\frpartial{\phi^j}
 \l(\ov{k+1}\V_{k+1}(\Phi)\r)
 =\flpartial{\phi^i}{\bar c}^i_lc^l(\phi)\ ,
 \label{vtom}
\end{equation*}
that comes from 
$\V_L^{+,ii_1}\omega_{i_1l}={\bar c}^i_l$ 
where $Q^+(\eb_l)=\eb_i{\bar c}^i_l$. 
Lemma \ref{lem:main}, that is, the equivalence of 
eq.(\ref{Smatrix}) and (\ref{tiV}) follows from 
concentrating the inverse of $t$ for each element in $G^{cyc}_n$. 
\qed
\begin{rem}
Corollary \ref{cor:main} indicates that 
the collection of on-shell correlation functions 
in Definition \ref{defn:Smatrix} 
forms the `minimal action' $\ti{S}(\ti\Phi^p)$ in eq.(\ref{Sp}), 
\begin{equation*}
 \ti{S}(\ti\Phi^p)=\sum_{k\ge 2}\ov{k+1}
 \omega(\ti\Phi^p,\ti{m}^p_k(\ti\Phi^p))\ .
\end{equation*}
The action $\ti{S}(\ti\Phi^p)$ is in fact an effective action 
in the following sense. 
$\ti{S}(\ti\Phi^p)$ is obtained by substituting 
$\Phi|_{gf}=\ti\cF^p(\ti\Phi^p)$ 
into $S(\Phi)$ as explained above. When we express 
$\Phi=\Phi^t+\Phi|_{gf}=\Phi^t+\Phi^p+\Phi^u$ where 
$\Phi^t$, $\Phi^p$ and $\Phi^u$ denotes 
the trivial, physical and unphysical modes of $\Phi$, 
respectively, the substitution means $\Phi^p=\ti\Phi^p$ and 
$\Phi^u=f(\ti\Phi^p)$. 
As seen from eq.(\ref{eom3}), 
the latter is nothing but the equation of motion for $\Phi^u$. 
Moreover, $\ti{S}(\ti\Phi^p)$ is related to $S(\Phi)$ by integrating 
$\Phi^u$ at tree level through a gauge fixing as 
\begin{equation}
 \int\mes\Phi^u\ e^{-S(\Phi)}|_{gf}=e^{-\ti{S}(\ti\Phi^p)} \ ,
\end{equation}
under an appropriate definition of the integration. 
In this sense the action $\ti{S}(\ti\Phi^p)$ is an effective action. 
The pair of gauge fixing $|_{gf}$, which extract 
$\Phi^t$, and 
the extraction of $\Phi^u$ by substituting is an analogue of 
symplectic reduction, 
where the extraction of $\Phi^u$ can also be regarded as the 
restriction $\Phi|_{S(\Phi)=0}$ (see also a comment 
in Remark \ref{rem:marginal}). 
Also, it is clear that our arguments are applicable to constructions 
of any other tree level effective actions. 
 \label{rem:eff}
\end{rem}
\begin{rem}[Independence of the choice of gauge fixing]
We encoded the dependence of gauge fixing $\Psi$ in the following two; 
a cyclic $A_\infty$-isomorphism from 
the original cyclic $A_\infty$-algebra $(\cH,\omega, S)$, 
and propagator $\V_L^+$. 
The higher terms than the quadratic term of $\Psi$ are absorbed into 
the cyclic $A_\infty$-isomorphism. Whereas, 
a propagator $\V_L^+$ is derived by a well-defined gauge fixing. 
Note that $\V_L^+$ determines the decomposition 
$\cH=\cH^t\oplus\cH^p\oplus\cH^u$. 
Any $\V_L^+$ derived by a well-defined gauge fixing 
picks up the same (or isomorphic) minimal part $\cH^p$. 

Suppose that we have two minimal cyclic $A_\infty$-algebra 
$(\cH^p,\ti\omega^p,\ti\m^p)$ and $(\cH^{p'},\ti\omega^{p'},\ti\m^{p'})$ 
that are obtained by two gauge fixing conditions. 
The uniqueness of minimal model (Corollary \ref{cor:uniqueness}) 
implies that minimal cyclic $A_\infty$-algebras obtained by
perturbative expansion are independent of 
the choice of the gauge fixing at least up to 
isomorphisms on them. 

In general in field theory, it is known that 
S-matrices are `invariant' under 
(certain class of) field redefinition, 
which is called an {\em equivalence theorem} \cite{KOrS}. 
Physically, the fact stated above might be thought of as a version of 
this theorem. 
 \label{rem:gf-indep} 
\end{rem}

 \subsection{Equivalence of classical open string field theories}
\label{ssec:main}

In this subsection we shall apply the relation between 
the minimal model theorem and Feynman graph in the previous subsection 
to classical open string field theories constructed as in 
subsection \ref{ssec:osft}. 
For the family of the well-defined string field theories, 
Lemma \ref{lem:main} gives us 
a classification of string field theories (Theorem \ref{thm:main})
described below. 
The field transformations induce 
one-to-one correspondence of moduli spaces of classical solutions 
between such string field theories in the context of deformation theory. 
We shall explain that the classical solutions are regarded as those 
corresponding to marginal deformations.

A classical open string field theory is defined 
on a fixed {\em conformal background} (see the beginning of 
subsection \ref{ssec:osft}).  
For a fixed conformal background, 
a $\Z$-graded vector space $\cH$ is given canonically. 
$\cH$ is called an open string Hilbert space, where 
each base $\eb_i$ is called an open string state and 
the associated dual coordinate $\phi^i$ is a field 
in the sense of field theory. 
The superfield $\Phi:=\eb_i\phi^i\in\cH\otimes\cH^*$ 
is then called the {\em string field}. 
Moreover, a degree one coboundary operator $Q:\cH\to\cH$ and 
an odd constant symplectic structure $\omega(\ ,\ ):\cH\otimes\cH\to\C$ 
are defined canonically. 
$Q$ and $\omega$ are called 
the BRST-operator \cite{KO}
\footnote{Here BRST-operator indicates the operator that induces the
BRST transformation in the sense of string world sheet theory. 
Do not confuse it with the (BV-)BRST transformation in string field theory 
explained in subsection \ref{ssec:gf}. }
and the BPZ-inner product \cite{BPZ}, 
respectively, on the fixed conformal background. 
The information with which open string theory provides us is the collection 
of {\em on-shell open string scattering amplitudes}. 
As a subset of the collection of the on-shell open string 
scattering amplitudes, we have 
{\em on-shell tree open string scattering amplitudes} 
which are cyclic multilinear maps $(\cH^p)^{\otimes k}\to\C$ 
for $k\ge 3$. 

Though string field theories are usually constructed by 
decomposing the moduli spaces of Riemann surfaces into cell 
as explained in subsection \ref{ssec:osft}, we shall give a purely 
algebraic definition for classical open string field theories as follows. 
\begin{defn}[Classical open string field theory : axiom]
A {\it classical open string field theory} is an action $S(\Phi)\in C(\phi)_c$ 
satisfying the following properties: 
\begin{itemize}
 \item[(a)] $\Phi$ is the string field, $\omega:\cH\otimes\cH\to\C$ is the 
BPZ-inner product 
and $Q:\cH\to\cH$ is the BRST operator on a fixed conformal background\ , 

 \item[(b)] The action $S(\Phi)$ is of the form 
${\displaystyle S(\Phi)=\half\omega(\Phi,Q\Phi)
+\sum_{k\ge 2}\ov{k+1}\omega(\Phi,m_k(\Phi,\dots,\Phi))}$ 
where $S(\Phi)\in C(\phi)_c$\ ,

 \item[(c)] the action $S(\Phi)$ satisfies 
the classical BV-master equation $(S, S)=0$ and is proper\ ,

 \item[(d)] the on-shell scattering amplitudes 
(Definition \ref{defn:Smatrix}) of the action $S(\Phi)$ 
by perturbative expansion reproduces the tree on-shell 
open string scattering amplitudes 
on the fixed conformal background. 
\end{itemize}
 \label{defn:axiom}
\end{defn}
The fact that $(\cH,\omega,S)$ is a cyclic $A_\infty$-algebra 
follows from the condition (b) and (c).

\begin{thm}[cf. \cite{Ka}]
 All the well-defined classical open string field theories 
which are constructed 
on a fixed conformal background are cyclic $A_\infty$-isomorphic 
to each other. 
 \label{thm:main}
\end{thm}
\begin{pf} 
Let $(\cH,\omega, S)$ be a cyclic
$A_\infty$-algebra describing a classical 
open string field theory on a fixed conformal background. 
Lemma \ref{lem:main} states that the collection of the 
on-shell scattering amplitudes 
for the action $S$ forms a minimal cyclic $A_\infty$-algebra 
just in the way given in subsection \ref{ssec:cycminimal}. 
We denote it by $(\cH^p,\omega^p,\ti{S}^p)$. 
As stated in the end of subsection \ref{ssec:cycminimal}, 
Theorem \ref{thm:cycMandC} further implies that 
$(\cH,\omega, S)$ are cyclic $A_\infty$-isomorphic to 
the decomposed cyclic $A_\infty$-algebra 
$(\cH^p,\omega^p,\ti{S}^p)\oplus (\cH^t\oplus\cH^p,Q)$, 
where $(\cH^t\oplus\cH^u, Q)$ is the linear contractible 
cyclic $A_\infty$-algebra 
such that 
$\cH=\cH^p\oplus(\cH^t\oplus\cH^u)$ and $Q$ is the restriction
of the original differential $Q:\cH\to\cH$ in $(\cH,\omega, S)$ onto 
$\cH^t\oplus\cH^u$. 
Suppose $(\cH,\omega,S')$ be another cyclic 
$A_\infty$-algebra of a classical 
open string field theory on the same conformal background. 
Then, it is cyclic $A_\infty$-isomorphic to 
the decomposed cyclic $A_\infty$-algebra 
$(\cH^p,\omega^p,\ti{S'}^p)\oplus (\cH^t\oplus\cH^p,Q)$. 
On the other hand, by Definition \ref{defn:axiom} (d), 
both of these two minimal cyclic $A_\infty$-algebras 
$(\cH^p,\omega^p,\ti{S}^p)$ and $(\cH^p,\omega^p,\ti{S'}^p)$ 
give the same minimal cyclic $A_\infty$-algebra consisting of the 
on-shell {\em open string} scattering amplitudes. 
This implies that two decomposed cyclic $A_\infty$-algebras 
$(\cH^p,\omega^p,\ti{S}^p)\oplus (\cH^t\oplus\cH^p,Q)$ 
and $(\cH^p,\omega^p,\ti{S'}^p)\oplus (\cH^t\oplus\cH^p,Q)$ 
are cyclic $A_\infty$-isomorphic to each other. 
Since the composition of any two cyclic $A_\infty$-isomorphisms 
forms a cyclic $A_\infty$-isomorphism, 
one can conclude that $(\cH,\omega, S)$ and $(\cH,\omega,S')$
are cyclic $A_\infty$-isomorphic to each other. 
\qed
\end{pf}

In physical terms, this theorem means 
any two classical string field theories on a conformal background 
are transformed to each other by a field redefinition (preserving 
classical BV-structures). 
This may also be thought of as a converse statement of the equivalence 
theorem (`S-matrices are invariant under field
redefinitions') \cite{KOrS} 
for field theories equipped with classical BV-structures. 
In contrast two cyclic $A_\infty$-algebras are not connected by 
a field definition if their minimal cyclic 
$A_\infty$-algebras are different from each other. 
This fact is also clear and follows from Theorem \ref{thm:cycMandC}. 
Usually a string field theory is constructed by decomposing Riemann 
surfaces as in subsection \ref{ssec:osft}. 
However, the actions which are transformed to each other 
by field redefinitions with preserving BV-symplectic forms 
are regarded to be equivalent \cite{SZ1,GZ}. 
Thus, Theorem \ref{thm:main} implies also the sufficiency of the axiom of 
string field theory in Definition \ref{defn:axiom}. 
\begin{rem}
Since there exists a cyclic $A_\infty$-isomorphism between 
any two classical open string field theories 
on a fixed conformal background, 
one may expect that it transforms all the equations of motions in one side 
to those in another side. 
However, the problem depends on convergence problem of the
transformation, that depends on the models we consider, 
and the answer may be `no' in general. 
On the other hand, 
in the context of deformation theory, the solutions for equations of 
motions, \ie the Maurer-Cartan equations, 
are assumed to be of the form 
$\Phi=\hbar\ti\Phi^p+\cO(\hbar^2)$ for $\hbar$ a `small' 
deformation parameter. The argument in subsection \ref{ssec:minimal} 
is just the case where the small parameter is thought to be included in 
$\ti\Phi^p$. Such solutions are expressed as $\Phi=\ti\cF^p_*(\ti\Phi^p)$ 
where $\ti\Phi^p$ is a solution for Maurer-Cartan equation 
$\ti\m^p_*(e^{\ti\Phi^p})=0$. Namely, they are the solutions which 
corresponds to the continuous deformations from the origin $\Phi=0$. 
They are regarded as the marginal deformation 
in the terminology of conformal field theory. 
The solutions satisfy not only the equation of motion but $S(\Phi)=0$. 
Now the problem of the convergence still remains, though 
we expect in string field theories 
they are valid at least in some neighborhood of the origin. 
However, if we treat $\hbar$ as a formal deformation parameter, 
then all the solutions are valid. 
Namely, all classical open string field theories on a 
fixed conformal background have isomorphic moduli spaces of 
`formal marginal deformations'. It follows from Theorem \ref{thm:main} 
and Theorem \ref{thm:quasiisomoduli} in subsection \ref{ssec:homgauge}. 
 \label{rem:marginal}
\end{rem}
The fact that the 
collection of open string 
tree scattering amplitudes 
possesses a minimal cyclic $A_\infty$-structure 
has essentially appeared in some literatures. 
It is described in \cite{WZ} that 
the $S^2$ tree amplitudes for closed strings has a 
$L_\infty$-structure, where the external states are restricted to 
physical states and therefore it has vanishing $Q$. This implies that 
the tree level closed string free energy satisfies the classical 
BV-master equation. The result is extended to 
quantum closed string, and it is shown that 
the free energy which consists of the closed string loop amplitudes 
satisfies the quantum BV-master equation \cite{V}.
These structures are derived in the context of 2D-string theory, 
\ie the dimension of the target space is two. 
However they are in fact the general structures of the string world
sheet, 
and reinterpreted in \cite{Z1} from a world sheet viewpoint explained below. 
The open string version of this $L_\infty$-structure is nothing but 
the $A_\infty$-structure $\ti\m^p$ above. 

{}From string world sheet point of view, 
a string field theory action is usually constructed by 
decomposing $\{\cM_n\}_{n\ge 3}$, moduli spaces of disk with $n$ punctures 
on the boundary, into cell (see subsection \ref{ssec:osft}). 
The vertex map $\V:=\omega\circ(\1\otimes m_{n-1}):\cH^{\otimes n}\to\C$ 
is determined by a certain integral over the cell
$\cM^0_n\subset\cM_n$. 
$\{\cM^0_k\}_{k\ge 3}$ satisfy some consistency condition and 
then forms a topological operad. 
The string field theory constructed at the limit 
$\{\cM^0_k\}_{k\ge 3}\to \{\cM_k\}_{k\ge 3}$ 
then turns out to be the on-shell open string 
scattering amplitudes. 
This also implies that the on-shell open string 
scattering amplitudes 
form a minimal cyclic $A_\infty$-algebra. 
This interpretation is just the one given in \cite{Z1} for closed
string case. 
Note that, for a $\{\cM^0_k\}_{k\ge 3}$ constructed
consistently, 
one can in fact take a continuous deformation from 
$\{\cM^0_k\}_{k\ge 3}$ to $\{\cM_k\}_{k\ge 3}$ 
with preserving the cyclicity. 
This implies that all $\{\cM^0_k\}_{k\ge 3}$ that are
constructed consistently are homotopic to each other as topological
operads \cite{Ka2}.

The continuous deformation from 
$\{\cM^0_k\}_{k\ge 3}$ to $\{\cM_k\}_{k\ge 3}$ is related to 
renormalization group flow \cite{P} of string world sheet theory 
in the sense of \cite{BM,HLP}, etc. 
Note that, for a given conformal background, 
the field $\phi^i$ is the coupling constant for open string state
$\eb_i$ and also thought of as `source' for $\eb_i$ 
from world sheet theory viewpoints. 
Let $(\cH,\omega,\ti{S})$ be the decomposed cyclic $A_\infty$-algebra 
whose minimal part is $(\cH^p,\omega^p,\ti{S}^p)$. 
It is an algebra over operad $\{\cM_k\}_{k\ge 3}$ and 
is essentially the generating function (or free energy) of tree open
string world sheet theory. 
The action $S(\Phi)$ over $\{\cM^0_k\}_{k\ge 3}$ then defines some
kind of deformed free energy of world sheet theory. 
The renormalization group flow is the flow (solution of differential 
equations) of the coupling constants such that 
the free energy $S(\Phi)$ over $\{\cM^0_k\}_{k\ge 3}$ is preserved under 
the continuous deformation. 
Now, by Theorem \ref{thm:main}, there 
exists continuous field redefinition that preserves the free energy 
$S(\Phi)$ under the continuous deformation of $\{\cM^0_k\}_{k\ge 3}$. 
Namely, the field redefinition of string field theory is just the 
renormalization group flow of string world sheet. 
On the other hand, string field theory is a field theory on the target 
space where strings are mapped. 
One can also consider the renormalization group as a target space
field theory \cite{BdA}. 
The correspondence between the world sheet renormalization group and 
the target space one is then discussed in \cite{BR,N}. 
Namely, the world sheet renormalization group, which is target space 
field redefinition, can also be regarded as the target space
renormalization group, though our theory is classical as a target
space field theory. 

Though the (abstract) way of the construction is given, 
physically string field theory has had mainly two general puzzles 
that should be resolved. 
One is the relation between different string field theory
actions on the same conformal background. The another one is then 
the background independence of the action; namely, if 
a string field theory action does give an nonperturbative definition 
of string theory, 
any classical solution should correspond to a vacuum (conformal
background) of the string theory and the action expanded around the 
classical solution should describe a string field theory on the 
corresponding conformal background. 
Both problems have ever resolved only at infinitesimal level. 
In \cite{HZ} (for quantum closed string field theory) 
it is shown that on a fixed conformal background 
any infinitesimal variation of the way of decomposition of 
moduli space is absorbed by an infinitesimal field redefinition. 
The background independence for (infinitesimal) marginal deformation 
is discussed and shown in \cite{Sbg1,Sbg2,Sbg3,RSZ,SZ1,SZ2,SZ3,KZ}. 
Theorem \ref{thm:main} then gives the answer for the first problem 
for classical string field theory. 
We believe that similar result holds for quantum case 
\footnote{Some (homotopy) algebraic structures for quantum closed 
SFT are discussed in \cite{SZ2,SZ3,Markl}. }
and 
it could be the starting point for the problem of the background 
independence.

%section7

 \section{Homotopy equivalence, gauge equivalence and moduli spaces} 
\label{sec:homotopy}

In this section we shall discuss 
homotopy theoretical aspects of $A_\infty$-algebras. 
In subsection \ref{ssec:hom} we shall define homotopy 
between an $A_\infty$-morphism, and show that two $A_\infty$-algebras 
are homotopy equivalent 
to each other iff there exists an
$A_\infty$-quasi-isomorphism between them (Theorem \ref{thm:inversehom}). 
Similar results are obtained in the framework of twisting cochains 
\cite{kadei3,kadei4}, in terms of operads \cite{Ma1}, 
in a purely algebraic way \cite{FOOO}, more recently 
in \cite{Le-Ha} 
in terms of closed model categories \cite{Le-Ha}, and so on. 
In subsection \ref{ssec:homgauge} we shall define the notion of 
gauge equivalence, and define the moduli space of an $A_\infty$-algebra 
as the solution space of the Maurer-Cartan equation 
for the $A_\infty$-algebra modulo gauge equivalence. 
The homotopy invariance of the moduli spaces are then shown in 
Theorem \ref{thm:quasiisomoduli}. 
A characteristic of our approach is to apply the decomposition theorem 
(Theorem \ref{thm:MandC}), 
which enables us to show both Theorem \ref{thm:inversehom} 
and Theorem \ref{thm:quasiisomoduli}. 
See also \cite{KaTe}.

 \subsection{Homotopy equivalence}
\label{ssec:hom}

\begin{defn}[Homotopy of $A_\infty$-morphisms]
Given two $A_\infty$-algebras $(\cH,\m)$ and $(\cH',\m')$ 
and two (weak) $A_\infty$-morphisms $\cF, \cG : (\cH,\m)\to (\cH',\m')$, 
we say that $\cG$ is 
{\it homotopic} to $\cF$ 
when there exists degree minus one operator 
$H : C(\cH)\to C(\cH')$ such that 
\begin{equation}
 \cG-\cF=\m' H + H\m\ .
 \label{morphomotopy}
\end{equation}
We call $H$ the homotopy operator between $A_\infty$-morphisms 
$\cG$ and $\cF$. 
 \label{defn:morphomotopy}
\end{defn}
Note that if we forget the cohomomorphism structures of $\cF$, $G$ and 
regard $(C(\cH), \m)$, $(C(\cH'),\m')$ as complexes of groups such as 
deRham complex, then the definition is the homotopy operator $H$ is 
just the usual one. 
\begin{rem}[A condition for $H$]
Note that $\cF$ and $\cG$ are cohomomorphisms. 
Therefore $H$ has a condition defined by $(\cF\otimes\cF)\tri=\tri\cF$, 
$(\cG\otimes\cG)\tri=\tri\cG$ and eq.(\ref{morphomotopy}). 
One of the most simplest solutions for the condition is given by 
\cite{GuMu,smirnov}
\begin{equation}
 \tri H=(\cF\otimes H+H\otimes\cG)\tri\ .
 \label{simpleH}
\end{equation}
Under this condition, $H$ is determined if the collection 
of degree minus one map 
$h_n :\cH^{\otimes n}\to\cH'$ for $n\ge 0$ is given as follows: 
\begin{equation}
 H=\sum_{n\ge 0}\cF\otimes h_n\otimes\cG\ .
 \label{homconstruct}
\end{equation}
One can see that if the image of $H$ in $C(\cH')$ is restricted to $\cH'$, 
$H$ just reduces to $\sum_n h_n$. 
General solutions are then written of the form 
\begin{equation*}
 \tri H=\l((\cF\otimes H+H\otimes\cG)+A\r)\tri
\end{equation*}
for a degree minus one element $A\in C(\cH)\otimes C(\cH)$. 
In this paper we treat all such $H$ as homotopies 
without fixing a condition for $\tri H$. 
 \label{rem:H}
\end{rem}
\begin{lem}
The homotopy in Definition \ref{defn:morphomotopy} 
actually defines an equivalence relation. 
 \label{lem:homotopy}
\end{lem}
\begin{pf}
It is easy to confirm the fact. 
\end{pf}
\begin{lem}
Let $(\cH_{dc},\m_{dc})$ be a direct sum of a minimal $A_\infty$-algebra 
and a linear contractible $A_\infty$-algebra, 
and $(\cH_{dc}^p,\m_{dc}^p)$ the corresponding minimal $A_\infty$-algebra. 
We have $A_\infty$-quasi-isomorphisms 
$\pi: (\cH_{dc},\m_{dc})\to (\cH_{dc}^p,\m_{dc}^p)$ and 
$\iota: (\cH_{dc}^p,\m_{dc}^p)\to(\cH_{dc},\m_{dc})$ 
and the composition 
$P=\iota\circ\pi: (\cH_{dc},\m_{dc})\to (\cH_{dc},\m_{dc})$ is also an 
$A_\infty$-quasi-isomorphism. 
In this situation, 
$A_\infty$-quasi-isomorphisms $\Id:(\cH_{dc},\m_{dc})\to
(\cH_{dc},\m_{dc})$ and 
$P: (\cH_{dc},\m_{dc})\to (\cH_{dc},\m_{dc})$ 
are homotopic to each other. 
 \label{lem:homIdP}
\end{lem}
\begin{pf}
The homotopy $H$ satisfying 
$\Id-P=\m H+H\m$ 
is obtained explicitly by setting $\cF=\Id$, 
$\cG=P$, $h_1=Q^+$ and $h_n=0$ for $n\ne 1$ in the condition (\ref{simpleH}),
that is, 
\begin{equation*}
 H=\Id\otimes Q^+\otimes P\ . 
\end{equation*}
Note that here we denote $\Id$ and 
$P$ as cohomomorphisms on 
$C(\cH_{dc})$. 
\qed
\end{pf}

The following theorem is known, which follows from the results in 
\cite{kadei3,kadei4} in the framework of {\it twisting cochains}, 
where homotopy is defined as in Definition \ref{defn:morphomotopy} 
but with the additional condition (\ref{simpleH}). 
Now we can show this theorem by applying the decomposition theorem 
(Theorem \ref{thm:MandC}). 
\begin{thm}
For two $A_\infty$-algebras $(\cH,\m)$ and $(\cH',\m')$, 
suppose that there exists an $A_\infty$-quasi-isomorphism 
$\cF : (\cH,\m) \to (\cH',\m')$. 
Then there exists an inverse $A_\infty$-quasi-isomorphism 
$\cF^{-1} : (\cH',\m')\to (\cH,\m)$ such that 
$\cF^{-1}\circ\cF$ is homotopic to the identity $\Id$ on $(\cH,\m)$ 
and $\cF\circ\cF^{-1}$ is homotopic to the identity $\Id$ on $(\cH',\m')$. 
 \label{thm:inversehom}
\end{thm}
\begin{pf}
The inverse $A_\infty$-quasi-isomorphism is just given by 
Theorem \ref{thm:inverse} combined with Corollary \ref{cor:inverse}. 
It is of the form 
\begin{equation*}
 \cF^{-1}
=\cF_{dc}\circ\iota\circ (\cF^p)^{-1}\circ\pi'\circ ({\cF'}_{dc})^{-1}\ ,
\end{equation*}
where we have the following commutative diagram 
\begin{equation*}
\xymatrix{
 \ \ (\cH,\m)\ \  \ar[r]^{\cF} \ar@<0.5ex>[d]^{(\cF_{dc})^{-1}} & 
 \ \ (\cH',\m')\ \ \, \ar@<0.5ex>[d]^{(\cF_{dc}')^{-1}} \\
 \ \ (\cH_{dc},\m_{dc})\ \ \ar@<0.5ex>[u]^{\cF_{dc}} 
 \ar@<0.5ex>[d]^{\pi} & 
 \ \ (\cH'_{dc},\m'_{dc})\ \ \ar@<0.5ex>[u]^{{\cF'}_{dc}} 
 \ar@<0.5ex>[d]^{\pi'}\ \ \,  \\
 \ \ (\cH_{dc}^p,\m_{dc}^p)\ \ \ar@<0.5ex>[u]^{\iota} 
 \ar@<0.5ex>[r]^{\cF^p} & 
 \ \ ({\cH'}_{dc}^p,{\m'}_{dc}^p) \ar@<0.5ex>[u]^{\iota'} 
 \ar@<0.5ex>[l]^{(\cF^p)^{-1}}\ \ .
}
\end{equation*}
Since one already knows the existence of $\cF^{-1}$, 
one can see that 
\begin{equation*}
 (\cF_{dc})^{-1}\circ(\cF^{-1}\circ\cF)\circ\cF_{dc}
 = P\ .
\end{equation*}
From Lemma \ref{lem:homIdP}, there exists a homotopy operator $H_o$ 
such as 
\begin{equation*}
 \Id- P=\m_{dc} H_o +H_o\m_{dc}
\end{equation*}
on $(\cH_{dc},\m_{dc})$. 
Acting $\cF_{dc}$ from left and $(\cF_{dc})^{-1}$ from right in both sides 
of the equation above, one can obtain 
\begin{equation*}
 \Id-(\cF)^{-1}\circ\cF=\m H+H\m
\end{equation*}
where $H:=\cF_{dc}\circ H_o\circ(\cF_{dc})^{-1}$ and 
we used the fact that 
$\m=\cF_{dc}\circ\m_{dc}\circ(\cF_{dc})^{-1}$. 
Thus it is shown that $\Id$ and $\cF^{-1}\circ\cF$ are 
homotopic to each other with homotopy operator $H$. 
The fact that $\cF\circ\cF^{-1}$ is homotopic to $\Id$ can also be 
shown in a similar way. 
\qed
\end{pf}

Note that an $A_\infty$-morphism $(\cH,\m)\to (\cH,\m)$ that is not an 
$A_\infty$-quasi-isomorphism cannot be homotopic to the identity. 
However, the converse is not true. 
Namely, an $A_\infty$-quasi-isomorphism 
is not necessarily homotopic to the identity. 
For instance, in general there exist $A_\infty$-isomorphisms 
associated to discrete $A_\infty$-automorphisms on $(\cH,\m)$.

Given two $A_\infty$-algebras $(\cH,\m)$, $(\cH',\m')$, 
an $A_\infty$-morphism $\cF: (\cH,\m)\to (\cH',\m')$ 
is called 
a homotopy equivalence between the two $A_\infty$-algebras iff 
there exists an $A_\infty$-morphism 
$\cF^{-1}: (\cH',\m')\to (\cH,\m)$ such that 
$\cF^{-1}\circ\cF$ and $\cF\circ\cF^{-1}$ are homotopic to the
identity (see \cite{Fukaya2}). 
Then Theorem \ref{thm:inversehom} implies that the notion of 
the homotopy equivalence of $A_\infty$-algebras 
is equivalent to the existence of $A_\infty$-quasi-isomorphisms 
between the $A_\infty$-algebras.

We end this subsection with some byproducts from 
the decomposition theorem (Theorem \ref{thm:MandC}) 
and the notion of homotopy equivalence of $A_\infty$-algebras. 
There is the notion of homotopy invariant algebraic structures 
by Boardman and Vogt \cite{BoVo} and the notion is translated into 
homotopy algebraic language in \cite{Ma1}(see also \cite{MSS}). 
Namely, homotopy algebras should have the following three properties 
(Corollary \ref{cor:move1}, Corollary \ref{cor:move2},
Corollary \ref{cor:move3}) 
so that they actually define homotopy invariant algebraic structures. 
They can be shown due to our arguments until now. 
\begin{cor}
For an $A_\infty$-algebra $(\cH,\m)$, a chain complex $(\cH',m_1')$ 
and a chain complex $f_1 :(\cH, m_1)\to(\cH',m_1')$, 
there exists an $A_\infty$-structure on $\cH'$ and 
an $A_\infty$-morphism $\cF$ whose leading term is $f_1$. 
 \label{cor:move1}
\end{cor}
\begin{cor}
Suppose two $A_\infty$-algebras $(\cH,\m)$, $(\cH',\m')$ and 
an $A_\infty$-morphism $\cF: (\cH,\m)\to (\cH',\m')$ are given. 
Moreover suppose there exists a chain map $g_1:(\cH,m_1)\to (\cH',m_1')$ 
that is chain homotopic to $f_1$. 
Then there exists an $A_\infty$-morphism $\cG:(\cH,\m)\to (\cH',\m')$ 
whose leading term is $g_1$. 
 \label{cor:move2}
\end{cor}
\begin{cor}
Given an $A_\infty$-quasi-isomorphism $\cF:(\cH,\m)\to (\cH',\m')$ 
and $g_1: (\cH',m_1')\to (\cH,m_1)$, 
a chain homotopy inverse of $f_1$, 
there exists an $A_\infty$-quasi-isomorphism 
$\cG: (\cH',\m')\to (\cH,\m)$. 
 \label{cor:move3}
\end{cor}
\begin{pf}
Corollary \ref{cor:move1} can be shown due to the decomposition theorem. 
For an $A_\infty$-algebra $(\cH,\m)$ one can consider a 
decomposed $A_\infty$-algebra $(\cH_{dc},\m_{dc})$ and 
it is clear that $f_1$ is naturally extended to 
an $A_\infty$-quasi-isomorphism and it induces a decomposed 
$A_\infty$-structure on $\cH'$. 
$\cF$ is then obtained by the composition of the 
$A_\infty$-quasi-isomorphism with 
$(\cF_{dc})^{-1}:(\cH,\m)\to (\cH_{dc},\m_{dc})$. 

Corollary \ref{cor:move2} follows from the notion of homotopy 
of $A_\infty$-morphisms. 
That is, $\cG: (\cH,\m)\to (\cH',\m')$ is obtained by 
\begin{equation}
 \cG=\cF+\m'H+H\m
 \label{move2}
\end{equation} 
if a homotopy operator $H$ exists. 
Here the problem is that whether there exists 
the $H$ such that the restriction of eq.(\ref{move2}) 
to the term $\cH\to\cH'$ 
is $g_1=f_1+m_1'h_1+h_1 m_1$ for $h_1$ the chain homotopy. 
Such $H$ does exist. One such is constructed by 
eq.(\ref{homconstruct}) with $h_0=h_2=h_3=\cdots=0$. 
 
Corollary \ref{cor:move3} is almost the same as Theorem \ref{thm:inversehom} 
and also follows from the decomposition theorem. 
 \qed
\end{pf}

 \subsection{Gauge equivalence and the moduli space of an $A_\infty$-algebra}
\label{ssec:homgauge}

Roughly speaking, the moduli space of an $A_\infty$-algebra is  
degree zero solution space of the Maurer-Cartan equation over 
gauge equivalence. 
Therefore we should define the gauge equivalence. 
Two elements in $\cH$ are called gauge equivalent 
if they are transformed to each other by a gauge transformation. 
Then the gauge transformation is usually defined 
as integrals of some infinitesimal transformations along paths in $\cH$. 
There exist mainly two candidates of the infinitesimal
transformation. 
In terms of formal noncommutative supermanifolds, they are the
followings. 
\begin{itemize}
 \item[(a)] The infinitesimal gauge transformation is defined by 
\begin{equation*}
 \delta_\alpha=\flpartial{\phi^i}\l(c^i(\phi)\flpartial{\phi^j}\alpha^j\r)
\end{equation*}
where $\alpha^j\in\C$ for $\deg(\eb_j)=0$ and 
$\alpha^j=0$ for $\deg(\eb_j)\ne 0$. 
 \item[(b)] The infinitesimal gauge transformation is defined by 
\begin{equation*}
 \delta_\alpha=[\delta, \alpha(\phi)]\ ,\qquad 
 \alpha(\phi):=\flpartial{\phi^i}\alpha^i(\phi)
=\flpartial{\phi^i}\sum_{k\ge 0}\alpha^i_{j_1\cdots j_k}
\phi^{j_k}\cdots\phi^{j_1}\ ,\quad \alpha^i_{j_1\cdots j_k}\in\C
\end{equation*}
where $\delta$ is the $A_\infty$-odd vector field and 
$\alpha(\phi)$ is a degree minus one vector field (see \cite{Ko1}). 
\end{itemize}
The infinitesimal transformation (a) is just the gauge transformation
in cyclic field theory (\ref{gaugedual}) explained 
in subsection \ref{ssec:BV}. 
More precisely, restricting the graded gauge parameter $\alpha^j$, 
$\deg(\eb_j)\ne 0$, to zero leads to the definition (a) above. 
Note that the gauge transformation (\ref{gaugedual}) 
just reduces to transformation (a) in the degree-zero subvector
space of $\cH$. 
The infinitesimal transformation (a) forms a Lie algebra only 
{\it on-shell}, where on-shell means the 
solution space of the Maurer-Cartan equation. 
Namely, the Lie bracket closes modulo the Maurer-Cartan equation. 
This implies the infinitesimal transformations are integrable on-shell 
and so the definition (a) is sufficient to define the moduli space of 
$A_\infty$-algebras. 
The integral is given by the usual exponential map (see below).

On the other hand, infinitesimal transformation (b) forms 
a Lie algebra on the whole space $\cH$. This fact can also be confirmed by 
a direct calculation. Note that the transformation (b) reduces to 
transformation (a) if $\alpha^i_{j_1\cdots j_k}=0$ for $k\ge 1$. 
One can see that extending $\alpha$ to $\phi$-dependent one leads to 
a Lie algebra which closes even off-shell. 

A natural extension of gauge transformation in differential graded Lie 
algebra case \cite{GoMi, SchSta} 
to $A_\infty$-algebras leads to the choice (a) 
(see \cite{Fukaya2}). 
Thus, (a) is also natural and in fact sufficient for our purpose. 
However, in this subsection we shall use (b) 
as the definition of the gauge transformation. 
By employing the arguments in section \ref{sec:dual}, 
we translate the above into coalgebra language and 
define precisely the gauge transformation. 
\begin{defn}[Gauge transformation]
Given an $A_\infty$-algebra $(\cH,\m)$, 
let us consider a piecewise smooth path of weak $A_\infty$-automorphisms 
$U_{\alpha[0,t]}: (\cH,\m)\to (\cH,\m)$, $0\le t\le 1$ defined as 
a coalgebra homomorphism $U_{\alpha[0,t]}:C(\cH)\to C(\cH)$ 
satisfying the following differential equation
\begin{equation*}
 \fd{t}U_{\alpha[0,t]}
 =\l([\m,\alpha(t)]U_{\alpha[0,t]}\r)\ ,\qquad U_{\alpha[0,0]}=\Id\ .
\end{equation*}
Here 
$\alpha(t) : C(\cH)\to C(\cH)$ is a degree minus one coderivation 
which consists of 
$\alpha(t)=\alpha_0(t)+\alpha_1(t)+\alpha_2(t)+\cdots$, 
where each coderivation $\alpha_k(t)$ is defined by the lift of a 
multilinear map $\cH^{\otimes k}\to\cH$ 
like as $\m=\m_0+\m_1+\m_2+\cdots$, 
and then 
$[\m,\alpha(t)]:=\m\alpha(t)+\alpha(t)\m: C(\cH)\to C(\cH)$. 
We call a weak $A_\infty$-automorphism 
$U_{\alpha[0,1]}:(\cH,\m)\to (\cH,\m)$ a {\it gauge transformation}. 
 \label{defn:gauge}
\end{defn}
The integral along the path $[0,s]$ is represented explicitly 
by an iterated integral, 
\begin{equation}
 U_{\alpha[0,s]}:=\cP e^{\int_0^s dt [\m,\alpha(t)]}=
\1+\int_0^s dt[\m,\alpha(t)]+\int_{s>t>t'>0}dtdt'
[\m,\alpha(t)][\m,\alpha(t')]+\cdots\ ,
 \label{gaugeauto}
\end{equation}
where $\cP$ is the one which is what is called the path-ordering. 
Note that if $\alpha(t)$ is constant with respect to $t$ 
the equation above reduces to the ordinary exponential map 
(in terms of coalgebra side). 
\begin{equation*}
 U_{\alpha[0,1]}=\cP e^{\int_0^1 dt [\m,\alpha]}=e^{[\m,\alpha]}:=
\1+[\m,\alpha]+\ov{2!}([\m,\alpha])^2+\ov{3!}([\m,\alpha])^3
+\cdots\ .
\end{equation*}
\begin{lem}
$U_{\alpha[0,s]}$ in eq.(\ref{gaugeauto}) is actually 
a weak $A_\infty$-automorphism on $(\cH,\m)$. 
Namely, the gauge transformations are weak $A_\infty$-automorphisms. 
 \label{lem:gauge}
\end{lem}
\begin{pf}
The fact that $U_{\alpha[0,s]}\m=\m U_{\alpha[0,s]}$ 
can be shown directly. 
In this calculation, one may use $[\m,[\m,\alpha]]=0$, 
which is just the infinitesimal description of the statement of 
this Lemma and follows from the Jacobi identity. 
 \qed
\end{pf}
\begin{lem}
The gauge transformation $U_{\alpha[0,1]}$ is homotopic to the identity 
on $A_\infty$-algebra $(\cH,\m)$. 
 \label{lem:gaugehom}
\end{lem}
\begin{pf}
First we have 
\begin{equation*}
 U_{\alpha[0,1]}-\Id=\int_0^1 ds \fd{s}U_{\alpha[0,s]}\ ,
\end{equation*}
where by definition 
$\fd{s}U_{\alpha[0,s]}=[\m,\alpha(s)]U_{\alpha[0,s]}$. 
Moreover Lemma \ref{lem:gauge} 
and the fact that $\cF$ is an $A_\infty$-morphism 
imply that the equation above is rewritten as 
\begin{equation*}
 U_{\alpha[0,1]}-\Id=\m H+H\m\ ,\qquad 
 H=\int_0^1 ds\, \alpha(s)U_{\alpha[0,s]}\ . 
\end{equation*}
The arguments above are similar to those of 
Weinstein-Darboux theorem \cite{We}, 
though the usage is different. 
\qed
\end{pf}

Note that the converse of Lemma \ref{lem:gaugehom} is not true. 
As seen in Lemma \ref{lem:homIdP} the retraction for the linear contractible 
direction is homotopic to the identity, which is not a gauge transformation. 
\begin{rem}
The composition of gauge transformations is a gauge transformation. 
In fact, for given two gauge transformation 
$U_{\alpha[0,1]}$ and $U_{\alpha'[0,1]}$, the composition is given by 
\begin{equation*}
 U_{\alpha'[0,1]}\circ U_{\alpha[0,1]}=U_{\alpha''[0,1]}\ , \qquad
 \alpha''(t)=
 \begin{cases}
  \alpha(2t) & 0\le t\le \half \\
  \alpha'(2t-1) & \half \le t\le 1
 \end{cases}\ .
\end{equation*}
Of course it is clear that one can take other $\alpha''$ that come 
from a reparametrization of $t$. 
 \label{rem:gauge}
\end{rem}
We shall below consider the moduli space of an $A_\infty$-algebra 
in $\cH^0$, the degree zero part of $\cH$, 
though our formalism may be appropriate to extend it to 
the whole degree (cf. subsection \ref{ssec:superfield}, 
subsection \ref{ssec:minimal} and \cite{BK}). 
In order to avoid the problem of convergence, 
one often considers in general the tensor product 
of maximal ideal of an Artin algebra 
with $\cH$ (see \cite{Fukaya2}). 
{}For simplicity, we shall introduce a formal parameter $\hbar$ and 
deal with $\cH\otimes \hbar\C[[\hbar]]$ from now on.  
\begin{defn}[Moduli space of an $A_\infty$-algebra]
For an $A_\infty$-algebra $(\cH,\m)$, 
let $\Phi$ be a degree zero element in 
$\hbar\cH[[\hbar]]:=\cH\otimes\C\hbar[[\hbar]]$, 
where $\hbar$ is a formal deformation parameter. 
We define the solution space of Maurer-Cartan equation for 
$A_\infty$-algebra $(\cH,\m)$ by 
\begin{equation*}
 \cMC(\cH,\m)=\{\Phi\in\hbar\cH[[\hbar]]\ |\ \m_*(e^\Phi)=0,\ \deg(\Phi)=0\}\ .
\end{equation*}
By definition, gauge transformations (Definition \ref{defn:gauge}) 
preserve the space $\cMC(\cH,\m)$. 
We call $\Phi_1, \Phi_2\in\cMC(\cH,\m)$ are {\it gauge equivalent} to
each other when they are transformed to each other by a gauge
transformation, that is, 
$e^{\Phi_2}=U_{\alpha[0,1]}e^{\Phi_1}$ or equivalently 
$\Phi_2=(U_{\alpha[0,1]})_*(\Phi_1)$ for some $\alpha(t)$, $0\le t\le 1$.  
Remark \ref{rem:gauge} guarantees that this 
actually defines an equivalence relation. 
The {\it moduli space} of $A_\infty$-algebra $(\cH,\m)$ is then defined by 
$\cMC(\cH,\m)$ modulo the gauge equivalence $\sim$, 
\begin{equation*}
 \cM(\cH,\m):=\cMC(\cH,\m)/\sim\ .
\end{equation*}
 \label{defn:moduli} 
\end{defn}
\begin{lem}
Given two $A_\infty$-algebras $(\cH,\m)$, $(\cH',\m')$ and 
an $A_\infty$-morphism $\cF: (\cH,\m)\to (\cH',\m')$, 
one can define the gauge transformation on the image of $\cF$ 
so that it is compatible with the gauge transformation in $(\cH,\m)$ 
when restricted to the solution spaces of the Maurer-Cartan equations. 
 \label{lem:compati}
\end{lem}
\begin{pf}
It is sufficient to show it infinitesimally. 
What should be shown is then that 
there exists a degree minus one coderivation $\alpha'$ on $\cH'$ 
such that 
\begin{equation}
 [\m',\alpha']\cF\l(e^\Phi\r)=\cF[\m,\alpha]\l(e^\Phi\r)
 \label{gaugecompati}
\end{equation}
holds if $\m\l(e^\Phi\r)=0$. 

Let us define 
a degree minus one coderivation $\alpha'$ on $\cH'$ by 
\begin{equation*}
 \alpha'(\Phi')
=\left.\cF (\alpha (e^\Phi))\right|_{(\cH')^{\otimes 1}}
 =f_1(\alpha(\Phi))+f_2(\alpha(\Phi),\Phi)+
f_2(\Phi,\alpha(\Phi))+\cdots\  
\end{equation*}
where $\alpha(\Phi)=\alpha_0+\alpha_1(\Phi)+\alpha_2(\Phi,\Phi)+\cdots$. 
This $\alpha'$ actually satisfies eq.(\ref{gaugecompati}). 
One can see that it follows from the condition of 
$A_\infty$-morphisms (Definition \ref{defn:amorp})
\begin{equation*}
 (\m'\cF-\cF\m)\alpha\l(e^{\Phi}\r)=0\ .
\end{equation*}
 \qed
\end{pf}

Note that this Lemma, combined with the fact that an
$A_\infty$-morphism preserves the solution spaces of the Maurer-Cartan
equations as shown in subsection \ref{ssec:MCeq}, 
implies that an $A_\infty$-morphism 
$\cF:(\cH,\m)\to (\cH',\m')$ actually induces a well-defined map from 
$\cM(\cH,\m)$ to $\cM(\cH',\m')$. 
In particular, if $\cF$ is an $A_\infty$-isomorphism, it induces an 
isomorphism on the moduli spaces. 
\begin{lem}
For a direct sum $(\cH_{dc},\m_{dc})$ of a minimal $A_\infty$-algebra 
$(\cH_{dc}^p,\m_{dc}^p)$ 
and a linear contractible $A_\infty$-algebra, 
the projection $P:\cH_{dc}\to\cH_{dc}$ induces the identity map $\Id$ 
on $\cM(\cH_{dc},\m_{dc})$. 
 \label{lem:quasiisoMC}
\end{lem}
\begin{pf}
For the Maurer-Cartan equation 
of the decomposed $A_\infty$-algebra $(\cH_{dc},\m_{dc})$, 
\begin{equation*}
 m_{dc,1}\Phi+\sum_{k\ge 2}m_{dc,k}(\Phi)=0\ ,
\end{equation*}
the first term and the second term should be zero independently, 
and the second term is nothing but the Maurer-Cartan equation 
of the minimal part $(\cH_{dc}^p,\m_{dc}^p)$. 
Here, for the solutions of $m_{dc,1}(\Phi)=0$, 
the $m_{dc,1}$-exact part is obviously gauge equivalent to zero; 
for instance, one may take $\alpha(t)=\alpha:\C\to\cH_{dc}^{-1}$ 
and then $m_{dc,1}(\alpha)\in\cH_{dc}^0$ is generated by a 
gauge transformation. 
Thus, any gauge equivalent class in
$\cMC(\cH_{dc},\m_{dc})$ can be represented by an element in $\cH^p$ 
and from which 
the statement of this Lemma follows. 
\qed
\end{pf}

Then we obtain the following, the $A_\infty$ version of the theorem 
in differential graded Lie algebra case \cite{GoMi,SchSta}. 
\begin{thm}
Given two $A_\infty$-algebras $(\cH,\m)$, $(\cH',\m')$ and 
suppose that there exists an $A_\infty$-quasi-isomorphism 
$\cF:(\cH,\m)\to(\cH',\m')$. 
Then the moduli spaces of these two $A_\infty$-algebras 
are isomorphic to each other: 
\begin{equation*}
 \cM(\cH,\m)\simeq \cM(\cH',\m')\ .
\end{equation*}
 \label{thm:quasiisomoduli}
\end{thm}
\begin{pf}
By Theorem \ref{thm:inversehom}, 
we have an inverse $A_\infty$-quasi-isomorphism 
$\cF^{-1}: (\cH',\m')\to (\cH,\m)$ 
such that $\cF^{-1}\circ\cF$ and $\cF\circ\cF^{-1}$ are homotopic to the 
identity $\Id$. 
Denote by $\cF_{\sim}:\cM(\cH,\m)\to\cM(\cH',\m')$ 
and 
$(\cF)^{-1}_{\sim}:\cM(\cH',\m')\to\cM(\cH,\m)$ 
the corresponding maps between the moduli spaces. 
One may show the following two equations 
\begin{equation}
(\cF)^{-1}_{\sim}\circ\cF_{\sim}=\Id\ ,\qquad 
\cF_{\sim}\circ (\cF)^{-1}_{\sim}=\Id\ . 
\label{cFsimId}
\end{equation}
On the other hand, in the proof of Theorem \ref{thm:inversehom} 
we know that 
\begin{equation*}
 (\cF_{dc})^{-1}\circ(\cF^{-1}\circ\cF)\circ\cF_{dc} =
 P\ ,\qquad 
 (\cF'_{dc})^{-1}\circ((\cF')^{-1}\circ\cF')\circ\cF'_{dc} =
 P'\  
\end{equation*}
hold, and furthermore Lemma \ref{lem:quasiisoMC} implies that 
$P:\cH_{dc}\to\cH_{dc}$ and 
$P':\cH'_{dc}\to\cH'_{dc}$ induce the identities 
$\Id$ on $\cM(\cH_{dc},\m_{dc})$ and 
$\cM(\cH'_{dc},\m'_{dc})$, respectively. 
This implies eq.(\ref{cFsimId}). 
 \qed
\end{pf}

The $L_\infty$-version of this theorem is 
used in the proof of Kontsevich's deformation quantization \cite{Ko1}.

We ends with leaving some comments. 
Though we defined homotopy between $A_\infty$-morphisms 
as in Definition \ref{defn:morphomotopy}, 
one can consider another definition of homotopy based on 
interpolating two $A_\infty$-morphisms with one parameter family of 
$A_\infty$-morphisms. Such one is used in \cite{Fukaya2, KaTe}. 
This homotopy is included in the homotopy in Definition
\ref{defn:morphomotopy} in the sense that 
integrating the homotopy along the one parameter path yields a
homotopy operator $H$ in Definition \ref{defn:morphomotopy}. 
However the converse is not true. 
Actually, the homotopy obtained in such a way satisfies a condition 
in Remark \ref{rem:H} but the condition is different from, for instance,
eq.(\ref{simpleH}). 
This another definition of homotopy is more compatible with the gauge
transformation, where the gauge transformation (a) is in fact more
natural than (b) used in this paper. 
Moreover, all of the theorems and corollaries in this section, 
for instance, remain true even if we replace Definition
\ref{defn:morphomotopy} to this another one \cite{KaTe}. 

All the arguments in this section hold in a similar way 
for cyclic $A_\infty$-algebras. 
Note that, in the cyclic case, the gauge transformation preserves the form 
of the action $S$ but 
does not preserve the form of the odd constant symplectic form. 
Moreover, all the arguments above hold true if 
an $A_\infty$-algebra is replaced to an $L_\infty$-algebra, too. 
In the cyclic case, these properties are closely related 
to the BV-formalism \cite{BV1,BV2,GPS,HT}. 
To exploring these relations should be very interesting 
(see for instance \cite{Sta2,Sta3}).

%\newpage

\begin{center}
\noindent{\large \textbf{Acknowledgments}}
\end{center}

I am very grateful to A.~Kato for helpful discussions, precise advice 
and encouragement. 
Also, I would like to thank T.~Asakawa, H.~Hata, Y.~Kazama, T.~Kugo, 
H.~Kunitomo and T.~Nakatsu 
for valuable discussions in physics, and 
M.~Akaho, K.~Fukaya, T.~Gocho, S.~Hosono, T.~Kimura, 
J.~Stasheff, Y.~Terashima and A.~Voronov 
for valuable discussions in mathematics. 
Especially, I am very grateful to S.~Hosono for careful reading 
this manuscript and 
J.~Stasheff for reading this manuscript, 
many valuable comments and encouragement. 
The author is supported by JSPS Research Fellowships for Young
Scientists.

%%%%%%%%%%%%%%%%%%%%%%%%%%%%%%%%%%%%%%%%%%%%%%%
%%%%%%%%%%%%%%%%%%%%%%%%%%%%%%%%%%%%%%%%%%%%%%%

\end{document}

%% file: A1.tex
%WinTpicVersion3.08
\unitlength 0.1in
\begin{picture}( 14.5669,  8.2185)( 20.8661,-15.7480)
% CIRCLE 2 0 3 0
% 4 3200 1200 2900 1200 2900 1200 2900 1200
% 
\special{pn 8}%
\special{ar 3150 1182 296 296  0.0000000 6.2831853}%
% STR 2 0 3 0
% 3 3200 1100 3200 1200 5 0
% $a$
\put(31.4961,-11.8110){\makebox(0,0){$a$}}%
% DOT 1 0 3 0
% 2 2900 1200 2900 1200
% 
\special{pn 13}%
\special{sh 1}%
\special{ar 2855 1182 10 10 0  6.28318530717959E+0000}%
\special{sh 1}%
\special{ar 2855 1182 10 10 0  6.28318530717959E+0000}%
% DOT 1 0 3 0
% 2 3050 940 3050 940
% 
\special{pn 13}%
\special{sh 1}%
\special{ar 3002 926 10 10 0  6.28318530717959E+0000}%
\special{sh 1}%
\special{ar 3002 926 10 10 0  6.28318530717959E+0000}%
% DOT 1 0 3 0
% 2 2940 1050 2940 1050
% 
\special{pn 13}%
\special{sh 1}%
\special{ar 2894 1034 10 10 0  6.28318530717959E+0000}%
\special{sh 1}%
\special{ar 2894 1034 10 10 0  6.28318530717959E+0000}%
% DOT 1 0 3 0
% 2 2940 1350 2940 1350
% 
\special{pn 13}%
\special{sh 1}%
\special{ar 2894 1329 10 10 0  6.28318530717959E+0000}%
\special{sh 1}%
\special{ar 2894 1329 10 10 0  6.28318530717959E+0000}%
% STR 2 0 3 0
% 3 2800 1100 2800 1200 5 0
% $\phi^{i_n}$
\put(27.5591,-11.8110){\makebox(0,0){$\phi^{i_n}$}}%
% STR 2 0 3 0
% 3 2840 890 2840 990 5 0
% $\phi^{i_{n-1}}$
\put(27.9528,-9.7441){\makebox(0,0){$\phi^{i_{n-1}}$}}%
% STR 2 0 3 0
% 3 3000 750 3000 850 5 0
% $\phi^{i_{n-2}}$
\put(29.5276,-8.3661){\makebox(0,0){$\phi^{i_{n-2}}$}}%
% STR 2 0 3 0
% 3 2840 1300 2840 1400 5 0
% $\phi^{i_1}$
\put(27.9528,-13.7795){\makebox(0,0){$\phi^{i_1}$}}%
% LINE 2 0 3 0
% 2 3000 1250 2600 1360
% 
\special{pn 8}%
\special{pa 2953 1231}%
\special{pa 2560 1339}%
\special{fp}%
% STR 2 0 3 0
% 3 2400 1300 2400 1400 5 0
% cut
\put(23.6220,-13.7795){\makebox(0,0){cut}}%
% CIRCLE 2 2 3 0
% 4 3200 1200 3200 800 2070 3160 3200 400
% 
\special{pn 8}%
\special{ar 3150 1182 394 394  4.7123890 4.7423890}%
\special{ar 3150 1182 394 394  4.8323890 4.8623890}%
\special{ar 3150 1182 394 394  4.9523890 4.9823890}%
\special{ar 3150 1182 394 394  5.0723890 5.1023890}%
\special{ar 3150 1182 394 394  5.1923890 5.2223890}%
\special{ar 3150 1182 394 394  5.3123890 5.3423890}%
\special{ar 3150 1182 394 394  5.4323890 5.4623890}%
\special{ar 3150 1182 394 394  5.5523890 5.5823890}%
\special{ar 3150 1182 394 394  5.6723890 5.7023890}%
\special{ar 3150 1182 394 394  5.7923890 5.8223890}%
\special{ar 3150 1182 394 394  5.9123890 5.9423890}%
\special{ar 3150 1182 394 394  6.0323890 6.0623890}%
\special{ar 3150 1182 394 394  6.1523890 6.1823890}%
\special{ar 3150 1182 394 394  6.2723890 6.3023890}%
\special{ar 3150 1182 394 394  6.3923890 6.4223890}%
\special{ar 3150 1182 394 394  6.5123890 6.5423890}%
\special{ar 3150 1182 394 394  6.6323890 6.6623890}%
\special{ar 3150 1182 394 394  6.7523890 6.7823890}%
\special{ar 3150 1182 394 394  6.8723890 6.9023890}%
\special{ar 3150 1182 394 394  6.9923890 7.0223890}%
\special{ar 3150 1182 394 394  7.1123890 7.1423890}%
\special{ar 3150 1182 394 394  7.2323890 7.2623890}%
\special{ar 3150 1182 394 394  7.3523890 7.3823890}%
\special{ar 3150 1182 394 394  7.4723890 7.5023890}%
\special{ar 3150 1182 394 394  7.5923890 7.6223890}%
\special{ar 3150 1182 394 394  7.7123890 7.7423890}%
\special{ar 3150 1182 394 394  7.8323890 7.8623890}%
\special{ar 3150 1182 394 394  7.9523890 7.9823890}%
\special{ar 3150 1182 394 394  8.0723890 8.1023890}%
\special{ar 3150 1182 394 394  8.1923890 8.2223890}%
\special{ar 3150 1182 394 394  8.3123890 8.3423890}%
\end{picture}%

%% file: AB.tex
%WinTpicVersion3.08
\unitlength 0.1in
\begin{picture}( 16.7323,  7.8740)( 21.6535,-15.7480)
% CIRCLE 2 0 3 0
% 4 2600 1200 2900 1200 2900 1200 2900 1200
% 
\special{pn 8}%
\special{ar 2560 1182 296 296  0.0000000 6.2831853}%
% STR 2 0 3 0
% 3 2600 1100 2600 1200 5 0
% $a$
\put(25.5906,-11.8110){\makebox(0,0){$a$}}%
% DOT 1 0 3 0
% 2 2900 1200 2900 1200
% 
\special{pn 13}%
\special{sh 1}%
\special{ar 2855 1182 10 10 0  6.28318530717959E+0000}%
\special{sh 1}%
\special{ar 2855 1182 10 10 0  6.28318530717959E+0000}%
% DOT 1 0 3 0
% 2 2750 940 2750 940
% 
\special{pn 13}%
\special{sh 1}%
\special{ar 2707 926 10 10 0  6.28318530717959E+0000}%
\special{sh 1}%
\special{ar 2707 926 10 10 0  6.28318530717959E+0000}%
% DOT 1 0 3 0
% 2 2860 1050 2860 1050
% 
\special{pn 13}%
\special{sh 1}%
\special{ar 2815 1034 10 10 0  6.28318530717959E+0000}%
\special{sh 1}%
\special{ar 2815 1034 10 10 0  6.28318530717959E+0000}%
% DOT 1 0 3 0
% 2 2860 1350 2860 1350
% 
\special{pn 13}%
\special{sh 1}%
\special{ar 2815 1329 10 10 0  6.28318530717959E+0000}%
\special{sh 1}%
\special{ar 2815 1329 10 10 0  6.28318530717959E+0000}%
% CIRCLE 2 2 3 0
% 4 2600 1200 2600 800 2600 400 3730 3160
% 
\special{pn 8}%
\special{ar 2560 1182 394 394  1.0478125 1.0778125}%
\special{ar 2560 1182 394 394  1.1678125 1.1978125}%
\special{ar 2560 1182 394 394  1.2878125 1.3178125}%
\special{ar 2560 1182 394 394  1.4078125 1.4378125}%
\special{ar 2560 1182 394 394  1.5278125 1.5578125}%
\special{ar 2560 1182 394 394  1.6478125 1.6778125}%
\special{ar 2560 1182 394 394  1.7678125 1.7978125}%
\special{ar 2560 1182 394 394  1.8878125 1.9178125}%
\special{ar 2560 1182 394 394  2.0078125 2.0378125}%
\special{ar 2560 1182 394 394  2.1278125 2.1578125}%
\special{ar 2560 1182 394 394  2.2478125 2.2778125}%
\special{ar 2560 1182 394 394  2.3678125 2.3978125}%
\special{ar 2560 1182 394 394  2.4878125 2.5178125}%
\special{ar 2560 1182 394 394  2.6078125 2.6378125}%
\special{ar 2560 1182 394 394  2.7278125 2.7578125}%
\special{ar 2560 1182 394 394  2.8478125 2.8778125}%
\special{ar 2560 1182 394 394  2.9678125 2.9978125}%
\special{ar 2560 1182 394 394  3.0878125 3.1178125}%
\special{ar 2560 1182 394 394  3.2078125 3.2378125}%
\special{ar 2560 1182 394 394  3.3278125 3.3578125}%
\special{ar 2560 1182 394 394  3.4478125 3.4778125}%
\special{ar 2560 1182 394 394  3.5678125 3.5978125}%
\special{ar 2560 1182 394 394  3.6878125 3.7178125}%
\special{ar 2560 1182 394 394  3.8078125 3.8378125}%
\special{ar 2560 1182 394 394  3.9278125 3.9578125}%
\special{ar 2560 1182 394 394  4.0478125 4.0778125}%
\special{ar 2560 1182 394 394  4.1678125 4.1978125}%
\special{ar 2560 1182 394 394  4.2878125 4.3178125}%
\special{ar 2560 1182 394 394  4.4078125 4.4378125}%
\special{ar 2560 1182 394 394  4.5278125 4.5578125}%
\special{ar 2560 1182 394 394  4.6478125 4.6778125}%
% CIRCLE 2 0 3 0
% 4 3500 1200 3200 1200 3200 1200 3200 1200
% 
\special{pn 8}%
\special{ar 3445 1182 296 296  0.0000000 6.2831853}%
% STR 2 0 3 0
% 3 3500 1100 3500 1200 5 0
% $b$
\put(34.4488,-11.8110){\makebox(0,0){$b$}}%
% DOT 1 0 3 0
% 2 3200 1200 3200 1200
% 
\special{pn 13}%
\special{sh 1}%
\special{ar 3150 1182 10 10 0  6.28318530717959E+0000}%
\special{sh 1}%
\special{ar 3150 1182 10 10 0  6.28318530717959E+0000}%
% DOT 1 0 3 0
% 2 3350 940 3350 940
% 
\special{pn 13}%
\special{sh 1}%
\special{ar 3298 926 10 10 0  6.28318530717959E+0000}%
\special{sh 1}%
\special{ar 3298 926 10 10 0  6.28318530717959E+0000}%
% DOT 1 0 3 0
% 2 3240 1050 3240 1050
% 
\special{pn 13}%
\special{sh 1}%
\special{ar 3189 1034 10 10 0  6.28318530717959E+0000}%
\special{sh 1}%
\special{ar 3189 1034 10 10 0  6.28318530717959E+0000}%
% DOT 1 0 3 0
% 2 3240 1350 3240 1350
% 
\special{pn 13}%
\special{sh 1}%
\special{ar 3189 1329 10 10 0  6.28318530717959E+0000}%
\special{sh 1}%
\special{ar 3189 1329 10 10 0  6.28318530717959E+0000}%
% CIRCLE 2 2 3 0
% 4 3500 1200 3500 800 2370 3160 3500 400
% 
\special{pn 8}%
\special{ar 3445 1182 394 394  4.7123890 4.7423890}%
\special{ar 3445 1182 394 394  4.8323890 4.8623890}%
\special{ar 3445 1182 394 394  4.9523890 4.9823890}%
\special{ar 3445 1182 394 394  5.0723890 5.1023890}%
\special{ar 3445 1182 394 394  5.1923890 5.2223890}%
\special{ar 3445 1182 394 394  5.3123890 5.3423890}%
\special{ar 3445 1182 394 394  5.4323890 5.4623890}%
\special{ar 3445 1182 394 394  5.5523890 5.5823890}%
\special{ar 3445 1182 394 394  5.6723890 5.7023890}%
\special{ar 3445 1182 394 394  5.7923890 5.8223890}%
\special{ar 3445 1182 394 394  5.9123890 5.9423890}%
\special{ar 3445 1182 394 394  6.0323890 6.0623890}%
\special{ar 3445 1182 394 394  6.1523890 6.1823890}%
\special{ar 3445 1182 394 394  6.2723890 6.3023890}%
\special{ar 3445 1182 394 394  6.3923890 6.4223890}%
\special{ar 3445 1182 394 394  6.5123890 6.5423890}%
\special{ar 3445 1182 394 394  6.6323890 6.6623890}%
\special{ar 3445 1182 394 394  6.7523890 6.7823890}%
\special{ar 3445 1182 394 394  6.8723890 6.9023890}%
\special{ar 3445 1182 394 394  6.9923890 7.0223890}%
\special{ar 3445 1182 394 394  7.1123890 7.1423890}%
\special{ar 3445 1182 394 394  7.2323890 7.2623890}%
\special{ar 3445 1182 394 394  7.3523890 7.3823890}%
\special{ar 3445 1182 394 394  7.4723890 7.5023890}%
\special{ar 3445 1182 394 394  7.5923890 7.6223890}%
\special{ar 3445 1182 394 394  7.7123890 7.7423890}%
\special{ar 3445 1182 394 394  7.8323890 7.8623890}%
\special{ar 3445 1182 394 394  7.9523890 7.9823890}%
\special{ar 3445 1182 394 394  8.0723890 8.1023890}%
\special{ar 3445 1182 394 394  8.1923890 8.2223890}%
\special{ar 3445 1182 394 394  8.3123890 8.3423890}%
% BOX 2 5 2 0
% 2 2800 1150 3300 1250
% 
\special{pn 8}%
\special{sh 0}%
\special{pa 2756 1132}%
\special{pa 3249 1132}%
\special{pa 3249 1231}%
\special{pa 2756 1231}%
\special{pa 2756 1132}%
\special{ip}%
% LINE 2 0 3 0
% 2 2890 1150 3200 1150
% 
\special{pn 8}%
\special{pa 2845 1132}%
\special{pa 3150 1132}%
\special{fp}%
% LINE 2 0 3 0
% 2 2900 1250 3200 1250
% 
\special{pn 8}%
\special{pa 2855 1231}%
\special{pa 3150 1231}%
\special{fp}%
% CIRCLE 2 0 3 0
% 4 1800 2000 3100 1640 3030 1430 2850 1190
% 
\special{pn 8}%
\special{ar 1772 1969 1328 1328  5.6261104 5.8492319}%
% STR 2 0 3 0
% 3 3060 1460 3060 1560 5 0
% cut
\put(30.1181,-15.3543){\makebox(0,0){cut}}%
\end{picture}%

%% file: AB2.tex
%WinTpicVersion3.08
\unitlength 0.1in
\begin{picture}( 16.7323,  7.8740)( 21.6535,-15.7480)
% CIRCLE 2 0 3 0
% 4 2600 1200 2900 1200 2900 1200 2900 1200
% 
\special{pn 8}%
\special{ar 2560 1182 296 296  0.0000000 6.2831853}%
% STR 2 0 3 0
% 3 2600 1100 2600 1200 5 0
% $a$
\put(25.5906,-11.8110){\makebox(0,0){$a$}}%
% DOT 1 0 3 0
% 2 2900 1200 2900 1200
% 
\special{pn 13}%
\special{sh 1}%
\special{ar 2855 1182 10 10 0  6.28318530717959E+0000}%
\special{sh 1}%
\special{ar 2855 1182 10 10 0  6.28318530717959E+0000}%
% DOT 1 0 3 0
% 2 2750 940 2750 940
% 
\special{pn 13}%
\special{sh 1}%
\special{ar 2707 926 10 10 0  6.28318530717959E+0000}%
\special{sh 1}%
\special{ar 2707 926 10 10 0  6.28318530717959E+0000}%
% DOT 1 0 3 0
% 2 2860 1050 2860 1050
% 
\special{pn 13}%
\special{sh 1}%
\special{ar 2815 1034 10 10 0  6.28318530717959E+0000}%
\special{sh 1}%
\special{ar 2815 1034 10 10 0  6.28318530717959E+0000}%
% DOT 1 0 3 0
% 2 2860 1350 2860 1350
% 
\special{pn 13}%
\special{sh 1}%
\special{ar 2815 1329 10 10 0  6.28318530717959E+0000}%
\special{sh 1}%
\special{ar 2815 1329 10 10 0  6.28318530717959E+0000}%
% CIRCLE 2 2 3 0
% 4 2600 1200 2600 800 2600 400 3730 3160
% 
\special{pn 8}%
\special{ar 2560 1182 394 394  1.0478125 1.0778125}%
\special{ar 2560 1182 394 394  1.1678125 1.1978125}%
\special{ar 2560 1182 394 394  1.2878125 1.3178125}%
\special{ar 2560 1182 394 394  1.4078125 1.4378125}%
\special{ar 2560 1182 394 394  1.5278125 1.5578125}%
\special{ar 2560 1182 394 394  1.6478125 1.6778125}%
\special{ar 2560 1182 394 394  1.7678125 1.7978125}%
\special{ar 2560 1182 394 394  1.8878125 1.9178125}%
\special{ar 2560 1182 394 394  2.0078125 2.0378125}%
\special{ar 2560 1182 394 394  2.1278125 2.1578125}%
\special{ar 2560 1182 394 394  2.2478125 2.2778125}%
\special{ar 2560 1182 394 394  2.3678125 2.3978125}%
\special{ar 2560 1182 394 394  2.4878125 2.5178125}%
\special{ar 2560 1182 394 394  2.6078125 2.6378125}%
\special{ar 2560 1182 394 394  2.7278125 2.7578125}%
\special{ar 2560 1182 394 394  2.8478125 2.8778125}%
\special{ar 2560 1182 394 394  2.9678125 2.9978125}%
\special{ar 2560 1182 394 394  3.0878125 3.1178125}%
\special{ar 2560 1182 394 394  3.2078125 3.2378125}%
\special{ar 2560 1182 394 394  3.3278125 3.3578125}%
\special{ar 2560 1182 394 394  3.4478125 3.4778125}%
\special{ar 2560 1182 394 394  3.5678125 3.5978125}%
\special{ar 2560 1182 394 394  3.6878125 3.7178125}%
\special{ar 2560 1182 394 394  3.8078125 3.8378125}%
\special{ar 2560 1182 394 394  3.9278125 3.9578125}%
\special{ar 2560 1182 394 394  4.0478125 4.0778125}%
\special{ar 2560 1182 394 394  4.1678125 4.1978125}%
\special{ar 2560 1182 394 394  4.2878125 4.3178125}%
\special{ar 2560 1182 394 394  4.4078125 4.4378125}%
\special{ar 2560 1182 394 394  4.5278125 4.5578125}%
\special{ar 2560 1182 394 394  4.6478125 4.6778125}%
% CIRCLE 2 0 3 0
% 4 3500 1200 3200 1200 3200 1200 3200 1200
% 
\special{pn 8}%
\special{ar 3445 1182 296 296  0.0000000 6.2831853}%
% STR 2 0 3 0
% 3 3500 1100 3500 1200 5 0
% $b$
\put(34.4488,-11.8110){\makebox(0,0){$b$}}%
% DOT 1 0 3 0
% 2 3200 1200 3200 1200
% 
\special{pn 13}%
\special{sh 1}%
\special{ar 3150 1182 10 10 0  6.28318530717959E+0000}%
\special{sh 1}%
\special{ar 3150 1182 10 10 0  6.28318530717959E+0000}%
% DOT 1 0 3 0
% 2 3350 940 3350 940
% 
\special{pn 13}%
\special{sh 1}%
\special{ar 3298 926 10 10 0  6.28318530717959E+0000}%
\special{sh 1}%
\special{ar 3298 926 10 10 0  6.28318530717959E+0000}%
% DOT 1 0 3 0
% 2 3240 1050 3240 1050
% 
\special{pn 13}%
\special{sh 1}%
\special{ar 3189 1034 10 10 0  6.28318530717959E+0000}%
\special{sh 1}%
\special{ar 3189 1034 10 10 0  6.28318530717959E+0000}%
% DOT 1 0 3 0
% 2 3240 1350 3240 1350
% 
\special{pn 13}%
\special{sh 1}%
\special{ar 3189 1329 10 10 0  6.28318530717959E+0000}%
\special{sh 1}%
\special{ar 3189 1329 10 10 0  6.28318530717959E+0000}%
% CIRCLE 2 2 3 0
% 4 3500 1200 3500 800 2370 3160 3500 400
% 
\special{pn 8}%
\special{ar 3445 1182 394 394  4.7123890 4.7423890}%
\special{ar 3445 1182 394 394  4.8323890 4.8623890}%
\special{ar 3445 1182 394 394  4.9523890 4.9823890}%
\special{ar 3445 1182 394 394  5.0723890 5.1023890}%
\special{ar 3445 1182 394 394  5.1923890 5.2223890}%
\special{ar 3445 1182 394 394  5.3123890 5.3423890}%
\special{ar 3445 1182 394 394  5.4323890 5.4623890}%
\special{ar 3445 1182 394 394  5.5523890 5.5823890}%
\special{ar 3445 1182 394 394  5.6723890 5.7023890}%
\special{ar 3445 1182 394 394  5.7923890 5.8223890}%
\special{ar 3445 1182 394 394  5.9123890 5.9423890}%
\special{ar 3445 1182 394 394  6.0323890 6.0623890}%
\special{ar 3445 1182 394 394  6.1523890 6.1823890}%
\special{ar 3445 1182 394 394  6.2723890 6.3023890}%
\special{ar 3445 1182 394 394  6.3923890 6.4223890}%
\special{ar 3445 1182 394 394  6.5123890 6.5423890}%
\special{ar 3445 1182 394 394  6.6323890 6.6623890}%
\special{ar 3445 1182 394 394  6.7523890 6.7823890}%
\special{ar 3445 1182 394 394  6.8723890 6.9023890}%
\special{ar 3445 1182 394 394  6.9923890 7.0223890}%
\special{ar 3445 1182 394 394  7.1123890 7.1423890}%
\special{ar 3445 1182 394 394  7.2323890 7.2623890}%
\special{ar 3445 1182 394 394  7.3523890 7.3823890}%
\special{ar 3445 1182 394 394  7.4723890 7.5023890}%
\special{ar 3445 1182 394 394  7.5923890 7.6223890}%
\special{ar 3445 1182 394 394  7.7123890 7.7423890}%
\special{ar 3445 1182 394 394  7.8323890 7.8623890}%
\special{ar 3445 1182 394 394  7.9523890 7.9823890}%
\special{ar 3445 1182 394 394  8.0723890 8.1023890}%
\special{ar 3445 1182 394 394  8.1923890 8.2223890}%
\special{ar 3445 1182 394 394  8.3123890 8.3423890}%
% BOX 2 5 2 0
% 2 2800 1150 3300 1250
% 
\special{pn 8}%
\special{sh 0}%
\special{pa 2756 1132}%
\special{pa 3249 1132}%
\special{pa 3249 1231}%
\special{pa 2756 1231}%
\special{pa 2756 1132}%
\special{ip}%
% LINE 2 0 3 0
% 2 2890 1150 3200 1150
% 
\special{pn 8}%
\special{pa 2845 1132}%
\special{pa 3150 1132}%
\special{fp}%
% LINE 2 0 3 0
% 2 2900 1250 3200 1250
% 
\special{pn 8}%
\special{pa 2855 1231}%
\special{pa 3150 1231}%
\special{fp}%
% STR 2 0 3 0
% 3 2920 1500 2920 1600 5 0
% cut
\put(28.7402,-15.7480){\makebox(0,0){cut}}%
% DOT 1 0 3 0
% 2 3050 1150 3050 1150
% 
\special{pn 13}%
\special{sh 1}%
\special{ar 3002 1132 10 10 0  6.28318530717959E+0000}%
\special{sh 1}%
\special{ar 3002 1132 10 10 0  6.28318530717959E+0000}%
% DOT 1 0 3 0
% 2 3050 1250 3050 1250
% 
\special{pn 13}%
\special{sh 1}%
\special{ar 3002 1231 10 10 0  6.28318530717959E+0000}%
\special{sh 1}%
\special{ar 3002 1231 10 10 0  6.28318530717959E+0000}%
% DOT 1 0 3 0
% 2 3020 1150 3020 1150
% 
\special{pn 13}%
\special{sh 1}%
\special{ar 2973 1132 10 10 0  6.28318530717959E+0000}%
\special{sh 1}%
\special{ar 2973 1132 10 10 0  6.28318530717959E+0000}%
% DOT 1 0 3 0
% 2 3080 1150 3080 1150
% 
\special{pn 13}%
\special{sh 1}%
\special{ar 3032 1132 10 10 0  6.28318530717959E+0000}%
\special{sh 1}%
\special{ar 3032 1132 10 10 0  6.28318530717959E+0000}%
% DOT 1 0 3 0
% 2 3080 1250 3080 1250
% 
\special{pn 13}%
\special{sh 1}%
\special{ar 3032 1231 10 10 0  6.28318530717959E+0000}%
\special{sh 1}%
\special{ar 3032 1231 10 10 0  6.28318530717959E+0000}%
% DOT 1 0 3 0
% 2 3020 1250 3020 1250
% 
\special{pn 13}%
\special{sh 1}%
\special{ar 2973 1231 10 10 0  6.28318530717959E+0000}%
\special{sh 1}%
\special{ar 2973 1231 10 10 0  6.28318530717959E+0000}%
% STR 2 0 3 0
% 3 3050 950 3050 1050 5 0
% $\phi^{I}$
\put(30.0197,-10.3346){\makebox(0,0){$\phi^{I}$}}%
% STR 2 0 3 0
% 3 3050 1250 3050 1350 5 0
% $\phi^{J}$
\put(30.0197,-13.2874){\makebox(0,0){$\phi^{J}$}}%
% CIRCLE 2 0 3 0
% 4 1800 1390 2920 1510 2920 1510 2910 1220
% 
\special{pn 8}%
\special{ar 1772 1369 1109 1109  6.1312130 6.2831853}%
\special{ar 1772 1369 1109 1109  0.0000000 0.1067357}%
\end{picture}%

%% file: Ader.tex
%WinTpicVersion3.08
\unitlength 0.1in
\begin{picture}( 10.8268,  7.8740)( 31.4961,-21.6535)
% CIRCLE 2 0 3 0
% 4 3900 1800 3600 1800 3600 1800 3600 1800
% 
\special{pn 8}%
\special{ar 3839 1772 296 296  0.0000000 6.2831853}%
% STR 2 0 3 0
% 3 3900 1700 3900 1800 5 0
% $A$
\put(38.3858,-17.7165){\makebox(0,0){$A$}}%
% DOT 1 0 3 0
% 2 3600 1800 3600 1800
% 
\special{pn 13}%
\special{sh 1}%
\special{ar 3544 1772 10 10 0  6.28318530717959E+0000}%
\special{sh 1}%
\special{ar 3544 1772 10 10 0  6.28318530717959E+0000}%
% DOT 1 0 3 0
% 2 3750 1540 3750 1540
% 
\special{pn 13}%
\special{sh 1}%
\special{ar 3691 1516 10 10 0  6.28318530717959E+0000}%
\special{sh 1}%
\special{ar 3691 1516 10 10 0  6.28318530717959E+0000}%
% DOT 1 0 3 0
% 2 3640 1650 3640 1650
% 
\special{pn 13}%
\special{sh 1}%
\special{ar 3583 1625 10 10 0  6.28318530717959E+0000}%
\special{sh 1}%
\special{ar 3583 1625 10 10 0  6.28318530717959E+0000}%
% DOT 1 0 3 0
% 2 3640 1950 3640 1950
% 
\special{pn 13}%
\special{sh 1}%
\special{ar 3583 1920 10 10 0  6.28318530717959E+0000}%
\special{sh 1}%
\special{ar 3583 1920 10 10 0  6.28318530717959E+0000}%
% CIRCLE 2 2 3 0
% 4 3900 1800 3900 1400 2770 3760 3900 1000
% 
\special{pn 8}%
\special{ar 3839 1772 394 394  4.7123890 4.7423890}%
\special{ar 3839 1772 394 394  4.8323890 4.8623890}%
\special{ar 3839 1772 394 394  4.9523890 4.9823890}%
\special{ar 3839 1772 394 394  5.0723890 5.1023890}%
\special{ar 3839 1772 394 394  5.1923890 5.2223890}%
\special{ar 3839 1772 394 394  5.3123890 5.3423890}%
\special{ar 3839 1772 394 394  5.4323890 5.4623890}%
\special{ar 3839 1772 394 394  5.5523890 5.5823890}%
\special{ar 3839 1772 394 394  5.6723890 5.7023890}%
\special{ar 3839 1772 394 394  5.7923890 5.8223890}%
\special{ar 3839 1772 394 394  5.9123890 5.9423890}%
\special{ar 3839 1772 394 394  6.0323890 6.0623890}%
\special{ar 3839 1772 394 394  6.1523890 6.1823890}%
\special{ar 3839 1772 394 394  6.2723890 6.3023890}%
\special{ar 3839 1772 394 394  6.3923890 6.4223890}%
\special{ar 3839 1772 394 394  6.5123890 6.5423890}%
\special{ar 3839 1772 394 394  6.6323890 6.6623890}%
\special{ar 3839 1772 394 394  6.7523890 6.7823890}%
\special{ar 3839 1772 394 394  6.8723890 6.9023890}%
\special{ar 3839 1772 394 394  6.9923890 7.0223890}%
\special{ar 3839 1772 394 394  7.1123890 7.1423890}%
\special{ar 3839 1772 394 394  7.2323890 7.2623890}%
\special{ar 3839 1772 394 394  7.3523890 7.3823890}%
\special{ar 3839 1772 394 394  7.4723890 7.5023890}%
\special{ar 3839 1772 394 394  7.5923890 7.6223890}%
\special{ar 3839 1772 394 394  7.7123890 7.7423890}%
\special{ar 3839 1772 394 394  7.8323890 7.8623890}%
\special{ar 3839 1772 394 394  7.9523890 7.9823890}%
\special{ar 3839 1772 394 394  8.0723890 8.1023890}%
\special{ar 3839 1772 394 394  8.1923890 8.2223890}%
\special{ar 3839 1772 394 394  8.3123890 8.3423890}%
% BOX 2 5 2 0
% 2 3200 1750 3700 1850
% 
\special{pn 8}%
\special{sh 0}%
\special{pa 3150 1723}%
\special{pa 3642 1723}%
\special{pa 3642 1821}%
\special{pa 3150 1821}%
\special{pa 3150 1723}%
\special{ip}%
% LINE 2 0 3 0
% 2 3300 1750 3610 1750
% 
\special{pn 8}%
\special{pa 3249 1723}%
\special{pa 3554 1723}%
\special{fp}%
% LINE 2 0 3 0
% 2 3300 1850 3600 1850
% 
\special{pn 8}%
\special{pa 3249 1821}%
\special{pa 3544 1821}%
\special{fp}%
\end{picture}%

%% file: Ader-l.tex
%WinTpicVersion3.08
\unitlength 0.1in
\begin{picture}( 29.8720,  7.8740)( 12.4508,-21.6535)
% CIRCLE 2 0 3 0
% 4 3900 1800 3600 1800 3600 1800 3600 1800
% 
\special{pn 8}%
\special{ar 3839 1772 296 296  0.0000000 6.2831853}%
% STR 2 0 3 0
% 3 3900 1700 3900 1800 5 0
% $A$
\put(38.3858,-17.7165){\makebox(0,0){$A$}}%
% DOT 1 0 3 0
% 2 3600 1800 3600 1800
% 
\special{pn 13}%
\special{sh 1}%
\special{ar 3544 1772 10 10 0  6.28318530717959E+0000}%
\special{sh 1}%
\special{ar 3544 1772 10 10 0  6.28318530717959E+0000}%
% DOT 1 0 3 0
% 2 3750 1540 3750 1540
% 
\special{pn 13}%
\special{sh 1}%
\special{ar 3691 1516 10 10 0  6.28318530717959E+0000}%
\special{sh 1}%
\special{ar 3691 1516 10 10 0  6.28318530717959E+0000}%
% DOT 1 0 3 0
% 2 3640 1650 3640 1650
% 
\special{pn 13}%
\special{sh 1}%
\special{ar 3583 1625 10 10 0  6.28318530717959E+0000}%
\special{sh 1}%
\special{ar 3583 1625 10 10 0  6.28318530717959E+0000}%
% DOT 1 0 3 0
% 2 3640 1950 3640 1950
% 
\special{pn 13}%
\special{sh 1}%
\special{ar 3583 1920 10 10 0  6.28318530717959E+0000}%
\special{sh 1}%
\special{ar 3583 1920 10 10 0  6.28318530717959E+0000}%
% CIRCLE 2 2 3 0
% 4 3900 1800 3900 1400 2770 3760 3900 1000
% 
\special{pn 8}%
\special{ar 3839 1772 394 394  4.7123890 4.7423890}%
\special{ar 3839 1772 394 394  4.8323890 4.8623890}%
\special{ar 3839 1772 394 394  4.9523890 4.9823890}%
\special{ar 3839 1772 394 394  5.0723890 5.1023890}%
\special{ar 3839 1772 394 394  5.1923890 5.2223890}%
\special{ar 3839 1772 394 394  5.3123890 5.3423890}%
\special{ar 3839 1772 394 394  5.4323890 5.4623890}%
\special{ar 3839 1772 394 394  5.5523890 5.5823890}%
\special{ar 3839 1772 394 394  5.6723890 5.7023890}%
\special{ar 3839 1772 394 394  5.7923890 5.8223890}%
\special{ar 3839 1772 394 394  5.9123890 5.9423890}%
\special{ar 3839 1772 394 394  6.0323890 6.0623890}%
\special{ar 3839 1772 394 394  6.1523890 6.1823890}%
\special{ar 3839 1772 394 394  6.2723890 6.3023890}%
\special{ar 3839 1772 394 394  6.3923890 6.4223890}%
\special{ar 3839 1772 394 394  6.5123890 6.5423890}%
\special{ar 3839 1772 394 394  6.6323890 6.6623890}%
\special{ar 3839 1772 394 394  6.7523890 6.7823890}%
\special{ar 3839 1772 394 394  6.8723890 6.9023890}%
\special{ar 3839 1772 394 394  6.9923890 7.0223890}%
\special{ar 3839 1772 394 394  7.1123890 7.1423890}%
\special{ar 3839 1772 394 394  7.2323890 7.2623890}%
\special{ar 3839 1772 394 394  7.3523890 7.3823890}%
\special{ar 3839 1772 394 394  7.4723890 7.5023890}%
\special{ar 3839 1772 394 394  7.5923890 7.6223890}%
\special{ar 3839 1772 394 394  7.7123890 7.7423890}%
\special{ar 3839 1772 394 394  7.8323890 7.8623890}%
\special{ar 3839 1772 394 394  7.9523890 7.9823890}%
\special{ar 3839 1772 394 394  8.0723890 8.1023890}%
\special{ar 3839 1772 394 394  8.1923890 8.2223890}%
\special{ar 3839 1772 394 394  8.3123890 8.3423890}%
% BOX 2 5 2 0
% 2 3200 1750 3700 1850
% 
\special{pn 8}%
\special{sh 0}%
\special{pa 3150 1723}%
\special{pa 3642 1723}%
\special{pa 3642 1821}%
\special{pa 3150 1821}%
\special{pa 3150 1723}%
\special{ip}%
% LINE 2 0 3 0
% 2 3290 1750 3600 1750
% 
\special{pn 8}%
\special{pa 3239 1723}%
\special{pa 3544 1723}%
\special{fp}%
% LINE 2 0 3 0
% 2 3300 1850 3600 1850
% 
\special{pn 8}%
\special{pa 3249 1821}%
\special{pa 3544 1821}%
\special{fp}%
% LINE 2 0 3 0
% 2 1400 1800 2800 1800
% 
\special{pn 8}%
\special{pa 1378 1772}%
\special{pa 2756 1772}%
\special{fp}%
% DOT 1 0 3 0
% 6 1400 1800 1550 1800 2150 1800 2650 1800 2800 1800 2800 1800
% 
\special{pn 13}%
\special{sh 1}%
\special{ar 1378 1772 10 10 0  6.28318530717959E+0000}%
\special{sh 1}%
\special{ar 1526 1772 10 10 0  6.28318530717959E+0000}%
\special{sh 1}%
\special{ar 2117 1772 10 10 0  6.28318530717959E+0000}%
\special{sh 1}%
\special{ar 2609 1772 10 10 0  6.28318530717959E+0000}%
\special{sh 1}%
\special{ar 2756 1772 10 10 0  6.28318530717959E+0000}%
\special{sh 1}%
\special{ar 2756 1772 10 10 0  6.28318530717959E+0000}%
% STR 2 0 3 0
% 3 1400 1800 1400 1900 5 0
% $1$
\put(13.7795,-18.7008){\makebox(0,0){$1$}}%
% STR 2 0 3 0
% 3 1550 1800 1550 1900 5 0
% $2$
\put(15.2559,-18.7008){\makebox(0,0){$2$}}%
% STR 2 0 3 0
% 3 2150 1800 2150 1900 5 0
% $k$
\put(21.1614,-18.7008){\makebox(0,0){$k$}}%
% STR 2 0 3 0
% 3 2800 1800 2800 1900 5 0
% $n$
\put(27.5591,-18.7008){\makebox(0,0){$n$}}%
% CIRCLE 2 0 3 0
% 4 2720 2300 3260 1760 3290 1800 2200 1730
% 
\special{pn 8}%
\special{ar 2678 2264 752 752  3.9728302 5.5631146}%
% SARROW 2 0 3 1
% 2 2216 1726 2205 1736
% 
\special{pn 8}%
\special{pa 2182 1699}%
\special{pa 2171 1709}%
\special{fp}%
\special{sh 1}%
\special{pa 2171 1709}%
\special{pa 2233 1680}%
\special{pa 2209 1674}%
\special{pa 2206 1650}%
\special{pa 2171 1709}%
\special{fp}%
% STR 2 0 3 0
% 3 2100 1400 2100 1500 5 0
% $B$
\put(20.6693,-14.7638){\makebox(0,0){$B$}}%
% LINE 2 2 3 0
% 2 1700 1900 2000 1900
% 
\special{pn 8}%
\special{pa 1674 1871}%
\special{pa 1969 1871}%
\special{dt 0.045}%
% LINE 2 2 3 0
% 2 2300 1900 2500 1900
% 
\special{pn 8}%
\special{pa 2264 1871}%
\special{pa 2461 1871}%
\special{dt 0.045}%
\end{picture}%

%% file: Ader-r.tex
%WinTpicVersion3.08
\unitlength 0.1in
\begin{picture}( 15.1083, 10.8268)( 12.4508,-18.7008)
% CIRCLE 2 0 3 0
% 4 2150 1200 2150 1500 2150 1500 2150 1500
% 
\special{pn 8}%
\special{ar 2117 1182 296 296  0.0000000 6.2831853}%
% STR 2 0 3 0
% 3 2050 1200 2150 1200 5 0
% $A$
\put(21.1614,-11.8110){\special{rt 0 0  4.7124}\makebox(0,0){$A$}}%
\special{rt 0 0 0}%
% DOT 1 0 3 0
% 2 2150 1500 2150 1500
% 
\special{pn 13}%
\special{sh 1}%
\special{ar 2117 1477 10 10 0  6.28318530717959E+0000}%
\special{sh 1}%
\special{ar 2117 1477 10 10 0  6.28318530717959E+0000}%
% DOT 1 0 3 0
% 2 1890 1350 1890 1350
% 
\special{pn 13}%
\special{sh 1}%
\special{ar 1861 1329 10 10 0  6.28318530717959E+0000}%
\special{sh 1}%
\special{ar 1861 1329 10 10 0  6.28318530717959E+0000}%
% DOT 1 0 3 0
% 2 2000 1460 2000 1460
% 
\special{pn 13}%
\special{sh 1}%
\special{ar 1969 1438 10 10 0  6.28318530717959E+0000}%
\special{sh 1}%
\special{ar 1969 1438 10 10 0  6.28318530717959E+0000}%
% DOT 1 0 3 0
% 2 2300 1460 2300 1460
% 
\special{pn 13}%
\special{sh 1}%
\special{ar 2264 1438 10 10 0  6.28318530717959E+0000}%
\special{sh 1}%
\special{ar 2264 1438 10 10 0  6.28318530717959E+0000}%
% CIRCLE 2 2 3 0
% 4 2150 1200 1750 1200 4110 2330 1350 1200
% 
\special{pn 8}%
\special{ar 2117 1182 394 394  3.1415927 3.1715927}%
\special{ar 2117 1182 394 394  3.2615927 3.2915927}%
\special{ar 2117 1182 394 394  3.3815927 3.4115927}%
\special{ar 2117 1182 394 394  3.5015927 3.5315927}%
\special{ar 2117 1182 394 394  3.6215927 3.6515927}%
\special{ar 2117 1182 394 394  3.7415927 3.7715927}%
\special{ar 2117 1182 394 394  3.8615927 3.8915927}%
\special{ar 2117 1182 394 394  3.9815927 4.0115927}%
\special{ar 2117 1182 394 394  4.1015927 4.1315927}%
\special{ar 2117 1182 394 394  4.2215927 4.2515927}%
\special{ar 2117 1182 394 394  4.3415927 4.3715927}%
\special{ar 2117 1182 394 394  4.4615927 4.4915927}%
\special{ar 2117 1182 394 394  4.5815927 4.6115927}%
\special{ar 2117 1182 394 394  4.7015927 4.7315927}%
\special{ar 2117 1182 394 394  4.8215927 4.8515927}%
\special{ar 2117 1182 394 394  4.9415927 4.9715927}%
\special{ar 2117 1182 394 394  5.0615927 5.0915927}%
\special{ar 2117 1182 394 394  5.1815927 5.2115927}%
\special{ar 2117 1182 394 394  5.3015927 5.3315927}%
\special{ar 2117 1182 394 394  5.4215927 5.4515927}%
\special{ar 2117 1182 394 394  5.5415927 5.5715927}%
\special{ar 2117 1182 394 394  5.6615927 5.6915927}%
\special{ar 2117 1182 394 394  5.7815927 5.8115927}%
\special{ar 2117 1182 394 394  5.9015927 5.9315927}%
\special{ar 2117 1182 394 394  6.0215927 6.0515927}%
\special{ar 2117 1182 394 394  6.1415927 6.1715927}%
\special{ar 2117 1182 394 394  6.2615927 6.2915927}%
\special{ar 2117 1182 394 394  6.3815927 6.4115927}%
\special{ar 2117 1182 394 394  6.5015927 6.5315927}%
\special{ar 2117 1182 394 394  6.6215927 6.6515927}%
\special{ar 2117 1182 394 394  6.7415927 6.7715927}%
% LINE 2 0 3 0
% 2 2200 1800 2200 1500
% 
\special{pn 8}%
\special{pa 2166 1772}%
\special{pa 2166 1477}%
\special{fp}%
% LINE 2 0 3 0
% 2 1400 1800 2800 1800
% 
\special{pn 8}%
\special{pa 1378 1772}%
\special{pa 2756 1772}%
\special{fp}%
% DOT 1 0 3 0
% 6 1400 1800 1550 1800 2150 1800 2650 1800 2800 1800 2800 1800
% 
\special{pn 13}%
\special{sh 1}%
\special{ar 1378 1772 10 10 0  6.28318530717959E+0000}%
\special{sh 1}%
\special{ar 1526 1772 10 10 0  6.28318530717959E+0000}%
\special{sh 1}%
\special{ar 2117 1772 10 10 0  6.28318530717959E+0000}%
\special{sh 1}%
\special{ar 2609 1772 10 10 0  6.28318530717959E+0000}%
\special{sh 1}%
\special{ar 2756 1772 10 10 0  6.28318530717959E+0000}%
\special{sh 1}%
\special{ar 2756 1772 10 10 0  6.28318530717959E+0000}%
% STR 2 0 3 0
% 3 1400 1800 1400 1900 5 0
% $1$
\put(13.7795,-18.7008){\makebox(0,0){$1$}}%
% STR 2 0 3 0
% 3 1550 1800 1550 1900 5 0
% $2$
\put(15.2559,-18.7008){\makebox(0,0){$2$}}%
% STR 2 0 3 0
% 3 2800 1800 2800 1900 5 0
% $n$
\put(27.5591,-18.7008){\makebox(0,0){$n$}}%
% POLYGON 2 5 2 0
% 5 2100 1900 2200 1900 2200 1400 2100 1400 2100 1900
% 
\special{pn 8}%
\special{sh 0}%
\special{pa 2067 1871}%
\special{pa 2166 1871}%
\special{pa 2166 1378}%
\special{pa 2067 1378}%
\special{pa 2067 1871}%
\special{ip}%
% LINE 2 0 3 0
% 2 2100 1800 2100 1490
% 
\special{pn 8}%
\special{pa 2067 1772}%
\special{pa 2067 1467}%
\special{fp}%
% LINE 2 2 3 0
% 2 1700 1900 2000 1900
% 
\special{pn 8}%
\special{pa 1674 1871}%
\special{pa 1969 1871}%
\special{dt 0.045}%
% LINE 2 2 3 0
% 2 2300 1900 2500 1900
% 
\special{pn 8}%
\special{pa 2264 1871}%
\special{pa 2461 1871}%
\special{dt 0.045}%
% STR 2 0 3 0
% 3 2150 1800 2150 1900 5 0
% $k$
\put(21.1614,-18.7008){\makebox(0,0){$k$}}%
% STR 2 0 3 0
% 3 1750 1550 1750 1650 5 0
% $B$
\put(17.2244,-16.2402){\makebox(0,0){$B$}}%
\end{picture}%